%% file: main.tex
\documentclass[11pt,letterpaper]{article}
\usepackage[T1]{fontenc}
\usepackage[utf8]{inputenc}
\usepackage{newunicodechar}
\newunicodechar{’}{\textquoteright}
\usepackage[margin=1in]{geometry}
\usepackage{amsmath,amssymb,amsthm,mathtools,bm}
\usepackage{graphicx,booktabs,tabularx,array,enumitem,needspace,xcolor}
\usepackage{aliascnt}
\usepackage[colorlinks=true,allcolors=blue,pagebackref=true]{hyperref}
\renewcommand*{\backref}[1]{}
\renewcommand*{\backrefalt}[4]{\ifcase#1\or Cited on page~#2.\else Cited on pages~#2.\fi
}
\usepackage[nameinlink,noabbrev]{cleveref}

\renewenvironment{abstract}{\begin{center}\bfseries Abstract\end{center}\begin{quotation}\noindent
}{\end{quotation}}

\newtheorem{theorem}{Theorem}[section]
\newaliascnt{lemma}{theorem}
\newtheorem{lemma}[lemma]{Lemma}
\aliascntresetthe{lemma}
\newaliascnt{proposition}{theorem}
\newtheorem{proposition}[proposition]{Proposition}
\aliascntresetthe{proposition}
\newaliascnt{corollary}{theorem}
\newtheorem{corollary}[corollary]{Corollary}
\aliascntresetthe{corollary}
\theoremstyle{definition}
\newaliascnt{definition}{theorem}
\newtheorem{definition}[definition]{Definition}
\aliascntresetthe{definition}
\newaliascnt{assumption}{theorem}

\aliascntresetthe{assumption}
\theoremstyle{remark}
\newaliascnt{remark}{theorem}

\aliascntresetthe{remark}

\newcommand{\R}{\mathbb{R}}
\newcommand{\N}{\mathbb{N}}

\DeclarePairedDelimiter{\norm}{\lVert}{\rVert}

\newcommand{\D}{\nabla}
\newcommand{\ip}[2]{\left\langle#1,#2\right\rangle}
\newcommand{\Breg}[3]{D_{#1}\left(#2,#3\right)}
\newcommand{\resultpart}[2]{\ref{#1}\,\ref{#2}}
\makeatletter
\newcommand{\appendixcontents}{\@starttoc{apx}}
\newcommand{\recordappendixcontents}{\let\paper@addtocontents\addtocontents
  \renewcommand{\addtocontents}[2]{\def\paper@contentsfile{##1}\def\paper@maincontentsfile{toc}\ifx\paper@contentsfile\paper@maincontentsfile
      \paper@addtocontents{apx}{##2}\else
      \paper@addtocontents{##1}{##2}\fi
  }}
\makeatother
\theoremstyle{definition}
\newaliascnt{condition}{theorem}
\newtheorem{condition}[condition]{Design condition}
\aliascntresetthe{condition}
\newlist{claims}{enumerate}{1}
\setlist[claims]{label=\textbf{(\alph*)},ref=(\alph*),leftmargin=1.7em,
  topsep=3pt,itemsep=2pt,parsep=0pt}
\numberwithin{equation}{section}
\allowdisplaybreaks[2]
\setlist{topsep=3pt,itemsep=2pt,parsep=0pt,leftmargin=*}

\title{\rule{\linewidth}{1.5pt}\\[3pt]
  \textbf{Tight Lower Bounds under Matched High-Order
  H\"older Smoothness and Uniform Convexity}\\[3pt]
  \rule[8pt]{\linewidth}{1pt}}
\author{\normalsize\textbf{Cedar Site Bai}\\
  \normalsize Department of Computer Science\\
  \normalsize Purdue University\\
  \normalsize\href{mailto:best990317@gmail.com}{best990317@gmail.com}}
\date{}
\hypersetup{
  pdftitle={Tight Lower Bounds under Matched High-Order
    H\"older Smoothness and Uniform Convexity},
  pdfauthor={Cedar Site Bai}
}

\begin{document}
\hypersetup{pageanchor=false}
\maketitle
\begin{abstract}
\input{abstract-author}
\end{abstract}
\thispagestyle{empty}
\clearpage 
\setcounter{page}{1}
\hypersetup{pageanchor=true}

\begingroup
\setlength{\parskip}{2pt}
\setlength{\abovedisplayskip}{3pt plus 2pt minus 1pt}
\setlength{\belowdisplayskip}{3pt plus 2pt minus 1pt}
\setlength{\abovedisplayshortskip}{0pt plus 2pt}
\setlength{\belowdisplayshortskip}{3pt plus 2pt minus 1pt}
\setlength{\textfloatsep}{10pt plus 2pt minus 2pt}
\setlength{\floatsep}{10pt plus 2pt minus 2pt}
\setlength{\intextsep}{8pt plus 2pt minus 2pt}
\setlength{\topsep}{4pt plus 2pt minus 1pt}
\setlength{\jot}{2pt}
\setlist{topsep=2pt,itemsep=1pt,parsep=0pt,leftmargin=*}
\setlist[claims]{label=\textbf{(\alph*)},ref=(\alph*),leftmargin=1.7em,
  topsep=2pt,itemsep=1pt,parsep=0pt}
\makeatletter
\renewcommand\section{\@startsection{section}{1}{\z@}{-2.75ex \@plus -.8ex \@minus -.2ex}{1.4ex \@plus .2ex}{\normalfont\Large\bfseries}}
\renewcommand\subsection{\@startsection{subsection}{2}{\z@}{-2.5ex \@plus -.8ex \@minus -.2ex}{1ex \@plus .2ex}{\normalfont\large\bfseries}}
\renewcommand\paragraph{\@startsection{paragraph}{4}{\z@}{1.75ex \@plus .5ex \@minus .2ex}{-1em}{\normalfont\normalsize\bfseries}}
\makeatother

\input{manuscript-introduction}

\input{manuscript-related-work}

\input{manuscript-preliminaries}

\section{Lower Bound for the \texorpdfstring{$q=p+\nu$}{q=p+nu} Case}
\label{paper:lower-bound}
We prove Theorem~\ref{paper:main-theorem} by constructing a hard
function whose successive phases each require additional
oracle calls. Following the separation of geometry and hardness
in Sec.~1.3, we first build a uniformly convex base and then
add convex phase functions. Combining these components requires
the smooth maximum introduced below.
\input{manuscript-smoothing}
\subsection{Step 1: Construct the Uniformly Convex Base}\label{paper:construction}
To force the condition-number cost at successive accuracy scales,
we construct $L\ge1$ phases. We first give the base $B$ the same
local uniform-convexity geometry at each phase.

\input{manuscript-base}
\input{manuscript-phases}

\input{manuscript-oracle-rate}
\section{Conclusion and Future Work}
We establish a tight deterministic high-order oracle lower bound
of $\Omega_{p,\nu}\left(
\left(\frac{H}{\sigma}\right)^{\frac{2}{3q-2}}
\log\left(\frac{\Delta_0}{\epsilon}\right)\right)$ calls
in the matched case $q=p+\nu$.
The construction combines Euclidean curvature at successive centers
with short hard chains, while preserving earlier oracle outputs
and keeping the global class parameters independent of the number
of phases. Together with the asymmetric results, it resolves the
remaining equality case, with the stated function-class distinctions.
Randomized algorithms and matched lower bounds in other norm
geometries remain separate questions.
\label{paper:main-end}
\endgroup
\clearpage
\begingroup
\input{ai-statement}
\label{paper:ai-end}
\endgroup
\clearpage
\appendix
\begingroup
\recordappendixcontents
\input{appendix-guide}

\input{manuscript-supporting-proofs}

\endgroup
\clearpage
\begingroup
\raggedright
\bibliographystyle{alpha}
\bibliography{manuscript-references}
\endgroup
\end{document}

%% file: abstract-author.tex
In this paper, we provide tight lower bounds for minimizing high-order Hölder-smooth and uniformly convex functions when their orders match. Specifically, we consider functions whose $p^{th}$-order derivatives are Hölder continuous with degree $\nu$ and parameter $H$, and that are uniformly convex with degree $q$ and parameter $\sigma$. When these orders match, namely, $q = p + \nu$, we show that any deterministic algorithm with access to derivatives up to order $p$ requires $\Omega\left(\left(\frac{H}{\sigma}\right)^{\frac{2}{3(p+\nu)-2}}\log\frac{\Delta_0}{\epsilon}\right)$ oracle calls to reduce an initial optimality gap of $\Delta_0$ to $\epsilon$. This matches the corresponding upper bound and establishes a higher-order uniformly convex counterpart to Nesterov’s classical first-order lower bound for strongly convex functions. Rather than directly generalizing Nesterov’s classical quadratic-chain construction, we develop a multiscale phase-wise construction that preserves previous oracle responses while maintaining global Hölder smoothness and uniform convexity. Our result resolves the matched case left open by existing lower bounds for the two asymmetric regimes.

%% file: manuscript-introduction.tex
\begingroup
\linespread{0.9705882352941176}\selectfont
\setlength{\parskip}{2pt}
\section{Introduction}

Determining the optimal complexity of minimizing uniformly convex
functions is a fundamental problem in optimization.
Uniform convexity generalizes strong convexity through curvature
bounds involving higher powers of the distance between points, with
applications in high-order
optimization~\cite{Nesterov2021Tensor,Bullins2020}
and machine learning~\cite{ThomsenDoikov2024,Bullins2020,BaiBullins2025}.
The optimal complexity depends jointly on the order of derivative
oracle access, the smoothness of the function, and its
uniform-convexity degree.

We study this problem in a general setting, with $p^{th}$-order
derivative oracle access to functions $F\left(\mathbf{x}\right)$ whose $p^{th}$-order derivatives $\nabla^p F\left(\mathbf{x}\right)$
are H\"older continuous with exponent $\nu$ and parameter $H$,
and that are uniformly convex with degree $q$ and modulus $\sigma$. We focus on the matched regime
$q=p+\nu$, which includes the classical smooth strongly convex case when
$p=1,\nu=1,q=2$, in which $\frac{H}{\sigma}$ is the condition number.
In this case, Nesterov's accelerated gradient
method~\cite{Nesterov1983AcceleratedGradient} converges linearly,
attaining an optimality gap $F\left(\bm x_T\right)-F\left(\bm x^\star\right)\le\epsilon$ in
$T\in\mathcal{O}\left(\left(\frac{H}{\sigma}\right)^{\frac{1}{2}}\log\left(\frac{\Delta_0}{\epsilon}\right)\right)$ gradient oracle
calls~\cite{Nesterov2018,Bubeck2015ConvexOptimization}, where
$\bm x^\star$ is the minimizer and
$\Delta_0:=F\left(\bm x_0\right)-F\left(\bm x^\star\right)$ denotes the initial optimality gap.

On the upper-bound side, \cite{SongJiangMa2019Preprint} establishes linear convergence in the
matched higher-order regime $q=p+\nu>2$, with iteration complexity
\phantomsection\label{paper:upper-rate}$O_{p,\nu}\!\left(
 \left(\frac{H}{\sigma}\right)^{\frac{2}{3q-2}}
 \log\left(\frac{\Delta_0}{\epsilon}\right)\right)$.
These methods use derivatives through order $p$, and
$\frac{H}{\sigma}$ is the high-order condition
number~\cite{GasnikovEtAl2018Full,Bullins2020,DoikovNesterov2021}.
The logarithmic accuracy dependence follows from a restarting
mechanism that reapplies an accelerated convex method at smaller
accuracy scales, reducing the optimality gap by a constant factor
at each stage~\cite{GasnikovEtAl2018Full,SongJiangMa2019Preprint}.

On the lower-bound side, Nesterov's classical first-order
construction~\cite{Nesterov2018} establishes the optimality of the
first-order rate above, while a tight characterization for the
matched higher-order setting has remained open.

For general combinations of
$p,q,\nu$, \cite{BaiBullins2025} establishes tight
deterministic lower bounds for both asymmetric regimes, with
polynomial dependence on $\frac{1}{\epsilon}$ when $q>p+\nu$ and
doubly logarithmic dependence on $\frac{1}{\epsilon}$ when $q<p+\nu$,
matching the corresponding upper bounds~\cite{SongJiangMa2019Preprint}.
These complementary results leave $q=p+\nu$ as the only missing
case in the general complexity characterization, as summarized
in Table~\ref{paper:complexity-comparison}, under the respective
function-class assumptions.

\subsection{Our Result}
In this paper, we resolve this open problem by establishing a
lower bound matching the \hyperref[paper:upper-rate]{upper bound} for
deterministic algorithms with oracle access to all derivatives
up to order $p$.
Our result establishes a higher-order uniformly convex counterpart
to Nesterov's classical first-order lower bound~\cite{Nesterov2018} and complements
the asymmetric results of~\cite{BaiBullins2025}.

\noindent\begin{minipage}{\linewidth}
\setlength{\parskip}{1.5pt}
\setlength{\topsep}{2pt}
\makeatletter
\def\thm@space@setup{\thm@preskip=2pt\thm@postskip=2pt}
\makeatother
\begin{theorem}[Matched-order lower bound]\label{paper:main-theorem}
\label{main:theorem}
Fix an integer $p\ge2$, $0<\nu\le1$, and $q=p+\nu$.
There are constants $c_{p,\nu},K_{p,\nu}>0$ such that, for
$H,\sigma,\Delta_0>0$ with sufficiently large $\frac{H}{\sigma}$ and
$0<\epsilon\le c_{p,\nu}\Delta_0$, and every dimension
$d\ge K_{p,\nu}\left(\frac{H}{\sigma}\right)^{\frac{2}{3q-2}}
\log\left(\frac{\Delta_0}{\epsilon}\right)$,
every deterministic $p^{th}$-order oracle algorithm
initialized at $\bm x_0\in\R^d$ admits a
$p$ times continuously differentiable function $F:\R^d\to\R$
with minimizer $\bm x^\star$ and
$F\left(\bm x_0\right)-F\left(\bm x^\star\right)=\Delta_0$,
such that, for all
$\bm x,\bm y\in\R^d$:
\begin{claims}[topsep=1pt,itemsep=0pt,parsep=0pt]
\item\label{main:theorem:part:smoothness}
\textbf{$p^{th}$-order H\"older smoothness.}
$\norm{\D^pF\left(\bm y\right)-\D^pF\left(\bm x\right)}
\le H\norm{\bm y-\bm x}^{\nu}$.
\item\label{main:theorem:part:uc}
\textbf{Uniform convexity.}
$F\left(\bm y\right)\ge F\left(\bm x\right)
+\ip{\D F\left(\bm x\right)}{\bm y-\bm x}
+\frac{\sigma}{q}\norm{\bm y-\bm x}^{q}$.
\end{claims}
For a candidate $\widehat{\bm x}_T$ after $T$ oracle calls,
$F\left(\widehat{\bm x}_T\right)-F\left(\bm x^\star\right)\le\epsilon$ requires
$T\in\Omega_{p,\nu}\left(
 \left(\frac{H}{\sigma}\right)^{\frac{2}{3q-2}}
 \log\left(\frac{\Delta_0}{\epsilon}\right)\right)$.
\end{theorem}
\end{minipage}
\endgroup

\input{figures/multiphase-overview-figure}

\subsection{Motivation}
To prove an oracle lower bound, we construct, for any given
deterministic algorithm, a hard function for which reaching the
target optimality gap requires the number of oracle calls stated
in Theorem~\ref{paper:main-theorem}.
The construction limits the information that each oracle output
provides about the minimizer, so that fewer oracle calls cannot
identify an accurate solution.

\textbf{Classical Construction.}\quad
Nesterov's classical first-order lower-bound
construction~\cite{Nesterov2018} uses a quadratic
chain: a quadratic function whose terms couple neighboring
coordinates, such as $\left(x_i-x_{i+1}\right)^2$.
The coordinate directions can be chosen so that gradient oracle outputs
reveal the minimizer gradually. Its coordinates decay geometrically,
so the unrevealed tail prevents the algorithm from reaching the
target accuracy until sufficiently many directions have been
determined.
This construction does not directly extend to higher-order
oracles, because one Hessian oracle call reveals the entire quadratic
form.

\phantomsection\label{paper:norm-mismatch-motivation}
\textbf{Norm Mismatch.}\quad
A natural extension is therefore to replace the quadratic
terms with absolute $q^{th}$ powers, such as
$\left|x_i-x_{i+1}\right|^q$.
However, a long chain of scalar $q$-power residuals naturally
controls the $\ell_q$ norm, while our uniform-convexity condition
uses the Euclidean norm. For $q>2$,
$\norm{\bm x}_q^q\ge d^{\,1-\frac{q}{2}}\norm{\bm x}_2^q$.
This norm conversion introduces the dimension-dependent factor
$d^{\,1-\frac{q}{2}}$ into the Euclidean uniform-convexity modulus,
leading to a suboptimal lower-bound rate for this naive extension.\footnote{Appx.~\ref{app:paper-barrier} verifies this dimension loss for
the corresponding geometric $q$-power chain.}

\textbf{Successive Phases.}\quad
We therefore need a different construction to recover both
factors in the \hyperref[paper:upper-rate]{upper bound}.
To identify the condition-number cost that this construction
should preserve, we examine the asymmetric lower bound
in~\cite[Theorem~1]{BaiBullins2025}. For $q>p+\nu$, their rate is
$\Omega\left(
\left(\frac{H}{\sigma}\right)^{\frac{2}{3\left(p+\nu\right)-2}}
\left(\frac{\sigma}{\epsilon}\right)^{
\frac{2\left(q-p-\nu\right)}{q\left(3\left(p+\nu\right)-2\right)}}\right)$.
Although their theorem treats only $q>p+\nu$, formally substituting
$q=p+\nu$ in this expression gives
$\left(\frac{H}{\sigma}\right)^{\frac{2}{3q-2}}$.
This expression has the desired condition-number dependence but
is independent of $\epsilon$, and therefore misses the logarithmic
factor in the \hyperref[paper:upper-rate]{upper bound}.
In light of the restarting mechanism in the upper-bound
analyses~\cite{GasnikovEtAl2018Full,SongJiangMa2019Preprint},
we seek one hard function that forces this condition-number cost
at each constant-factor reduction of the optimality gap.
The successive-scale constructions used in earlier logarithmic
lower bounds~\cite{Nemirovski1994ComplexityNotes,
BraunGuzmanPokutta2017,DornEtAl2026GradientFree}
motivate arranging these repeated hard instances into phases
at decreasing accuracy scales.
Each phase must require a further
$\Omega_{p,\nu}\left(
\left(\frac{H}{\sigma}\right)^{\frac{2}{3q-2}}\right)$
calls even after the algorithm has observed all oracle outputs from
the preceding phases. To obtain the sharp product in
the \hyperref[paper:upper-rate]{upper bound}, we must realize these phases in
one objective without increasing $H$ or decreasing $\sigma$
with their number.

\subsection{Techniques and Proof Strategy}
Our technical novelty is a multi-center base $B$ that allows us
to separate Euclidean uniform convexity from the chains supplying
oracle hardness, avoiding the
\hyperref[paper:norm-mismatch-motivation]{norm mismatch}.
Figure~\ref{fig:multiphase-overview} illustrates this separation.
We construct $L$ successive phases, each with a short power
chain of length $m$, where
\begin{equation}\label{paper:phase-scales}
 m\asymp_{p,\nu}\left(\frac{H}{\sigma}\right)^{\frac{2}{3q-2}},\qquad
 L\asymp_{p,\nu}\log\left(\frac{\Delta_0}{\epsilon}\right).
\end{equation}
\textbf{Multi-phase Uniform Convexity.}\quad
For geometry, $B$ must reproduce the same local
$\smash{q^{th}}$-power behavior at each scale.
A fixed-center regularizer $\left\|\bm x\right\|_2^q$ has positive
quadratic curvature at nonzero centers, so it cannot provide
this behavior there.
We therefore construct $B$ around successive reference centers
$\bm c_0,\ldots,\bm c_L$, whose consecutive distances decrease by
a fixed factor to match the decreasing accuracy scales.
Near every center, the base has the exact form
\[
 B\left(\bm x\right)=B\left(\bm c_j\right)
 +\ip{\D B\left(\bm c_j\right)}{\bm x-\bm c_j}
 +\frac{\mu}{q}\left\|\bm x-\bm c_j\right\|_2^q
\]
with the same $\mu>0$.
Its higher-order derivatives through order $p$ consequently vanish
at every center, while $B$ remains globally high-order H\"older
smooth and Euclidean uniformly convex, with bounds independent of $L$.

\textbf{Hardness Construction.}\quad
For hardness, we add the short power chains and an initial
linear term giving a nonzero starting gradient.
To initiate the next phase, we join each chain to an affine
function near later centers. Its slope makes the complete objective's gradient at
the next center point along the next chain's starting direction,
while the full oracle output at that center does not reveal
the remaining chain directions.
Reducing the gap further therefore requires resolving another
short chain. Convexity preserves the base's uniform-convexity bound,
while orthogonality of nonadjacent phases prevents smoothness costs
from accumulating with $L$.
The local base replacements and affine continuations must also
retain their prescribed formulas and derivatives exactly.
An ordinary maximum preserves convexity but need not be smooth.
We therefore use a smooth maximum~\cite[Lemma~1]{Azagra2013} that
remains convex and equals the larger input whenever it exceeds
the other by at least the smoothing width.

\textbf{Oracle Consistency.}\quad
We enforce these restrictions for arbitrary deterministic
oracle algorithms by choosing unrevealed chain directions
as the algorithm makes its oracle calls, while preserving all earlier
oracle outputs. The final choices define one fixed objective
requiring $\Omega\left(m\right)$ calls per phase, hence
$\Omega\left(mL\right)$ calls in total.

%% file: figures/multiphase-overview-figure.tex
\begin{figure}[!t]
  \centering
  \includegraphics[width=\linewidth]{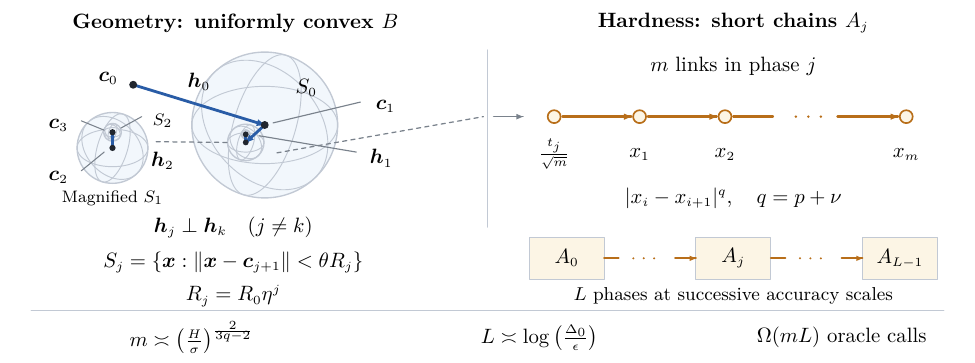}
  \setlength{\abovecaptionskip}{3pt}
  \caption{Geometry and hardness at successive accuracy scales (schematic).}
  \label{fig:multiphase-overview}
\end{figure}

%% file: manuscript-related-work.tex
\section{Related Work}

\paragraph{Upper Bounds.}
In the matched regime, \cite{GasnikovEtAl2018Full}
obtains the \hyperref[paper:upper-rate]{upper bound} in tensor iterations for
$\nu=1$, and \cite{SongJiangMa2019Preprint}
obtains the general H\"older-smooth iteration rate using a restarting
mechanism for accelerated high-order methods.
For $p=2$ and $q=2+\nu$, \cite{DoikovNesterov2021} establishes linear convergence
with the larger condition-number exponent $\frac{1}{1+\nu}$,
compared with the sharp $\frac{2}{4+3\nu}$.
\cite{ContrerasEtAl2024Full}
develops general non-Euclidean high-order oracle bounds.
Their Euclidean convex rate, combined with the restarting mechanism,
also gives the \hyperref[paper:upper-rate]{upper bound} in full $p^{th}$-order oracle calls.

\paragraph{Lower Bounds.}
Before this work, a tight lower bound for the matched regime
was known only in Nesterov's classical first-order setting
$p=1,\nu=1,q=2$~\cite{Nesterov2018}.
Beyond extending this classical result, our bound completes the
general $p,q,\nu$ characterization in
Table~\ref{paper:complexity-comparison}, complementing the
asymmetric results of~\cite{BaiBullins2025}.
Smaller entries record earlier special cases.

\begin{table}[!htbp]
\centering
\setlength{\tabcolsep}{0pt}
\setlength{\parskip}{0pt}
\renewcommand{\arraystretch}{0.92}
\caption{Complexity overview.
$\kappa:=\frac{H}{\sigma}$; $\Delta_0$: initial optimality gap;
$\epsilon$: target accuracy.}
\label{paper:complexity-comparison}
\newsavebox{\paperComparisonMeasure}
\newsavebox{\paperMatchedRateMeasure}
\newsavebox{\paperMatchedSpecialMeasure}
\newcommand{\paperSpecialCase}[2]{{\fontsize{9}{10}\selectfont
\renewcommand{\arraystretch}{0.85}\begin{tabular}{@{}c@{}}#1\\[-0.4pt]$#2$\end{tabular}}}
\newcommand{\paperMatchedUpperCases}{\paperSpecialCase{\cite{GasnikovEtAl2018Full}}{p\ge2,\nu=1}\hspace{8pt}\paperSpecialCase{\cite{Nesterov2018}}{p=1,\nu=1,q=2}}
\newcommand{\paperMatchedLowerCases}{\paperSpecialCase{\cite{Nesterov2018}}{p=1,\nu=1,q=2}}
\newcommand{\paperComparisonRows}[1]{\begin{tabular}{@{}l@{\hspace{8pt}}c@{\hspace{8pt}}c@{\hspace{8pt}}c@{\hspace{8pt}}c@{\hspace{8pt}}c@{}}
\toprule
\textbf{Regime} & \textbf{Tight rate} $\Theta\left(\cdot\right)$ &
\textbf{Upper} & \textbf{Special cases} &
\textbf{Lower} & \textbf{Special cases}\\
\cmidrule(lr){1-2}\cmidrule(lr){3-4}\cmidrule(lr){5-6}
$q>p+\nu$ &
$\displaystyle\kappa^{\frac{2}{3\left(p+\nu\right)-2}}
 \left(\frac{\sigma}{\epsilon}\right)^{
 \frac{2\left(q-p-\nu\right)}{q\left(3\left(p+\nu\right)-2\right)}}$ &
\cite{SongJiangMa2019Preprint} &
\paperSpecialCase{\cite{RouletDAspremont2017}}{p=1,\nu=1} &
\cite{BaiBullins2025} &
\paperSpecialCase{\cite{ThomsenDoikov2024}}{p=1,\nu=1,q>2}\\
\midrule
#1{\boldmath$q=p+\nu$} &
{\boldmath$\displaystyle
\kappa^{\frac{2}{3q-2}}\log\left(\frac{\Delta_0}{\epsilon}\right)$} &
\cite{SongJiangMa2019Preprint} & \paperMatchedUpperCases &
\textbf{This work} & \paperMatchedLowerCases\\
\midrule
$q<p+\nu$ &
$\displaystyle
 \kappa^{\frac{2}{3\left(p+\nu\right)-2}}
 +\log\log\left(
 \left(\frac{\sigma^{p+\nu}}{H^q}\right)^{\frac{1}{p+\nu-q}}
 \frac{1}{\epsilon}\right)$ &
\cite{SongJiangMa2019Preprint} &
\paperSpecialCase{$\widetilde O$~\cite{ArjevaniShamirShiff2019}}{p=2,\nu=1,q=2}\hspace{8pt}\paperSpecialCase{$\widetilde O$~\cite{KornowskiShamir2020}}{p\ge2,\nu=1,q=2} &
\cite{BaiBullins2025} &
\paperSpecialCase{\cite{ArjevaniShamirShiff2019}}{p=2,\nu=1,q=2}\hspace{8pt}\paperSpecialCase{\cite{KornowskiShamir2020}}{p\ge2,\nu=1,q=2}\\
\bottomrule
\end{tabular}}
\sbox{\paperComparisonMeasure}{\paperComparisonRows{}}
\sbox{\paperMatchedRateMeasure}{\strut{\boldmath$\displaystyle
\kappa^{\frac{2}{3q-2}}\log\left(\frac{\Delta_0}{\epsilon}\right)$}}
\sbox{\paperMatchedSpecialMeasure}{\paperMatchedUpperCases}
\ifdim\ht\paperMatchedSpecialMeasure>\ht\paperMatchedRateMeasure
  \ht\paperMatchedRateMeasure=\ht\paperMatchedSpecialMeasure
\fi
\ifdim\dp\paperMatchedSpecialMeasure>\dp\paperMatchedRateMeasure
  \dp\paperMatchedRateMeasure=\dp\paperMatchedSpecialMeasure
\fi
\resizebox{\linewidth}{!}{\paperComparisonRows{\rlap{\raisebox{-\dimexpr\dp\paperMatchedRateMeasure+0.6pt\relax}[0pt][0pt]{\color{gray!10}\rule{\wd\paperComparisonMeasure}{\dimexpr\ht\paperMatchedRateMeasure+\dp\paperMatchedRateMeasure+1.2pt\relax}}}}}
\end{table}

%% file: manuscript-preliminaries.tex
\section{Notations}\label{paper:settings}

We use bold lowercase letters for vectors (e.g., $\bm x,\bm y$),
with subscripts for vectors at different iterations
(e.g., $\bm x_T$). We use regular lowercase letters for scalars,
and $x_i$ for the $i^{th}$ coordinate of $\bm x$.
We use $\norm{\cdot}$ to denote the Euclidean norm for vectors
and the induced operator norm for multilinear forms, and
$\ip{\cdot}{\cdot}$ for the Euclidean inner product.
For a $k$-linear form $\mathcal T$ on $\R^d$, its operator norm is
$\norm{\mathcal T}:=\sup_{\norm{\bm z_i}\le1,\,1\le i\le k}
\left|\mathcal T\left[\bm z_1,\ldots,\bm z_k\right]\right|$.
We use $\D F\left(\bm x\right)$ for the gradient vector and
$\D^kF\left(\bm x\right)$ for the $k^{th}$ derivative of $F$ at
$\bm x$, a symmetric $k$-linear form for $k\ge2$. We set $\D^0F:=F$.
For a differentiable $F:\R^d\to\R$, we denote its Bregman divergence by
$D_F\left(\bm y,\bm x\right)
:=F\left(\bm y\right)-F\left(\bm x\right)
-\ip{\D F\left(\bm x\right)}{\bm y-\bm x}$
for $\bm x,\bm y\in\R^d$.
We use $C^p\left(\R^d\right)$ for the $p$ times continuously
differentiable functions and $C^{p,\nu}\left(\R^d\right)$ for those
whose $p^{th}$ derivatives are globally $\nu$-H\"older continuous.
We use $I$ for the identity matrix and $P_E$ for the orthogonal
projection onto a linear subspace $E$.
For pairwise orthogonal subspaces $E_1,\ldots,E_n$, i.e.,
$\ip{\bm u_i}{\bm u_j}=0$ whenever
$\bm u_i\in E_i$, $\bm u_j\in E_j$, and $i\ne j$,
we denote their orthogonal direct sum by
$\bigoplus_{i=1}^nE_i
:=\left\{\sum_{i=1}^n\bm u_i:\bm u_i\in E_i\right\}$, with
$\norm{\sum_{i=1}^n\bm u_i}^2=\sum_{i=1}^n\norm{\bm u_i}^2$.
We denote the closure of a set $S$ by $\overline S$.
For a nonempty set $S$, we denote
$\operatorname{dist}\left(\bm x,S\right):=\inf_{\bm z\in S}\norm{\bm x-\bm z}$.

%% file: manuscript-smoothing.tex
\subsection{Smooth Maximum}\label{paper:smoothing}

A standard lower-bound construction takes the pointwise maximum
of convex functions to restrict the information revealed by each
oracle output. Since this maximum need not be differentiable,
a smoothing operator supplies the required derivative regularity
while preserving that restriction~\cite{GuzmanNemirovski2013Preprint,
AgarwalHazan2018,GargEtAl2021,BaiBullins2025}.
Our construction additionally requires exact agreement with the
base and chain formulas on specified regions, to preserve their
derivatives and earlier oracle outputs. Averaging over nearby
points can alter nonlinear norm-power formulas even where one
input dominates. We therefore use the smooth maximum~\cite{Azagra2013}, which
smoothly compares input values and equals the larger
input whenever their difference exceeds the smoothing width.

\begin{definition}[Smooth maximum]\label{paper:smooth-maximum-definition}
We use the even smooth probability density
$\varphi\left(u\right)=c_\varphi\exp\left(-\frac{1}{1-u^2}\right)$
for $|u|<1$ and $\varphi\left(u\right)=0$ otherwise, with
$c_\varphi>0$ chosen so that
$\int_{-1}^1\varphi\left(u\right)\,du=1$. For
$t_+:=\max\left\{t,0\right\}$ and smoothing width $\delta>0$, define
\begin{equation}\label{paper:smooth-maximum}
 \psi_\delta\left(t\right)
 :=\int_{-1}^1\left(t-\delta u\right)_+\varphi\left(u\right)\,du,
 \qquad
 M_\delta\left(a,b\right):=a+\psi_\delta\left(b-a\right)
 \quad\left(a,b\in\R\right).
\end{equation}
\end{definition}
\begin{lemma}[{\cite[Lemma~1]{Azagra2013}}]
\label{paper:smooth-maximum-properties}
For every $\delta>0$, the function $M_\delta$ in
Definition~\ref{paper:smooth-maximum-definition} is $C^\infty$,
convex, and nondecreasing in each argument. For all $a,b\in\R$,
\[
 \max\left\{a,b\right\}
 \le M_\delta\left(a,b\right)
 \le\max\left\{a,b\right\}+\frac{\delta}{2},
 \qquad
 M_\delta\left(a,b\right)=
 \begin{cases}
 a,&a-b\ge\delta,\\
 b,&b-a\ge\delta.
 \end{cases}
\]
\end{lemma}

%% file: manuscript-base.tex
\textbf{Step 1.1: Choose the phase centers.}
To give each phase the same local uniform-convexity geometry,
we anchor phase $j$ at a center $\bm c_j$, with endpoint
$\bm c_{j+1}$ and $\bm c_0=\bm 0$. We define its phase displacement by
$\bm h_j=\bm c_{j+1}-\bm c_j$ for $0\le j<L$.
Later chain directions must be chosen without altering earlier
chains, so we choose mutually orthogonal phase subspaces, with
each $\bm h_j$ in its corresponding subspace.
We will specify the coordinates of $\bm h_j$ in Step~2.

To impose Euclidean uniform convexity, we begin with
$B_0\left(\bm x\right)=\frac{\mu}{q}\norm{\bm x}^q$, with
$\mu>0$, the simplest norm-power function.
We then reproduce its geometry at every phase center through
successive local replacements. Let $B_j$ denote the base after $j$
replacements, and define $B=B_L$. We require
\begin{equation}\label{paper:base-local}
 B_j\left(\bm x\right)=B_j\left(\bm c_j\right)
 +\ip{\D B_j\left(\bm c_j\right)}{\bm x-\bm c_j}
 +\frac{\mu}{q}\norm{\bm x-\bm c_j}^q
\end{equation}
on a ball containing $\bm c_{j+1}$.
Here growth above the tangent is measured by the Euclidean
$q^{th}$ power, avoiding the norm mismatch. Since $q>p$, the
higher-order derivatives satisfy
$\D^kB_j\left(\bm c_j\right)=0$ for $2\le k\le p$.

\textbf{Step 1.2: Preserve the completed phase's gradient.}
At the end of phase $j$, the objective's gradient should point
only in the next phase's directions. The objective's gradient at $\bm c_j$,
the increment from the base, and the gradient of the hardness
function $A_j$, which we will construct in Step~2, must therefore
cancel in the completed phase's directions.
Define $R_j=\norm{\bm h_j}$. Differentiating
Eq.~\eqref{paper:base-local} at
$\bm c_{j+1}=\bm c_j+\bm h_j$ identifies the base increment as
$\D B_j\left(\bm c_{j+1}\right)-\D B_j\left(\bm c_j\right)
=\mu R_j^{q-2}\bm h_j$.
We will make $A_j$ equal an affine function near later centers,
so its gradient there is fixed. For the combined gradient to
remain canceled along $\bm h_j$, the base gradient component
along $\bm h_j$ must also stay unchanged.

To check whether this cancellation survives movement in the
next phase, we evaluate the gradient given by
Eq.~\eqref{paper:base-local} at
$\bm c_{j+1}+\bm z$ near the next center, with
$\bm z\ne\bm 0$ in the next phase's directions and hence
$\bm z\perp\bm h_j$. Where Eq.~\eqref{paper:base-local} holds,
Eq.~\eqref{base:old-gradient-increment} in Appx.~\ref{app:paper-base}
gives
\[
 \ip{\D B_j\left(\bm c_{j+1}+\bm z\right)
       -\D B_j\left(\bm c_{j+1}\right)}{\bm h_j}
 =\mu R_j^2
 \left(\left(R_j^2+\norm{\bm z}^2\right)^{\frac{q-2}{2}}
       -R_j^{q-2}\right)>0.
\]
The positive inner product shows that movement in the next
phase changes the base gradient component along $\bm h_j$.
Therefore, we need to replace Eq.~\eqref{paper:base-local}
with an alternative formula $N_j^B$, centered at $\bm c_{j+1}$,
which we define below in Eq.~\eqref{paper:base-next-local}.
Its gradient increment at $\bm c_{j+1}+\bm z$ should be
$\mu\norm{\bm z}^{q-2}\bm z$, orthogonal to $\bm h_j$ as required.
To make the smooth maximum equal $N_j^B$ near $\bm c_{j+1}$,
we give $N_j^B$ value $B_j\left(\bm c_{j+1}\right)+C_j$
there, with $C_j>0$, while retaining the endpoint gradient.
We therefore define
\begin{equation}\label{paper:base-next-local}
 N_j^B\left(\bm x\right)
 =B_j\left(\bm c_{j+1}\right)+C_j
   +\ip{\D B_j\left(\bm c_{j+1}\right)}{\bm x-\bm c_{j+1}}
   +\frac{\mu}{q}\norm{\bm x-\bm c_{j+1}}^q.
\end{equation}

\textbf{Step 1.3: Join the formulas near the next center.}
To preserve earlier values and derivatives, we confine
the change to a small ball
$S_j=\left\{\bm x:\norm{\bm x-\bm c_{j+1}}<\theta R_j\right\}$,
with radius parameter $\theta>0$. We denote its closure by
$\overline S_j=\left\{\bm x:\norm{\bm x-\bm c_{j+1}}\le\theta R_j\right\}$
and define
\begin{equation}\label{paper:base-reset}
 \begin{aligned}
 B_{j+1}\left(\bm x\right)
 &=
 \begin{cases}
 M_{\delta_j^B}\left(B_j\left(\bm x\right),N_j^B\left(\bm x\right)\right),
       &\norm{\bm x-\bm c_{j+1}}<\theta R_j,\\
 B_j\left(\bm x\right),&\norm{\bm x-\bm c_{j+1}}\ge\theta R_j,
 \end{cases}
 \end{aligned}
\end{equation}
where $M_\delta$ is defined in Eq.~\eqref{paper:smooth-maximum}.
Smoothness and convexity across the boundary require the
smooth maximum to equal $B_j$ there. The next phase's geometry
instead requires equality with $N_j^B$ near the center.
Let $\gamma>0$ denote the inner radius parameter.
By Lemma~\ref{paper:smooth-maximum-properties}, these equalities follow
from
\begin{equation}\label{paper:base-agreement-conditions}
 \begin{aligned}
 B_j\left(\bm x\right)-N_j^B\left(\bm x\right)&\ge\delta_j^B
 &&\text{if }\frac{\theta R_j}{2}\le\norm{\bm x-\bm c_{j+1}}\le\theta R_j,\\
 N_j^B\left(\bm x\right)-B_j\left(\bm x\right)&\ge\delta_j^B
 &&\text{if }\norm{\bm x-\bm c_{j+1}}\le\gamma R_j.
 \end{aligned}
\end{equation}
Above $B_j$'s tangent at $\bm c_{j+1}$, $B_j$ grows
quadratically with distance, whereas $N_j^B$ adds
$C_j+\frac{\mu}{q}\norm{\bm x-\bm c_{j+1}}^q$.
Since $q>2$, choosing the value increase and radii together
allows $B_j$ to be larger near the boundary while $N_j^B$
is larger near the center. We choose
\begin{equation}\label{paper:base-constants}
 \begin{gathered}
 a_0=\frac{1}{2}\left(\frac{3}{4}\right)^{q-2},\qquad
 a_1=\frac{q-1}{2}\left(\frac{5}{4}\right)^{q-2},\qquad
 a_{\rm lift}=\frac{a_0}{8},\\
 \theta=\min\left\{\frac{1}{16},
      \left(\frac{3qa_0}{32}\right)^{\frac{1}{q-2}}\right\},\quad
 \gamma=\theta\sqrt{\frac{3a_0}{32a_1}},\quad
 C_j=a_{\rm lift}\mu\theta^2R_j^q,\quad
 \delta_j^B=\frac{C_j}{4}.
 \end{gathered}
\end{equation}
Lemma~\ref{base:reset-lemma} in Appx.~\ref{app:paper-base}
derives these choices and proves Eq.~\eqref{paper:base-agreement-conditions}.
Consequently, $B_{j+1}=N_j^B$ for
$\norm{\bm x-\bm c_{j+1}}\le\gamma R_j$, and $B_{j+1}=B_j$
for $\norm{\bm x-\bm c_{j+1}}\ge\frac{\theta R_j}{2}$.
Thus the next phase has its required geometry and starting
gradient, while the earlier formula remains unchanged outside
the replacement ball.

\textbf{Step 1.4: Repeat at smaller accuracy scales.}
With the same chain and normalization parameters across
phases, the assembled function $F$ will satisfy
$F\left(\bm x\right)-\min F\ge\mu R_j^q$
on the phase-$j$ point sets specified in Lemma~\ref{paper:assembly}.
Here $\mu R_j^q=\mu\norm{\bm c_{j+1}-\bm c_j}^q$.
A fixed-factor decrease in these quantities therefore supplies
logarithmically many phases before the gap reaches $\epsilon$.
We obtain this decrease by defining
\[
 R_j=R_0\eta^j,\qquad
 \mu R_j^q=\mu R_0^q\eta^{qj},\qquad
 \frac{\mu R_{j+1}^q}{\mu R_j^q}=\eta^q.
\]
Each subsequent replacement must lie inside the ball where
$B_{j+1}=N_j^B$, so that its starting formula has the geometry
of Eq.~\eqref{paper:base-next-local}.
The triangle inequality gives this containment if
$\left(1+\theta\right)R_{j+1}\le\frac{\gamma R_j}{2}$.
We consequently require
$0<\eta\le\frac{\gamma}{2\left(1+\theta\right)}$.
The next lemma establishes the resulting local formulas,
retains earlier values and derivatives, and verifies global
smoothness and Euclidean uniform convexity independently of $L$.

\Needspace{5\baselineskip}
\begin{lemma}[Uniformly convex base]\label{paper:base}
For the constants chosen in Step~1.3, there is
$C_B>0$, depending only on $p,\nu$, such that the following holds.
We have $\theta\le\frac{1}{16}$ and $0<\gamma<\frac{\theta}{2}$.
For any $L\ge1$, $\mu,R_0>0$,
$0<\eta\le\frac{\gamma}{2\left(1+\theta\right)}$, and mutually orthogonal
$\bm h_0,\ldots,\bm h_{L-1}$ with $\norm{\bm h_j}=R_0\eta^j$,
Eq.~\eqref{paper:base-reset} produces convex $B_j\in C^{p,\nu}\left(\R^d\right)$
and $B=B_L$ satisfying:
\begin{claims}
\item\label{paper:base:class}
\textbf{Smoothness and uniform convexity.}
For $c_B=\frac{1}{2q\,16^{q-2}}$, all $\bm x,\bm y\in\R^d$,
and $0\le j\le L$,
\[
 \norm{\D^pB_j\left(\bm y\right)-\D^pB_j\left(\bm x\right)}
 \le C_B\mu\norm{\bm y-\bm x}^\nu,\qquad
 \Breg B{\bm y}{\bm x}\ge c_B\mu\norm{\bm y-\bm x}^q,
\]
\item\label{paper:base:centers}
\textbf{Values and gradients at the centers.}
$B\left(\bm 0\right)=0$, $\D B\left(\bm 0\right)=\bm 0$, and, for $0\le j<L$,
\begin{equation}\label{paper:base-data}
 \D B\left(\bm c_{j+1}\right)-\D B\left(\bm c_j\right)=\mu R_j^{q-2}\bm h_j,
 \qquad
 B\left(\bm c_{j+1}\right)-B\left(\bm c_j\right)=\beta_q\mu R_j^q,
\end{equation}
where $\beta_q=\frac{1}{q}+a_{\rm lift}\theta^2$.
\item\label{paper:base:locality}
\textbf{Exact local formulas.}
For $0\le j<L$,
$\displaystyle
 B_{j+1}\left(\bm x\right)=
 \begin{cases}
 N_j^B\left(\bm x\right),&
       \norm{\bm x-\bm c_{j+1}}\le\gamma R_j,\\
 B_j\left(\bm x\right),&
       \norm{\bm x-\bm c_{j+1}}\ge\frac{\theta R_j}{2}.
 \end{cases}
$
Eq.~\eqref{paper:base-local} holds globally for $j=0$ and on
$\norm{\bm x-\bm c_j}\le\gamma R_{j-1}$ for $1\le j\le L$.
\item\label{paper:base:preservation}
\textbf{Preservation of earlier formulas.}
For the replacement ball $S_j$ defined in Step~1.3 and
$0\le j<L$,
$B\left(\bm x\right)=B_j\left(\bm x\right)$ and
$\D^kB\left(\bm x\right)=\D^kB_j\left(\bm x\right)$
for $\bm x\notin\overline S_j$ and $1\le k\le p$.
Also, $B=B_j$ on a neighborhood of $\bm c_j$ for
$0\le j\le L$. The partial base $B_j$ depends only on
the previously chosen centers $\bm c_0,\ldots,\bm c_j$.
\end{claims}
\end{lemma}
\begin{proof}[Proof sketch]
\textbf{\ref{paper:base:class}}
Uniform convexity requires curvature along every segment, even though
our local geometry has zero Hessian at the centers. For
$\mathcal C=\{\bm c_0,\ldots,\bm c_L\}$,
$\bm d=\bm y-\bm x\ne\bm0$, and $\ell=\norm{\bm d}$,
Lemma~\ref{base:hessian-lemma} gives
$\D^2B\left(\bm z\right)[\bm d,\bm d]
\ge\mu\ell^2\operatorname{dist}\left(\bm z,\mathcal C\right)^{q-2}$.
Lemma~\ref{base:porosity} guarantees that at least half the segment
stays at distance $\frac{\ell}{16}$ from every center.
Integrating the curvature over this part gives gradient growth.
A further integration yields
$\Breg B{\bm y}{\bm x}\ge c_B\mu\norm{\bm y-\bm x}^q$.
For smoothness, each correction $B_{j+1}-B_j$ has
$p^{th}$-derivative magnitude bounded by $C\mu R_j^\nu$.
The geometric decrease of $R_j$ bounds their sum independently of $L$;
the exact local formulas handle nearby points.
The complete proof is given by
Corollary~\ref{base:complete-base-corollary} in Appx.~\ref{app:paper-base}.
\end{proof}

%% file: manuscript-phases.tex
\subsection{Step 2: Add Hardness to the Uniformly Convex Base}
\label{paper:hardness-construction}
The base $B$ supplies the required uniform convexity and smoothness,
but does not yet force oracle calls.
In particular, the initial point $\bm c_0=\bm 0$ already minimizes $B$,
since
\begin{equation}\label{paper:base-starts-optimal}
 B\left(\bm x\right)
 \ge B\left(\bm 0\right)+\ip{\D B\left(\bm 0\right)}{\bm x}
                +c_B\mu\norm{\bm x}^{q}
 =c_B\mu\norm{\bm x}^{q},
\end{equation}
where the inequality follows from
Lemma~\resultpart{paper:base}{paper:base:class} and the equality from
Lemma~\resultpart{paper:base}{paper:base:centers}.
Thus returning $\bm 0$ attains zero gap without an oracle call.
We therefore add convex chains that restrict the directions revealed
by each oracle output and leave a positive gap at points using only
the first few chain coordinates,
while preserving the base's uniform convexity.

\textbf{Step 2.1: Define the short chains.}
Since each phase should force the same condition-number cost,
we use a common length $m\ge1$, independent of $L$.
To give each phase fresh chain directions and keep the center
increments orthogonal, we choose mutually orthogonal blocks for
$0\le j<L$,
$E_j=\operatorname{span}\left\{\bm v_j,\bm u_{j,1},\ldots,\bm u_{j,m}\right\}$,
with orthonormal bases and starting directions $\bm v_j$.
We refer to the subspace
$U_{j,b}:=\operatorname{span}\left\{\bm v_j,\bm u_{j,1},\ldots,\bm u_{j,b}\right\}$,
for $b\in\{0,\ldots,m\}$,
as a \emph{chain prefix}, and call it a \emph{short prefix} when
$b\le\left\lfloor\frac{m}{4}\right\rfloor$.
Following the power chain in~\cite[Eq.~(5.4)]{GrapigliaNesterov2020},
we adapt its scaling to our phase geometry.
For $\bm w\in E_j$, we denote its coordinates in this orthonormal basis by
$t=\ip{\bm v_j}{\bm w}$ and $x_i=\ip{\bm u_{j,i}}{\bm w}$,
and define, for $\lambda>0$,
\begin{equation}\label{paper:chain}
 Q_j\left(\bm w\right)=\frac{\lambda}{q}
 \left(\left|\frac{t}{\sqrt{m}}-x_1\right|^q
 +\sum_{i=1}^{m-1}|x_i-x_{i+1}|^q+|x_m|^q\right).
\end{equation}
\emph{Choose the phase endpoint.}
To pass from one phase to the next, we need an endpoint at which
the chain's gradient can be balanced with the base's gradient increment.
We choose one with no chain gradient in the internal directions
$\bm u_{j,i}$.
We refer to the differences inside Eq.~\eqref{paper:chain} as
\emph{chain residuals}.
Since each $x_i$ appears with opposite signs in consecutive residuals,
we cancel the internal gradient components by requiring
\[
 \frac{t}{\sqrt m}-x_1=x_1-x_2=\cdots=x_{m-1}-x_m=x_m=\tau_j>0.
\]
The middle differences are omitted when $m=1$.
This follows the minimizer calculation
of~\mbox{\cite[Lemma~5.2]{GrapigliaNesterov2020}}.
Solving these equalities gives the center increment
\begin{equation}\label{paper:ramp}
 t_j^\star=\sqrt m\left(m+1\right)\tau_j,\qquad
 \bm h_j=t_j^\star\bm v_j+
          \sum_{i=1}^m\left(m-i+1\right)\tau_j\bm u_{j,i},\quad
 R_j=\norm{\bm h_j}.
\end{equation}
The remaining chain gradient is
$\D Q_j\left(\bm h_j\right)=\Gamma_j\bm v_j$, with
$\Gamma_j=\frac{\lambda\tau_j^{q-1}}{\sqrt{m}}$.
The base adds $\mu R_j^{q-2}\bm h_j$ by
Eq.~\eqref{paper:base-data}.
To bound the base-to-chain gradient ratio in Step~3 and keep
earlier evaluation points outside later replacement balls in
Step~5, we need the endpoint length in terms of $\tau_j$ and $m$.
Lemma~\ref{collar:chain-lemma} in Appx.~\ref{app:paper-phases}
proves the chain gradient above and gives
$R_j=\beta_m\tau_jm^{\frac{3}{2}}$ with $1\le\beta_m\le\sqrt5$.

\emph{Limit the revealed directions.}
The endpoint involves all $m$ internal directions, so the oracle
must not reveal all of them in a single oracle output.
For an evaluation point $\bm w\in U_{j,b}$, we have $x_i=0$ for $i>b$.
For $b<m$, its last possibly nonzero residual uses
$\bm u_{j,b+1}$. Since $q>p$, the $|\cdot|^q$ terms have zero
derivatives of orders $1,\ldots,p$ at zero. Thus the value and
derivatives through order $p$ depend only on this prefix and
next direction, regardless of $\norm{\bm w}$.

\textbf{Step 2.2: Create a nonzero starting gradient.}
The same vanishing derivatives give $\D^kQ_j\left(\bm 0\right)=0$
for $1\le k\le p$, so the chain alone creates no starting slope.
Our phase functions will retain these derivative values at their
local zero coordinates.
Since $\bm c_j\in\bigoplus_{k<j}E_k$ and $E_k\perp E_\ell$ for $k\ne\ell$, we have
$P_{E_j}\left(\bm 0-\bm c_j\right)=\bm 0$ and
$\ip{\bm v_{j+1}}{\bm 0-\bm c_j}=0$ when $j<L-1$,
so adding the phase functions alone would still leave $\bm 0$ optimal.
We therefore include
$\ell_{\mathrm{init}}\left(\bm x\right)=-s_0\ip{\bm v_0}{\bm x}$,
with $s_0>0$ to be chosen later, to give starting gradient $-s_0\bm v_0$.
By Lemma~\ref{paper:bregman-additions} in Appx.~\ref{app:paper-elementary},
this affine term preserves the Bregman divergence and, since
$p\ge2$, leaves $\D^pB$ unchanged.

\textbf{Step 2.3: Keep the gradient of each completed phase fixed.}
To repeat the same chain hardness in each phase, we require
the starting gradient at $\bm c_j$ to be $-s_j\bm v_j$,
with no component in any earlier block $E_k$, $k<j$.
At $\bm c_{j+1}$, the base adds
$\bm b_j=\mu R_j^{q-2}\bm h_j$ to its gradient.
To cancel this increment and transfer the starting slope from
$\bm v_j$ to $\bm v_{j+1}$, the phase function must contribute
a gradient $\bm g_j$ satisfying
$-s_j\bm v_j+\bm b_j+\bm g_j=-s_{j+1}\bm v_{j+1}$.
We match $\ip{\bm g_j}{\bm v_j}$ to the chain's endpoint
gradient component $\Gamma_j$, which determines $s_j$:
\begin{equation}\label{paper:affine-slope}
 \bm g_j=s_j\bm v_j-s_{j+1}\bm v_{j+1}-\bm b_j,\qquad
 s_j=\Gamma_j+\frac{G_jt_j^\star}{R_j},\quad G_j=\norm{\bm b_j}.
\end{equation}
For $j<L-1$, we denote the coordinate along $\bm v_{j+1}$ by
$t_{j+1}$, distinguishing this variable from the fixed endpoint
coordinate $t_{j+1}^\star$.
Relative to $\bm c_j$, the center $\bm c_{j+1}$ has coordinates
$\left(\bm h_j,0\right)$, and all subsequent centers have
$\left(\bm h_j,t_{j+1}^\star\right)$.
Since the completed phase must contribute $\bm g_j$ at each of
them, we define the affine function
\begin{subequations}\label{paper:phase-definition}
\begin{equation}\label{paper:phase-affine}
 I_j\left(\bm w,t_{j+1}\right)
 =\ip{\bm g_j}{\bm w-\bm h_j+t_{j+1}\bm v_{j+1}}.
\end{equation}
To bound the derivatives of the smooth maximum, we need to control
the chain's values where it meets the affine function.
However, $Q_j-I_j$ contains the term $s_{j+1}t_{j+1}$, which
is unbounded below as $t_{j+1}\to-\infty$.
To control this direction as well, we add a convex regularizer,
\begin{equation}\label{paper:phase-regularized-chain}
 P_j\left(\bm w,t_{j+1}\right)
 =Q_j\left(\bm w\right)+\frac{\alpha}{q}|t_{j+1}|^q,\qquad
 \alpha=\lambda m^{1-\frac{3q}{2}}.
\end{equation}
This choice makes the added term comparable to
$\lambda m\tau_j^q$ when $|t_{j+1}|=\tau_jm^{\frac{3}{2}}$.
To make the affine formula larger at $\left(\bm h_j,0\right)$,
we choose $\delta_j=c_1\lambda m\tau_j^q$, with
$c_0=c_1=\frac{2^{-q-3}}{q}$, and shift the regularized chain down:
\begin{equation}\label{paper:phase-shift}
 O_j\left(\bm w,t_{j+1}\right)
 =P_j\left(\bm w,t_{j+1}\right)-P_j\left(\bm h_j,0\right)-4\delta_j.
\end{equation}
Now $I_j\left(\bm h_j,0\right)=0$ and
$O_j\left(\bm h_j,0\right)=-4\delta_j$.
The smooth maximum therefore agrees with $I_j$ near this point,
while retaining $O_j$ wherever $O_j-I_j>\delta_j$.
We define the phase function by
\begin{equation}\label{paper:phase-smooth-maximum}
 A_j\left(\bm w,t_{j+1}\right)
 =M_{\delta_j}\left(I_j,O_j\right).
\end{equation}
\end{subequations}
For $j=L-1$, we omit $t_L,\bm v_L$ and set $s_L=t_L^\star=0$,
with $s_L\bm v_L:=\bm 0$.
These functions are expressed relative to the center $\bm c_j$.
To add them to $B$, we extend each $T\in\{P,I,O,A\}$ to $\R^d$
using the projection onto $E_j$ and the next starting direction:
\begin{equation}\label{paper:ambient-phase}
 T_j\left(\bm x\right):=
 T_j\left(P_{E_j}\left(\bm x-\bm c_j\right),
          \ip{\bm v_{j+1}}{\bm x-\bm c_j}\right),
\end{equation}
with the second input omitted when $j=L-1$.
It remains to verify that the smooth maximum uses the chain at short
prefixes and the affine function near completed centers, while
preserving global smoothness.
For $\rho_j^A=\left(\frac{\delta_j}{\lambda}\right)^{\frac{1}{q}}$,
the critical relation $q=p+\nu$ gives
$\frac{\delta_j}{\left(\rho_j^A\right)^{p+\nu}}=\lambda$, so the smoothing
normalization is independent of $\tau_j$.
The following lemma verifies these formulas and starting values,
with a global smoothness bound independent of $m$ and dimension.

\Needspace{5\baselineskip}
\begin{lemma}[The phase function]\label{paper:phase}
There is a constant $C_A>0$,
depending only on $p,\nu$, with the following property.
In Eqs.~\eqref{paper:chain}--\eqref{paper:phase-definition},
take $m\ge1$, $\lambda,\tau_j>0$, and
$G_j,s_{j+1},t_{j+1}^\star\ge0$.
Define $\bm b_j=\frac{G_j}{R_j}\bm h_j$ and the slope by
Eq.~\eqref{paper:affine-slope}, using the scalar parameters above.
Assume
\begin{equation}\label{paper:phase-conditions}
 G_j\le c_0\Gamma_j,\quad
 \alpha^{-\frac{1}{q-1}}s_{j+1}^{\frac{q}{q-1}}
 \le c_0\lambda m\tau_j^q,\quad
 s_{j+1}t_{j+1}^\star+\frac{\alpha}{q}\left(t_{j+1}^\star\right)^q
 \le c_0\lambda m\tau_j^q.
\end{equation}
\begin{claims}
\item\label{paper:phase:class} \textbf{Convexity and smoothness.}
$A_j$ is convex and $C^{p,\nu}$, with
$\norm{\D^pA_j\left(\bm x\right)-\D^pA_j\left(\bm y\right)}
\le C_A\lambda\norm{\bm x-\bm y}^\nu$
for all $\bm x,\bm y\in\R^d$,
where $C_A$ is independent of $m$ and dimension.
\item\label{paper:phase:prefix} \textbf{Chain formulas and oracle outputs.}
For every $\bm w\in U_{j,b}$ and
$0\le b\le\left\lfloor \frac{m}{4}\right\rfloor$,
$O_j\left(\bm w,0\right)-I_j\left(\bm w,0\right)\ge7\delta_j$,
and the functions on $\R^d$ defined by
Eq.~\eqref{paper:ambient-phase} satisfy $A_j=O_j$ near $\bm c_j+\bm w$.
In particular,
$A_j\left(\bm c_j+\bm w\right)
=Q_j\left(\bm w\right)-\frac{\lambda\left(m+1\right)\tau_j^q}{q}-4\delta_j$.
For $1\le k\le p$ and all $\bm a_1,\ldots,\bm a_k\in\R^d$,
the derivative tensors on $\R^d$ satisfy
\[
 \D^kA_j\left(\bm c_j+\bm w\right)\left[\bm a_1,\ldots,\bm a_k\right]
 =\D^kQ_j\left(\bm w\right)
 \left[P_{U_{j,b+1}}\bm a_1,\ldots,P_{U_{j,b+1}}\bm a_k\right].
\]
The value and these tensors are unchanged by varying the remaining
chain directions and, when present, $\bm v_{j+1}$ in an orthonormal completion.
\item\label{paper:phase:completed} \textbf{Gradient after completing the phase.}
$A_j=I_j$ near
$\left\{\left(\bm h_j,t_{j+1}\right):0\le t_{j+1}\le t_{j+1}^\star\right\}$.
On a neighborhood of
$\left\{\bm c_j+\bm h_j+t\bm v_{j+1}:0\le t\le t_{j+1}^\star\right\}$
in $\R^d$ (with the last term omitted when $j=L-1$),
$\D A_j=\bm g_j$ and $\D^kA_j=0$ for $2\le k\le p$.
Its values at the two endpoints are $0$ and
$-s_{j+1}t_{j+1}^\star$.
\item\label{paper:phase:starting} \textbf{Starting and terminal data.}
$A_j\left(\bm 0,0\right)=-\frac{\lambda\left(m+1\right)\tau_j^q}{q}-4\delta_j$ and
$\D^kA_j\left(\bm 0,0\right)=0$ for $1\le k\le p$.
The terminal phase satisfies the same conclusions with the
next coordinate omitted and the completed segment reduced to $\{\bm c_L\}$.
\end{claims}
\end{lemma}
The proof is given by Theorem~\ref{collar:theorem}
in Appx.~\ref{app:paper-phases}.

\input{manuscript-assembly}

%% file: manuscript-assembly.tex
\subsection{Step 3: Choose Compatible Parameters and Assemble the Function}
\label{paper:assembly-section}
Selecting $O_j$ on short chain prefixes requires
$G_j\le c_0\Gamma_j$ in Eq.~\eqref{paper:phase-conditions}.
Since $R_j=\beta_m\tau_jm^{\frac{3}{2}}$, the ratio of the base
increment $G_j$ to $\Gamma_j$ is independent of $\tau_j$.
We choose $\tau_0>0$ and
\begin{equation}\label{paper:parameters}
 a_q=\frac{3q-2}{2},\quad
 \lambda=K\mu m^{a_q},\quad \alpha=K\mu,\quad
 \tau_j=\eta^j\tau_0,\quad \bm c_j=\sum_{k<j}\bm h_k.
\end{equation}
Then
\begin{equation}\label{paper:balance}
 \frac{G_j}{\Gamma_j}
 =\frac{\mu R_j^{q-1}\sqrt m}{\lambda\tau_j^{q-1}}
 =\frac{\mu}{\lambda}\beta_m^{q-1}m^{a_q}
 =\frac{\beta_m^{q-1}}{K}.
\end{equation}
Thus $K=\frac{5^{\frac{q-1}{2}}}{c_0}$ gives $G_j\le c_0\Gamma_j$
for every phase.
Choosing $\eta=\frac{\gamma}{2\left(1+\theta\right)}$ also
satisfies the next-phase slope, endpoint, and base-nesting
conditions by Proposition~\ref{global:parameters} in
Appx.~\ref{app:paper-assembly}.
All constants depend only on $p,\nu$.

We then assemble
\begin{equation}\label{paper:objective}
 F\left(\bm x\right)=B\left(\bm x\right)+\ell_{\mathrm{init}}\left(\bm x\right)
                      +\sum_{j=0}^{L-1}A_j\left(\bm x\right).
\end{equation}
An algorithm may retain completed-phase information, so we include
those subspaces in the current-prefix space
$W_{j,b}:=\left(\bigoplus_{k<j}E_k\right)\oplus U_{j,b}$,
for $0\le j<L$ and $0\le b\le m$.
The sum over $k<j$ is $\{\bm 0\}$ when $j=0$.

\Needspace{6\baselineskip}
\begin{lemma}[The complete objective]\label{paper:assembly}
Let $m,L\ge1$ and $\mu,\tau_0>0$. Choose the parameters in
Eq.~\eqref{paper:parameters} with
$K=\frac{5^{\frac{q-1}{2}}}{c_0}$ and
$\eta=\frac{\gamma}{2\left(1+\theta\right)}$, and construct $B,A_j$
by Eqs.~\eqref{paper:base-reset} and \eqref{paper:phase-definition}.
Then $F$ in Eq.~\eqref{paper:objective} satisfies:
\begin{claims}
\item\label{paper:assembly:class}\textbf{Smoothness and uniform convexity.}
For all $\bm x,\bm y\in\R^d$,
\[
 \Breg F{\bm y}{\bm x}\ge c_B\mu\norm{\bm y-\bm x}^q,\qquad
 \norm{\D^pF\left(\bm x\right)-\D^pF\left(\bm y\right)}
 \le C_H\mu m^{a_q}\norm{\bm x-\bm y}^\nu,
\]
where
$C_H:=C_B+2C_AK$ depends only on $p,\nu$.
\item\label{paper:assembly:centers}\textbf{Center gradients and minimizer.}
$\D F\left(\bm c_j\right)=-s_j\bm v_j$ for $j<L$,
$\D F\left(\bm c_L\right)=\bm 0$, and $\bm c_L$ is the unique minimizer.
\item\label{paper:assembly:prefix}\textbf{Optimality gap on chain prefixes.}
Every $\bm x\in W_{j,b}$ with $b\le\left\lfloor \frac{m}{4}\right\rfloor$
satisfies
$F\left(\bm x\right)-F\left(\bm c_L\right)\ge7\delta_j$, where
$\delta_j=c_1\lambda m\tau_j^q$ and $c_1=\frac{2^{-q-3}}{q}$.
\item\label{paper:assembly:initial}\textbf{Objective decay and initial gap.}
Define the decrease between consecutive centers by
$\chi_j:=F\left(\bm c_j\right)-F\left(\bm c_{j+1}\right)$.
The initial optimality gap $D_{m,L}:=F\left(\bm c_0\right)-F\left(\bm c_L\right)$ satisfies
\begin{equation}\label{paper:initial-gap}
 \chi_j=\eta^{qj}\chi_0,\qquad
 D_{m,L}=\chi_0\frac{1-\eta^{qL}}{1-\eta^q}.
\end{equation}
Here $0<\chi_0\le2\left(1+c_0\right)\lambda m\tau_0^q$.
The value of $D_{m,L}$ is independent of the orthonormal phase directions.
\end{claims}
\end{lemma}
\begin{proof}[Proof sketch]
\textbf{\ref{paper:assembly:class}}
Adding convex phase functions preserves the uniform convexity of $B$.
To avoid accumulating a smoothness cost proportional to $L$, we group
phases by parity. Each $A_j$ depends only on
$V_j=E_j\oplus\operatorname{span}\{\bm v_{j+1}\}$, with
$V_{L-1}=E_{L-1}$, and these subspaces are orthogonal within each group.
The norm of each group's derivative sum is therefore its largest
component norm, giving a total phase H\"older bound $2C_A\lambda$.
\textbf{\ref{paper:assembly:centers}}
The center gradients follow by telescoping
$\bm b_j+\bm g_j=s_j\bm v_j-s_{j+1}\bm v_{j+1}$.
\textbf{\ref{paper:assembly:prefix}}
To keep the other components from offsetting the chain gap, consider
the convex function $\Psi_j=F-A_j+I_j$, which equals $F$ near
$\bm c_{j+1}$. For $\bm x\in W_{j,b}$, the center gradient is
perpendicular to $\bm x-\bm c_{j+1}$; it is zero in the terminal phase.
Convexity therefore gives $\Psi_j\left(\bm x\right)\ge F\left(\bm c_{j+1}\right)$.
Adding $A_j\left(\bm x\right)-I_j\left(\bm x\right)\ge7\delta_j$
gives $F\left(\bm x\right)\ge\min F+7\delta_j$.
The full proof follows this order in Appx.~\ref{app:assembly-full-proof}.
\end{proof}

%% file: manuscript-oracle-rate.tex
\subsection{Steps 4--5: Generate the Trajectory and Prove Oracle Indistinguishability}
\label{paper:trajectory}

\textbf{Step 4: Generate the trajectory.}\quad
To apply Lemma~\resultpart{paper:assembly}{paper:assembly:prefix}
to an arbitrary deterministic algorithm $\mathcal A$, we must place its candidate
solutions in short prefixes $W_{j,b}$. We therefore choose the unfixed
chain directions as the algorithm runs, while fixing the scalar
parameters of Step~3 and $d\ge L\left(m+1\right)$ beforehand.
For $m\ge32$, define $n_m=\left\lfloor\frac{m}{16}\right\rfloor$.
Group $j$ consists of oracle calls $jn_m+1,\ldots,\left(j+1\right)n_m$.

For oracle call $s$, denote the evaluation point by $\bm x_s$,
the returned value and full derivative tensors by
$\mathcal O_s=\left(a_{s,0},\ldots,a_{s,p}\right)$,
and the candidate solution after this call by $\widehat{\bm x}_s$.
We impose no norm bound on these points, and a candidate solution
need not be an evaluation point.
Lemma~\ref{oracle:point} in Appx.~\ref{app:paper-oracle} allows any
point to be included in $W_{j,b}$ by choosing at most one additional
chain direction. We begin by choosing $\bm v_0$ and including
$\widehat{\bm x}_0$, using at most one direction.
Within group $j$, let $b$ count the chosen directions
$\bm u_{j,1},\ldots,\bm u_{j,b}$.
Each call consists of three ordered selections, with $b$ increased
whenever a direction is chosen.
\begin{enumerate}
\item The algorithm determines $\bm x_{s+1}$ from the earlier
oracle outputs. We include it in $W_{j,b}$ using at most one direction.
\item The current chain can involve $\bm u_{j,b+1}$ in its
full derivative tensors, so we fix one further unit direction perpendicular
to $W_{j,b}$. Any orthonormal completion gives an objective
$\widetilde F$, from which we return
$a_{s+1,k}=\D^k\widetilde F\left(\bm x_{s+1}\right)$ for
$0\le k\le p$.
Step~5 proves independence from the remaining completion.
\item The new oracle output determines $\widehat{\bm x}_{s+1}$.
We include this candidate solution in $W_{j,b}$ using at most one
further direction.
\end{enumerate}
Thus each oracle call fixes at most three chain directions.
After the last call in group $j$, we complete its block $E_j$.
If $j<L-1$, we choose $\bm v_{j+1}$ perpendicular to the completed
blocks and begin the next group with $b=0$.

\textbf{Step 5: Prove oracle indistinguishability and complete one instance.}\quad
We prove that one fixed objective reproduces the adaptive
trajectory in Step~4. Let $\mathcal F_s$ be the family of objectives in
Eq.~\eqref{paper:objective} obtained by completing the directions
fixed after $s$ oracle calls, including those chosen for
$\widehat{\bm x}_s$ and preparation of the next group.
Equality of the recorded oracle outputs across these completions
will make the deterministic algorithm reproduce the
trajectory, allowing us to apply the prefix gaps from Step~3.

\begin{proposition}[Oracle indistinguishability and one fixed instance]
\label{paper:oracle}
Fix the scalar parameters specified in Step~3, $m\ge32$,
$L\ge1$, and $d\ge L\left(m+1\right)$. Put $n_m=\left\lfloor \frac{m}{16}\right\rfloor$.
For every deterministic $p^{th}$-order oracle algorithm, the ordered
orthonormal directions in Step~4 can be chosen so that:
\begin{claims}
\item\label{paper:oracle:consistency}\textbf{Oracle indistinguishability.}
For every $0\le s\le n_mL$,
\begin{equation}\label{paper:oracle-consistency}
 \mathcal F_s\ne\varnothing,\qquad
 \D^kF\left(\bm x_\ell\right)=a_{\ell,k}
 \quad\left(F\in\mathcal F_s,\ 1\le\ell\le s,\ 0\le k\le p\right).
\end{equation}
\item\label{paper:oracle:fixed}\textbf{One fixed objective and candidate gaps.}
Every $F\in\mathcal F_{n_mL}$ reproduces the recorded evaluation
points, oracle outputs, and candidate solutions when the algorithm
runs on it through $n_mL$ oracle calls, and satisfies, for
$jn_m<s\le\left(j+1\right)n_m$ and $0\le j<L$,
\begin{equation}\label{paper:oracle-gaps}
 F\left(\widehat{\bm x}_s\right)-\min F\ge7\delta_j,\qquad
 F\left(\widehat{\bm x}_0\right)-\min F\ge7\delta_0.
\end{equation}
\end{claims}
\end{proposition}
\Needspace{3\baselineskip}
\begin{proof}[Proof sketch]
\textbf{\ref{paper:oracle:consistency}}
\textbf{Step 5.1: Available directions.}\quad
The selections in Step~4 leave candidate solutions in short prefixes,
with unused directions for completion.
\textbf{Step 5.2: Identical oracle outputs.}\quad
Since
$\operatorname{dist}\left(\bm c_{j+1},W_{j,b}\right)>
\theta R_j$, prefix points miss all later replacement balls, so
the base's values and derivatives depend only on fixed centers.
The current phase uses only the prefix and the additional direction
fixed in Step~4, while future and completed phases give fixed values
and derivative tensors.
Thus every completion gives the same full oracle output, and further
selections preserve earlier outputs by restricting the completion family.
\textbf{\ref{paper:oracle:fixed}}
\textbf{Step 5.3: One fixed objective.}\quad
Determinism reproduces the trajectory on any final completion, and
Lemma~\resultpart{paper:assembly}{paper:assembly:prefix} gives the stated gaps.
The full construction and proof are given in
Proposition~\ref{oracle:realization} in Appx.~\ref{app:paper-oracle}.
\end{proof}

\subsection{Step 6: Choose the Parameters and Prove the Main Theorem}\label{paper:rate}

We now formally prove Theorem~\ref{paper:main-theorem} by combining
Lemma~\ref{paper:assembly} with Proposition~\ref{paper:oracle}.
We choose the chain length to meet $H,\sigma$, the phase count to
keep the candidate gaps above $\epsilon$, and rescale to the
prescribed initial gap $\Delta_0$.

\begin{proof}
\setlength{\abovedisplayskip}{5pt plus 2pt minus 1pt}
\setlength{\belowdisplayskip}{5pt plus 2pt minus 1pt}
\setlength{\belowdisplayshortskip}{5pt plus 2pt minus 1pt}
\textbf{High-Order H\"older Smoothness and Uniform Convexity.}\quad
With $\tau_0=1$ and $\kappa=\frac{H}{\sigma}$,
Lemma~\resultpart{paper:assembly}{paper:assembly:class}
gives the prescribed class parameters provided
$qc_B\mu\ge\sigma$ and $C_H\mu m^{a_q}\le H$.
Taking $\mu=\frac{\sigma}{qc_B}$ leaves
$m^{a_q}\le\frac{qc_B}{C_H}\kappa$ as the smoothness requirement.
We therefore choose the admissible integer
$m\asymp_{p,\nu}\kappa^{\frac{2}{3q-2}}$ given by
Lemma~\resultpart{final:calibration}{final:calibration:part:chain}
in Appx.~\ref{final:parameters}.
For sufficiently large $\kappa$, this choice also gives
$m\ge32$ and $n_m\asymp m$, fixing the condition-number cost
of each phase.

\textbf{Phase Count.}\quad
Because these class bounds are independent of $L$, we can repeat
this cost at successively smaller gaps.
Comparing $\delta_j=c_1\lambda m\eta^{qj}$ with the initial gap
in Lemma~\resultpart{paper:assembly}{paper:assembly:initial} gives,
for a constant $c_*\in\left(0,1\right)$ depending only on $p,\nu$,
\begin{equation}\label{paper:relative-protected-gap}
 \frac{7\delta_j}{D_{m,L}}\ge c_*\eta^{qj}.
\end{equation}
Normalizing the initial gap to $\Delta_0$ therefore gives a
phase-$j$ candidate gap of at least $c_*\Delta_0\eta^{qj}$.
We choose the largest integer $L\ge1$ such that
$c_*\Delta_0\eta^{q\left(L-1\right)}\ge2\epsilon$.
Since $0<\eta<1$, this keeps every earlier phase gap above
$2\epsilon$ as well.
Lemma~\ref{final:calibration}\,\ref{final:calibration:part:phases}--\ref{final:calibration:part:dimension}
verifies $L\asymp_{p,\nu}\log\left(\frac{\Delta_0}{\epsilon}\right)$
and $L\left(m+1\right)\le d$ under the theorem's dimension bound.

\textbf{Initial Gap.}\quad
To apply this comparison to the prescribed initial gap, we take
$a=\left(\frac{\Delta_0}{D_{m,L}}\right)^{\frac{1}{q}}$
and define
$\widehat F\left(\bm x\right)
=a^qF\left(\frac{\bm x-\bm x_0}{a}\right)$.
Lemma~\resultpart{paper:assembly}{paper:assembly:initial}
makes $D_{m,L}$ independent of the chain directions, so this
normalization is fixed before the oracle trajectory.
Since $q=p+\nu$, Lemma~\ref{final:rescaling} preserves both
class parameters and gives
$\widehat F\left(\bm x_0\right)-\min\widehat F
=\Delta_0$.

\textbf{Oracle Lower Bound.}\quad
Applying Proposition~\ref{paper:oracle} to the algorithm in
coordinates $\frac{\bm x-\bm x_0}{a}$, with the derivative
transformation in Lemma~\ref{final:rescaling}, gives one fixed
$\widehat F$ reproducing its oracle trajectory.
For a candidate solution $\widehat{\bm x}_s$ in group $j$,
Proposition~\resultpart{paper:oracle}{paper:oracle:fixed},
Eq.~\eqref{paper:relative-protected-gap},
and our choice of $L$ give
\[
 \widehat F\left(\widehat{\bm x}_s\right)-\min\widehat F
 \ge7a^q\delta_j
 \ge c_*\Delta_0\eta^{qj}
 \ge2\epsilon>\epsilon.
\]
The zero-call candidate solution satisfies the same bound.
Thus reaching the target gap after $T$ oracle calls requires
\[
 T>n_mL=\Omega_{p,\nu}\left(
 \left(\frac{H}{\sigma}\right)^{\frac{2}{3q-2}}
 \log\left(\frac{\Delta_0}{\epsilon}\right)\right)
\]
Choosing $c_{p,\nu}>0$ sufficiently small therefore proves both
the theorem's accuracy range and lower bound.
\end{proof}

%% file: ai-statement.tex
\Needspace{5\baselineskip}
\makeatletter
\renewcommand\section{\@startsection{section}{1}{\z@}{-1.4ex \@plus -.4ex \@minus -.2ex}{.6ex \@plus .2ex}{\normalfont\Large\bfseries}}
\makeatother
\section*{AI Statement}

We used OpenAI Codex 6 Astra and 6.1 Sol to assist with literature
searches, identifying suitable constants, verifying proof arguments,
and refining the presentation of proofs. This assistance materially
influenced the main text, appendices, and figures. The author independently
verified all mathematical arguments and takes full responsibility
for the content of the paper.

%% file: appendix-guide.tex
\section*{Contents of the Appendices}

Appx.~\ref{app:paper-introduction} gives the norm-mismatch witness
from the introduction.
Appx.~\ref{app:paper-elementary} supplies the preliminary analytic
tools. Appx.~\ref{app:paper-base}--\ref{app:paper-assembly} give
the proofs for the base, phase functions, and complete objective
in Steps~1--3. Appx.~\ref{app:paper-oracle} proves the oracle
realization in Sec.~\ref{paper:trajectory}, and
Appx.~\ref{app:paper-calibration} gives the parameter calibration
and rescaling for Sec.~\ref{paper:rate}.
Imported results are restated with their source citations beside the
corresponding arguments.

\begingroup
\setcounter{tocdepth}{2}
\appendixcontents
\endgroup
\clearpage

%% file: manuscript-supporting-proofs.tex
\input{appendix-introduction}
\input{appendix-preliminaries}
\input{appendix-base}
\input{appendix-phases}

\input{appendix-assembly}

\input{appendix-oracle}
\input{appendix-rate}

%% file: appendix-introduction.tex
\section{Norm Mismatch in a Long Residual Chain}
\label{app:paper-introduction}\label{app:paper-barrier}

To quantify the \hyperref[paper:norm-mismatch-motivation]{norm mismatch},
we consider a single long residual chain. Let $q=p+\nu>2$ and
let $S\left(x_1,\ldots,x_n\right)=\left(0,x_1,\ldots,x_{n-1}\right)$
be the lower shift. For $\alpha,a>0$ and $0<\xi<1$, define
\begin{equation}\label{paper:barrier-chain}
 A_\xi=I-\xi S,\qquad
 G\left(\bm x\right)=\frac{\alpha}{q}
          \norm{A_\xi\bm x-a\bm e_1}_q^q.
\end{equation}
Here $\norm{\bm z}_q^q=\sum_i|z_i|^q$ and $\bm e_1$ is the first
standard basis vector. The minimizer of $G$ is
$\bm x^\star=a\left(1,\xi,\ldots,\xi^{n-1}\right)$, since
$A_\xi\bm x^\star=a\bm e_1$ and $A_\xi$ is invertible.
In particular, $G\left(\bm x^\star\right)=0$ and
$\D G\left(\bm x^\star\right)=\bm 0$ by
Lemma~\resultpart{elem:power}{elem:power:part:flatness} in
Appx.~\ref{app:paper-elementary}.

For $\phi\left(t\right)=\frac{|t|^q}{q}$, the scalar divergence bound in
Lemma~\resultpart{elem:power}{elem:power:part:bregman} is expressed
in the residual $\ell_q$ norm. Converting it to a Euclidean bound
uses
\[
 \begin{aligned}
 \norm{\bm z}_2^2
 &=\sum_{i=1}^n |z_i|^2\cdot1\\
 &\le\left(\sum_{i=1}^n|z_i|^q\right)^{\frac{2}{q}}
       n^{\frac{q-2}{q}}
 &&\text{by H\"older's inequality},\\
 \norm{\bm z}_q^q
 &\ge n^{1-\frac{q}{2}}\norm{\bm z}_2^q
 &&\text{by raising to the power }\frac{q}{2}.
 \end{aligned}
\]
The shift satisfies $\norm{S\bm z}_2\le\norm{\bm z}_2$, so
the reverse triangle inequality gives
$\norm{A_\xi\bm z}_2
\ge\norm{\bm z}_2-\xi\norm{S\bm z}_2
\ge\left(1-\xi\right)\norm{\bm z}_2$.
For $\bm x,\bm y\in\R^n$, define their residual vectors
$\bm u=A_\xi\bm x-a\bm e_1$ and $\bm v=A_\xi\bm y-a\bm e_1$.
Applying the two norm bounds to $\bm v-\bm u=A_\xi\left(\bm y-\bm x\right)$ gives
\[
 \begin{aligned}
 D_G\left(\bm y,\bm x\right)
 &=\alpha\sum_{i=1}^n D_\phi\left(v_i,u_i\right)
 &&\text{by the chain rule and Bregman definition},\\
 &\ge\frac{2^{2-q}\alpha}{q}
       \norm{A_\xi\left(\bm y-\bm x\right)}_q^q
 &&\text{by Lemma~\resultpart{elem:power}{elem:power:part:bregman}},\\
 &\ge\frac{2^{2-q}\alpha}{q}n^{1-\frac{q}{2}}
       \norm{A_\xi\left(\bm y-\bm x\right)}_2^q
 &&\text{by the preceding norm comparison},\\
 &\ge\frac{2^{2-q}\alpha}{q}n^{1-\frac{q}{2}}
       \left(1-\xi\right)^q\norm{\bm y-\bm x}_2^q
 &&\text{by the lower bound for }A_\xi.
 \end{aligned}
\]
Thus $G$ is uniformly convex with modulus
\begin{equation}\label{paper:barrier-certificates}
 \sigma_{\rm cert}
 =2^{2-q}\alpha n^{1-\frac{q}{2}}\left(1-\xi\right)^q.
\end{equation}
The following calculation shows that the factor $n^{1-\frac{q}{2}}$ is
also present in the best Euclidean modulus of this family.

\begin{lemma}[A direction witnessing the Euclidean modulus loss]
\label{paper:barrier-witness}
If $n\ge2\left\lceil\left(1-\xi\right)^{-1}\right\rceil$, every valid degree-$q$
Euclidean uniform-convexity modulus $\sigma_G$ of
$G$ in Eq.~\eqref{paper:barrier-chain} satisfies
\begin{equation}\label{paper:barrier-modulus}
 \sigma_G\le
 \frac{2^{\frac{q}{2}}}{\left(1-e^{-1}\right)^q}
 \alpha n^{1-\frac{q}{2}}\left(1-\xi\right)^q.
\end{equation}
\end{lemma}
\begin{proof}
To bound the modulus from above, we denote $\delta=1-\xi$ and
choose $\bm w=\left(w_1,\ldots,w_n\right)$ with $w_i=1-\xi^i$.
This choice makes all residual coordinates equal to $\delta$, since
\[
 \begin{aligned}
 \left(A_\xi\bm w\right)_1&=1-\xi=\delta,\\
 \left(A_\xi\bm w\right)_i
 &=\left(1-\xi^i\right)-\xi\left(1-\xi^{i-1}\right)\\
 &=1-\xi^i-\xi+\xi^i=\delta
 \qquad\left(2\le i\le n\right).
 \end{aligned}
\]
For $i\ge\left\lceil \frac{n}{2}\right\rceil$, the length assumption gives
$i\delta\ge \frac{n\delta}{2}\ge1$. Consequently,
\[
 \begin{aligned}
 w_i&=1-\left(1-\delta\right)^i\\
 &\ge1-e^{-i\delta}
 &&\text{since }1-\delta\le e^{-\delta},\\
 &\ge1-e^{-1}
 &&\text{since }i\delta\ge1,\\
 \norm{\bm w}_2^2
 &=\sum_{i=1}^n w_i^2\\
 &\ge\frac{n}{2}\left(1-e^{-1}\right)^2
 &&\text{since at least }\frac{n}{2}\text{ coordinates satisfy this bound}.
 \end{aligned}
\]
Testing uniform convexity at
$\left(\bm x,\bm y\right)=\left(\bm x^\star,\bm x^\star+\bm w\right)$ gives
\[
 \begin{aligned}
 \frac{\sigma_G}{q}\norm{\bm w}_2^q
 &\le\Breg G{\bm x^\star+\bm w}{\bm x^\star}
 &&\text{by uniform convexity},\\
 &=G\left(\bm x^\star+\bm w\right)
 &&\text{since }G\left(\bm x^\star\right)=0,\ \D G\left(\bm x^\star\right)=\bm 0,\\
 &=\frac{\alpha}{q}\norm{A_\xi\bm w}_q^q
 &&\text{since }A_\xi\bm x^\star=a\bm e_1,\\
 &=\frac{\alpha}{q}n\delta^q
 &&\text{since }A_\xi\bm w=\delta\bm 1,\\
 \sigma_G
 &\le\frac{\alpha n\delta^q}{\norm{\bm w}_2^q}\\
 &\le\frac{\alpha n\delta^q}{\left(\frac{n}{2}\right)^{\frac{q}{2}}\left(1-e^{-1}\right)^q}
 &&\text{by the norm lower bound},\\
 &=\frac{2^{\frac{q}{2}}}{\left(1-e^{-1}\right)^q}
       \alpha n^{1-\frac{q}{2}}\left(1-\xi\right)^q.
 \end{aligned}
\]
The last line is Eq.~\eqref{paper:barrier-modulus}.
\end{proof}

This dimension dependence cannot be canceled by a dimension-dependent
improvement in the smoothness estimate. Since $A_\xi\bm e_n=\bm e_n$,
$G\left(\bm x^\star+t\bm e_n\right)=\frac{\alpha|t|^q}{q}$.
For $S_{q,p}:=\prod_{\ell=1}^{p-1}\left(q-\ell\right)$ and $t>0$,
Lemma~\resultpart{elem:power}{elem:power:part:derivatives} therefore gives
\[
 \begin{aligned}
 \norm{\D^pG\left(\bm x^\star+t\bm e_n\right)
               -\D^pG\left(\bm x^\star\right)}_{\rm op}
 &\ge
 \left|\left(\D^pG\left(\bm x^\star+t\bm e_n\right)
               -\D^pG\left(\bm x^\star\right)\right)
       \left[\bm e_n,\ldots,\bm e_n\right]\right|\\
 &=\alpha S_{q,p}t^\nu.
 \end{aligned}
\]
Every valid H\"older parameter $H_G$ consequently satisfies
$H_G\ge\alpha S_{q,p}$. Together with
Eq.~\eqref{paper:barrier-modulus}, this gives
\[
 \frac{H_G}{\sigma_G}
 \ge
 \frac{S_{q,p}\left(1-e^{-1}\right)^q}
      {2^{\frac{q}{2}}}
 n^{\frac{q-2}{2}}\left(1-\xi\right)^{-q}.
\]
The obstruction concerns this residual family;
it does not assert that every lower-bound construction must use
multiple phases.

%% file: appendix-preliminaries.tex
\section{Definitions and Elementary Tools}\label{app:paper-elementary}

\subsection{Additional notation and oracle model}\label{def:conventions}
We use the notation and Bregman divergence of
Sec.~\ref{paper:settings}, and the class assumptions of
Theorem~\ref{paper:main-theorem}.
For derivative evaluations, this gives
\[
 \D f\left(\bm x\right)\left[\bm z\right]=\ip{\D f\left(\bm x\right)}{\bm z},\qquad
 \D^2 f\left(\bm x\right)\left[\bm u,\bm v\right]=\bm u^\top \D^2 f\left(\bm x\right)\bm v.
\]
For symmetric matrices, we use $A\succeq B$ to mean that
$\bm z^\top\left(A-B\right)\bm z\ge0$ for every $\bm z$.
For a normed-space-valued map $g$ on $U$, we denote
$\norm g_{\infty;U}:=\sup_{\bm x\in U}\norm{g\left(\bm x\right)}$,
omitting $U$ when the domain is fixed.

If $E$ has orthonormal basis $\bm e_1,\ldots,\bm e_k$, its projection is
\[
 P_E\bm x=\sum_{i=1}^k\ip{\bm e_i}{\bm x}\bm e_i.
\]
The open and closed balls are
\[
 B\left(\bm c,R\right)=\{\bm x:\norm{\bm x-\bm c}<R\},\qquad
 \overline{B\left(\bm c,R\right)}=\{\bm x:\norm{\bm x-\bm c}\le R\}.
\]
For any $\bm x,\bm y\in\R^d$, the closed line segment joining them is
\[
 \left[\bm x,\bm y\right]:=\{\left(1-t\right)\bm x+t\bm y:0\le t\le1\}.
\]
Its length is $\norm{\bm y-\bm x}$.
The factor $\frac{1}{q}$ in Theorem~\ref{paper:main-theorem}'s
uniform-convexity convention means that a bound
$\Breg f{\bm y}{\bm x}\ge c\mu\norm{\bm y-\bm x}^q$
gives modulus $qc\mu$.

\begin{definition}[Deterministic algorithm with $p^{th}$-order oracle access]\label{def:oracle}
An oracle call at $\bm x$ returns the oracle output
$\left(F\left(\bm x\right),\D F\left(\bm x\right),\ldots,\D^pF\left(\bm x\right)\right)$.
Each evaluation point is a deterministic function of the prescribed
parameters, initialization, and preceding evaluation points and oracle outputs.
After $s$ calls, a candidate solution $\widehat{\bm x}_s$ is any deterministic
function of this information; it need not be an evaluation point.
A candidate solution $\widehat{\bm x}_0$ before the first oracle call is also allowed.
\end{definition}
Evaluation points and candidate solutions can have arbitrary norms.
In particular, no restriction to the span of previous gradients is assumed.

\subsection{Elementary functions used in the construction}\label{elem:section}
Throughout, $p\ge2$, $0<\nu\le1$, and $q=p+\nu>2$; vector norms are Euclidean and derivative norms are operator norms.

\subsubsection{Scalar and Euclidean functions}\label{elem:power-motivation}
We define
\[
 \phi\left(t\right)=\frac{|t|^q}{q},\qquad
 f\left(\bm x\right)=\frac{\norm{\bm x}^q}{q}.
\]
The factor $\frac{1}{q}$ simplifies the first derivative. Differentiation $p$ times reduces the degree
to $q-p=\nu$. The norm-power theorems cited in Lemma~\ref{elem:power} justify this scaling globally, including at zero.

For scalar arguments, $D_\phi\left(a,b\right)$ denotes the divergence
defined in Sec.~\ref{paper:settings} with $d=1$; the affine approximation
is taken at the second argument $b$.

\subsubsection{The fixed smooth-maximum profile}\label{elem:smooth-max}
We use the smooth maximum in Definition~\ref{paper:smooth-maximum-definition},
whose scalar convexity and equality properties are given in
Lemma~\ref{paper:smooth-maximum-properties}.
With its fixed density $\varphi$, we define
$\Theta\left(v\right):=\int_{-1}^1|v-u|\varphi\left(u\right)\,du$.
This function is smooth, convex, even, and $1$-Lipschitz, with
$\Theta\left(v\right)=|v|$ if and only if $|v|\ge1$, as in
\cite[Lemma~1 and its proof]{Azagra2013}.
It is fixed independently of the dimension and all phase parameters.
For $\delta>0$, the definitions in Eq.~\eqref{paper:smooth-maximum} equivalently give
\begin{equation}\label{elem:psi-M-def}
 \psi_\delta\left(v\right):=\frac{v+\delta\Theta\left(\frac{v}{\delta}\right)}{2},
 \qquad M_\delta\left(s,t\right):=s+\psi_\delta\left(t-s\right).
\end{equation}
Thus $\psi_\delta\left(v\right)=\delta\psi_1\left(\frac{v}{\delta}\right)$, and evenness gives
\[
 M_\delta\left(s,t\right)=\frac{s+t+\delta\Theta\left(\frac{s-t}{\delta}\right)}{2}.
\]
Lemma~\ref{elem:max} records the positive-part derivative and translation
properties used below; Lemma~\ref{elem:composition-lemma} gives the Hessian
and joining rules for its composition with convex functions.
All smoothing constants refer to this one fixed choice of $\Theta$.

\subsection{Elementary analytic tools}\label{app:elementary-support}

\subsubsection{Taylor's formula}

\begin{theorem}[Taylor's formula {\cite[Theorem~5.1]{ConradTaylor}}]
\label{sources:taylor}
Let $V,W$ be finite-dimensional real normed spaces, let $U\subset V$
be open, and let $f:U\to W$ be $C^k$, where $k\ge1$ is an integer.
Fix $\bm x\in U$ and $r>0$ such that $B\left(\bm x,r\right)\subset U$.
For $\norm{\bm h}<r$,
\[
\begin{aligned}
 f\left(\bm x+\bm h\right)
 &=\sum_{j=0}^{k}\frac{1}{j!}\D^jf\left(\bm x\right)\left[\bm h,\ldots,\bm h\right]
   +\mathcal E_k\left(\bm x,\bm h\right),\\
 \mathcal E_k\left(\bm x,\bm h\right)
 &=\frac{1}{\left(k-1\right)!}\int_0^1\left(1-t\right)^{k-1}
   \left(\D^kf\left(\bm x+t\bm h\right)-\D^kf\left(\bm x\right)\right)
       \left[\bm h,\ldots,\bm h\right]\,dt.
\end{aligned}
\]
The $j=0$ summand means $f\left(\bm x\right)$. With
\[
 C_{k,\bm h,\bm x}
 :=\frac{1}{k!}\sup_{0\le t\le1}
       \norm{\D^kf\left(\bm x+t\bm h\right)-\D^kf\left(\bm x\right)},
\]
one has $\norm{\mathcal E_k\left(\bm x,\bm h\right)}\le
C_{k,\bm h,\bm x}\norm{\bm h}^k$ and
$C_{k,\bm h,\bm x}\to0$ as $\bm h\to0$.
This convergence is uniform when $\bm x$ ranges over a compact subset of $U$.
Equivalently,
\begin{equation}\label{sources:taylor-integral}
 f\left(\bm x+\bm h\right)
 =\sum_{j=0}^{k-1}\frac{1}{j!}\D^jf\left(\bm x\right)\left[\bm h,\ldots,\bm h\right]
 +\frac{1}{\left(k-1\right)!}\int_0^1\left(1-t\right)^{k-1}
       \D^kf\left(\bm x+t\bm h\right)\left[\bm h,\ldots,\bm h\right]\,dt.
\end{equation}
\end{theorem}

For $f\in C^2$ and $\bm d=\bm y-\bm x$, Theorem~\ref{sources:taylor} gives
\begin{equation}\label{def:bregman-integral}
 \Breg f{\bm y}{\bm x}
 =\int_0^1\left(1-t\right)\D^2f\left(\bm x+t\bm d\right)\left[\bm d,\bm d\right]\,dt.
\end{equation}

\subsubsection{Mean-value estimates}

\begin{theorem}[Mean-value inequality {\cite[Theorem~3.1.3, p.~39]{Cartan1971}}]
\label{sources:mean-value}
Let $a<b$ be real. If $v:\left[a,b\right]\to X$ and $g:\left[a,b\right]\to\R$ are continuous,
$X$ is a Banach space, and outside a countable exceptional set their
right derivatives exist on $[a,b)$ and satisfy
$\norm{v'_+\left(t\right)}\le g'_+\left(t\right)$, then
\[
 \norm{v\left(b\right)-v\left(a\right)}\le g\left(b\right)-g\left(a\right).
\]
\end{theorem}

\begin{lemma}[Mean-value bound on a convex domain]\label{elem:mean-value}
Let $U\subset\R^d$ be open and convex, and let $f:U\to\R^a$
be differentiable. If $\norm{\D f\left(\bm z\right)}\le M$ for every
$\bm z\in U$, where $M\ge0$, then
\[
 \norm{f\left(\bm y\right)-f\left(\bm x\right)}
 \le M\norm{\bm y-\bm x}\qquad\left(\bm x,\bm y\in U\right).
\]
In particular, $\D f=0$ on $U$ implies that $f$ is constant on $U$.
\end{lemma}
\begin{proof}
Fix $\bm x,\bm y\in U$ and define, for $t\in[0,1]$,
$v\left(t\right)=f\left(\bm x+t\left(\bm y-\bm x\right)\right)$ and
$g\left(t\right)=Mt\norm{\bm y-\bm x}$. Convexity of $U$ places this segment
in the domain, and differentiability gives
\[
 \begin{aligned}
 \norm{v'_+\left(t\right)}
 &=\norm{\D f\left(\bm x+t\left(\bm y-\bm x\right)\right)\left[\bm y-\bm x\right]}\\
 &\le M\norm{\bm y-\bm x}
 &&\text{by the derivative bound},\\
 &=g'_+\left(t\right).
 \end{aligned}
\]
The continuous maps $v,g$ therefore satisfy the hypotheses of
Theorem~\ref{sources:mean-value}. It gives
\[
 \begin{aligned}
 \norm{f\left(\bm y\right)-f\left(\bm x\right)}
 &=\norm{v\left(1\right)-v\left(0\right)}\\
 &\le g\left(1\right)-g\left(0\right)
 &&\text{by Theorem~\ref{sources:mean-value}},\\
 &=M\norm{\bm y-\bm x}.
 \end{aligned}
\]
Taking $M=0$ proves the final assertion.
\end{proof}

\subsubsection{A quantitative H\"older composition theorem}\label{sources:holder-composition-section}

\begin{definition}[Lipschitz jets {\cite[Definition~1.1]{Boutaib2020}}]
\label{sources:lipschitz-jet}
Let $p\ge0$ be an integer, $0<\nu\le1$, and $q=p+\nu$.
Let $X,Y$ be normed spaces and $U\subset X$. A jet
$\mathfrak f=\left(f_0,f_1,\ldots,f_p\right)$ consists of $f_0:U\to Y$ and symmetric
$k$-linear maps $f_k\left(\bm x\right):X^k\to Y$ for $1\le k\le p$.
Regard each multilinear map as its linear map on the tensor power.
For $0\le k\le p$, $\bm x,\bm y\in U$, and every
$u\in X^{\otimes k}$, define the remainder by
\[
 f_k\left(\bm x\right)\left[u\right]
 =\sum_{a=0}^{p-k}\frac{1}{a!}f_{k+a}\left(\bm y\right)
 \left[u\otimes\left(\bm x-\bm y\right)^{\otimes a}\right]
 +\mathcal R_k\left(\bm x,\bm y\right)\left[u\right].
\]
For $u=\bm z_1\otimes\cdots\otimes\bm z_k$, the bracket is the original
multilinear evaluation, with $a$ additional copies of $\bm x-\bm y$;
for $k=0$ use $X^{\otimes0}=\R$.
The norm
$\norm{\mathfrak f}_{\operatorname{Lip}\left(q\right)}$ is the least $C$ such that,
for all $0\le k\le p$, $\bm x,\bm y\in U$, and
$\bm z_1,\ldots,\bm z_k\in X$,
\[
 \begin{aligned}
 \norm{f_k\left(\bm x\right)\left[\bm z_1,\ldots,\bm z_k\right]}
 &\le C\norm{\bm z_1\otimes\cdots\otimes\bm z_k},\\
 \norm{\mathcal R_k\left(\bm x,\bm y\right)\left[\bm z_1,\ldots,\bm z_k\right]}
 &\le C\norm{\bm x-\bm y}^{q-k}
           \norm{\bm z_1\otimes\cdots\otimes\bm z_k}.
 \end{aligned}
\]
The tensor norms have the \emph{projective property}
$\norm{a\otimes b}\le\norm a\norm b$ whenever the tensor orders match.
No openness or convexity of $U$ is required.
\end{definition}

\begin{theorem}[Characterization of Lipschitz jets {\cite[Theorem~1.20]{Boutaib2020}}]
\label{sources:holder-characterization}
Let $p\ge0$ be an integer, $0<\nu\le1$, $C\ge0$, let $X,Y$ be
normed spaces, and let $U\subset X$. Suppose the tensor norms on $X^{\otimes k}$, $1\le k\le p$,
are the canonical projective tensor norms. Given a jet $\mathfrak f=\left(f_0,\ldots,f_p\right)$ as in
Definition~\ref{sources:lipschitz-jet},
consider the assertions
\begin{enumerate}[label=\textnormal{(A\arabic*)}]
\item $\mathfrak f$ is $\operatorname{Lip}\left(p+\nu\right)$ and
$\norm{\mathfrak f}_{\operatorname{Lip}\left(p+\nu\right)}\le C$.
\item $f_0$ is $p$ times differentiable, its successive derivatives are
$f_1,\ldots,f_p$, all $\norm{f_k}_\infty\le C$, and
$\norm{f_p\left(\bm x\right)-f_p\left(\bm y\right)}\le C\norm{\bm x-\bm y}^{\nu}$ for
$\bm x,\bm y\in U$.
\end{enumerate}
If $U$ is open, (A1) implies (A2). If $U$ is also convex, the two assertions
are equivalent.
\end{theorem}

\begin{theorem}[Composition of Lipschitz jets {\cite[Theorem~1.48]{Boutaib2020}}]
\label{sources:holder-composition}
Let $X,Y,Z$ be normed spaces, $U\subset X$, $V\subset Y$, and $q>1$.
Define $p=\left\lceil q\right\rceil-1$ and $\nu=q-p$, so $p\ge1$ and $0<\nu\le1$.
Equip the tensor powers of $X$ and $Y$ with the canonical projective
tensor norms, and use the Lipschitz-jet convention of
Definition~\ref{sources:lipschitz-jet}. If $f:U\to Y$ and $h:V\to Z$
are $\operatorname{Lip}\left(q\right)$ maps
with $f\left(U\right)\subset V$, then the jet of $h\circ f$ obtained by formal
successive chain rules is $\operatorname{Lip}\left(q\right)$, and
\[
 \norm{h\circ f}_{\operatorname{Lip}\left(q\right)}
 \le C_q^{\rm comp}\norm{h}_{\operatorname{Lip}\left(q\right)}
          \max\left\{\norm{f}_{\operatorname{Lip}\left(q\right)}^q,1\right\}.
\]
The constant depends only on $q$; neither $U,V$ nor the dimensions enter it.
\end{theorem}

\begin{lemma}[Bounded H\"older norms under scalar composition]
\label{elem:holder-composition}
Let $p\ge1$ be an integer, $0<\nu\le1$, and $q=p+\nu$. For a $C^p$ scalar function
on an open set $W$, define the bounded norm
\[
 \norm{f}_{p,\nu;W}:=\max\left\{
   \max_{0\le k\le p}\norm{\D^kf}_{\infty;W},
   \sup_{\substack{\bm x,\bm y\in W\\\bm x\ne\bm y}}
   \frac{\norm{\D^pf\left(\bm x\right)-\D^pf\left(\bm y\right)}}
        {\norm{\bm x-\bm y}^{\nu}}\right\}.
\]
The domain subscript restricts both suprema to $W$. The first term
bounds the function and all its derivatives through order $p$; the second
bounds variations of the $p^{th}$ derivative.
Let $U\subset\R^d$ be convex and open, let $J\subset\R$ be an open
interval, and suppose $f\in C^p\left(U\right)$, $h\in C^p\left(J\right)$, $f\left(U\right)\subset J$,
and both bounded norms below are finite. Then
\[
 \norm{h\circ f}_{p,\nu;U}
 \le C_q^{\rm comp}\norm{h}_{p,\nu;J}
       \max\left\{1,\norm{f}_{p,\nu;U}^{q}\right\},
\]
where $C_q^{\rm comp}$ depends only on $q$, not on the dimension or
the domains. This includes $\nu=1$ and uses derivatives only through $p$.
\end{lemma}
\begin{proof}
Equip each tensor power of $\R^d$ with the projective norm
\[
 \norm{u}_\pi
 :=\inf_{u=\sum_\ell
       \bm z_{\ell,1}\otimes\cdots\otimes\bm z_{\ell,k}}
       \sum_\ell\prod_{i=1}^k\norm{\bm z_{\ell,i}}.
\]
These norms satisfy
$\norm{u\otimes\bm z}_\pi=\norm u_\pi\norm{\bm z}$.
Indeed, the projective property gives the upper bound. For the reverse
bound, suppose $\bm z\ne0$ and define
$\ell\left(\bm w\right)=\ip{\frac{\bm z}{\norm{\bm z}}}{\bm w}$.
Applying $I\otimes\ell$ to any pure-tensor decomposition and using
$\norm\ell=1$ gives $\norm{I\otimes\ell}\le1$ in the projective norms.
Therefore
\[
 \begin{aligned}
 \norm u_\pi\norm{\bm z}
 &=\norm{\left(I\otimes\ell\right)\left(u\otimes\bm z\right)}_\pi
 &&\text{since }\ell\left(\bm z\right)=\norm{\bm z},\\
 &\le\norm{u\otimes\bm z}_\pi
 &&\text{since }\norm{I\otimes\ell}\le1.
 \end{aligned}
\]
The case $\bm z=0$ is immediate.
For a $k$-linear form $\mathcal Q$, linearity gives, for every such
decomposition,
\[
 \begin{aligned}
 |\mathcal Q\left[u\right]|
 &=\left|\sum_\ell
       \mathcal Q\left[\bm z_{\ell,1},\ldots,\bm z_{\ell,k}\right]\right|
 &&\text{by linearity},\\
 &\le\sum_\ell
       |\mathcal Q\left[\bm z_{\ell,1},\ldots,\bm z_{\ell,k}\right]|
 &&\text{by the triangle inequality},\\
 &\le\norm{\mathcal Q}\sum_\ell
       \prod_{i=1}^k\norm{\bm z_{\ell,i}}
 &&\text{by the operator norm}.
 \end{aligned}
\]
Taking the infimum gives
$|\mathcal Q\left[u\right]|\le\norm{\mathcal Q}\norm{u}_\pi$.
Conversely, testing the associated linear map on pure tensors and using
$\norm{\bm z_1\otimes\cdots\otimes\bm z_k}_\pi
=\prod_i\norm{\bm z_i}$ recovers the full multilinear operator norm.
Thus the pure-tensor bounds in Definition~\ref{sources:lipschitz-jet} extend to
arbitrary tensors with exactly the same constant.

Define $\mathfrak f=\left(f,\D f,\ldots,\D^pf\right)$ and
$\mathfrak h=\left(h,h',\ldots,h^{\left(p\right)}\right)$.
Since $U,J$ are convex and open,
Theorem~\ref{sources:holder-characterization} gives
\[
 \norm{\mathfrak f}_{\operatorname{Lip}\left(q\right)}
 =\norm f_{p,\nu;U},\qquad
 \norm{\mathfrak h}_{\operatorname{Lip}\left(q\right)}
 =\norm h_{p,\nu;J}.
 \]
The assumption $f\left(U\right)\subset J$ verifies the range hypothesis of
Theorem~\ref{sources:holder-composition}. Its composition jet therefore
satisfies
\[
 \begin{aligned}
 \norm{h\circ f}_{p,\nu;U}
 &=\norm{\mathfrak h\circ\mathfrak f}_{\operatorname{Lip}\left(q\right)}
 &&\text{by Theorem~\ref{sources:holder-characterization}},\\
 &\le C_q^{\rm comp}\norm{\mathfrak h}_{\operatorname{Lip}\left(q\right)}
       \max\left\{1,\norm{\mathfrak f}_{\operatorname{Lip}\left(q\right)}^q\right\}
 &&\text{by Theorem~\ref{sources:holder-composition}},\\
 &=C_q^{\rm comp}\norm h_{p,\nu;J}
       \max\left\{1,\norm f_{p,\nu;U}^q\right\}
 &&\text{by the two norm identities}.
 \end{aligned}
\]
When $\nu=1$, we still have $p=q-1$, so
Theorems~\ref{sources:holder-characterization} and
\ref{sources:holder-composition} require derivatives only through $p$.
\end{proof}

\subsubsection{Orthogonal multilinear forms}\label{sources:direct-sum-section}

\begin{lemma}[Tensor direct sums {\cite[Definition~2.2 and Lemma~3.1, Eq.~(20)]{KongLiWang2018}}]
\label{sources:direct-sum}
Let $k\ge3$ and $n\ge1$ be integers, and let $d_j^{\left(i\right)}$ be positive
integers for $1\le j\le k$ and $1\le i\le n$. Let
$\mathcal Q_i\in\R^{d_1^{\left(i\right)}\times\cdots\times d_k^{\left(i\right)}}$ be real tensors.
For their tensor direct sum
\[
 \mathcal Q=\mathcal Q_1\oplus\cdots\oplus\mathcal Q_n
 \in\R^{\left(\sum_i d_1^{\left(i\right)}\right)\times\cdots\times\left(\sum_i d_k^{\left(i\right)}\right)},
\]
the operator norm satisfies
\[
 \norm{\mathcal Q}=\max_{1\le i\le n}\norm{\mathcal Q_i}.
\]
The spectral norm is defined by
\[
 \norm{\mathcal Q}=\max_{\norm{\bm z_1}=\cdots=\norm{\bm z_k}=1}
       |\mathcal Q\left[\bm z_1,\ldots,\bm z_k\right]|.
\]
The test vectors in the $k$ modes are independent. The direct sum places
each tensor in its
matching diagonal block and makes all mixed-block entries zero.
\end{lemma}

\begin{lemma}[Orthogonal sums of multilinear forms]
\label{elem:orthogonal-sum}\label{sources:orthogonal-forms}
Let $k\ge2$ be an integer, let $V_1,\ldots,V_n$ be mutually
orthogonal Euclidean subspaces, and let $\mathcal Q_i$ be a
$k$-linear form on $V_i$. Define the ambient form
\[
 \mathcal Q\left[\bm z_1,\ldots,\bm z_k\right]
 =\sum_{i=1}^n
 \mathcal Q_i\left[P_{V_i}\bm z_1,\ldots,P_{V_i}\bm z_k\right].
\]
Then its full multilinear operator norm satisfies
\[
 \norm{\mathcal Q}=\max_{1\le i\le n}\norm{\mathcal Q_i}.
\]
\end{lemma}
\begin{proof}
Zero blocks contribute nothing. If only one nonzero block is present,
restricting all test vectors to that block proves the identity.
For $k\ge3$ and at least two nonzero blocks, choose orthonormal coordinates
on each nonzero $V_i$.
The ambient form is the tensor direct sum of the coordinate
representations of $\mathcal Q_i$, with all mixed-block entries zero.
If the ambient space contains an orthogonal complement of
$\bigoplus_i V_i$, include its zero tensor block in this direct sum.
Lemma~\ref{sources:direct-sum} therefore gives the identity.

For $k=2$, define $M=\max_i\norm{\mathcal Q_i}$. For arbitrary unit
vectors $\bm z_1,\bm z_2$,
\[
 \begin{aligned}
 |\mathcal Q\left[\bm z_1,\bm z_2\right]|
 &\le M\sum_{i=1}^n
       \norm{P_{V_i}\bm z_1}\norm{P_{V_i}\bm z_2}
 &&\text{by the block operator norms},\\
 &\le M\left(\sum_i\norm{P_{V_i}\bm z_1}^2\right)^{\frac{1}{2}}
         \left(\sum_i\norm{P_{V_i}\bm z_2}^2\right)^{\frac{1}{2}}
 &&\text{by Cauchy--Schwarz},\\
 &\le M
 &&\text{by orthogonality and }\norm{\bm z_1}=\norm{\bm z_2}=1.
 \end{aligned}
\]
Taking the supremum gives $\norm{\mathcal Q}\le M$.
Conversely, choose an index $i$ attaining the largest block norm and
restrict both test vectors to $V_i$. The other blocks vanish, so
$\norm{\mathcal Q}\ge\norm{\mathcal Q_i}=M$. This proves equality.
\end{proof}

\subsubsection{Estimates for \texorpdfstring{$\frac{|t|^q}{q}$ and $\frac{\norm{\bm x}^q}{q}$}{the scalar and Euclidean functions}}\label{sources:power}

The scalar divergence $D_\phi\left(a,b\right)$ is defined in Appx.~\ref{elem:power-motivation}.

\begin{definition}[Derivative polynomials {\cite[Definition~3.1]{RodomanovNesterov2020}}]
\label{sources:power-polynomials}
For real $q$ and integers $k\ge0$, define the polynomials by
\[
 g_{0,q}\left(u\right)=1,\qquad
 g_{k,q}\left(u\right)=\left(1-u^2\right)g_{k-1,q}'\left(u\right)+\left(q-k+1\right)u g_{k-1,q}\left(u\right)
 \quad\left(k\ge1\right).
\]
\end{definition}

\begin{lemma}[Symmetric derivative norms and norm-power derivatives
{\cite[Sec.~2, Eq.~(1), Theorem~3.1, and Proposition~4.5]{RodomanovNesterov2020}}]
\label{sources:power-derivatives}
Equip $\R^d$ with its Euclidean norm. For $q>0$, define
$f_q\left(\bm x\right)=\norm{\bm x}^q$, and let $g_{k,q}$ be the polynomials
in Definition~\ref{sources:power-polynomials}.
\begin{claims}
\item\label{sources:power-derivatives:norm}
\textbf{Symmetric operator norm.}
For a symmetric $k$-linear form $\mathcal A$ on $\R^d$, $k\ge1$,
\[
 \norm{\mathcal A}
 =\max_{\norm{\bm z}=1}|\mathcal A\left[\bm z,\ldots,\bm z\right]|.
\]
\item\label{sources:power-derivatives:formula}
\textbf{Derivative formula.}
For real $q>0$ and every integer $0\le k<q$, $f_q$ is $k$ times
differentiable, and for every unit vector $\bm z$,
\[
 \D^kf_q\left(\bm x\right)\left[\bm z,\ldots,\bm z\right]
 =\norm{\bm x}^{q-k}g_{k,q}\left(\tau_{\bm z}\left(\bm x\right)\right),\qquad
 \tau_{\bm z}\left(\bm x\right):=
 \begin{cases}
 \frac{\ip{\bm x}{\bm z}}{\norm{\bm x}},&\bm x\ne0,\\
 0,&\bm x=0.
 \end{cases}
\]
At $k=0$ the displayed derivative means $f_q$ itself.
\item\label{sources:power-derivatives:polynomial}
\textbf{Polynomial bound.}
For every integer $k\ge0$ and real $q\ge k$,
\[
 \max_{|u|\le1}|g_{k,q}\left(u\right)|
 =\prod_{\ell=0}^{k-1}\left(q-\ell\right).
\]
An empty product equals one.
\end{claims}
\end{lemma}

\begin{lemma}[H\"older continuity of norm-power derivatives
{\cite[Theorems~6.2--6.3 and Eq.~(24)]{RodomanovNesterov2020}}]
\label{sources:power-holder}
Let $p\ge0$ be an integer, $0<\nu\le1$, $q=p+\nu$, and
$f_q\left(\bm x\right)=\norm{\bm x}^q$ on Euclidean $\R^d$. Let $g_{k,q}$ be
the polynomials in Definition~\ref{sources:power-polynomials}.
\begin{claims}
\item\label{sources:power-holder:line}
\textbf{Restriction to a line through the origin.}
The restriction of $\D^pf_q$ to any line through the origin is
$\nu$-H\"older continuous with constant
\[
 C_{p,\nu}:=
 \begin{cases}
 \displaystyle\prod_{i=1}^{p}\left(\nu+i\right),&p\text{ even},\\[3pt]
 \displaystyle2^{1-\nu}\prod_{i=1}^{p}\left(\nu+i\right),&p\text{ odd}.
 \end{cases}
\]
\item\label{sources:power-holder:global}
\textbf{Global H\"older continuity.}
The map $\D^pf_q$ is globally $\nu$-H\"older continuous with constant
\[
 A_{p,\nu}:=
 \begin{cases}
 \displaystyle\left(p-1\right)!!\prod_{i=1}^{\frac{p}{2}}\left(\nu+2i\right)+H_{p,\nu},
       &p\text{ even},\\[3pt]
 \displaystyle2^{1-\nu}\prod_{i=1}^{p}\left(\nu+i\right),&p\text{ odd},
 \end{cases}
\]
where $H_{p,\nu}$ denotes the $\nu$-H\"older constant of
$g_{p,p+\nu}$ on $[0,1]$:
\[
 H_{p,\nu}:=\sup_{0\le u<v\le1}
 \frac{|g_{p,p+\nu}\left(v\right)-g_{p,p+\nu}\left(u\right)|}{|v-u|^\nu}.
\]
In particular, a valid global constant is
\[
 \widetilde A_{p,\nu}:=
 \begin{cases}
 \displaystyle\left(p-1\right)!!\prod_{i=1}^{\frac{p}{2}}\left(\nu+2i\right)
                 +\prod_{i=1}^{p}\left(\nu+i\right),
       &p\text{ even},\\[3pt]
 \displaystyle2^{1-\nu}\prod_{i=1}^{p}\left(\nu+i\right),&p\text{ odd}.
 \end{cases}
\]
For even $p$, $\left(p-1\right)!!=\prod_{i=1}^{\frac{p}{2}}\left(2i-1\right)$, including
$\left(-1\right)!!=1$. Global H\"older continuity here means
$\norm{\D^pf_q\left(\bm x\right)-\D^pf_q\left(\bm y\right)}
\le A_{p,\nu}\norm{\bm x-\bm y}^\nu$ for all $\bm x,\bm y$,
and likewise with $\widetilde A_{p,\nu}$.
\end{claims}
\end{lemma}

\begin{lemma}[Uniform monotonicity implies uniform convexity
{\cite[Lemma~2.1]{DoikovNesterov2021}}]
\label{sources:uniform-monotonicity}
Let $C$ be a convex subset of the open domain of a differentiable
function $f$. If, for some $\sigma>0$ and real $q\ge2$,
\[
 \ip{\D f\left(\bm y\right)-\D f\left(\bm x\right)}{\bm y-\bm x}
 \ge \sigma\norm{\bm y-\bm x}^q\qquad\left(\bm x,\bm y\in C\right),
\]
then $f$ is uniformly convex of degree $q$ on $C$ with parameter
$\sigma$, i.e.,
\[
 f\left(\bm y\right)\ge f\left(\bm x\right)+\ip{\D f\left(\bm x\right)}{\bm y-\bm x}
             +\frac{\sigma}{q}\norm{\bm y-\bm x}^q
 \qquad\left(\bm x,\bm y\in C\right).
\]
\end{lemma}

\begin{lemma}[Uniform monotonicity of the Euclidean norm power
{\cite[Lemma~2.5(1), Eq.~(16)]{DoikovNesterov2021}}]
\label{sources:power-monotonicity}
For a fixed real $q\ge2$, the function
$f\left(\bm x\right)=\frac{\norm{\bm x}^q}{q}$, $\bm x\in\R^d$, is uniformly
convex of degree $q$ and satisfies
\[
 \ip{\D f\left(\bm y\right)-\D f\left(\bm x\right)}{\bm y-\bm x}
 \ge 2^{2-q}\norm{\bm y-\bm x}^q\qquad\left(\bm x,\bm y\in\R^d\right).
\]
\end{lemma}

\begin{lemma}[Scalar and Euclidean norm estimates]\label{elem:power}
For every $d\geq1$, let $\phi\left(t\right)=\frac{|t|^q}{q}$ and $f\left(\bm x\right)=\frac{\norm{\bm x}^q}{q}$.
The constants below depend only on $q$ and the indicated derivative order.
\begin{claims}
\item\label{elem:power:part:derivatives} \textbf{Derivatives away from zero.}
The function $f$ is smooth on $\R^d\setminus\{\bm0\}$. With $C_{q,k}$ defined in
Eq.~\eqref{elem:Cqk}, every integer $1\leq k\leq p$ satisfies
\begin{equation}\label{elem:derivative-bound}
 \norm{\D^k f\left(\bm x\right)}\leq C_{q,k}\norm{\bm x}^{q-k}
 \qquad\left(\bm x\ne\bm0\right).
\end{equation}
For $1\leq k\leq p$, the scalar derivatives are
\begin{equation}\label{elem:scalar-derivatives}
 \begin{aligned}
 S_{q,k}&:=\prod_{\ell=1}^{k-1}\left(q-\ell\right),\qquad S_{q,1}=1,\\
 \phi^{\left(k\right)}\left(a\right)&=S_{q,k}\operatorname{sign}\left(a\right)^k|a|^{q-k}\qquad\left(a\ne0\right).
 \end{aligned}
\end{equation}
\item\label{elem:power:part:flatness} \textbf{Derivatives at zero.}
Both functions belong to $C^p$, and all their derivatives of orders $1,\ldots,p$ vanish at
zero. For $1\leq k\leq p$, Eq.~\eqref{elem:derivative-bound} also holds at $\bm x=\bm0$.
For every linear map $r:\R^d\to\R$, all derivatives of $\phi\circ r$ through order $p$
vanish at any point where $r\left(\bm x\right)=0$.
\item\label{elem:power:part:hessian} \textbf{Hessian lower bound.}
For every $\bm x\in\R^d$,
\begin{equation}\label{elem:hessian-lower}
 \D^2f\left(\bm x\right)\succeq\norm{\bm x}^{q-2}I.
\end{equation}
In particular, the Hessian is positive definite when $\bm x\ne\bm0$.
\item\label{elem:power:part:holder} \textbf{Global H\"older estimates.}
With $H_{\rm pow}$ defined in Eq.~\eqref{elem:holder-constant} and
$L_\phi:=2^{1-\nu}S_{q,p}$, for every $\bm x,\bm y\in\R^d$ and $a,b\in\R$,
\begin{align}
 \norm{\D^p f\left(\bm x\right)-\D^p f\left(\bm y\right)}
 &\leq H_{\rm pow}\norm{\bm x-\bm y}^{\nu},\label{elem:holder-bound}\\
 |\phi^{\left(p\right)}\left(a\right)-\phi^{\left(p\right)}\left(b\right)|
 &\leq L_\phi|a-b|^\nu.\label{elem:scalar-holder}
\end{align}
\item\label{elem:power:part:bregman} \textbf{Scalar divergence lower bound.}
With $c_q:=\frac{2^{2-q}}{q}>0$,
\begin{equation}\label{elem:bregman-lower}
 D_\phi\left(a,b\right)\geq c_q|a-b|^q
 \qquad\left(a,b\in\R\right).
\end{equation}
\end{claims}
\end{lemma}
\begin{proof}
\begin{itemize}
\item \textbf{Proof of \ref{elem:power:part:derivatives}.}
Away from zero, $f=h\circ s$ with $s\left(\bm x\right)=\ip{\bm x}{\bm x}$ and $h\left(v\right)=\frac{v^{\frac{q}{2}}}{q}$ for $v>0$, so $f$ is smooth. Direct differentiation gives $\D f\left(\bm x\right)=r^{q-2}\bm x$, where $r=\norm{\bm x}$, and
\begin{equation}\label{elem:hessian-formula}
 \D^2 f\left(\bm x\right)=r^{q-2}I+\left(q-2\right)r^{q-4}\bm x\bm x^{\mathsf T}.
\end{equation}
For $1\le k\le p<q$,
Lemma~\ref{sources:power-derivatives} gives the derivative norm bound with the constants
\begin{equation}\label{elem:Cqk}
 C_{q,k}:=S_{q,k}=\prod_{\ell=1}^{k-1}\left(q-\ell\right)
 \qquad\left(1\le k\le p\right),
\end{equation}
where the empty product is $1$. Lemma~\resultpart{sources:power-derivatives}{sources:power-derivatives:norm}
identifies the diagonal norm of a symmetric form with its full operator norm, giving
\begin{equation}\label{elem:nonzero-derivative-bound}
 \norm{\D^k f\left(\bm x\right)}\leq C_{q,k}\norm{\bm x}^{q-k}
 \qquad\left(\bm x\ne\bm0,\ 1\le k\le p\right).
\end{equation}
Dividing the derivative product by $q$ gives $C_{q,k}$.
Ordinary differentiation on the two half-lines gives
Eq.~\eqref{elem:scalar-derivatives}.

\item \textbf{Proof of \ref{elem:power:part:flatness}.}
Lemma~\resultpart{sources:power-derivatives}{sources:power-derivatives:formula} gives derivatives of $\norm{\bm x}^q$ for every integer $k<q$, with those of positive order zero at the origin. Applying this formula for $1\le k\le p<q$ and dividing by $q$ proves the vanishing derivatives of $f$. Eq.~\eqref{elem:nonzero-derivative-bound}, with $q-k>0$, gives continuity there and extends the bound to zero. The scalar case is $d=1$. For a linear $r$, the chain rule gives
\[
 \D^k\left(\phi\circ r\right)\left(\bm x\right)\left[\bm z_1,\ldots,\bm z_k\right]
 =\phi^{\left(k\right)}\left(r\left(\bm x\right)\right)\prod_{i=1}^k r\left(\bm z_i\right),
\]
which is zero whenever $r\left(\bm x\right)=0$.

\item \textbf{Proof of \ref{elem:power:part:hessian}.}
For $\bm x\ne0$ and any $\bm z\in\R^d$,
Eq.~\eqref{elem:hessian-formula} gives
\[
 \begin{aligned}
 \D^2f\left(\bm x\right)\left[\bm z,\bm z\right]
 &=\norm{\bm x}^{q-2}\norm{\bm z}^2
   +\left(q-2\right)\norm{\bm x}^{q-4}\ip{\bm x}{\bm z}^2\\
 &\ge\norm{\bm x}^{q-2}\norm{\bm z}^2
 &&\text{since }q>2.
 \end{aligned}
\]
This is the matrix inequality Eq.~\eqref{elem:hessian-lower}.
At zero both sides vanish by part~\ref{elem:power:part:flatness}.

\item \textbf{Proof of \ref{elem:power:part:holder}.}
Lemma~\resultpart{sources:power-holder}{sources:power-holder:global}
gives the dimension-independent constant $\widetilde A_{p,\nu}$.
For even $p$, compare each odd factor separately:
\[
 \begin{aligned}
 \left(p-1\right)!!\prod_{i=1}^{\frac{p}{2}}\left(\nu+2i\right)
 &=\prod_{i=1}^{\frac{p}{2}}\left(2i-1\right)\left(\nu+2i\right)\\
 &\le\prod_{i=1}^{\frac{p}{2}}\left(\nu+2i-1\right)\left(\nu+2i\right)
 &&\text{since }\nu>0,\\
 &=\prod_{i=1}^{p}\left(\nu+i\right).
 \end{aligned}
\]
and hence $\widetilde A_{p,\nu}\le2\prod_{i=1}^p\left(\nu+i\right)$.
For odd $p$, the same inequality follows from $2^{1-\nu}\le2$.
Thus the $\nu$-H\"older constant of $\D^p\norm{\bm x}^{p+\nu}$ is at most
$2\prod_{i=1}^p\left(\nu+i\right)=2qS_{q,p}$ for $0<\nu\le1$.
Dividing by $q$ gives Eq.~\eqref{elem:holder-bound} with
\begin{equation}\label{elem:holder-constant}
 H_{\rm pow}:=2S_{q,p}.
\end{equation}
For the scalar estimate, Lemma~\resultpart{sources:power-holder}{sources:power-holder:line}, restricted to $d=1$ and divided by $q$, gives the constant $S_{q,p}$ for even $p$ and $2^{1-\nu}S_{q,p}$ for odd $p$. Both are at most the stated $L_\phi=2^{1-\nu}S_{q,p}$, proving Eq.~\eqref{elem:scalar-holder}.

\item \textbf{Proof of \ref{elem:power:part:bregman}.}
Lemma~\ref{sources:power-monotonicity}, applied on $\R$, gives
$\left(\phi'\left(a\right)-\phi'\left(b\right)\right)\left(a-b\right)\ge2^{2-q}|a-b|^q$. Lemma~\ref{sources:uniform-monotonicity} then gives
$D_\phi\left(a,b\right)\ge\left(\frac{2^{2-q}}{q}\right)|a-b|^q$, which is Eq.~\eqref{elem:bregman-lower}.
\end{itemize}
\end{proof}

\subsubsection{Positive-part derivatives and common translation}
We use the scaled functions in Eq.~\eqref{elem:psi-M-def}.

\begin{lemma}[Positive-part derivatives and common translation]\label{elem:max}
For $\delta>0$, let $\psi_\delta$ and $M_\delta$ be defined by Eq.~\eqref{elem:psi-M-def}.
\begin{claims}
\item\label{elem:max:part:psi-exact} \textbf{Exact formulas for $\psi_\delta$.}
\begin{align}
 \psi_\delta\left(v\right)&=0 &&\left(v\leq-\delta\right),\label{elem:psi-left}\\
 \psi_\delta\left(v\right)&=v &&\left(v\geq\delta\right).\label{elem:psi-right}
\end{align}
\item\label{elem:max:part:psi-derivatives} \textbf{Smoothness and derivatives of $\psi_\delta$.}
The function $\psi_\delta$ belongs to $C^\infty\left(\R\right)$ and is convex and nondecreasing, with
\begin{equation}\label{elem:psi-estimates}
 0\leq\psi_\delta'\left(v\right)\leq1.
\end{equation}
\item\label{elem:max:part:translation} \textbf{Common translation.}
For every $a\in\R$,
\begin{equation}\label{elem:max-shift}
 M_\delta\left(s+a,t+a\right)=M_\delta\left(s,t\right)+a.
\end{equation}
\end{claims}
\end{lemma}
\begin{proof}
\begin{itemize}
\item \textbf{Proof of \ref{elem:max:part:psi-exact}.}
Substitute $\Theta\left(\frac{v}{\delta}\right)=\frac{|v|}{\delta}$ into Eq.~\eqref{elem:psi-M-def}
when $|v|\ge\delta$. Then
\[
 \begin{aligned}
 \psi_\delta\left(v\right)
 &=\frac{v+\delta\left(\frac{|v|}{\delta}\right)}{2}
 &&\text{by Eq.~\eqref{elem:psi-M-def}},\\
 &=\frac{v+|v|}{2}
 =\begin{cases}0,&v\le-\delta,\\v,&v\ge\delta.\end{cases}
 \end{aligned}
\]

\item \textbf{Proof of \ref{elem:max:part:psi-derivatives}.}
Differentiate the fixed smooth profile:
\begin{equation}\label{elem:psi-derivatives}
 \psi_\delta'\left(v\right)=\frac{1+\Theta'\left(\frac{v}{\delta}\right)}{2},\qquad
 \psi_\delta''\left(v\right)=\frac{\Theta''\left(\frac{v}{\delta}\right)}{2\delta}.
\end{equation}
The bounds $|\Theta'|\le1$ and $\Theta''\ge0$ give
Eq.~\eqref{elem:psi-estimates}, monotonicity, and convexity.

\item \textbf{Proof of \ref{elem:max:part:translation}.}
By Eq.~\eqref{elem:psi-M-def},
\[
 \begin{aligned}
 M_\delta\left(s+a,t+a\right)
 &=\left(s+a\right)+\psi_\delta\left(\left(t+a\right)-\left(s+a\right)\right)\\
 &=s+\psi_\delta\left(t-s\right)+a
 =M_\delta\left(s,t\right)+a.
 \end{aligned}
\]
\end{itemize}
\end{proof}

\subsubsection{Composition and joining}\label{elem:composition}\label{sources:convexity}

\begin{corollary}[Smooth convex replacement]
\label{paper:smooth-convex-replacement}
Let $U\subset\R^d$ be open and convex, let $f,g\in C^p\left(U\right)$
be convex, and fix $\delta>0$. Define
$h\left(\bm x\right):=M_\delta\left(f\left(\bm x\right),
g\left(\bm x\right)\right)$ for $\bm x\in U$. Then:
\begin{claims}
\item\label{paper:smooth-convex-replacement:regularity}
$h$ is convex and belongs to $C^p\left(U\right)$.
\item\label{paper:smooth-convex-replacement:selection}
If $\left|f\left(\bm x\right)-g\left(\bm x\right)\right|>\delta$
at $\bm x\in U$, there is an open neighborhood $V\subset U$
of $\bm x$ such that, for every $\bm y\in V$ and $0\le k\le p$,
$\displaystyle
 \D^kh\left(\bm y\right)=
 \begin{cases}
 \D^kf\left(\bm y\right),
   &f\left(\bm x\right)>g\left(\bm x\right)+\delta,\\
 \D^kg\left(\bm y\right),
   &g\left(\bm x\right)>f\left(\bm x\right)+\delta.
 \end{cases}
$
\end{claims}
\end{corollary}
\begin{proof}[Proof of Corollary~\ref{paper:smooth-convex-replacement}]
\textbf{\ref{paper:smooth-convex-replacement:regularity}}
Convexity follows from convexity and coordinatewise monotonicity of
$M_\delta$; $C^p$ regularity follows by the chain rule.

\textbf{\ref{paper:smooth-convex-replacement:selection}}
Continuity of $f-g$ preserves whichever strict inequality holds
at $\bm x$ on an open neighborhood $V\subset U$.
Lemma~\ref{paper:smooth-maximum-properties} then gives $h=f$
or $h=g$ throughout $V$, so their derivatives through order $p$
agree there.
\end{proof}

\begin{lemma}[Convexity criteria and operations
{\cite[Secs.~3.1.3--3.1.4 and~3.2.1--3.2.2]{BoydVandenberghe2004}}]
\label{sources:convex-rules}
For real-valued functions on Euclidean spaces, the following statements hold.
\begin{claims}
\item\label{sources:convex-rules:first}
Suppose $f$ is differentiable on its open domain $U$.
Then $f$ is convex if and only if $U$ is convex and
\[
 f\left(\bm y\right)\ge f\left(\bm x\right)+\ip{\D f\left(\bm x\right)}{\bm y-\bm x}
 \qquad\left(\bm x,\bm y\in U\right).
\]
Equivalently,
$\Breg f{\bm y}{\bm x}\ge0$.
For convex $f$, $\bm x,\bm y\in U$, and $0\le t\le1$,
the equivalent value inequality is
\[
 f\left(\left(1-t\right)\bm x+t\bm y\right)
 \le\left(1-t\right)f\left(\bm x\right)+t f\left(\bm y\right).
\]
\item\label{sources:convex-rules:second}
Suppose $f$ is twice differentiable on its open domain $U$.
Then $f$ is convex if and only if $U$ is convex and
$\D^2f\left(\bm x\right)\succeq0$ for all $\bm x\in U$.
\item\label{sources:convex-rules:sum}
A nonnegative weighted sum $\sum_{i=1}^n a_i f_i$ of convex functions
is convex, and a nonnegative weighted sum of concave functions is concave.
\item\label{sources:convex-rules:affine}
For a linear map $\mathcal L:\R^d\to\R^n$ and $\bm b\in\R^n$,
set $g\left(\bm x\right)=f\left(\mathcal L\bm x+\bm b\right)$ with
$\operatorname{dom}g=\{\bm x:\mathcal L\bm x+\bm b\in\operatorname{dom}f\}$.
If $f$ is convex, then $g$ is convex; if $f$ is concave, then $g$ is concave.
\end{claims}
\end{lemma}

\begin{lemma}[Convex additions and affine terms]\label{paper:bregman-additions}
Let $F,G:\R^d\to\R$ be differentiable, with $G$ convex, and let
$\ell:\R^d\to\R$ be affine. For all $\bm x,\bm y\in\R^d$,
\[
 D_{F+G+\ell}\left(\bm y,\bm x\right)
 =D_F\left(\bm y,\bm x\right)+D_G\left(\bm y,\bm x\right)
 \ge D_F\left(\bm y,\bm x\right).
\]
\end{lemma}
\begin{proof}
Expanding the Bregman definition in Sec.~\ref{paper:settings} gives additivity and
$D_\ell=0$. Lemma~\resultpart{sources:convex-rules}{sources:convex-rules:first}
gives $D_G\ge0$.
\end{proof}

\begin{lemma}[Composition and exact joining]\label{elem:composition-lemma}
Let $U\subseteq\R^d$ be open and convex, let $r\geq2$ be an integer, and let
$f_1,f_2\in C^r\left(U\right)$ be convex. For $\delta>0$, define
\[
 F\left(\bm x\right):=M_\delta\left(f_1\left(\bm x\right),f_2\left(\bm x\right)\right),\qquad
 g:=f_2-f_1,\qquad a\left(\bm x\right):=\psi_\delta'\left(g\left(\bm x\right)\right).
\]
\begin{claims}
\item\label{elem:composition-lemma:part:convexity} \textbf{Smooth convex composition.}
The function $F$ belongs to $C^r\left(U\right)$ and is convex. Its gradient and Hessian satisfy
Eq.~\eqref{elem:composition-gradient} and \eqref{elem:composition-hessian}.
\item\label{elem:composition-lemma:part:hessian} \textbf{Preservation of a common Hessian bound.}
If $\D^2f_i\left(\bm x\right)\succeq K\left(\bm x\right)$ for $i=1,2$, where $K\left(\bm x\right)$ is symmetric, then
$\D^2F\left(\bm x\right)\succeq K\left(\bm x\right)$.
\item\label{elem:composition-lemma:part:exact} \textbf{Agreement with the larger input.}
On any open region where $f_i-f_{3-i}\geq\delta$, the equality $F=f_i$ holds,
with $\D^kF=\D^kf_i$ for $0\le k\le r$.
\item\label{elem:composition-lemma:part:joining} \textbf{Joining to the unchanged function.}
Suppose $\overline{B\left(\bm c,R\right)}\subset U$ and $f_1-f_2\geq\delta$ on an open
neighborhood of the sphere $\{\bm x:\norm{\bm x-\bm c}=R\}$. Then
\[
 \widetilde F\left(\bm x\right)=
 \begin{cases}
 M_\delta\left(f_1\left(\bm x\right),f_2\left(\bm x\right)\right),&\bm x\in B\left(\bm c,R\right),\\
 f_1\left(\bm x\right),&\bm x\in U\setminus B\left(\bm c,R\right)
 \end{cases}
\]
is convex and belongs to $C^r\left(U\right)$.
\end{claims}
\end{lemma}
\begin{proof}
\begin{itemize}
\item \textbf{Proof of \ref{elem:composition-lemma:part:convexity}.}
Corollary~\ref{paper:smooth-convex-replacement}, applied with $p=r$,
gives convexity and $C^r$ regularity of
$F=M_\delta\left(f_1,f_2\right)$ on the convex open set $U$.
Equivalently, $F=f_1+\psi_\delta\left(g\right)$, with $g=f_2-f_1$.
The chain rule gives
\begin{equation}\label{elem:composition-gradient}
 \begin{aligned}
 \D F
 &=\D f_1+\psi_\delta'\left(g\right)\D g\\
 &=\D f_1+a\left(\D f_2-\D f_1\right)
 &&\text{by the definitions of }a,g,\\
 &=\left(1-a\right)\D f_1+a\D f_2.
 \end{aligned}
\end{equation}
\begin{align}
 \D^2F
 &=\D^2f_1+\psi_\delta'\left(g\right)\D^2g
   +\psi_\delta''\left(g\right)\D g\D g^{\mathsf T}
 &&\text{by the chain rule},\notag\\
 &=\D^2f_1+a\left(\D^2f_2-\D^2f_1\right)
   +\psi_\delta''\left(g\right)\D g\D g^{\mathsf T}
 &&\text{by the definitions of }a,g,\notag\\
 &=\left(1-a\right)\D^2f_1+a\D^2f_2
   +\psi_\delta''\left(g\right)\D g\D g^{\mathsf T}.
 \label{elem:composition-hessian}
\end{align}

\item \textbf{Proof of \ref{elem:composition-lemma:part:hessian}.}
For any $\bm z\in\R^d$, Eq.~\eqref{elem:composition-hessian} gives
\[
 \begin{aligned}
 \bm z^{\mathsf T}\left(\D^2F-K\right)\bm z
 &=\left(1-a\right)\bm z^{\mathsf T}\left(\D^2f_1-K\right)\bm z
   +a\bm z^{\mathsf T}\left(\D^2f_2-K\right)\bm z\\
 &\qquad+\psi_\delta''\left(g\right)\ip{\D g}{\bm z}^2\\
 &\ge0.
 \end{aligned}
\]
The last line uses the two assumed Hessian bounds and
$0\le a\le1$, $\psi_\delta''\ge0$ from
Lemma~\resultpart{elem:max}{elem:max:part:psi-derivatives}.
Thus $\D^2F\succeq K$.

\item \textbf{Proof of \ref{elem:composition-lemma:part:exact}.}
Lemma~\ref{paper:smooth-maximum-properties} gives $F=f_i$ throughout
the stated open region, since $f_i-f_{3-i}\ge\delta$ there.
Differentiating this equality gives agreement through order $r$.

\item \textbf{Proof of \ref{elem:composition-lemma:part:joining}.}
Part~\ref{elem:composition-lemma:part:exact} makes the two definitions
agree on an open neighborhood of the joining sphere. Thus
$\widetilde F\in C^r\left(U\right)$ and its Hessian is positive semidefinite everywhere,
using part~\ref{elem:composition-lemma:part:convexity} inside and convexity
of $f_1$ outside. Lemma~\resultpart{sources:convex-rules}{sources:convex-rules:second},
applied on the convex set $U$, proves convexity.
By continuity, the strict inequality $f_1-f_2>\delta$ on the sphere
suffices for the neighborhood hypothesis.
\end{itemize}
\end{proof}

%% file: appendix-base.tex
\section{Supporting Proofs for the Uniformly Convex Base}\label{app:paper-base}

\subsection{Constructing and verifying the uniformly convex component}
\label{base:section}
To give each of $L$ phases the required Euclidean $q^{th}$-power local formula,
we anchor phase $j$ at a center $\bm c_j$, as in Step~1. We use
$q=p+\nu>2$, $\mu>0$, and consecutive centers satisfying
\[
 \bm c_0=0,\qquad \bm c_{j+1}=\bm c_j+\bm h_j,
 \qquad R_j:=\norm{\bm h_j}>0,
\]
where the phase displacements $\bm h_j$ are mutually orthogonal, and their lengths
will decrease geometrically. The first center is the origin, so we start
with $B_0\left(\bm x\right)=\frac{\mu\norm{\bm x}^q}{q}$.
Let $B_j$ denote the partial base after adding the local formula at
$\bm c_j$, and set $B=B_L$.

\subsubsection{Requirements for the base}
\label{base:motivation}\label{base:consistency}\label{base:gradient-design}
For the cancellation Eq.~\eqref{design:gradient-cancellation}, we choose the base increment below. We also require uniform convexity and H\"older smoothness:
\begin{equation}\label{base:requirements}
\begin{aligned}
 \D B\left(\bm c_{j+1}\right)-\D B\left(\bm c_j\right)
   &=\mu R_j^{q-2}\bm h_j &&\left(0\le j<L\right),\\
 \Breg B{\bm y}{\bm x}
   &\ge c_B\mu\norm{\bm y-\bm x}^q &&\left(\bm x,\bm y\in\R^d\right),\\
 \norm{\D^pB\left(\bm y\right)-\D^pB\left(\bm x\right)}
   &\le C_B\mu\norm{\bm y-\bm x}^\nu
   &&\left(\bm x,\bm y\in\R^d\right).
\end{aligned}
\end{equation}
The positive constants $c_B,C_B$, independent of dimension and the number of phases, will be derived in the proofs.
Lemma~\ref{base:nesting-lemma} establishes the common construction. We then verify H\"older smoothness, preservation of oracle outputs, and uniform convexity in separate subsections. 

Preserving earlier oracle outputs is a separate requirement. If the contribution of $B_j$ has already been used at an evaluation point $\bm x_s$, we require
\begin{equation}\label{base:past-answers}
 B_k\left(\bm x_s\right)=B_j\left(\bm x_s\right),\qquad
 \D^rB_k\left(\bm x_s\right)=\D^rB_j\left(\bm x_s\right)
 \quad\left(1\le r\le p, k>j\right).
\end{equation}
We will change $B_j$ only where $\norm{\bm x-\bm c_{j+1}}<\theta R_j$, with $\theta>0$ to be chosen later so that this ball misses the previous evaluation points $\bm x_s$. All later changes will lie inside that same ball. Lemma~\resultpart{base:prefix-lemma}{base:prefix-lemma:part:neighborhood} proves that this gives neighborhood agreement at every earlier evaluation point. Centers index phases and need not be evaluation points. Thus Eq.~\eqref{base:past-answers} does not prescribe the gradient at a new center that has not been evaluated.

\subsubsection{Choose a local formula that can be repeated at successive centers}
\label{base:obstruction}\label{design:local-model}
To obtain the first line of Eq.~\eqref{base:requirements}, we impose the following local gradient rule on a ball about $\bm c_j$:
\begin{equation}\label{base:gradient-profile}
 \D B_j\left(\bm x\right)-\D B_j\left(\bm c_j\right)
 =\mu\norm{\bm x-\bm c_j}^{q-2}\left(\bm x-\bm c_j\right)
\end{equation}
Lemma~\resultpart{base:profile-lemma}{base:profile-lemma:part:equivalence}, proved in Appx.~\ref{app:base-support}, shows that this is equivalent to
\begin{equation}\label{base:old-local-form}
 B_j\left(\bm x\right)=B_j\left(\bm c_j\right)
 +\ip{\D B_j\left(\bm c_j\right)}{\bm x-\bm c_j}
 +\frac{\mu}{q}\norm{\bm x-\bm c_j}^q.
\end{equation}
When this ball contains $\bm c_{j+1}$, evaluating the gradient there gives the desired increment $\mu R_j^{q-2}\left(\bm c_{j+1}-\bm c_j\right)$.
Keeping this same expression for the next phase would also change the preceding block's gradient component. Indeed, when $\bm x$ and $\bm c_{j+1}$ lie in the region where Eq.~\eqref{base:old-local-form} holds and $\bm x-\bm c_{j+1}\perp\bm c_{j+1}-\bm c_j$, differentiation gives
\begin{equation}\label{base:old-gradient-increment}
\begin{split}
 &\D B_j\left(\bm x\right)-\D B_j\left(\bm c_{j+1}\right)\\
 ={}&\mu\left(\norm{\bm c_{j+1}-\bm c_j}^2
                 +\norm{\bm x-\bm c_{j+1}}^2\right)^{\frac{q-2}{2}}
       \left(\bm x-\bm c_{j+1}\right)\\
 &+\mu\left[
 \left(\norm{\bm c_{j+1}-\bm c_j}^2
       +\norm{\bm x-\bm c_{j+1}}^2\right)^{\frac{q-2}{2}}
       -\norm{\bm c_{j+1}-\bm c_j}^{q-2}\right]
       \left(\bm c_{j+1}-\bm c_j\right).
\end{split}
\end{equation}
The coefficient on the last line is positive when $\bm x-\bm c_{j+1}\ne\bm 0$, because $q>2$. Thus movement away from $\bm c_{j+1}$ changes the gradient component along $\bm h_j=\bm c_{j+1}-\bm c_j$. To keep this component fixed during the next phase, we require
\begin{equation}\label{base:new-own-form}
 B_{j+1}\left(\bm x\right)=B_{j+1}\left(\bm c_{j+1}\right)
 +\ip{\D B_{j+1}\left(\bm c_{j+1}\right)}{\bm x-\bm c_{j+1}}
 +\frac{\mu}{q}\norm{\bm x-\bm c_{j+1}}^q.
\end{equation}
To retain the gradient increment just obtained, we set
\begin{equation}\label{base:same-center-gradient}
 \D B_{j+1}\left(\bm c_{j+1}\right)=\D B_j\left(\bm c_{j+1}\right).
\end{equation}
We also set $B_{j+1}\left(\bm c_{j+1}\right)=B_j\left(\bm c_{j+1}\right)+C_j$, where $C_j>0$ will allow the two expressions to join. These choices give
\begin{equation}\label{base:next-local-choice}
 B_{j+1}\left(\bm x\right)=B_j\left(\bm c_{j+1}\right)+C_j
 +\ip{\D B_j\left(\bm c_{j+1}\right)}{\bm x-\bm c_{j+1}}
 +\frac{\mu}{q}\norm{\bm x-\bm c_{j+1}}^q.
\end{equation}
Lemma~\ref{base:recenter-lemma} and Corollaries~\ref{base:replacement-effects}--\ref{base:derivative-effects} verify the resulting changes in value and derivatives. We next construct a function that has this formula near $\bm c_{j+1}$ while leaving $B_j$ unchanged outside a small ball.

\subsubsection{Define the replacement and its two radii}
\label{base:local}\label{base:localregularity}
For each $j$, we extend the desired expression to a function on $\R^d$:
\[
 N_j^B\left(\bm x\right):=B_j\left(\bm c_{j+1}\right)+C_j
 +\ip{\D B_j\left(\bm c_{j+1}\right)}{\bm x-\bm c_{j+1}}
 +\frac{\mu}{q}\norm{\bm x-\bm c_{j+1}}^q.
\]
Using the smooth maximum $M_{\delta_j^B}$ from Eq.~\eqref{elem:psi-M-def}, define
\begin{equation}\label{base:reset}
 B_{j+1}\left(\bm x\right):=
 \begin{cases}
 M_{\delta_j^B}\left(B_j\left(\bm x\right),N_j^B\left(\bm x\right)\right),
       &\norm{\bm x-\bm c_{j+1}}<\rho_j,\\
 B_j\left(\bm x\right),&\norm{\bm x-\bm c_{j+1}}\ge\rho_j.
 \end{cases}
\end{equation}
The replacement must equal $N_j^B$ near $\bm c_{j+1}$ and
$B_j$ near the sphere $\norm{\bm x-\bm c_{j+1}}=\rho_j$,
so the inside and outside formulas join smoothly.
We therefore scale the support radius, value increase, and
smoothing width by
\begin{equation}\label{base:scaledconstants}
 \rho_j:=\theta R_j,\qquad
 C_j:=a_{\mathrm{lift}}\theta^2\mu R_j^q,\qquad
 \delta_j^B:=\frac{C_j}{4}
\end{equation}
The value comparisons derived in Lemma~\ref{base:reset-lemma}
in Appx.~\ref{app:base-support} determine the constants
needed for these two equalities. They give the following
explicit choices, all depending only on $q$:
\begin{equation}\label{base:quadraticconstants}
 a_0:=\frac{1}{2}\left(\frac{3}{4}\right)^{q-2},\qquad
 a_1:=\frac{q-1}{2}\left(\frac{5}{4}\right)^{q-2},
\end{equation}
\begin{equation}\label{base:liftwidth}
 a_{\mathrm{lift}}:=\frac{a_0}{8},\qquad d_*:=\frac{a_0}{32}=\frac{a_{\mathrm{lift}}}{4}.
\end{equation}
\begin{equation}\label{base:theta}
 \theta:=\min\left\{\frac{1}{16},\left(\frac{3qa_0}{32}\right)^{\frac{1}{q-2}}\right\}.
\end{equation}
\begin{equation}\label{base:gamma}
 \gamma:=\theta\sqrt{\frac{3a_0}{32a_1}}.
\end{equation}
Define
\[
\begin{aligned}
 S_j&:=B\left(\bm c_{j+1},\theta R_j\right)
     =\left\{\bm x:\norm{\bm x-\bm c_{j+1}}<\theta R_j\right\},\\
 U_j&:=B\left(\bm c_{j+1},\gamma R_j\right)
     =\left\{\bm x:\norm{\bm x-\bm c_{j+1}}<\gamma R_j\right\},\\
 \Delta_j^B&:=B_{j+1}-B_j.
\end{aligned}
\]
For one replacement, assume that $B_j$ is globally convex and $C^p$ and that Eq.~\eqref{base:old-local-form} holds on an open set containing $\overline S_j$. Lemma~\resultpart{base:nesting-lemma}{base:nesting-lemma:part:construction} below verifies these hypotheses at every step.

The support ball $S_j$ is the region where Eq.~\eqref{base:reset} uses $M_{\delta_j^B}$. The inner ball $U_j$ is the smaller region where $B_{j+1}=N_j^B$. Lemma~\resultpart{base:reset-lemma}{base:reset-lemma:part:exterior}--\ref{base:reset-lemma:part:interior} gives
\[
 B_{j+1}=N_j^B\quad\text{on }\overline U_j,
 \qquad
 B_{j+1}=B_j\quad\text{where }\norm{\bm x-\bm c_{j+1}}\ge\frac{\theta R_j}{2}.
\]
Thus the two definitions agree on an open region containing the joining sphere. Lemmas~\ref{base:branch-regularity-lemma} and \ref{base:regularity-lemma} bound the derivatives and H\"older constant of $\Delta_j^B$. The following main results verify that the replacements can be repeated and that the completed base has all the properties in Eq.~\eqref{base:requirements}--\eqref{base:past-answers}.

\subsubsection{Repeat the construction and obtain the center recurrences}
\label{base:nesting}\label{base:centers}
The next replacement must lie where the current base already has its
exact local formula. We therefore place its support ball inside $U_j$
and, to control later corrections for nearby pairs, require the smaller
containment below.

To shrink the accuracy scales geometrically, we choose
$R_j=R_0\eta^j$ and require
\[
 \overline S_{j+1}\subseteq\overline B\left(\bm c_{j+1},\frac{\gamma R_j}{2}\right)\subset U_j.
\]
A pair less than $\frac{\gamma R_j}{2}$ apart with one point outside $U_j$ then misses this smaller ball. For $\bm x\in\overline S_{j+1}$, the triangle inequality gives $\norm{\bm x-\bm c_{j+1}}\le\left(1+\theta\right)R_{j+1}=\left(1+\theta\right)\eta R_j$. Thus the required containment follows from
\begin{equation}\label{base:etaconstraint}
 0<\eta\le\frac{\gamma}{2\left(1+\theta\right)}.
\end{equation}
Lemma~\ref{base:reset-lemma} gives $\gamma<\frac{\theta}{2}$ and $\theta\le\frac{1}{16}$, so this already implies $\eta<\frac{1}{64}<\frac{1}{8}$, as needed for geometric tails and center spacing. No additional minimum is required.

\begin{lemma}[Nested replacements preserve earlier centers and give explicit recurrences]
\label{base:nesting-lemma}
Assume $R_j=R_0\eta^j$ and Eq.~\eqref{base:etaconstraint}, and start from $B_0\left(\bm x\right)=\frac{\mu\norm{\bm x}^q}{q}$. For the functions defined recursively by Eq.~\eqref{base:reset}, define
\[
 B:=B_L,\qquad \bm d_j:=\D B\left(\bm c_j\right),\qquad
 \omega_j:=B\left(\bm c_j\right),\qquad \beta_q:=\frac{1}{q}+a_{\mathrm{lift}}\theta^2,
\]
\begin{claims}
\item \label{base:nesting-lemma:part:construction}
\begin{enumerate}[label=\textnormal{(\roman*)},leftmargin=*]
\item \textbf{Local formula and containment of the next support.}
For $1\le j\le L$,
\[
\begin{aligned}
 B_j\left(\bm x\right)
 &=B_j\left(\bm c_j\right)+\ip{\D B_j\left(\bm c_j\right)}{\bm x-\bm c_j}
   +\frac{\mu}{q}\norm{\bm x-\bm c_j}^q,\\
 &\hspace{2em}\text{whenever }\norm{\bm x-\bm c_j}\le\gamma R_{j-1}.
\end{aligned}
\]
For $j=0$, the same function identity holds on all of $\R^d$. For $1\le j<L$,
\[
 \begin{aligned}
 \overline S_j
 &=\left\{\bm x:\norm{\bm x-\bm c_{j+1}}\le\theta R_j\right\}\\
 &\subseteq\overline B\left(\bm c_j,\frac{\gamma R_{j-1}}{2}\right)\\
 &\subset\left\{\bm x:\norm{\bm x-\bm c_j}<\gamma R_{j-1}\right\}=U_{j-1}.
 \end{aligned}
\]
Thus the displayed formula for $B_j$ holds on an open set containing $\overline S_j$, as assumed in Lemma~\ref{base:reset-lemma}.
\item \textbf{Global $C^p$ regularity and convexity.}
\[
 B_j\in C^p\left(\R^d\right),\qquad B_j\text{ is convex}
 \qquad\left(0\le j\le L\right).
\]
The H\"older bound is proved in Lemma~\ref{base:global-holder-lemma}.
\item \textbf{Preservation around every earlier center.}
For every $0\le i<L$, the functions agree on the open ball with center $\bm c_i$ and radius $\left(1-\theta\right)R_i>0$:
\[
 B_k\left(\bm x\right)=B_i\left(\bm x\right)
 \qquad\left(i\le k\le L,\ \norm{\bm x-\bm c_i}<\left(1-\theta\right)R_i\right).
\]
 \item \textbf{Preservation outside the replacement ball.}
 For $0\le i<L$ and $i\le k\le L$,
 \[
  B_k\left(\bm x\right)=B_i\left(\bm x\right),\qquad
  \D^rB_k\left(\bm x\right)=\D^rB_i\left(\bm x\right)
  \quad\left(\bm x\notin\overline S_i,\ 1\le r\le p\right).
 \]
 The partial base $B_j$ depends only on $\bm c_0,\ldots,\bm c_j$.
\end{enumerate}
\item \label{base:nesting-lemma:part:recurrences} The center recurrences are
\begin{equation}\label{base:centerrecurrences}
 \bm d_{j+1}-\bm d_j=\mu R_j^{q-2}\bm h_j,
 \qquad
 \omega_{j+1}-\omega_j=\ip{\bm d_j}{\bm h_j}+\beta_q\mu R_j^q.
\end{equation}

\item \label{base:nesting-lemma:part:center-formulas} The center values and derivatives are
\begin{equation}\label{base:centerformulas}
 \bm d_j=\mu\sum_{i<j}R_i^{q-2}\bm h_i,
 \qquad \omega_j=\beta_q\mu\sum_{i<j}R_i^q,
 \qquad \D^kB\left(\bm c_j\right)=0\ \left(2\le k\le p\right).
\end{equation}

\item \label{base:nesting-lemma:part:local-formula} The final function satisfies
\begin{equation}\label{base:local-invariant}
 B\left(\bm x\right)=B\left(\bm c_j\right)+\ip{\D B\left(\bm c_j\right)}{\bm x-\bm c_j}
                        +\frac{\mu}{q}\norm{\bm x-\bm c_j}^q.
\end{equation}
For $L\ge1$, this identity holds on the following explicit balls:
\[
 \begin{cases}
  \norm{\bm x-\bm c_j}<\left(1-\theta\right)R_j,&0\le j<L,\\
  \norm{\bm x-\bm c_L}\le\gamma R_{L-1},&j=L.
 \end{cases}
\]
For $L=0$, it holds on all of $\R^d$.
\end{claims}
\end{lemma}
\begin{proof}
\begin{itemize}
\item \textbf{Proof of \ref{base:nesting-lemma:part:construction}.}
\textbf{Proof of (i) and (ii).}
Lemma~\ref{base:reset-lemma} gives the two exact formulas for one
replacement. To repeat the construction, we place each subsequent
replacement inside the region where $B_{j+1}=N_j^B$.
This containment depends only on the centers, radii, and
Eq.~\eqref{base:etaconstraint}, so we can verify it before
constructing the functions. We then prove the local formulas
and global convexity and $C^p$ regularity together by induction.

\textbf{Geometric preparation for (i).}
To apply Lemma~\ref{base:reset-lemma} to a replacement $B_{j+1}\mapsto B_{j+2}$, the local formula for $B_{j+1}$ must hold on an open set containing the next closed support $\overline S_{j+1}$. The formula in (i), with $j+1$ in place of $j$, is required on $U_j$. Thus the geometric condition we need is $\overline S_{j+1}\subset U_j$.
\begin{enumerate}[label=\textbf{Step \arabic*.},leftmargin=*]
\item \textbf{Identify the two balls and their centers.}
Fix $0\le j<L-1$. By the definitions of $S_{j+1}$ and $U_j$,
\[
\begin{aligned}
 \overline S_{j+1}
 &=\left\{\bm x:\norm{\bm x-\bm c_{j+2}}\le\theta R_{j+1}\right\},
 &&\text{center }\bm c_{j+2},\ \text{radius }\theta R_{j+1},\\
 U_j
 &=\left\{\bm x:\norm{\bm x-\bm c_{j+1}}<\gamma R_j\right\},
 &&\text{center }\bm c_{j+1},\ \text{radius }\gamma R_j.
\end{aligned}
\]
Their centers are separated by
\[
 \norm{\bm c_{j+2}-\bm c_{j+1}}=R_{j+1}=\eta R_j.
\]
For an arbitrary $\bm x\in\overline S_{j+1}$, we know its distance from $\bm c_{j+2}$ and need an upper bound on its distance from $\bm c_{j+1}$:
\[
 \underbrace{\norm{\bm x-\bm c_{j+2}}\le\theta R_{j+1}}_{\text{known from membership in }\overline S_{j+1}}
 \quad\Longrightarrow\quad
 \underbrace{\norm{\bm x-\bm c_{j+1}}\le\frac{\gamma R_j}{2}}_{\text{bound to prove}}.
\]

\item \textbf{Bound the distance to $\bm c_{j+1}$ and obtain the containment.}
The vector identity
\[
 \bm x-\bm c_{j+1}
 =\left(\bm x-\bm c_{j+2}\right)+\left(\bm c_{j+2}-\bm c_{j+1}\right)
\]
and the triangle inequality give
\[
\begin{aligned}
 \norm{\bm x-\bm c_{j+1}}
 &\le
 \underbrace{\norm{\bm x-\bm c_{j+2}}}_{\le\theta R_{j+1}}
 +\underbrace{\norm{\bm c_{j+2}-\bm c_{j+1}}}_{=R_{j+1}}\\
 &\le\theta R_{j+1}+R_{j+1}\\
 &=\left(1+\theta\right)R_{j+1}\\
 &=\left(1+\theta\right)\eta R_j\\
 &\le\frac{\gamma R_j}{2},
\end{aligned}
\]
where the last line uses Eq.~\eqref{base:etaconstraint}. Since $\bm x$ was arbitrary, the whole ball centered at $\bm c_{j+2}$ lies in the displayed ball centered at $\bm c_{j+1}$. The remaining inclusions compare radii around the same center, using $\frac{\gamma R_j}{2}<\gamma R_j<\theta R_j$:
\begin{equation}\label{base:containment}
\begin{aligned}
 \overline S_{j+1}
 &=\overline B\left(\bm c_{j+2},\theta R_{j+1}\right)\\
 &\subseteq\overline B\left(\bm c_{j+1},\frac{\gamma R_j}{2}\right)\\
 &\subset B\left(\bm c_{j+1},\gamma R_j\right)=U_j\\
 &\subset B\left(\bm c_{j+1},\theta R_j\right)=S_j.
\end{aligned}
\end{equation}
Replacing $j+1$ by $j$ gives $\overline S_j\subset U_{j-1}$ for $1\le j<L$, as asserted in (i).

\end{enumerate}

\textbf{Joint induction for the local formula in (i) and the properties in (ii).}
\begin{enumerate}[label=\textbf{\arabic*.},leftmargin=*]
\item \textbf{Base case: $j=0$.}
By definition, $B_0\left(\bm x\right)=\frac{\mu\norm{\bm x}^q}{q}$ and $\bm c_0=0$. Lemma~\resultpart{elem:power}{elem:power:part:flatness} in Appx.~\ref{app:paper-elementary} gives $B_0\left(\bm c_0\right)=0$ and $\D B_0\left(\bm c_0\right)=0$. Hence
\[
 B_0\left(\bm x\right)=B_0\left(\bm c_0\right)
 +\ip{\D B_0\left(\bm c_0\right)}{\bm x-\bm c_0}
 +\frac{\mu}{q}\norm{\bm x-\bm c_0}^q
 \qquad\left(\bm x\in\R^d\right),
\]
which proves the local-formula claim (i) at $j=0$. Parts~\ref{elem:power:part:flatness} and \ref{elem:power:part:hessian} of the same lemma give
\[
 B_0\in C^p\left(\R^d\right),\qquad B_0\text{ is convex},
\]
proving (ii) at $j=0$.

\item \textbf{Inductive hypothesis: fix $0\le j<L$ and assume both claims for $B_j$.}
For (i), assume
\[
 B_j\left(\bm x\right)=B_j\left(\bm c_j\right)
 +\ip{\D B_j\left(\bm c_j\right)}{\bm x-\bm c_j}
 +\frac{\mu}{q}\norm{\bm x-\bm c_j}^q,
\]
on all of $\R^d$ if $j=0$, and on $\norm{\bm x-\bm c_j}\le\gamma R_{j-1}$ if $j\ge1$. For (ii), assume
\[
 B_j\in C^p\left(\R^d\right),\qquad B_j\text{ is convex}.
\]
These are assumptions at index $j$. The induction step must establish both conclusions at index $j+1$.

\item \textbf{Induction step: construct $B_{j+1}$ and prove both claims.}
\begin{enumerate}[label=\textnormal{(\alph*)},leftmargin=*]
\item \textbf{Verify the assumptions of Lemma~\ref{base:reset-lemma}.}
The inductive hypothesis (ii) supplies global convexity and $C^p$ regularity of $B_j$. For $j\ge1$, the inductive hypothesis (i) supplies its local formula on the open set
\[
 U_{j-1}=\left\{\bm x:\norm{\bm x-\bm c_j}<\gamma R_{j-1}\right\}.
\]
The geometric containment Eq.~\eqref{base:containment}, with $j-1$ in place of $j$, gives
\[
 \overline S_j=\overline B\left(\bm c_{j+1},\theta R_j\right)
 \subset U_{j-1}.
\]
Thus the local formula holds on an open set containing the entire closed support used for $B_j\mapsto B_{j+1}$. For $j=0$, the global formula gives the same requirement. All assumptions of Lemma~\ref{base:reset-lemma} are therefore verified using the hypotheses at index $j$ and the independently proved geometric containment.

\item \textbf{Prove (i) for $B_{j+1}$.}
Define $B_{j+1}$ by Eq.~\eqref{base:reset}. Lemma~\resultpart{base:reset-lemma}{base:reset-lemma:part:interior} gives $B_{j+1}=N_j^B$ on $\norm{\bm x-\bm c_{j+1}}\le\gamma R_j$. Part~\ref{base:reset-lemma:part:center-data} of that lemma identifies its center value and gradient, yielding
\[
\begin{aligned}
 B_{j+1}\left(\bm x\right)=N_j^B\left(\bm x\right)
 ={}&B_{j+1}\left(\bm c_{j+1}\right)
 +\ip{\D B_{j+1}\left(\bm c_{j+1}\right)}{\bm x-\bm c_{j+1}}\\
 &+\frac{\mu}{q}\norm{\bm x-\bm c_{j+1}}^q
 \qquad\left(\norm{\bm x-\bm c_{j+1}}\le\gamma R_j\right).
\end{aligned}
\]
This is exactly the local formula in (i) at index $j+1$. When $j<L-1$, Eq.~\eqref{base:containment} also gives $\overline S_{j+1}\subset U_j$, so this formula is available on an open set containing the following closed support.

\item \textbf{Prove (ii) for $B_{j+1}$.}
Apply Lemma~\resultpart{base:reset-lemma}{base:reset-lemma:part:smoothness} and part~\ref{base:reset-lemma:part:convexity}, whose assumptions were verified in (a). They give
\[
 B_{j+1}\in C^p\left(\R^d\right),\qquad B_{j+1}\text{ is convex},
\]
which proves (ii) at index $j+1$.
\end{enumerate}
\end{enumerate}
The base case establishes both claims at $j=0$, and the induction step establishes both at $j+1$ from both at $j$. They therefore hold for every $0\le j\le L$. Together with the geometric containment already proved, this completes (i) and (ii).

\textbf{Proof of (iii) and (iv): preservation of earlier formulas.}
We fix $0\le i<L$. Since the later replacement balls are
contained in $S_i$ by Eq.~\eqref{base:containment}, every later
modification vanishes on $\R^d\setminus\overline S_i$, an open set.
For $\bm x$ in this set and $i\le k\le L$, the definition
Eq.~\eqref{base:reset} therefore gives
\[
 \begin{aligned}
 B_k\left(\bm x\right)-B_i\left(\bm x\right)
 &=\sum_{j=i}^{k-1}
       \left(B_{j+1}\left(\bm x\right)-B_j\left(\bm x\right)\right)\\
 &=0.
 \end{aligned}
\]
Equality of the functions on this open set gives equality of
their derivatives, since $B_i,B_k\in C^p$ by (ii):
\[
 \D^rB_k\left(\bm x\right)=\D^rB_i\left(\bm x\right)
 \qquad\left(\bm x\notin\overline S_i,\ 1\le r\le p,\ i\le k\le L\right).
\]
The same agreement leaves a neighborhood of every earlier
center unchanged. Indeed, if
$\norm{\bm x-\bm c_i}<\left(1-\theta\right)R_i$, then
\[
 \begin{aligned}
 \norm{\bm x-\bm c_{i+1}}
 &\ge R_i-\norm{\bm x-\bm c_i}
 &&\text{by the reverse triangle inequality},\\
 &>\theta R_i,
 \end{aligned}
\]
so $\bm x\notin\overline S_i$. This proves (iii) and (iv),
including the preservation in
Lemma~\resultpart{paper:base}{paper:base:preservation}.
The recursion uses only its previous base and next center,
so $B_j$ depends only on $\bm c_0,\ldots,\bm c_j$.

\item \textbf{Proof of \ref{base:nesting-lemma:part:recurrences}.}
Each replacement retains its endpoint gradient and adds $C_j$
to its value. To obtain the center recurrences, we evaluate
the local formula at $\bm c_{j+1}$ and then account for this
value increase. The preservation just proved identifies the
partial and final center data: for $0\le j<L$,
part~\ref{base:nesting-lemma:part:construction}(iii), with
$i=j$ and $k=L$, gives
\[
 B_j\left(\bm x\right)=B_L\left(\bm x\right)=B\left(\bm x\right)
 \qquad\text{whenever }\norm{\bm x-\bm c_j}<\left(1-\theta\right)R_j.
\]
Thus the equality holds throughout the open ball $B\left(\bm c_j,\left(1-\theta\right)R_j\right)$, whose radius is positive. Evaluating and differentiating this equality at its center $\bm c_j$ yields
\[
 \begin{aligned}
 B_j\left(\bm c_j\right)&=B\left(\bm c_j\right)=\omega_j,\\
 \D B_j\left(\bm c_j\right)&=\D B\left(\bm c_j\right)=\bm d_j.
 \end{aligned}
\]
These identities compare the partial function $B_j$ and the completed function $B=B_L$ at the same center $\bm c_j$. Substituting them into the local formula in part~\ref{base:nesting-lemma:part:construction}(i) gives, on its stated region,
\[
\begin{aligned}
 B_j\left(\bm x\right)
 &=B_j\left(\bm c_j\right)+\ip{\D B_j\left(\bm c_j\right)}{\bm x-\bm c_j}
   +\frac{\mu}{q}\norm{\bm x-\bm c_j}^q\\
 &=\omega_j+\ip{\bm d_j}{\bm x-\bm c_j}
   +\frac{\mu}{q}\norm{\bm x-\bm c_j}^q.
\end{aligned}
\]
The point $\bm c_{j+1}=\bm c_j+\bm h_j$ belongs to the local-formula region:
for $j\ge1$, Eq.~\eqref{base:etaconstraint} gives
$R_j=\eta R_{j-1}<\gamma R_{j-1}$, and for $j=0$ the formula is global.
Evaluating and differentiating it there gives
\[
\begin{aligned}
 B_j\left(\bm c_{j+1}\right)
 &=\omega_j+\ip{\bm d_j}{\bm h_j}+\frac{\mu}{q}R_j^q,\\
 \D B_j\left(\bm c_{j+1}\right)
 &=\bm d_j+\mu\norm{\bm h_j}^{q-2}\bm h_j
 =\bm d_j+\mu R_j^{q-2}\bm h_j.
\end{aligned}
\]
Lemma~\resultpart{base:reset-lemma}{base:reset-lemma:part:center-data}
and the neighborhood preservation in
part~\ref{base:nesting-lemma:part:construction}(iii) now give the following
identities. At $j+1=L$, use $B=B_L$ in place of that preservation property.
\[
\begin{aligned}
 \bm d_{j+1}
 &=\D B_{j+1}\left(\bm c_{j+1}\right)
 =\D B_j\left(\bm c_{j+1}\right)
 =\bm d_j+\mu R_j^{q-2}\bm h_j,\\
 \omega_{j+1}
 &=B_{j+1}\left(\bm c_{j+1}\right)
 =B_j\left(\bm c_{j+1}\right)+C_j\\
 &=\omega_j+\ip{\bm d_j}{\bm h_j}
   +\frac{\mu}{q}R_j^q+a_{\mathrm{lift}}\theta^2\mu R_j^q\\
 &=\omega_j+\ip{\bm d_j}{\bm h_j}
   +\left(\frac{1}{q}+a_{\mathrm{lift}}\theta^2\right)\mu R_j^q\\
 &=\omega_j+\ip{\bm d_j}{\bm h_j}+\beta_q\mu R_j^q.
\end{aligned}
\]
These are exactly the two recurrences in Eq.~\eqref{base:centerrecurrences}.

\item \textbf{Proof of \ref{base:nesting-lemma:part:center-formulas}.}
The initial value and gradient remain zero, so
$\bm d_0=0$ and $\omega_0=0$.
Since successive gradient increments are parallel to $\bm h_i$,
orthogonality gives $\ip{\bm d_j}{\bm h_j}=0$ for $0\le j<L$.
We obtain the explicit center data by telescoping the recurrences
in Eq.~\eqref{base:centerrecurrences}:
\[
\begin{aligned}
 \bm d_j
 &=\bm d_0+\sum_{i=0}^{j-1}\left(\bm d_{i+1}-\bm d_i\right)
 =\mu\sum_{i<j}R_i^{q-2}\bm h_i,\\
 \ip{\bm d_j}{\bm h_j}
 &=\mu\sum_{i<j}R_i^{q-2}\ip{\bm h_i}{\bm h_j}=0,\\
 \omega_j
 &=\omega_0+\sum_{i<j}\left(\omega_{i+1}-\omega_i\right)\\
 &=\sum_{i<j}\left(\ip{\bm d_i}{\bm h_i}
                         +\beta_q\mu R_i^q\right)
 =\beta_q\mu\sum_{i<j}R_i^q.
\end{aligned}
\]
For $1\le j\le L$, Lemma~\resultpart{base:reset-lemma}{base:reset-lemma:part:center-data}
gives $\D^kB_j\left(\bm c_j\right)=0$ for $2\le k\le p$. If $j<L$,
part~\ref{base:nesting-lemma:part:construction}(iii) gives $B=B_j$
on a neighborhood of $\bm c_j$, so
\[
 \D^kB\left(\bm c_j\right)=\D^kB_j\left(\bm c_j\right)=0
 \qquad\left(2\le k\le p\right).
\]
For $j=L$, the same equality follows from $B=B_L$.
At $\bm c_0$, use Lemma~\resultpart{elem:power}{elem:power:part:flatness} for $B_0$ and the same neighborhood
preservation. This proves every identity in Eq.~\eqref{base:centerformulas}.
\item \textbf{Proof of \ref{base:nesting-lemma:part:local-formula}.}
\begin{enumerate}[label=\textbf{\arabic*.},leftmargin=*]
\item \textbf{Earlier centers: $0\le j<L$.}
For $1\le j<L$, Eq.~\eqref{base:etaconstraint} gives
\[
 \left(1-\theta\right)R_j=\left(1-\theta\right)\eta R_{j-1}<\gamma R_{j-1}.
\]
Thus the formula in part~\ref{base:nesting-lemma:part:construction}(i) holds on $\norm{\bm x-\bm c_j}<\left(1-\theta\right)R_j$; for $j=0$, that formula holds globally. On this ball, part~\ref{base:nesting-lemma:part:construction}(iii) gives $B=B_j$, hence their center values and gradients agree by part~\ref{base:nesting-lemma:part:construction}(ii). Substitution into part~\ref{base:nesting-lemma:part:construction}(i) yields
\[
\begin{aligned}
 &B\left(\bm x\right)-B\left(\bm c_j\right)-\ip{\D B\left(\bm c_j\right)}{\bm x-\bm c_j}\\
 &\quad=B_j\left(\bm x\right)-B_j\left(\bm c_j\right)-\ip{\D B_j\left(\bm c_j\right)}{\bm x-\bm c_j}\\
 &\quad=\frac{\mu}{q}\norm{\bm x-\bm c_j}^q,
\end{aligned}
\]
which is Eq.~\eqref{base:local-invariant}.
\item \textbf{Last center: $j=L$.}
Since $B=B_L$, part~\ref{base:nesting-lemma:part:construction}(i) gives Eq.~\eqref{base:local-invariant} directly on $\norm{\bm x-\bm c_L}\le\gamma R_{L-1}$ if $L\ge1$, and on all of $\R^d$ if $L=0$.
\end{enumerate}
\end{itemize}
\end{proof}
By Lemma~\resultpart{base:nesting-lemma}{base:nesting-lemma:part:recurrences}, the original base increment and its magnitude are
\begin{equation}\label{base:gradient-drift}
 \bm b_j:=\D B\left(\bm c_{j+1}\right)-\D B\left(\bm c_j\right)
 =\mu R_j^{q-2}\bm h_j,
 \qquad G_j:=\norm{\bm b_j}=\mu R_j^{q-1}.
\end{equation}
Before the replacement, Eq.~\eqref{base:old-local-form} holds along the full segment from $\bm c_j$ to $\bm c_{j+1}$. The final neighborhood formula Eq.~\eqref{base:local-invariant} need not hold at $\bm c_{j+1}$. Thus across consecutive centers the final tangent remainder is $\beta_q\mu R_j^q=\frac{\mu R_j^q}{q}+C_j$, whereas within the final neighborhood of a center it is given by Eq.~\eqref{base:local-invariant}.

\Needspace{9\baselineskip}
\subsubsection{H\"older smoothness: compare a pair at the first relevant scale}
\label{base:holder}
For smoothness, each correction $B_{j+1}-B_j$ has
$p^{th}$-derivative magnitude at most a fixed multiple of
$\mu R_j^\nu$. Their geometric decrease bounds the sum beyond
any scale by its first term times a constant, avoiding a factor
of $L$. To turn this magnitude bound into the third line of
Eq.~\eqref{base:requirements}, we combine the correction sum
for distant pairs with the exact local norm-power formulas for
close pairs. The two local estimates from
Lemma~\ref{base:regularity-lemma}, with the second holding for all
$\bm x,\bm y\in\R^d$, are
\[
 \begin{aligned}
 \norm{\D^p\Delta_j^B}_\infty
 &\le M_p^\Delta\mu R_j^\nu,\\
 \norm{\D^p\Delta_j^B\left(\bm x\right)-\D^p\Delta_j^B\left(\bm y\right)}
 &\le H_\Delta\mu\norm{\bm x-\bm y}^\nu.
 \end{aligned}
\]
Their proofs, including the estimates for $N_j^B-B_j$ in Lemma~\ref{base:branch-regularity-lemma}, are grouped in Appx.~\ref{app:base-support}. Combined with the exact formulas and nesting in Lemma~\ref{base:nesting-lemma}, they give the following bound independent of the number of replacements.

\begin{lemma}[H\"older smoothness independent of the number of replacements]
\label{base:global-holder-lemma}
\leavevmode
\begin{claims}
\item \label{base:global-holder-lemma:part:holder} For all $\bm x,\bm y\in\R^d$, the final base satisfies
\[
 \norm{\D^pB\left(\bm x\right)-\D^pB\left(\bm y\right)}
 \le C_B\mu\norm{\bm x-\bm y}^\nu.
\]
The finite constant in Eq.~\eqref{base:holderconstant} is independent of dimension, $L$, $R_0$, and every allowed smaller value of $\eta$. For $0\le j\le L$, every partial base belongs to $C^{p,\nu}\left(\R^d\right)$ and satisfies, for all $\bm x,\bm y\in\R^d$,
\[
 \norm{\D^pB_j\left(\bm x\right)-\D^pB_j\left(\bm y\right)}
 \le C_B\mu\norm{\bm x-\bm y}^\nu.
\]
\end{claims}
\end{lemma}
\begin{proof}
\begin{itemize}
\item \textbf{Proof of \ref{base:global-holder-lemma:part:holder}.}
Fix the completed function $B=B_L$ and the constant $C_B$ in Eq.~\eqref{base:holderconstant}. We use backward induction on $j=L,L-1,\ldots,0$, since exact agreement in the next inner ball lets us apply the induction hypothesis to the same pair without an additional H\"older estimate.

\textbf{Induction claim $\mathcal P_j$ ($0\le j\le L$).}
For every pair of distinct $\bm x,\bm y\in\R^d$, put $d_0=\norm{\bm x-\bm y}$. We prove the following implication.

\emph{Assumption:} the local formula Eq.~\eqref{base:old-local-form} holds on a neighborhood of each point:
\[
\begin{aligned}
 B_j\left(\bm w\right)={}&B_j\left(\bm c_j\right)
 +\ip{\D B_j\left(\bm c_j\right)}{\bm w-\bm c_j}\\
 &+\frac{\mu}{q}\norm{\bm w-\bm c_j}^q\qquad
 \text{near both }\bm x,\bm y.
\end{aligned}
\]
\emph{Conclusion:} the final function $B$ satisfies
\[
 \norm{\D^pB\left(\bm x\right)-\D^pB\left(\bm y\right)}_{\mathrm{op}}\le C_B\mu d_0^\nu.
\]
The symbol $\mathcal P_j$ denotes this entire implication. The final function $B$ and constant $C_B$ are the same at every index.

\textbf{Bound the size of all remaining corrections.}
Lemma~\resultpart{base:regularity-lemma}{base:regularity-lemma:part:derivative-bounds} gives $\norm{\D^p\Delta_k^B}_\infty\le M_p^\Delta\mu R_k^\nu$.
For $0\le j<L$, let $k=j+a$ and use $R_k^\nu=R_j^\nu\eta^{a\nu}$:
\[
\begin{aligned}
 \left\|\sum_{k=j}^{L-1}\D^p\Delta_k^B\right\|_\infty
 &\le\sum_{k=j}^{L-1}\norm{\D^p\Delta_k^B}_\infty\\
 &\le M_p^\Delta\mu R_j^\nu
               \sum_{a=0}^{L-1-j}\eta^{a\nu}\\
 &\le M_p^\Delta\mu R_j^\nu
               \sum_{a=0}^{\infty}\eta^{a\nu}
 =\frac{M_p^\Delta\mu R_j^\nu}{1-\eta^\nu}.
\end{aligned}
\]
The first line is the triangle inequality; the second uses the displayed
derivative bound, and the last line uses $0<\eta^\nu<1$. Thus
\begin{equation}\label{base:tailbound}
 \left\|\sum_{k=j}^{L-1}\D^p\Delta_k^B\right\|_\infty
 \le\frac{M_p^\Delta\mu R_j^\nu}{1-\eta^\nu}.
\end{equation}
\textbf{An estimate from the local formula.}
Fix $0\le j\le L$ and suppose Eq.~\eqref{base:old-local-form} holds near both $\bm x$ and $\bm y$. Differentiating that formula at the two points, with $p\ge2$, removes its affine terms:
\begin{equation}\label{base:pairlocalidentity}
 \D^pB_j\left(\bm z\right)
 =\left.\mu\D_{\bm w}^p\!\left(\frac{\norm{\bm w-\bm c_j}^q}{q}\right)\right|_{\bm w=\bm z}
 \qquad\left(\bm z\in\left\{\bm x,\bm y\right\}\right).
\end{equation}
Lemma~\resultpart{elem:power}{elem:power:part:holder} in Appx.~\ref{app:paper-elementary} therefore gives
\begin{equation}\label{base:pairlocalbound}
 \norm{\D^pB_j\left(\bm x\right)-\D^pB_j\left(\bm y\right)}_{\mathrm{op}}
 \le H_{\rm pow}\mu d_0^\nu.
\end{equation}
\textbf{Base case: $j=L$.}
Take arbitrary distinct $\bm x,\bm y$ and assume Eq.~\eqref{base:old-local-form}, with $j=L$, holds near both points. Put $f\left(\bm u\right)=\frac{\norm{\bm u}^q}{q}$. Eq.~\eqref{base:pairlocalidentity}, obtained by differentiating this local formula, reads
\[
 \D^pB_L\left(\bm z\right)=\mu\D^pf\left(\bm z-\bm c_L\right)
 \qquad\left(\bm z\in\left\{\bm x,\bm y\right\}\right).
\]
Since $B=B_L$, Lemma~\resultpart{elem:power}{elem:power:part:holder} gives
\begin{align*}
 &\norm{\D^pB\left(\bm x\right)-\D^pB\left(\bm y\right)}_{\mathrm{op}}\\
 &\quad=\mu\norm{\D^pf\left(\bm x-\bm c_L\right)-\D^pf\left(\bm y-\bm c_L\right)}_{\mathrm{op}}\\
 &\quad\le H_{\rm pow}\mu\norm{\left(\bm x-\bm c_L\right)-\left(\bm y-\bm c_L\right)}^\nu\\
 &\quad=H_{\rm pow}\mu d_0^\nu\le C_B\mu d_0^\nu.
\end{align*}
We have derived the claimed bound from the assumed local formula for every such pair. This proves the implication $\mathcal P_L$.

\textbf{Inductive hypothesis.}
Fix $0\le j<L$ and assume the implication $\mathcal P_{j+1}$ has been proved for every distinct pair: if Eq.~\eqref{base:old-local-form}, with $j+1$ in place of $j$, holds near both points, then the displayed bound for $B$ follows.
The hypothesis is applied only in \hyperlink{base:holder-ih-use}{\textbf{Case~2 below}}; Cases~1 and~3 give direct estimates.

\textbf{Induction step: prove $\mathcal P_j$.}
Take arbitrary distinct $\bm x,\bm y$ and assume Eq.~\eqref{base:old-local-form} for $B_j$ holds near both points. Put $d_0=\norm{\bm x-\bm y}$. Eq.~\eqref{base:pairlocalbound} follows from this local formula. It remains to control the corrections in
\[
 B=B_j+\underbrace{\Delta_j^B}_{\text{current correction}}
       +\underbrace{\sum_{k=j+1}^{L-1}\Delta_k^B}_{\text{later corrections}}.
\]

\textbf{Why split at $d_0=\frac{\gamma R_j}{2}$?}
Lemma~\resultpart{base:nesting-lemma}{base:nesting-lemma:part:construction}(i) gives the exact formula for $B_{j+1}$ on $U_j=B\left(\bm c_{j+1},\gamma R_j\right)$. Its containment Eq.~\eqref{base:containment} confines all later corrections to the concentric closed ball of radius $\frac{\gamma R_j}{2}$. The gap between these radii is $\frac{\gamma R_j}{2}$. A pair closer than this gap either lies entirely in $U_j$ (Case~2), or has a point outside $U_j$, forcing both points outside the smaller ball (Case~3). For $d_0\ge\frac{\gamma R_j}{2}$, the size bound Eq.~\eqref{base:tailbound} suffices (Case~1).

\begin{enumerate}[label=\textbf{Case \arabic*.},leftmargin=*]
\item \textbf{Pairs at distance at least $\frac{\gamma R_j}{2}$.}
If $d_0\ge\frac{\gamma R_j}{2}$, then
\[
 R_j^\nu\le\left(\frac{2}{\gamma}\right)^\nu d_0^\nu.
\]
Use the exact decomposition $B=B_j+\sum_{k=j}^{L-1}\Delta_k^B$. Differentiating and subtracting at the two points gives
\[
\begin{aligned}
 \D^pB\left(\bm x\right)-\D^pB\left(\bm y\right)
 ={}&\D^pB_j\left(\bm x\right)-\D^pB_j\left(\bm y\right)\\
 &+\sum_{k=j}^{L-1}\D^p\Delta_k^B\left(\bm x\right)
  -\sum_{k=j}^{L-1}\D^p\Delta_k^B\left(\bm y\right).
\end{aligned}
\]
Apply the triangle inequality, then Eq.~\eqref{base:pairlocalbound} to the $B_j$ difference and Eq.~\eqref{base:tailbound} to each correction sum:
\begin{align*}
 &\norm{\D^pB\left(\bm x\right)-\D^pB\left(\bm y\right)}_{\mathrm{op}}\\
 &\quad\le\norm{\D^pB_j\left(\bm x\right)-\D^pB_j\left(\bm y\right)}_{\mathrm{op}}\\
 &\qquad+\left\|\sum_{k=j}^{L-1}\D^p\Delta_k^B\left(\bm x\right)\right\|_{\mathrm{op}}
          +\left\|\sum_{k=j}^{L-1}\D^p\Delta_k^B\left(\bm y\right)\right\|_{\mathrm{op}}\\
 &\quad\le\underbrace{H_{\rm pow}\mu d_0^\nu}_{Eq.~\eqref{base:pairlocalbound}}
       +\underbrace{\frac{M_p^\Delta\mu R_j^\nu}{1-\eta^\nu}}_{Eq.~\eqref{base:tailbound}\text{ at }\bm x}
       +\underbrace{\frac{M_p^\Delta\mu R_j^\nu}{1-\eta^\nu}}_{Eq.~\eqref{base:tailbound}\text{ at }\bm y}\\
 &\quad=H_{\rm pow}\mu d_0^\nu
       +\frac{2M_p^\Delta\mu R_j^\nu}{1-\eta^\nu}\\
 &\quad\le\left(H_{\rm pow}
       +\frac{2M_p^\Delta\left(\frac{2}{\gamma}\right)^\nu}{1-8^{-\nu}}\right)\mu d_0^\nu
 \le C_B\mu d_0^\nu.
\end{align*}
The two applications of Eq.~\eqref{base:tailbound} give the factor $2$; the last line uses $R_j^\nu\le\left(\frac{2}{\gamma}\right)^\nu d_0^\nu$ and $\eta\le\frac{1}{8}$. This proves $\mathcal P_j$ in Case~1. For $d_0<\frac{\gamma R_j}{2}$, the ratio $\frac{R_j^\nu}{d_0^\nu}$ can be arbitrarily large, so we use the locations of the points in the next two cases.

\item \textbf{Both points lie in $U_j$: use the exact formula for $B_{j+1}$.}
Assume $d_0<\frac{\gamma R_j}{2}$ and $\bm x,\bm y\in U_j$. Lemma~\resultpart{base:nesting-lemma}{base:nesting-lemma:part:construction}(i), applied at index $j+1$, gives
\[
 \begin{aligned}
 B_{j+1}\left(\bm w\right)
 ={}&B_{j+1}\left(\bm c_{j+1}\right)
 +\ip{\D B_{j+1}\left(\bm c_{j+1}\right)}{\bm w-\bm c_{j+1}}\\
 &+\frac{\mu}{q}\norm{\bm w-\bm c_{j+1}}^q
 \qquad\left(\norm{\bm w-\bm c_{j+1}}\le\gamma R_j\right).
 \end{aligned}
\]
For each $\bm z\in\left\{\bm x,\bm y\right\}$,
\[
 r_{\bm z}:=\gamma R_j-\norm{\bm z-\bm c_{j+1}}>0,
 \qquad B\left(\bm z,r_{\bm z}\right)\subset U_j.
\]
Thus the displayed formula can be differentiated at both points. Recall $f\left(\bm u\right)=\frac{\norm{\bm u}^q}{q}$. Since $p\ge2$, the affine terms have zero $p^{th}$ derivative, so
\[
 \D^pB_{j+1}\left(\bm z\right)=\mu\D^pf\left(\bm z-\bm c_{j+1}\right)
 \qquad\left(\bm z\in\left\{\bm x,\bm y\right\}\right).
\]
This is Eq.~\eqref{base:pairlocalidentity} at index $j+1$. Lemma~\resultpart{elem:power}{elem:power:part:holder} now gives
\begin{align*}
 &\norm{\D^pB_{j+1}\left(\bm x\right)-\D^pB_{j+1}\left(\bm y\right)}_{\mathrm{op}}\\
 &\quad=\mu\norm{\D^pf\left(\bm x-\bm c_{j+1}\right)-\D^pf\left(\bm y-\bm c_{j+1}\right)}_{\mathrm{op}}\\
 &\quad\le\mu H_{\rm pow}\norm{\left(\bm x-\bm c_{j+1}\right)-\left(\bm y-\bm c_{j+1}\right)}^{\nu}\\
 &\quad=H_{\rm pow}\mu\norm{\bm x-\bm y}^{\nu}
       =H_{\rm pow}\mu d_0^\nu.
\end{align*}
This proves Eq.~\eqref{base:pairlocalbound} at index $j+1$ with the same coefficient; no bound for $\Delta_j^B$ has been added.

\Needspace{7\baselineskip}
\hypertarget{base:holder-ih-use}{}\textbf{Apply the induction hypothesis $\mathcal P_{j+1}$.}
We have verified Eq.~\eqref{base:old-local-form} for $B_{j+1}$ near the same $\bm x,\bm y$. Applying $\mathcal P_{j+1}$ to this pair gives the bound for the final function $B$:
\[
 \norm{\D^pB\left(\bm x\right)-\D^pB\left(\bm y\right)}_{\mathrm{op}}
 \le C_B\mu d_0^\nu
 \qquad\text{by }\mathcal P_{j+1}.
\]
This proves $\mathcal P_j$ in Case~2, including $j=L-1$, where the hypothesis is the base case $\mathcal P_L$. No correction estimate is added to the induction bound.

\item \textbf{Remaining close pairs: at least one point lies outside $U_j$.}
Assume $d_0<\frac{\gamma R_j}{2}$ and, exchanging the points if needed, $\bm x\notin U_j$. The key identity is
\[
 B=B_{j+1},\qquad B-B_j=\Delta_j^B
 \quad\text{near both }\bm x\text{ and }\bm y.
\]
For $j=L-1$, it holds globally. For $j\le L-2$, apply Lemma~\resultpart{design:allowed-changes}{design:allowed-changes:part:neighborhood} in Appx.~\ref{app:paper-oracle} at starting index $j+1$. The required distance bounds are
\[
 \norm{\bm x-\bm c_{j+1}}\ge\gamma R_j,\qquad
 \norm{\bm y-\bm c_{j+1}}\ge\gamma R_j-d_0>\frac{\gamma R_j}{2}.
\]
Lemma~\resultpart{base:nesting-lemma}{base:nesting-lemma:part:construction}(i), through Eq.~\eqref{base:containment}, gives
\[
 \overline S_k\subseteq\overline S_{j+1}
 \subseteq\overline B\left(\bm c_{j+1},\frac{\gamma R_j}{2}\right)
 \qquad\left(j+1\le k\le L-1\right).
\]
Thus both points have positive distance from $\overline S_{j+1}$. Every correction $\Delta_k^B$ with $k\ge j+1$ vanishes outside $S_k$ by Eq.~\eqref{base:reset}; part~\ref{base:nesting-lemma:part:construction}(ii) supplies $C^p$ regularity. These verify the hypotheses of Lemma~\ref{design:allowed-changes} and prove the neighborhood identity.

Differentiating that identity gives
\[
 \D^p\left(B-B_j\right)\left(\bm z\right)=\D^p\Delta_j^B\left(\bm z\right)
 \qquad\left(\bm z\in\left\{\bm x,\bm y\right\}\right).
\]
Consequently, Eq.~\eqref{base:pairlocalbound} and Lemma~\resultpart{base:regularity-lemma}{base:regularity-lemma:part:holder} yield
\[
\begin{aligned}
 &\norm{\D^pB\left(\bm x\right)-\D^pB\left(\bm y\right)}_{\mathrm{op}}\\
 &\quad=\left\|\D^pB_j\left(\bm x\right)-\D^pB_j\left(\bm y\right)
       +\D^p\Delta_j^B\left(\bm x\right)-\D^p\Delta_j^B\left(\bm y\right)\right\|_{\mathrm{op}}\\
 &\quad\le\norm{\D^pB_j\left(\bm x\right)-\D^pB_j\left(\bm y\right)}_{\mathrm{op}}
       +\norm{\D^p\Delta_j^B\left(\bm x\right)-\D^p\Delta_j^B\left(\bm y\right)}_{\mathrm{op}}\\
 &\quad\le H_{\rm pow}\mu d_0^\nu+H_\Delta\mu d_0^\nu
 =\left(H_{\rm pow}+H_\Delta\right)\mu d_0^\nu\\
 &\quad\le C_B\mu d_0^\nu
 \qquad\text{by Eq.~\eqref{base:holderconstant}}.
\end{aligned}
\]
This proves $\mathcal P_j$ in Case~3.
\end{enumerate}
The three cases prove $\mathcal P_{j+1}\Longrightarrow\mathcal P_j$. The base case and all three estimates use the single constant
\begin{equation}\label{base:holderconstant}
 C_B:=H_{\rm pow}+\max\left\{H_\Delta,
   \frac{2M_p^\Delta\left(\frac{2}{\gamma}\right)^\nu}{1-8^{-\nu}}\right\}.
\end{equation}
These constants are finite and dimension independent by Lemmas~\resultpart{elem:power}{elem:power:part:holder} and \resultpart{base:regularity-lemma}{base:regularity-lemma:part:derivative-bounds}, \ref{base:regularity-lemma:part:holder}.

\textbf{Conclusion of the induction.}
Backward induction proves the implication $\mathcal P_j$ for every $0\le j\le L$. Since $B_0\left(\bm w\right)=\frac{\mu\norm{\bm w}^q}{q}$ globally, Eq.~\eqref{base:old-local-form} at $j=0$ holds near every pair of points. Applying the proved implication $\mathcal P_0$ to every distinct pair gives
\begin{equation}\label{base:holderbound}
 \norm{\D^pB\left(\bm x\right)-\D^pB\left(\bm y\right)}
 \le C_B\mu\norm{\bm x-\bm y}^\nu
 \qquad\left(\bm x,\bm y\in\R^d\right).
\end{equation}
Stopping the construction after any $j\le L$ gives the same bound for $B_j$. Together with Lemma~\resultpart{base:nesting-lemma}{base:nesting-lemma:part:construction}(ii), this proves $B_j\in C^{p,\nu}$.
\end{itemize}
\end{proof}

\subsubsection{Preserve every base oracle output on a separated subspace}
\label{base:prefix}
To verify Eq.~\eqref{base:past-answers}, we separate evaluation points $\bm x\in W$ from all future corrections $\Delta_k^B$, $k\ge j$.

\begin{lemma}[Future base changes are invisible on a separated subspace]
\label{base:prefix-lemma}
Fix $0\le j<L$, and let $W$ be a linear subspace. If
\begin{equation}\label{base:separationcondition}
 \operatorname{dist}\left(\bm c_{j+1},W\right)>\theta R_j,
\end{equation}
then:
\begin{claims}
\item \label{base:prefix-lemma:part:neighborhood}
\begin{equation}\label{base:prefixagreement}
 B=B_j\quad\text{on a neighborhood of every }\bm x\in W.
\end{equation}
\item \label{base:prefix-lemma:part:oracle} The values and all derivative tensors through order $p$ agree at every $\bm x\in W$.
\end{claims}
\end{lemma}
\begin{proof}
\begin{itemize}
\item \textbf{Proof of \ref{base:prefix-lemma:part:neighborhood}.} 
For every $\bm x\in W$, the separation assumption Eq.~\eqref{base:separationcondition} and the triangle inequality give
\[
 \operatorname{dist}\left(\bm x,\overline S_j\right)
 \ge\operatorname{dist}\left(\bm c_{j+1},W\right)-\theta R_j>0.
\]
Lemma~\resultpart{base:nesting-lemma}{base:nesting-lemma:part:construction} gives $C^p$ regularity of all partial bases and places every later support $S_k$ ($k\ge j$) inside $S_j=B\left(\bm c_{j+1},\theta R_j\right)$. Lemma~\resultpart{base:reset-lemma}{base:reset-lemma:part:exterior} gives agreement $B_{k+1}=B_k$ outside each support $S_k=B\left(\bm c_{k+1},\theta R_k\right)$. These verify the hypotheses of Lemma~\resultpart{design:allowed-changes}{design:allowed-changes:part:neighborhood} in Appx.~\ref{app:paper-oracle}, which proves Eq.~\eqref{base:prefixagreement}.
\item \textbf{Proof of \ref{base:prefix-lemma:part:oracle}.}
Fix $\bm x\in W$. Part~\ref{base:prefix-lemma:part:neighborhood} gives $B=B_j$ on an open set containing $\bm x$, and Lemma~\resultpart{base:nesting-lemma}{base:nesting-lemma:part:construction}(ii) gives $B,B_j\in C^p$. Applying Lemma~\resultpart{design:allowed-changes}{design:allowed-changes:part:oracle} to this equality yields
\[
 B\left(\bm x\right)=B_j\left(\bm x\right),\qquad
 \D^kB\left(\bm x\right)=\D^kB_j\left(\bm x\right)\quad\left(1\le k\le p\right).
\]
\end{itemize}
\end{proof}
The formula for $B_j$ uses only previously chosen centers, their values and gradients, and Euclidean distances. It is therefore unchanged by rotations of the remaining, unchosen directions. Some derivative tensors act on the orthogonal complement of the known span, but that action identifies no particular direction there. Lemma~\resultpart{base:prefix-lemma}{base:prefix-lemma:part:oracle} proves the base part of oracle consistency. The phase functions $A_j$ require their separate calculation on the chain-prefix subspaces $W_{j,b}$ from Step~3. Appx.~\ref{collar:base-separation} verifies Eq.~\eqref{base:separationcondition} for an incomplete phase using the endpoint defined in Appx.~\ref{collar:chain}.

\subsubsection{Prove uniform convexity for every pair of points}
\label{base:uniformconvexity}
Let $\mathcal C:=\left\{\bm c_0,\ldots,\bm c_L\right\}$. The distance from a point $\bm x$ to this finite set is
\[
 \operatorname{dist}\left(\bm x,\mathcal C\right)
 :=\min_{0\le i\le L}\norm{\bm x-\bm c_i},
\]
the Euclidean distance to the nearest center. To prove the second line of Eq.~\eqref{base:requirements}, we use the standard implication from uniform monotonicity of the gradient to uniform convexity. Lemmas~\ref{base:hessian-lemma} and \ref{base:porosity} in Appx.~\ref{app:base-curvature} verify its hypothesis for this base.

\begin{lemma}[Uniform convexity of the completed base]
\label{base:uc-lemma}
For the constant $c_B$ computed below:
\begin{claims}
\item \label{base:uc-lemma:part:remainder}
\begin{equation}\label{base:bregmanbound}
 \Breg B{\bm y}{\bm x}\ge c_B\mu\norm{\bm y-\bm x}^q
 \qquad\left(\bm x,\bm y\in\R^d\right).
\end{equation}
\item \label{base:uc-lemma:part:modulus} The uniform-convexity modulus, under the convention with factor $\frac{1}{q}$, is at least $q c_B\mu$.
\end{claims}
\end{lemma}
\begin{proof}
\begin{itemize}
\item \textbf{Proof of \ref{base:uc-lemma:part:remainder}.}
Uniform convexity requires sufficient curvature along every
segment from $\bm x$ to $\bm y$.
Here $B\in C^2$ by
Lemma~\resultpart{base:nesting-lemma}{base:nesting-lemma:part:construction}(ii),
since $p\ge2$.
For distinct $\bm x,\bm y\in\R^d$, define
$\bm d=\bm y-\bm x$ and $\ell=\norm{\bm d}$.
To identify how this curvature gives growth above the tangent,
we view $B$ along the segment as the curve
$g\left(s\right):=B\left(\bm x+s\bm u\right)$ for
$\bm u=\frac{\bm d}{\ell}$ and $0\le s\le\ell$.
The tangent at its starting point has height
$g\left(0\right)+sg'\left(0\right)$, and its gap from the curve
at the endpoint is exactly
\[
 g\left(\ell\right)-g\left(0\right)-\ell g'\left(0\right)
 =\Breg B{\bm y}{\bm x}.
\]
Thus we need a lower bound of order $\mu\ell^q$ on this gap,
uniformly over the starting point and direction.

Lemma~\ref{base:hessian-lemma} bounds directional curvature by
$\D^2B\left(\bm z\right)[\bm d,\bm d]
\ge\mu\ell^2\operatorname{dist}\left(\bm z,\mathcal C\right)^{q-2}$.
Equivalently, along the curve,
$g''\left(s\right)\ge
\mu\operatorname{dist}\left(\bm x+s\bm u,\mathcal C\right)^{q-2}$.
This can vanish at the centers, but Lemma~\ref{base:porosity}
guarantees that at least half the segment stays at distance
$\frac{\ell}{16}$ from every center, i.e., there is a measurable
$\Omega\subset[0,1]$ with $|\Omega|\ge\frac{1}{2}$ and
$\operatorname{dist}\left(\bm x+t\bm d,\mathcal C\right)
\ge\frac{\ell}{16}$ on $\Omega$.
Elsewhere the curvature is nonnegative, so the accumulated
increase in slope cannot be undone.
Integrating the curvature over this part gives
\begin{align*}
 \ip{\D B\left(\bm y\right)-\D B\left(\bm x\right)}{\bm d}
 &=\int_0^1\D^2B\left(\bm x+t\bm d\right)[\bm d,\bm d]\,dt
 &&\text{by Theorem~\ref{sources:taylor}},\\
 &\ge\mu\ell^2\int_0^1
       \operatorname{dist}\left(\bm x+t\bm d,\mathcal C\right)^{q-2}\,dt
 &&\text{by Lemma~\ref{base:hessian-lemma}},\\
 &\ge\mu\ell^2\int_\Omega
       \operatorname{dist}\left(\bm x+t\bm d,\mathcal C\right)^{q-2}\,dt
 &&\text{by nonnegativity},\\
 &\ge\mu\ell^2|\Omega|\left(\frac{\ell}{16}\right)^{q-2}
 &&\text{by Lemma~\ref{base:porosity}},\\
 &\ge\frac{\mu}{2\,16^{q-2}}\ell^q
 &&\text{since }|\Omega|\ge\frac{1}{2}.
\end{align*}
In terms of the curve, this gradient increase is
$\ell\left(g'\left(\ell\right)-g'\left(0\right)\right)$.
Dividing by $\ell$ gives the accumulated excess slope
\[
 g'\left(\ell\right)-g'\left(0\right)
 \ge\underbrace{\frac{\ell}{2}}_{\text{good length}}
     \underbrace{\mu\left(\frac{\ell}{16}\right)^{q-2}}_{\text{curvature there}}
 =\frac{\mu\ell^{q-1}}{2\cdot16^{q-2}}.
\]
The same estimate applies to every shorter subsegment
$[\bm x,\bm x+s\bm u]$, giving excess slope at least
$\frac{\mu s^{q-1}}{2\cdot16^{q-2}}$.
The height above the initial tangent accumulates this excess
slope. Formally, the uniform-monotonicity criterion
\cite[Lemma~2.1]{DoikovNesterov2021}, restated as
Lemma~\ref{sources:uniform-monotonicity}, gives
\[
 \Breg B{\bm y}{\bm x}\ge
 \frac{\mu}{2q\,16^{q-2}}\norm{\bm y-\bm x}^q.
\]
Thus the Hessian bound supplies curvature, and the geometric
lemma ensures that enough of it is encountered between every
pair of points. Eq.~\eqref{base:bregmanbound} holds with
\begin{equation}\label{base:ucconstant}
 c_B:=\frac{1}{2q\cdot16^{q-2}}.
\end{equation}
The case $\bm x=\bm y$ is immediate.
\item \textbf{Proof of \ref{base:uc-lemma:part:modulus}.}
Part~\ref{base:uc-lemma:part:remainder} gives
\[
 \Breg B{\bm y}{\bm x}
 \ge c_B\mu\norm{\bm y-\bm x}^q
 =\frac{q c_B\mu}{q}\norm{\bm y-\bm x}^q.
\]
Comparison with the uniform-convexity condition in
Theorem~\ref{paper:main-theorem}
therefore gives modulus at least $qc_B\mu$.
\end{itemize}
\end{proof}

\begin{corollary}[All properties of the uniformly convex base]
\label{base:complete-base-corollary}
Under the hypotheses and constant choices of
Lemma~\ref{paper:base}, the construction in
Eq.~\eqref{paper:base-reset} satisfies all conclusions of that
lemma. We may take $C_B$ from Eq.~\eqref{base:holderconstant}
and $c_B$ from Eq.~\eqref{base:ucconstant}.
\end{corollary}
\begin{proof}
\textbf{\ref{paper:base:class}}
Lemma~\resultpart{base:nesting-lemma}{base:nesting-lemma:part:construction}(ii)
gives convexity and $C^p$ regularity of every $B_j$.
Lemma~\ref{base:global-holder-lemma} gives their common global
H\"older bound, and Lemma~\ref{base:uc-lemma} gives the
Euclidean uniform-convexity bound for $B=B_L$.

\textbf{\ref{paper:base:centers}}
The gradient and value recurrences in
Lemma~\resultpart{base:nesting-lemma}{base:nesting-lemma:part:recurrences},
with the orthogonality computation in
part~\ref{base:nesting-lemma:part:center-formulas}, give
Eq.~\eqref{paper:base-data}. Their initial value and gradient are zero.

\textbf{\ref{paper:base:locality}}
The exterior and interior conclusions of Lemma~\ref{base:reset-lemma}
give the two exact formulas for $B_{j+1}$.
Lemma~\resultpart{base:nesting-lemma}{base:nesting-lemma:part:construction}(i)
verifies the local formula for each subsequent $B_j$.

\textbf{\ref{paper:base:preservation}}
Lemma~\resultpart{base:nesting-lemma}{base:nesting-lemma:part:construction}(iii)--(iv)
gives neighborhood agreement at the earlier centers, equality
of values and derivatives outside $\overline S_j$, and dependence
only on the previously chosen centers.
For $j=L$, the neighborhood agreement is the identity $B=B_L$.
\end{proof}

\subsection{Technical lemmas for the uniformly convex component}
\label{app:base-support}

\subsubsection{Local gradient formulas and recentering}
\begin{lemma}[A prescribed local gradient increment determines the local function]
\label{base:profile-lemma}
Let $B_j$ be differentiable on $B\left(\bm c_j,r\right)=\left\{\bm x:\norm{\bm x-\bm c_j}<r\right\}$ for some $r>0$. The following statements hold.
\begin{claims}
\item \label{base:profile-lemma:part:equivalence} On that ball, the gradient rule Eq.~\eqref{base:gradient-profile} is equivalent to the function identity Eq.~\eqref{base:old-local-form}.
\item \label{base:profile-lemma:part:endpoint} When either equivalent condition holds and $\norm{\bm c_{j+1}-\bm c_j}<r$, the gradient increment is in the direction $\bm c_{j+1}-\bm c_j$ and has magnitude $\mu\norm{\bm c_{j+1}-\bm c_j}^{q-1}$.
\end{claims}
\end{lemma}
\begin{proof}
\begin{itemize}
\item \textbf{Proof of \ref{base:profile-lemma:part:equivalence}.}
Define the difference
\[
 f_j\left(\bm x\right):=B_j\left(\bm x\right)-B_j\left(\bm c_j\right)
 -\ip{\D B_j\left(\bm c_j\right)}{\bm x-\bm c_j}
 -\frac{\mu}{q}\norm{\bm x-\bm c_j}^q.
\]
Lemma~\resultpart{elem:power}{elem:power:part:derivatives} in Appx.~\ref{app:paper-elementary}, including
its extension at zero in part~\ref{elem:power:part:flatness}, gives
\[
 \D f_j\left(\bm x\right)
 =\D B_j\left(\bm x\right)-\D B_j\left(\bm c_j\right)
 -\mu\norm{\bm x-\bm c_j}^{q-2}\left(\bm x-\bm c_j\right)
 =0
 \qquad\text{by Eq.~\eqref{base:gradient-profile}}.
\]
Fix $\bm x\in B\left(\bm c_j,r\right)$ and set
$v\left(t\right)=f_j\left(\bm c_j+t\left(\bm x-\bm c_j\right)\right)$ and $g\left(t\right)=0$ on $[0,1]$.
The segment stays in this convex ball, and
$|v'\left(t\right)|=|\D f_j\left(\bm c_j+t\left(\bm x-\bm c_j\right)\right)\left[\bm x-\bm c_j\right]|=0=g'\left(t\right)$.
Theorem~\ref{sources:mean-value} therefore gives
\[
 |f_j\left(\bm x\right)-f_j\left(\bm c_j\right)|=|v\left(1\right)-v\left(0\right)|\le g\left(1\right)-g\left(0\right)=0.
\]
Since $f_j\left(\bm c_j\right)=0$, this proves Eq.~\eqref{base:old-local-form}.
Conversely, differentiation of Eq.~\eqref{base:old-local-form} gives
\[
\begin{aligned}
 \D B_j\left(\bm x\right)
 &=\D B_j\left(\bm c_j\right)
    +\mu\norm{\bm x-\bm c_j}^{q-2}\left(\bm x-\bm c_j\right),\\
 \D B_j\left(\bm x\right)-\D B_j\left(\bm c_j\right)
 &=\mu\norm{\bm x-\bm c_j}^{q-2}\left(\bm x-\bm c_j\right),
\end{aligned}
\]
by the same norm-power derivative formula. This is
Eq.~\eqref{base:gradient-profile}, including at $\bm x=\bm c_j$.
\item \textbf{Proof of \ref{base:profile-lemma:part:endpoint}.}
Set $\bm h=\bm c_{j+1}-\bm c_j\ne0$ and evaluate
Eq.~\eqref{base:gradient-profile} at $\bm c_{j+1}$:
\[
\begin{aligned}
 \D B_j\left(\bm c_{j+1}\right)-\D B_j\left(\bm c_j\right)
 &=\mu\norm{\bm h}^{q-2}\bm h,\\
 \norm{\D B_j\left(\bm c_{j+1}\right)-\D B_j\left(\bm c_j\right)}
 &=\mu\norm{\bm h}^{q-2}\norm{\bm h}
 =\mu\norm{\bm h}^{q-1}.
\end{aligned}
\]
The coefficient $\mu\norm{\bm h}^{q-2}$ is positive, so the increment
has the direction of $\bm h$, as claimed.
\end{itemize}
\end{proof}

\begin{lemma}[Expressing $B_j$ relative to the next center]
\label{base:recenter-lemma}
Suppose Eq.~\eqref{base:old-local-form} holds on a neighborhood containing $\bm c_{j+1}$ and $\bm x$. Then:
\begin{claims}
\item \label{base:recenter-lemma:part:identity} The function $B_j$ satisfies
\begin{equation}\label{base:old-recentered}
\begin{split}
 B_j\left(\bm x\right)={}&B_j\left(\bm c_{j+1}\right)
 +\ip{\D B_j\left(\bm c_{j+1}\right)}{\bm x-\bm c_{j+1}}\\
 &+\frac{\mu}{q}\left[
 \norm{\left(\bm c_{j+1}-\bm c_j\right)+\left(\bm x-\bm c_{j+1}\right)}^q
 -\norm{\bm c_{j+1}-\bm c_j}^q\right]\\
 &-\mu\norm{\bm c_{j+1}-\bm c_j}^{q-2}
       \ip{\bm c_{j+1}-\bm c_j}{\bm x-\bm c_{j+1}}.
\end{split}
\end{equation}
Its last two lines generally differ from $\frac{\mu\norm{\bm x-\bm c_{j+1}}^q}{q}$.
\end{claims}
\end{lemma}
\begin{proof}
\begin{itemize}
\item \textbf{Proof of \ref{base:recenter-lemma:part:identity}.}
Set $\bm h=\bm c_{j+1}-\bm c_j$ and $\bm z=\bm x-\bm c_{j+1}$. Eq.~\eqref{base:old-local-form} and Lemma~\resultpart{base:profile-lemma}{base:profile-lemma:part:equivalence} give
\[
\begin{aligned}
 B_j\left(\bm x\right)-B_j\left(\bm c_{j+1}\right)
 &=\ip{\D B_j\left(\bm c_j\right)}{\bm z}
   +\frac{\mu}{q}\left(\norm{\bm h+\bm z}^q-\norm{\bm h}^q\right),\\
 \D B_j\left(\bm c_j\right)&=\D B_j\left(\bm c_{j+1}\right)-\mu\norm{\bm h}^{q-2}\bm h.
\end{aligned}
\]
Substituting the second equality into the first, and distributing the
inner product, gives
\[
\begin{aligned}
 B_j\left(\bm x\right)-B_j\left(\bm c_{j+1}\right)
 &=\ip{\D B_j\left(\bm c_{j+1}\right)-\mu\norm{\bm h}^{q-2}\bm h}{\bm z}
    +\frac{\mu}{q}\left(\norm{\bm h+\bm z}^q-\norm{\bm h}^q\right)\\
 &=\ip{\D B_j\left(\bm c_{j+1}\right)}{\bm z}
    -\mu\norm{\bm h}^{q-2}\ip{\bm h}{\bm z}\\
 &\quad+\frac{\mu}{q}\left(\norm{\bm h+\bm z}^q-\norm{\bm h}^q\right).
\end{aligned}
\]
Add $B_j\left(\bm c_{j+1}\right)$ and substitute the definitions of $\bm h,\bm z$.
The resulting expression is exactly Eq.~\eqref{base:old-recentered}.
\end{itemize}
\end{proof}

\begin{corollary}[The specified replacement changes the remainder and preserves the gradient]
\label{base:replacement-effects}
Suppose Eq.~\eqref{base:old-recentered} and Eq.~\eqref{base:next-local-choice} hold on the neighborhood under consideration. Then:
\begin{claims}
\item \label{base:replacement-effects:part:difference} Their complete difference is
\begin{equation}\label{base:full-difference}
\begin{split}
 B_{j+1}\left(\bm x\right)-B_j\left(\bm x\right)
 ={}&C_j+\mu\norm{\bm c_{j+1}-\bm c_j}^{q-2}
              \ip{\bm c_{j+1}-\bm c_j}{\bm x-\bm c_{j+1}}\\
 &+\frac{\mu}{q}\left[
 \norm{\bm x-\bm c_{j+1}}^q+\norm{\bm c_{j+1}-\bm c_j}^q
 -\norm{\left(\bm c_{j+1}-\bm c_j\right)+\left(\bm x-\bm c_{j+1}\right)}^q\right].
\end{split}
\end{equation}
\item \label{base:replacement-effects:part:next-gradient} The new local gradient rule is
\[
 \D B_{j+1}\left(\bm x\right)-\D B_{j+1}\left(\bm c_{j+1}\right)
 =\mu\norm{\bm x-\bm c_{j+1}}^{q-2}\left(\bm x-\bm c_{j+1}\right).
\]
\item \label{base:replacement-effects:part:preserved-increment} The preceding center-to-center increment is retained:
\begin{equation}\label{base:successive-center-gradient}
 \D B_{j+1}\left(\bm c_{j+1}\right)-\D B_j\left(\bm c_j\right)
 =\mu\norm{\bm c_{j+1}-\bm c_j}^{q-2}\left(\bm c_{j+1}-\bm c_j\right).
\end{equation}
\end{claims}
\end{corollary}
\begin{proof}
\begin{itemize}
\item \textbf{Proof of \ref{base:replacement-effects:part:difference}.}
Set $\bm h=\bm c_{j+1}-\bm c_j$ and $\bm z=\bm x-\bm c_{j+1}$.
Subtract Eq.~\eqref{base:old-recentered} from
Eq.~\eqref{base:next-local-choice}; both contain the same term
$B_j\left(\bm c_{j+1}\right)+\ip{\D B_j\left(\bm c_{j+1}\right)}{\bm z}$.
After canceling that term, the remaining expressions are
\[
\begin{aligned}
 B_{j+1}\left(\bm x\right)-B_j\left(\bm x\right)
 &=C_j+\frac{\mu}{q}\norm{\bm z}^q
   -\frac{\mu}{q}\left(\norm{\bm h+\bm z}^q-\norm{\bm h}^q\right)
   +\mu\norm{\bm h}^{q-2}\ip{\bm h}{\bm z}\\
 &=C_j+\mu\norm{\bm h}^{q-2}\ip{\bm h}{\bm z}
   +\frac{\mu}{q}\norm{\bm z}^q
   +\frac{\mu}{q}\norm{\bm h}^q-\frac{\mu}{q}\norm{\bm h+\bm z}^q\\
 &=C_j+\mu\norm{\bm h}^{q-2}\ip{\bm h}{\bm z}\\
 &\quad+\frac{\mu}{q}
       \left(\norm{\bm z}^q+\norm{\bm h}^q-\norm{\bm h+\bm z}^q\right).
\end{aligned}
\]
Substituting $\bm h,\bm z$ gives exactly Eq.~\eqref{base:full-difference}.
\item \textbf{Proof of \ref{base:replacement-effects:part:next-gradient}.}
Differentiate Eq.~\eqref{base:next-local-choice} using
Lemma~\resultpart{elem:power}{elem:power:part:derivatives}, with its extension at the center in
part~\ref{elem:power:part:flatness}. With $\bm z=\bm x-\bm c_{j+1}$,
\[
\begin{aligned}
 \D B_{j+1}\left(\bm x\right)
 &=\D B_j\left(\bm c_{j+1}\right)+\mu\norm{\bm z}^{q-2}\bm z,\\
 \D B_{j+1}\left(\bm c_{j+1}\right)
 &=\D B_j\left(\bm c_{j+1}\right),\\
 \D B_{j+1}\left(\bm x\right)-\D B_{j+1}\left(\bm c_{j+1}\right)
 &=\mu\norm{\bm z}^{q-2}\bm z.
\end{aligned}
\]
The second line uses the vanishing norm-power gradient at $\bm z=0$.
\item \textbf{Proof of \ref{base:replacement-effects:part:preserved-increment}.}
Using Eq.~\eqref{base:same-center-gradient} and then
Lemma~\resultpart{base:profile-lemma}{base:profile-lemma:part:endpoint},
\[
\begin{aligned}
 \D B_{j+1}\left(\bm c_{j+1}\right)-\D B_j\left(\bm c_j\right)
 &=\D B_j\left(\bm c_{j+1}\right)-\D B_j\left(\bm c_j\right)\\
 &=\mu\norm{\bm c_{j+1}-\bm c_j}^{q-2}
             \left(\bm c_{j+1}-\bm c_j\right).
\end{aligned}
\]
This is Eq.~\eqref{base:successive-center-gradient}.
\end{itemize}
\end{proof}

\begin{corollary}[The replacement cancels old curvature at its center]
\label{base:derivative-effects}
Under the preceding local specifications:
\begin{claims}
\item \label{base:derivative-effects:part:flatness} The higher derivatives satisfy
\[
 \D^kB_{j+1}\left(\bm c_{j+1}\right)=0\qquad\left(2\le k\le p\right).
\]
\item \label{base:derivative-effects:part:correction-hessian} The correction satisfies
\[
 \D^2\left(B_{j+1}-B_j\right)\left(\bm c_{j+1}\right)\prec0,
\]
so it cannot be a convex function added to the old base.
\item \label{base:derivative-effects:part:remainder} For every $\bm x$ in the neighborhood where Eq.~\eqref{base:next-local-choice} holds,
\begin{equation}\label{base:tangent-remainder}
 B_{j+1}\left(\bm x\right)-B_{j+1}\left(\bm c_{j+1}\right)
 -\ip{\D B_{j+1}\left(\bm c_{j+1}\right)}{\bm x-\bm c_{j+1}}
 =\frac{\mu}{q}\norm{\bm x-\bm c_{j+1}}^q.
\end{equation}
\end{claims}
\end{corollary}
\begin{proof}
\begin{itemize}
\item \textbf{Proof of \ref{base:derivative-effects:part:flatness}.}  The affine terms in Eq.~\eqref{base:next-local-choice} have zero derivatives of order at least two. The derivatives of $\frac{\norm{\bm x-\bm c_{j+1}}^q}{q}$ through order $p$ vanish there by Lemma~\resultpart{elem:power}{elem:power:part:flatness}, since $q>p$.
\item \textbf{Proof of \ref{base:derivative-effects:part:correction-hessian}.} For the old function, differentiating Eq.~\eqref{base:old-local-form} with the Hessian formula in Lemma~\resultpart{elem:power}{elem:power:part:derivatives} gives
\begin{equation}\label{base:old-hessian}
\begin{split}
 \D^2B_j\left(\bm c_{j+1}\right)={}&\mu\norm{\bm c_{j+1}-\bm c_j}^{q-2}I\\
 &+\mu\left(q-2\right)\norm{\bm c_{j+1}-\bm c_j}^{q-4}
 \left(\bm c_{j+1}-\bm c_j\right)\left(\bm c_{j+1}-\bm c_j\right)^\top\succ0.
\end{split}
\end{equation}
Set $\bm h=\bm c_{j+1}-\bm c_j\ne0$. For every nonzero $\bm u$,
Eq.~\eqref{base:old-hessian} and
part~\ref{base:derivative-effects:part:flatness} give
\[
\begin{aligned}
 \D^2\left(B_{j+1}-B_j\right)\left(\bm c_{j+1}\right)[\bm u,\bm u]
 &=0-\D^2B_j\left(\bm c_{j+1}\right)[\bm u,\bm u]\\
 &=-\mu\norm{\bm h}^{q-2}\norm{\bm u}^2
   -\mu\left(q-2\right)\norm{\bm h}^{q-4}\ip{\bm h}{\bm u}^2\\
 &\le-\mu\norm{\bm h}^{q-2}\norm{\bm u}^2<0.
\end{aligned}
\]
The inequality uses $q>2$; its strictness uses $\mu>0$ and
$\bm h,\bm u\ne0$. A twice differentiable convex function has a
positive semidefinite Hessian, so the correction is not convex.
\item \textbf{Proof of \ref{base:derivative-effects:part:remainder}.}
In Eq.~\eqref{base:next-local-choice}, the center value is
$B_j\left(\bm c_{j+1}\right)+C_j$ and its gradient is $\D B_j\left(\bm c_{j+1}\right)$.
Thus, with $\bm z=\bm x-\bm c_{j+1}$,
\[
\begin{aligned}
 &B_{j+1}\left(\bm x\right)-B_{j+1}\left(\bm c_{j+1}\right)
       -\ip{\D B_{j+1}\left(\bm c_{j+1}\right)}{\bm z}\\
 &\quad=\left[B_j\left(\bm c_{j+1}\right)+C_j
        +\ip{\D B_j\left(\bm c_{j+1}\right)}{\bm z}
        +\frac{\mu}{q}\norm{\bm z}^q\right]\\
 &\qquad-\left[B_j\left(\bm c_{j+1}\right)+C_j\right]
        -\ip{\D B_j\left(\bm c_{j+1}\right)}{\bm z}\\
 &\quad=\frac{\mu}{q}\norm{\bm z}^q.
\end{aligned}
\]
This is Eq.~\eqref{base:tangent-remainder}.
\end{itemize}
\end{proof}

\subsubsection{Joining one replacement to the old function}
\begin{lemma}[Local quadratic remainder]\label{base:local-quadratic-remainder}
Let $V\subset\R^d$ be open and convex, let $f\in C^2\left(V\right)$, and suppose
that $m_-I\preceq\D^2f\left(\bm z\right)\preceq m_+I$ for all $\bm z\in V$,
where $m_-\le m_+$ are finite real constants. Then
\begin{equation}\label{sources:local-quadratic-remainder}
 \frac{m_-}{2}\norm{\bm y-\bm x}^2
 \le\Breg f{\bm y}{\bm x}
 \le\frac{m_+}{2}\norm{\bm y-\bm x}^2
 \qquad\left(\bm x,\bm y\in V\right).
\end{equation}
\end{lemma}
\begin{proof}
Fix $\bm x,\bm y\in V$ and put $\bm d=\bm y-\bm x$.
Convexity of $V$ places the entire segment $\bm x+t\bm d$ in $V$.
For $v\left(t\right)=f\left(\bm x+t\bm d\right)$, the chain rule gives
$v'\left(0\right)=\ip{\D f\left(\bm x\right)}{\bm d}$ and
$v''\left(t\right)=\D^2f\left(\bm x+t\bm d\right)[\bm d,\bm d]$.
Thus
\[
 \begin{aligned}
 \Breg f{\bm y}{\bm x}
 &=v\left(1\right)-v\left(0\right)-v'\left(0\right)
 &&\text{by the Bregman definition},\\
 &=\int_0^1\left(v'\left(t\right)-v'\left(0\right)\right)\,dt
 &&\text{by the fundamental theorem of calculus},\\
 &=\int_0^1\int_0^t v''\left(s\right)\,ds\,dt
 &&\text{by the same theorem},\\
 &=\int_0^1\left(1-s\right)\D^2f\left(\bm x+s\bm d\right)[\bm d,\bm d]\,ds
 &&\text{by interchanging the continuous integrals}.
 \end{aligned}
\]
The assumed Hessian bounds give
$m_-\norm{\bm d}^2\le\D^2f\left(\bm x+s\bm d\right)[\bm d,\bm d]
\le m_+\norm{\bm d}^2$ for every $s\in[0,1]$. Since $1-s\ge0$
and $\int_0^1\left(1-s\right)\,ds=\frac{1}{2}$, substitution yields
\[
 \frac{m_-}{2}\norm{\bm d}^2
 \le\Breg f{\bm y}{\bm x}
 \le\frac{m_+}{2}\norm{\bm d}^2.
\]
This is Eq.~\eqref{sources:local-quadratic-remainder}.
\end{proof}

Fix $j$ and $R_j=\norm{\bm c_{j+1}-\bm c_j}$. Assume that $B_j$ is globally convex and $C^p$ and that Eq.~\eqref{base:old-local-form} holds on an open set containing $\overline S_j=\overline B\left(\bm c_{j+1},\theta R_j\right)$. Use $N_j^B$ and $B_{j+1}$ from Appx.~\ref{base:local}, and define
\begin{align*}
 \ell_j^B\left(\bm x\right)&:=B_j\left(\bm c_{j+1}\right)
                  +\ip{\D B_j\left(\bm c_{j+1}\right)}{\bm x-\bm c_{j+1}},\\
 u_j^B\left(\bm x\right)&:=B_j\left(\bm x\right)-B_j\left(\bm c_{j+1}\right)
                  -\ip{\D B_j\left(\bm c_{j+1}\right)}{\bm x-\bm c_{j+1}}.
\end{align*}
\begin{equation}\label{base:newbranch}
 v_j^B\left(\bm x\right):=C_j+\frac{\mu}{q}\norm{\bm x-\bm c_{j+1}}^q,
 \qquad
 \ell_j^B+u_j^B=B_j,\quad \ell_j^B+v_j^B=N_j^B.
\end{equation}
The scalar function $v_j^B$ is distinct from the vector $\bm v_j$.

\begin{lemma}[One replacement gives the new local formula and preserves the exterior]
\label{base:reset-lemma}
Under the assumptions above, there are constants $\theta,\gamma,a_{\mathrm{lift}}>0$, depending only on $q$, with $0<\gamma<\frac{\theta}{2}$ and $\theta\le\frac{1}{16}$, such that the choices in Eq.~\eqref{base:scaledconstants} give the following conclusions. Their explicit values are Eq.~\eqref{base:quadraticconstants}, Eq.~\eqref{base:liftwidth}, Eq.~\eqref{base:theta}, and Eq.~\eqref{base:gamma}.
\begin{claims}
\item \label{base:reset-lemma:part:exterior} \textbf{Agreement with $B_j$ outside radius $\frac{\rho_j}{2}=\frac{\theta R_j}{2}$.}
\[
 B_{j+1}\left(\bm x\right)=B_j\left(\bm x\right)
 \quad\text{if }\norm{\bm x-\bm c_{j+1}}\ge\frac{\rho_j}{2}.
\]
In particular, the formulas defining $B_{j+1}$ on $\norm{\bm x-\bm c_{j+1}}<\rho_j$ and $\norm{\bm x-\bm c_{j+1}}\ge\rho_j$ agree near the sphere $\norm{\bm x-\bm c_{j+1}}=\rho_j=\theta R_j$.
\item \label{base:reset-lemma:part:interior} \textbf{The prescribed expression near the new center.}
\[
 B_{j+1}\left(\bm x\right)=N_j^B\left(\bm x\right)
 \quad\text{if }\norm{\bm x-\bm c_{j+1}}\le\gamma R_j.
\]
\item \label{base:reset-lemma:part:center-data} \textbf{Value and derivatives at the new center.}
\[
\begin{aligned}
 B_{j+1}\left(\bm c_{j+1}\right)&=B_j\left(\bm c_{j+1}\right)+C_j,\\
 \D B_{j+1}\left(\bm c_{j+1}\right)&=\D B_j\left(\bm c_{j+1}\right),\\
 \D^kB_{j+1}\left(\bm c_{j+1}\right)&=0\qquad\left(2\le k\le p\right).
\end{aligned}
\]
\item \label{base:reset-lemma:part:smoothness} \textbf{Global $C^p$ smoothness.} The complete function $B_{j+1}$ belongs to $C^p\left(\R^d\right)$.
\item \label{base:reset-lemma:part:convexity} \textbf{Global convexity.} The complete function $B_{j+1}$ is convex.
\end{claims}
The H\"older bound for the correction $\Delta_j^B=B_{j+1}-B_j$ is proved separately in Lemma~\ref{base:regularity-lemma}.
\end{lemma}
\begin{proof}
We now choose $C_j,\delta_j^B,\theta$, and $\gamma$ to obtain
the two agreement inequalities in
Eq.~\eqref{paper:base-agreement-conditions}.
Both formulas have gradient $\D B_j\left(\bm c_{j+1}\right)$
at the next center, so subtracting $B_j$'s tangent there
reduces their comparison to their values above $B_j$'s tangent.
For $\bm z=\bm x-\bm c_{j+1}$, we have
\begin{equation}\label{base:calibration-comparison}
 \begin{aligned}
 N_j^B\left(\bm c_{j+1}+\bm z\right)
 -B_j\left(\bm c_{j+1}+\bm z\right)
 &=C_j+\frac{\mu}{q}\norm{\bm z}^q
   -\Breg{B_j}{\bm c_{j+1}+\bm z}{\bm c_{j+1}}\\
 &=v_j^B\left(\bm c_{j+1}+\bm z\right)
   -u_j^B\left(\bm c_{j+1}+\bm z\right).
 \end{aligned}
\end{equation}
\begin{itemize}
\item \textbf{Proof of \ref{base:reset-lemma:part:exterior}: agreement with $B_j$ near the boundary of $S_j$.}
We first bound the Bregman growth of
Eq.~\eqref{base:old-local-form}, which determines how much
the constant $C_j$ and the new $q^{th}$-power term can add
while leaving $B_j$ larger near the boundary.
We initially require $\rho_j=\theta R_j$ with
$0<\theta\le\frac{1}{4}$ and define
$\bm h=\bm c_{j+1}-\bm c_j$.
For $\norm{\bm z}\le\rho_j$ and every $t\in[0,1]$, put
$\bm w_t=\bm c_{j+1}+t\bm z$. The entire segment lies in
$\overline S_j$, and
\[
\begin{aligned}
 \norm{\bm w_t-\bm c_j}
 &=\norm{\bm h+t\bm z}\\
 &\ge R_j-t\norm{\bm z}
 \ge\left(1-\theta\right)R_j\ge\frac{3R_j}{4},\\
 \norm{\bm w_t-\bm c_j}
 &\le R_j+t\norm{\bm z}
 \le\left(1+\theta\right)R_j\le\frac{5R_j}{4}.
\end{aligned}
\]
The bounds follow from the reverse triangle and triangle inequalities.
The local formula Eq.~\eqref{base:old-local-form} holds on a neighborhood
of this segment; its affine terms have zero Hessian.
For any $\bm u\in\R^d$, Eq.~\eqref{elem:hessian-formula} in
Appx.~\ref{app:paper-elementary} therefore gives
\[
\begin{aligned}
 \D^2B_j\left(\bm w_t\right)[\bm u,\bm u]
 &=\mu\norm{\bm w_t-\bm c_j}^{q-2}\norm{\bm u}^2\\
 &\quad+\mu\left(q-2\right)\norm{\bm w_t-\bm c_j}^{q-4}
                   \ip{\bm w_t-\bm c_j}{\bm u}^2,\\
 \mu\norm{\bm w_t-\bm c_j}^{q-2}\norm{\bm u}^2
 &\le\D^2B_j\left(\bm w_t\right)[\bm u,\bm u]\\
 &\le\mu\left(q-1\right)\norm{\bm w_t-\bm c_j}^{q-2}\norm{\bm u}^2.
\end{aligned}
\]
Here $q>2$ makes the second Hessian term nonnegative, and
Cauchy--Schwarz bounds it by
$\mu\left(q-2\right)\norm{\bm w_t-\bm c_j}^{q-2}\norm{\bm u}^2$.
The preceding segment bounds and Eq.~\eqref{base:quadraticconstants}
now imply
\[
 2a_0\mu R_j^{q-2}I\preceq\D^2B_j\left(\bm x\right)
 \preceq2a_1\mu R_j^{q-2}I,
\]
with $a_0,a_1$ from Eq.~\eqref{base:quadraticconstants}.
To bound growth above $B_j$'s tangent, we apply
Lemma~\ref{base:local-quadratic-remainder} with
$f=B_j$, $V=S_j$, $\bm x=\bm c_{j+1}$,
$\bm y=\bm c_{j+1}+\bm z$, and
$m_-=2a_0\mu R_j^{q-2}$, $m_+=2a_1\mu R_j^{q-2}$.
Since $u_j^B\left(\bm c_{j+1}+\bm z\right)=
\Breg{B_j}{\bm c_{j+1}+\bm z}{\bm c_{j+1}}$, and continuity extends
the bound to $\norm{\bm z}=\theta R_j$, we obtain
\begin{equation}\label{base:quadraticbounds}
 a_0\mu R_j^{q-2}\norm{\bm z}^2
 \le u_j^B\left(\bm c_{j+1}+\bm z\right)
 \le a_1\mu R_j^{q-2}\norm{\bm z}^2.
\end{equation}

Thus $B_j$ grows quadratically above its tangent at
$\bm c_{j+1}$, whereas $N_j^B$ adds
$C_j+\frac{\mu}{q}\norm{\bm z}^q$ above $B_j$'s tangent.
Since $q>2$, we can choose $C_j$ so that $B_j$ is larger near
the boundary of $S_j$ while $N_j^B$ is larger near its center.

We choose $\delta_j^B=\frac{C_j}{4}$. On the annulus
$\frac{\theta R_j}{2}\le\norm{\bm z}\le\theta R_j$,
Eq.~\eqref{base:quadraticbounds} shows that
$B_j-N_j^B\ge\delta_j^B$ follows if
\[
 \underbrace{\frac{a_0}{4}\mu\theta^2R_j^q}_{\text{growth bound for }B_j}
 \ge
 \underbrace{C_j+\frac{\mu}{q}\left(\theta R_j\right)^q}_{\text{growth bound for }N_j^B}
 +\underbrace{\frac{C_j}{4}}_{\text{smoothing width}}.
\]
We choose $C_j=a_{\mathrm{lift}}\mu\theta^2R_j^q$, with
$a_{\mathrm{lift}}=\frac{a_0}{8}$, and denote the common
positive scale by
\[
 K_j^B:=\mu R_j^{q-2}\rho_j^2=\mu\theta^2R_j^q,
 \qquad C_j=a_{\mathrm{lift}}K_j^B,\qquad
 \delta_j^B=d_*K_j^B,\quad d_*=\frac{a_0}{32}.
\]
Dividing the required inequality by
$K_j^B=\mu\theta^2R_j^q$ leaves
\begin{equation}\label{base:outerconstraint}
 a_{\mathrm{lift}}+\frac{\theta^{q-2}}{q}+d_*
 \le\frac{a_0}{4}.
 \end{equation}
Using the choices of $a_{\mathrm{lift}}$ and $d_*$, this leaves
\[
 \begin{aligned}
 \frac{\theta^{q-2}}{q}
 &\le\frac{a_0}{4}-a_{\mathrm{lift}}-d_*\\
 &=\frac{a_0}{4}-\frac{a_0}{8}-\frac{a_0}{32}
 &&\text{by Eq.~\eqref{base:liftwidth}},\\
 &=\frac{3a_0}{32}.
 \end{aligned}
\]
We therefore choose
$\theta=\min\left\{\frac{1}{16},
\left(\frac{3qa_0}{32}\right)^{\frac{1}{q-2}}\right\}$,
as in Eq.~\eqref{base:theta}.
This is positive because $q>2$.
The cap $\frac{1}{16}$ meets the initial requirement
$\theta\le\frac{1}{4}$ and will also keep evaluation points
outside the replacement ball in Step~5, where
$\theta<\frac{1}{\sqrt{120}}$ is sufficient.
Consequently, on $\frac{\rho_j}{2}\le\norm{\bm z}\le\rho_j$,
\[
 \frac{u_j^B-v_j^B}{K_j^B}
 \ge \frac{a_0}{4}-\frac{a_0}{8}-\frac{3a_0}{32}
 =\frac{a_0}{32}=\frac{\delta_j^B}{K_j^B}.
\]
Thus $B_j-N_j^B=u_j^B-v_j^B\ge\delta_j^B$ on the annulus $\frac{\rho_j}{2}\le\norm{\bm x-\bm c_{j+1}}\le\rho_j$, equivalently $\frac{\theta R_j}{2}\le\norm{\bm z}\le\theta R_j$. Lemma~\ref{paper:smooth-maximum-properties} gives $B_{j+1}=B_j$ there; beyond radius $\rho_j$, equality is part of definition Eq.~\eqref{base:reset}. Together these prove
\begin{equation}\label{base:outerexact}
 B_{j+1}\left(\bm x\right)=B_j\left(\bm x\right)
 \quad\text{if }\norm{\bm x-\bm c_{j+1}}\ge\frac{\rho_j}{2}.
\end{equation}
This includes equality near the joining sphere $\norm{\bm x-\bm c_{j+1}}=\rho_j=\theta R_j$, as claimed in part~\ref{base:reset-lemma:part:exterior}.

\item \textbf{Proof of \ref{base:reset-lemma:part:interior}: agreement with $N_j^B$ near $\bm c_{j+1}$.}
On an inner ball $\norm{\bm z}\le\gamma R_j$, the value
increase $C_j$ must exceed $B_j$'s growth above its tangent
by $\delta_j^B$.
The upper quadratic bound in Eq.~\eqref{base:quadraticbounds}
therefore gives the sufficient condition
\begin{equation}\label{base:innerconstraint}
 \begin{aligned}
 a_1\mu R_j^q\gamma^2
 &\le C_j-\delta_j^B=\frac{3C_j}{4}
 &&\text{since }\delta_j^B=\frac{C_j}{4},\\
 a_1\gamma^2&\le\frac{3a_0\theta^2}{32}
 &&\text{after division by }\mu R_j^q,\\
 a_1\left(\frac{\gamma}{\theta}\right)^2+d_*&\le a_{\mathrm{lift}}
 &&\text{by Eq.~\eqref{base:liftwidth}}.
 \end{aligned}
\end{equation}
We therefore define
$\gamma=\theta\sqrt{\frac{3a_0}{32a_1}}$, as in Eq.~\eqref{base:gamma}.
The required restriction $\gamma<\frac{\theta}{2}$ follows
without another minimum, since Eq.~\eqref{base:quadraticconstants} gives
\[
 \frac{a_0}{a_1}=\frac{\left(\frac{3}{5}\right)^{q-2}}{q-1}<1,
 \qquad \left(\frac{\gamma}{\theta}\right)^2<\frac{3}{32}<\frac{1}{4}.
\]
These choices give Eq.~\eqref{base:scaledconstants} and $0<\gamma<\frac{\theta}{2}$. For $\norm{\bm z}\le\gamma R_j$,
\[
 \frac{v_j^B-u_j^B}{K_j^B}
 \ge \frac{a_0}{8}-\frac{3a_0}{32}
 =\frac{a_0}{32}=\frac{\delta_j^B}{K_j^B}.
\]
Therefore $N_j^B-B_j\ge\delta_j^B$ throughout the closed inner ball $\norm{\bm x-\bm c_{j+1}}\le\gamma R_j$. Lemma~\ref{paper:smooth-maximum-properties} and definition Eq.~\eqref{base:reset} give
\begin{equation}\label{base:innerexact}
 B_{j+1}\left(\bm x\right)=N_j^B\left(\bm x\right)
 \quad\text{if }\norm{\bm x-\bm c_{j+1}}\le\gamma R_j,
\end{equation}
which proves part~\ref{base:reset-lemma:part:interior}.

\item \textbf{Proof of \ref{base:reset-lemma:part:center-data}: the value and derivatives at the center.}
By part~\ref{base:reset-lemma:part:interior}, $B_{j+1}$ and $N_j^B$ agree on a neighborhood of $\bm c_{j+1}$. Evaluate and differentiate the defining formula for $N_j^B$. Its constant term is $B_j\left(\bm c_{j+1}\right)+C_j$, its affine term has gradient $\D B_j\left(\bm c_{j+1}\right)$, and every derivative of $\frac{\mu\norm{\bm x-\bm c_{j+1}}^q}{q}$ of order $1,\ldots,p$ vanishes at the center by Lemma~\resultpart{elem:power}{elem:power:part:flatness}, since $q>p$. Hence
\begin{equation}\label{base:centerreset}
\begin{split}
 B_{j+1}\left(\bm c_{j+1}\right)&=B_j\left(\bm c_{j+1}\right)+C_j,\\
 \D B_{j+1}\left(\bm c_{j+1}\right)&=\D B_j\left(\bm c_{j+1}\right),\\
 \D^kB_{j+1}\left(\bm c_{j+1}\right)&=0\qquad\left(2\le k\le p\right),
\end{split}
\end{equation}
proving part~\ref{base:reset-lemma:part:center-data}.

\item \textbf{Proof of \ref{base:reset-lemma:part:smoothness}: global $C^p$ smoothness.}
The function $B_j$ is convex and $C^p$ by assumption. Lemma~\resultpart{elem:power}{elem:power:part:flatness}, together with part~\ref{elem:power:part:hessian}, gives the same properties for $N_j^B$, which is a translated norm expression plus an affine function. To verify the open-neighborhood hypothesis for joining, use Eq.~\eqref{base:quadraticbounds} at $\norm{\bm z}=\rho_j$ and Eq.~\eqref{base:theta}:
\[
 \frac{B_j-N_j^B}{K_j^B}
 =\frac{u_j^B-v_j^B}{K_j^B}
 \ge a_0-\frac{a_0}{8}-\frac{3a_0}{32}
 =\frac{25a_0}{32}>\frac{a_0}{32}
 =\frac{\delta_j^B}{K_j^B}.
\]
By continuity, this strict inequality holds on an open neighborhood of the joining sphere $\norm{\bm x-\bm c_{j+1}}=\rho_j$. Thus Lemma~\resultpart{elem:composition-lemma}{elem:composition-lemma:part:joining} applies with $U=\R^d$, $r=p$, $f_1=B_j$, $f_2=N_j^B$, and ball $B\left(\bm c_{j+1},\rho_j\right)$. Its $C^p$ conclusion proves this part.

\item \textbf{Proof of \ref{base:reset-lemma:part:convexity}: global convexity.}
Apply the convexity conclusion of Lemma~\resultpart{elem:composition-lemma}{elem:composition-lemma:part:joining}, with the assumptions verified in part~\ref{base:reset-lemma:part:smoothness}.
\end{itemize}
\end{proof}

\subsubsection{Derivative and H\"older estimates for \texorpdfstring{$\Delta_j^B$}{the base correction}}
We use $\Delta_j^B=B_{j+1}-B_j$, $S_j=B\left(\bm c_{j+1},\theta R_j\right)$, and $U_j=B\left(\bm c_{j+1},\gamma R_j\right)$ from Appx.~\ref{base:local}.
By Eq.~\eqref{base:reset}, the definition $M_\delta\left(s,t\right)=s+\psi_\delta\left(t-s\right)$, and Lemma~\resultpart{base:reset-lemma}{base:reset-lemma:part:exterior}--\ref{base:reset-lemma:part:interior},
\begin{equation}\label{base:correctionformula}
 \Delta_j^B\left(\bm x\right)=
 \begin{cases}
 v_j^B\left(\bm x\right)-u_j^B\left(\bm x\right),
      &\norm{\bm x-\bm c_{j+1}}\le\gamma R_j,\\
 \psi_{\delta_j^B}\left(v_j^B\left(\bm x\right)-u_j^B\left(\bm x\right)\right),
      &\gamma R_j<\norm{\bm x-\bm c_{j+1}}<\frac{\rho_j}{2},\\
 0,&\norm{\bm x-\bm c_{j+1}}\ge\frac{\rho_j}{2}.
 \end{cases}
\end{equation}
The composition $\psi_{\delta_j^B}\left(v_j^B-u_j^B\right)$ also equals $\Delta_j^B$ throughout the closed support ball $\overline S_j$, where $\norm{\bm x-\bm c_{j+1}}\le\theta R_j$. Subtracting the two expressions in Eq.~\eqref{base:newbranch} and inserting Eq.~\eqref{base:old-local-form} gives on that same set,
\begin{equation}\label{base:branchdifference}
\begin{aligned}
 v_j^B\left(\bm x\right)-u_j^B\left(\bm x\right)
 ={}&N_j^B\left(\bm x\right)-B_j\left(\bm x\right)\\
 ={}&C_j+\frac{\mu}{q}\norm{\bm x-\bm c_{j+1}}^q
       -\frac{\mu}{q}\left(\norm{\bm x-\bm c_j}^q-R_j^q\right)\\
 &+\mu R_j^{q-2}
   \ip{\bm c_{j+1}-\bm c_j}{\bm x-\bm c_{j+1}}.
\end{aligned}
\end{equation}
The inherited value $B_j\left(\bm c_j\right)$ and gradient $\D B_j\left(\bm c_j\right)$ have canceled.

\begin{lemma}[Distance, regularity, and derivative bounds for $v_j^B-u_j^B$]
\label{base:branch-regularity-lemma}
Under the assumptions and parameter choices of Lemma~\ref{base:reset-lemma}, the difference in Eq.~\eqref{base:branchdifference} has the following properties. The constants $L_k^B$ are defined explicitly in the corresponding proof parts and depend only on $p,\nu$.
\begin{claims}
\item \label{base:branch-regularity-lemma:part:distances}
\textbf{Distance from the old center.} On the closed support ball $\overline S_j=\left\{\bm x:\norm{\bm x-\bm c_{j+1}}\le\theta R_j\right\}$,
\[
 \left(1-\theta\right)R_j\le\norm{\bm x-\bm c_j}\le\left(1+\theta\right)R_j.
\]
\item \label{base:branch-regularity-lemma:part:smoothness}
\textbf{Regularity of the difference.} The function $v_j^B-u_j^B$ belongs to $C^p\left(\R^d\right)$ and to $C^\infty\left(S_j\setminus\left\{\bm c_{j+1}\right\}\right)$, where
\[
 S_j\setminus\left\{\bm c_{j+1}\right\}
 =\left\{\bm x:0<\norm{\bm x-\bm c_{j+1}}<\theta R_j\right\}.
\]
\item \label{base:branch-regularity-lemma:part:derivative-bounds}
\textbf{Bounds through order $p$ on the closed support ball.} For $0\le k\le p$,
\[
 \norm{\D^k\left(v_j^B-u_j^B\right)\left(\bm x\right)}\le L_k^B\mu R_j^{q-k}
 \quad\text{when }\norm{\bm x-\bm c_{j+1}}\le\theta R_j.
\]
At $k=0$, the norm denotes absolute value.
\end{claims}
\end{lemma}
\begin{proof}
Put $\bm z=\bm x-\bm c_{j+1}$ and $\bm h=\bm c_{j+1}-\bm c_j$, so
$\norm{\bm h}=R_j$. With $f\left(\bm u\right)=\frac{\norm{\bm u}^q}{q}$,
Eq.~\eqref{base:branchdifference} reads on $\overline S_j$
\[
 v_j^B-u_j^B
 =C_j+\mu\left(f\left(\bm z\right)-f\left(\bm h+\bm z\right)+f\left(\bm h\right)
                   +\ip{\D f\left(\bm h\right)}{\bm z}\right).
\]
\begin{itemize}
\item \textbf{Proof of \ref{base:branch-regularity-lemma:part:distances}.}
The triangle and reverse triangle inequalities give
$R_j-\norm{\bm z}\le\norm{\bm h+\bm z}\le R_j+\norm{\bm z}$;
use $\norm{\bm z}\le\theta R_j$.

\item \textbf{Proof of \ref{base:branch-regularity-lemma:part:smoothness}.}
Globally, $u_j^B=B_j-\ell_j^B$ and $v_j^B=C_j+\mu f\left(\bm z\right)$ are $C^p$,
by the hypothesis on $B_j$ and Lemma~\resultpart{elem:power}{elem:power:part:flatness} in Appx.~\ref{app:paper-elementary}.
On $S_j\setminus\left\{\bm c_{j+1}\right\}$, both $\bm z$ and $\bm h+\bm z$ are
nonzero by part~\ref{base:branch-regularity-lemma:part:distances}; the displayed
formula and Lemma~\resultpart{elem:power}{elem:power:part:derivatives} give $C^\infty$
regularity there.

\item \textbf{Proof of \ref{base:branch-regularity-lemma:part:derivative-bounds}.}
For $k=0$, the displayed formula and Cauchy--Schwarz give
\[
\begin{aligned}
 |v_j^B\left(\bm x\right)-u_j^B\left(\bm x\right)|
 &\le C_j+\frac{\mu}{q}
       \left(\norm{\bm z}^q+\norm{\bm h+\bm z}^q+\norm{\bm h}^q\right)
       +\mu\norm{\bm h}^{q-1}\norm{\bm z}\\
 &\le\mu R_j^q\left(
       a_{\mathrm{lift}}\theta^2
       +\frac{\theta^q+\left(1+\theta\right)^q+1}{q}+\theta\right).
\end{aligned}
\]
The last line uses $C_j=a_{\mathrm{lift}}\theta^2\mu R_j^q$ and the
distance bounds in part~\ref{base:branch-regularity-lemma:part:distances}.
For $k=1$, differentiate the same formula:
\[
\begin{aligned}
 \D\left(v_j^B-u_j^B\right)\left(\bm x\right)
 &=\mu\left(\D f\left(\bm z\right)-\D f\left(\bm h+\bm z\right)+\D f\left(\bm h\right)\right),\\
 \norm{\D\left(v_j^B-u_j^B\right)\left(\bm x\right)}
 &\le\mu\left(\norm{\bm z}^{q-1}
       +\norm{\bm h+\bm z}^{q-1}+\norm{\bm h}^{q-1}\right)\\
 &\le\mu R_j^{q-1}
          \left(\theta^{q-1}+\left(1+\theta\right)^{q-1}+1\right).
\end{aligned}
\]
The norm-power gradient formula is Lemma~\resultpart{elem:power}{elem:power:part:derivatives}, with its extension at zero.
For $2\le k<p$, the constant and affine terms disappear, so that lemma
and the triangle inequality give
\[
\begin{aligned}
 \D^k\left(v_j^B-u_j^B\right)\left(\bm x\right)
 &=\mu\left(\D^kf\left(\bm z\right)-\D^kf\left(\bm h+\bm z\right)\right),\\
 \norm{\D^k\left(v_j^B-u_j^B\right)\left(\bm x\right)}
 &\le\mu C_{q,k}
          \left(\norm{\bm z}^{q-k}+\norm{\bm h+\bm z}^{q-k}\right)\\
 &\le\mu R_j^{q-k}C_{q,k}
          \left(\theta^{q-k}+\left(1+\theta\right)^{q-k}\right).
\end{aligned}
\]
Finally, for $k=p$, use the H\"older estimate in
Lemma~\resultpart{elem:power}{elem:power:part:holder}:
\[
\begin{aligned}
 \norm{\D^p\left(v_j^B-u_j^B\right)\left(\bm x\right)}
 &=\mu\norm{\D^pf\left(\bm z\right)-\D^pf\left(\bm h+\bm z\right)}\\
 &\le\mu H_{\rm pow}\norm{\bm z-\left(\bm h+\bm z\right)}^\nu\\
 &=\mu H_{\rm pow}R_j^\nu
 =\mu H_{\rm pow}R_j^{q-p}.
\end{aligned}
\]
These four calculations give the stated bounds with
\[
 L_k^B:=
 \begin{cases}
 a_{\mathrm{lift}}\theta^2+
   \frac{\theta^q+\left(1+\theta\right)^q+1}{q}+\theta,&k=0,\\
 \theta^{q-1}+\left(1+\theta\right)^{q-1}+1,&k=1,\\
 C_{q,k}\left[\theta^{q-k}+\left(1+\theta\right)^{q-k}\right],&2\le k<p,\\
 H_{\rm pow},&k=p.
 \end{cases}
\]
The third range is empty for $p=2$. All bounds include the center and
boundary by $C^p$ continuity.\qedhere
\end{itemize}
\end{proof}

\begin{lemma}[Derivative and H\"older bounds for one replacement]
\label{base:regularity-lemma}
Under the assumptions and parameter choices of Lemma~\ref{base:reset-lemma}, there is a finite constant $C_{\rm base}$, depending only on $p,\nu$ and the fixed smoothing function, for which the following hold.
\begin{claims}
\item \label{base:regularity-lemma:part:derivative-bounds}
\textbf{Derivative bounds at scale $R_j$.} The function $\Delta_j^B$ is $C^p$, and
\begin{equation}\label{base:derivativescale}
 \norm{\D^k\Delta_j^B}_\infty\le M_k^\Delta\mu R_j^{q-k},
 \qquad M_k^\Delta:=C_{\rm base}\quad\left(0\le k\le p\right).
\end{equation}
For $k=0$, the norm is absolute value.
\item \label{base:regularity-lemma:part:holder}
\textbf{H\"older continuity of the $p^{th}$ derivative.} For all $\bm x,\bm y\in\R^d$,
\begin{equation}\label{base:singleholder}
 \norm{\D^p\Delta_j^B\left(\bm x\right)-\D^p\Delta_j^B\left(\bm y\right)}
 \le H_\Delta\mu\norm{\bm x-\bm y}^\nu,
 \qquad H_\Delta:=C_{\rm base}.
\end{equation}
If also $B_j\in C^{p,\nu}$, then $B_{j+1}=B_j+\Delta_j^B$ belongs to $C^{p,\nu}$.
\end{claims}
\end{lemma}
\begin{proof}
We prove parts~\ref{base:regularity-lemma:part:derivative-bounds} and
\ref{base:regularity-lemma:part:holder} together by one bounded composition estimate.
Rotation and rescaling by $R_j=\norm{\bm c_{j+1}-\bm c_j}$ and
$\mu R_j^q$ give a correction whose derivative bounds do not depend on the phase.
Let $\bm e_1$ be the first coordinate
vector and choose an orthogonal matrix $Q_j$ with
$Q_j\bm e_1=\frac{\bm c_{j+1}-\bm c_j}{R_j}$.
The normalized input difference and resulting correction are
\[
 \delta_*:=\frac{a_{\mathrm{lift}}\theta^2}{4},\qquad
 \mathcal G_d\left(\bm z\right):=a_{\mathrm{lift}}\theta^2
 +\frac{\norm{\bm z}^q-\norm{\bm e_1+\bm z}^q+1}{q}
 +\ip{\bm e_1}{\bm z},
\]
\[
 \mathcal D_d\left(\bm z\right):=
 \begin{cases}
 \psi_{\delta_*}\left(\mathcal G_d\left(\bm z\right)\right),&\norm{\bm z}<\theta,\\
 0,&\norm{\bm z}\ge\theta.
 \end{cases}
\]
Define $a=\mu R_j^q>0$. Orthogonality of $Q_j$, its defining identity
$Q_j\bm e_1=\frac{\bm c_{j+1}-\bm c_j}{R_j}$, and
Eq.~\eqref{base:branchdifference} give
\[
\begin{aligned}
 &\left(v_j^B-u_j^B\right)\left(\bm c_{j+1}+R_jQ_j\bm z\right)\\
 &\quad=a\,a_{\mathrm{lift}}\theta^2
       +\frac{a}{q}
          \left(\norm{Q_j\bm z}^q-\norm{Q_j\left(\bm e_1+\bm z\right)}^q+1\right)
       +a\ip{Q_j\bm e_1}{Q_j\bm z}\\
 &\quad=a\left(
       a_{\mathrm{lift}}\theta^2
       +\frac{\norm{\bm z}^q-\norm{\bm e_1+\bm z}^q+1}{q}
       +\ip{\bm e_1}{\bm z}\right)
 =a\mathcal G_d\left(\bm z\right).
\end{aligned}
\]
Also $\delta_j^B=a\delta_*$, by Eq.~\eqref{base:scaledconstants}.
The scalar definition Eq.~\eqref{elem:psi-M-def} in
Appx.~\ref{app:paper-elementary} gives
\[
 \psi_{a\delta_*}\left(at\right)
 =\frac{at+a\delta_*\Theta\left(\frac{at}{a\delta_*}\right)}{2}
 =a\,\frac{t+\delta_*\Theta\left(\frac{t}{\delta_*}\right)}{2}
 =a\psi_{\delta_*}\left(t\right).
\]
Substituting these two identities into Eq.~\eqref{base:correctionformula}
proves
\[
 \Delta_j^B\left(\bm c_{j+1}+R_jQ_j\bm z\right)
 =\mu R_j^q\mathcal D_d\left(\bm z\right).
\]
Lemma~\resultpart{base:reset-lemma}{base:reset-lemma:part:smoothness}
and part~\ref{base:reset-lemma:part:exterior} imply that
$\mathcal D_d\in C^p$ and $\mathcal D_d=0$ on $\norm{\bm z}\ge\frac{\theta}{2}$.

On $U:=B\left(0,\theta\right)$, Lemma~\resultpart{base:branch-regularity-lemma}{base:branch-regularity-lemma:part:derivative-bounds}
gives $\norm{\D^k\mathcal G_d}_{\infty;U}\le L_k^B$ for $0\le k\le p$.
The $p^{th}$ derivatives of the two norm powers in $\mathcal G_d$ are each
H\"older continuous with parameter $H_{\rm pow}$ by Lemma~\resultpart{elem:power}{elem:power:part:holder};
its affine terms contribute zero. More explicitly, for distinct
$\bm z,\bm w\in U$,
\[
\begin{aligned}
 &\norm{\D^p\mathcal G_d\left(\bm z\right)-\D^p\mathcal G_d\left(\bm w\right)}\\
 &\quad\le
   \norm{\D^pf\left(\bm z\right)-\D^pf\left(\bm w\right)}
   +\norm{\D^pf\left(\bm e_1+\bm z\right)-\D^pf\left(\bm e_1+\bm w\right)}\\
 &\quad\le H_{\rm pow}\norm{\bm z-\bm w}^{\nu}
        +H_{\rm pow}\norm{\left(\bm e_1+\bm z\right)-\left(\bm e_1+\bm w\right)}^{\nu}\\
 &\quad=2H_{\rm pow}\norm{\bm z-\bm w}^{\nu},
\end{aligned}
\]
where $f\left(\bm u\right)=\frac{\norm{\bm u}^q}{q}$. Together with the bounds through
order $p$, this gives
\[
 C_G:=\max\left\{1,L_0^B,\ldots,L_p^B,2H_{\rm pow}\right\},\qquad
 \norm{\mathcal G_d}_{p,\nu;U}\le C_G.
\]
Set $J:=\left(-L_0^B-1,L_0^B+1\right)$ and
$C_\psi:=\norm{\psi_{\delta_*}}_{p,\nu;J}<\infty$.
The interval $J$ contains $\mathcal G_d\left(U\right)$, so
Lemma~\ref{elem:holder-composition} on the convex open sets $U$ and $J$ gives
\[
 \norm{\mathcal D_d}_{p,\nu;U}\le C_{\rm base},\qquad
 C_{\rm base}:=C_q^{\rm comp}C_\psi C_G^q.
\]
The bounded interval and smooth outer function make $C_\psi$ finite;
all constants are independent of $d$ and the phase scale.

The same bound holds globally. Derivative suprema extend because
$\mathcal D_d=0$ outside $U$. For H\"older continuity, pairs in $U$
are already covered. For a pair with one point outside $U$, exchange the
names if needed so that $\bm y\notin U$. If $\norm{\bm x}<\frac{\theta}{2}$, choose the point
$\bm z\in[\bm x,\bm y]\cap\partial B\left(0,\frac{\theta}{2}\right)$.
Both $\D^p\mathcal D_d\left(\bm z\right)$ and $\D^p\mathcal D_d\left(\bm y\right)$ vanish,
including at the sphere by $C^p$ continuity. Since $\bm x,\bm z\in U$,
\[
\begin{aligned}
 \norm{\D^p\mathcal D_d\left(\bm x\right)-\D^p\mathcal D_d\left(\bm y\right)}
 &=\norm{\D^p\mathcal D_d\left(\bm x\right)-\D^p\mathcal D_d\left(\bm z\right)}\\
 &\le C_{\rm base}\norm{\bm x-\bm z}^\nu\\
 &\le C_{\rm base}\norm{\bm x-\bm y}^\nu.
\end{aligned}
\]
Every remaining pair has both tensors zero. Thus
$\norm{\mathcal D_d}_{p,\nu;\R^d}\le C_{\rm base}$.
Finally, put $\bm z=R_j^{-1}Q_j^\top\left(\bm x-\bm c_{j+1}\right)$.
The scaling identity and the chain rule imply, for arbitrary vectors
$\bm u_1,\ldots,\bm u_k$ and $1\le k\le p$,
\[
\begin{aligned}
 \D^k\Delta_j^B\left(\bm x\right)[\bm u_1,\ldots,\bm u_k]
 &=\mu R_j^q\D^k\mathcal D_d\left(\bm z\right)
       [R_j^{-1}Q_j^\top\bm u_1,\ldots,R_j^{-1}Q_j^\top\bm u_k]\\
 &=\mu R_j^{q-k}\D^k\mathcal D_d\left(\bm z\right)
       [Q_j^\top\bm u_1,\ldots,Q_j^\top\bm u_k].
\end{aligned}
\]
The map $Q_j^\top$ preserves norms and maps the unit sphere onto itself.
Taking the full multilinear operator norm therefore gives equality of
the two derivative norms. For $k=0$, use the value identity directly.
For $\bm y=\bm c_{j+1}+R_jQ_j\bm w$ with $\bm w\ne\bm z$,
the same change of coordinates gives
\[
\begin{aligned}
 &\frac{\norm{\D^p\Delta_j^B\left(\bm x\right)-\D^p\Delta_j^B\left(\bm y\right)}}
             {\norm{\bm x-\bm y}^{\nu}}\\
 &\quad=\mu R_j^{q-p}
   \frac{\norm{\D^p\mathcal D_d\left(\bm z\right)-\D^p\mathcal D_d\left(\bm w\right)}}
        {R_j^\nu\norm{\bm z-\bm w}^{\nu}}\\
 &\quad=\mu R_j^{q-p-\nu}
   \frac{\norm{\D^p\mathcal D_d\left(\bm z\right)-\D^p\mathcal D_d\left(\bm w\right)}}
        {\norm{\bm z-\bm w}^{\nu}}.
\end{aligned}
\]
Every point and every distinct pair occur under this invertible change of
coordinates. Taking suprema in the derivative identities and using $q=p+\nu$
in the pairwise ratio yields
\begin{equation}\label{base:chainscaling}
\begin{aligned}
 \norm{\D^k\Delta_j^B}_\infty
 &=\mu R_j^{q-k}\norm{\D^k\mathcal D_d}_\infty\quad\left(0\le k\le p\right),\\
 \frac{\norm{\D^p\Delta_j^B\left(\bm x\right)-\D^p\Delta_j^B\left(\bm y\right)}}
      {\norm{\bm x-\bm y}^{\nu}}
 &=\mu
   \frac{\norm{\D^p\mathcal D_d\left(\bm z\right)-\D^p\mathcal D_d\left(\bm w\right)}}
        {\norm{\bm z-\bm w}^{\nu}}
 \le\mu C_{\rm base}.
\end{aligned}
\end{equation}
These prove both claims with $M_k^\Delta=H_\Delta=C_{\rm base}$.
If also $B_j\in C^{p,\nu}$, the triangle inequality gives, for all $\bm x,\bm y\in\R^d$,
\[
 \norm{\D^pB_{j+1}\left(\bm x\right)-\D^pB_{j+1}\left(\bm y\right)}
 \le\norm{\D^pB_j\left(\bm x\right)-\D^pB_j\left(\bm y\right)}
    +H_\Delta\mu\norm{\bm x-\bm y}^\nu.
\]
Thus $B_{j+1}\in C^{p,\nu}$.
\end{proof}

\subsubsection{Hessian and segment-length bounds for uniform convexity}
\label{app:base-curvature}
Let $B=B_L$ be as in Lemma~\ref{base:nesting-lemma}, and set
$\mathcal C=\left\{\bm c_0,\ldots,\bm c_L\right\}$ and
$\operatorname{dist}\left(\bm x,\mathcal C\right)=\min_i\norm{\bm x-\bm c_i}$.

\begin{lemma}[The Hessian is bounded below by distance to the centers]
\label{base:hessian-lemma}
\begin{claims}
\item \label{base:hessian-lemma:part:lower-bound} The final base satisfies
\begin{equation}\label{base:hessianlower}
 \D^2B\left(\bm x\right)\succeq
 \mu\operatorname{dist}\left(\bm x,\mathcal C\right)^{q-2}I
 \qquad\left(\bm x\in\R^d\right).
\end{equation}
\end{claims}
\end{lemma}
\begin{proof}
\begin{itemize}
\item \textbf{Proof of \ref{base:hessian-lemma:part:lower-bound}.}
Keep the full set $\mathcal C$ fixed and put $K\left(\bm x\right)=\mu\operatorname{dist}\left(\bm x,\mathcal C\right)^{q-2}I$. For every center $\bm c_i$, Lemma~\resultpart{elem:power}{elem:power:part:hessian} in Appx.~\ref{app:paper-elementary} gives
\[
 \D^2\!\left(\frac{\mu}{q}\norm{\bm x-\bm c_i}^q\right)
 \succeq\mu\norm{\bm x-\bm c_i}^{q-2}I\succeq K\left(\bm x\right).
\]
In particular $\D^2B_0\succeq K$. Assume $\D^2B_j\succeq K$; the displayed
bound also gives $\D^2N_j^B\succeq K$, since the remaining terms of $N_j^B$
are affine. On $S_j$, set
$s=N_j^B\left(\bm x\right)-B_j\left(\bm x\right)$,
$\vartheta=\psi'_{\delta_j^B}\left(s\right)$, and
$\bm a=\D N_j^B\left(\bm x\right)-\D B_j\left(\bm x\right)$.
Lemma~\resultpart{elem:max}{elem:max:part:psi-derivatives} gives $0\le\vartheta\le1$ and
$\psi''_{\delta_j^B}\left(s\right)\ge0$.
Both inputs have lower curvature $K$, and the smooth maximum's
additional positive semidefinite term can only increase it.
The Hessian formula in
Lemma~\resultpart{elem:composition-lemma}{elem:composition-lemma:part:hessian}
in that appendix gives
\[
\begin{aligned}
 \D^2B_{j+1}\left(\bm x\right)
 &=\left(1-\vartheta\right)\D^2B_j\left(\bm x\right)
       +\vartheta\D^2N_j^B\left(\bm x\right)
       +\psi''_{\delta_j^B}\left(s\right)\bm a\bm a^\top,\\
 \D^2B_{j+1}\left(\bm x\right)-K\left(\bm x\right)
 &=\left(1-\vartheta\right)\left(\D^2B_j\left(\bm x\right)-K\left(\bm x\right)\right)\\
 &\quad+\vartheta\left(\D^2N_j^B\left(\bm x\right)-K\left(\bm x\right)\right)
       +\psi''_{\delta_j^B}\left(s\right)\bm a\bm a^\top
 \succeq0.
\end{aligned}
\]
Each term on the last two lines is positive semidefinite. Hence
\[
 \D^2B_{j+1}=\D^2M_{\delta_j^B}\left(B_j,N_j^B\right)\succeq K.
\]
Outside $S_j$, including its boundary, Lemma~\resultpart{base:reset-lemma}{base:reset-lemma:part:exterior} gives neighborhood equality $B_{j+1}=B_j$, so the same bound holds. Induction and $B=B_L$ prove Eq.~\eqref{base:hessianlower}.
\end{itemize}
\end{proof}

\begin{lemma}[A segment spends a fixed fraction away from all centers]
\label{base:porosity}
\begin{claims}
\item \label{base:porosity:part:length} For every pair of distinct points $\bm x,\bm y\in\R^d$, set $\ell=\norm{\bm y-\bm x}$. The segment $[\bm x,\bm y]$ has a measurable subset of length at least $\frac{\ell}{2}$ such that every point $\bm z$ in that subset satisfies
\[
 \operatorname{dist}\left(\bm z,\mathcal C\right)\ge\frac{\ell}{16}.
\]
\end{claims}
\end{lemma}

\begin{proof}
\begin{itemize}
\item \textbf{Proof of \ref{base:porosity:part:length}.}
Fix distinct $\bm x,\bm y\in\R^d$, set $\ell=\norm{\bm y-\bm x}$ and $\rho_\ell=\frac{\ell}{16}$, and let $\bm u=\frac{\bm y-\bm x}{\ell}$. Parameterize the supporting line by $\bm x+t\bm u$; its speed is one, and the segment corresponds to $0\le t\le\ell$.
\begin{enumerate}[label=\textbf{Step \arabic*.},leftmargin=*]
\item \textbf{Bound the separation of three center projections.}
Define
\[
 \alpha_i:=\ip{\bm c_i-\bm x}{\bm u},\qquad
 \bm z_i:=\bm c_i-\bm x-\alpha_i\bm u,\qquad
 \mathcal I:=\left\{i:\norm{\bm z_i}\le\rho_\ell\right\}.
\]
\input{figures/segment-length-figure}
The vector $\bm z_i$ points from the projection $\bm x+\alpha_i\bm u$ to the
center $\bm c_i$, as illustrated in Figure~\ref{fig:segment-length-geometry}.
Thus $\bm z_i\perp\bm u$, and only centers indexed by $\mathcal I$ can have their radius-$\rho_\ell$ balls meet the line. For $i<j<k$ in $\mathcal I$, orthogonality of the phase displacements $\bm h_a$ gives
\[
 \ip{\bm c_i-\bm c_j}{\bm c_k-\bm c_j}
 =\left\langle-\sum_{a=i}^{j-1}\bm h_a,
                    \sum_{b=j}^{k-1}\bm h_b\right\rangle=0.
\]
Substitute
$\bm c_i-\bm c_j=\left(\alpha_i-\alpha_j\right)\bm u+\left(\bm z_i-\bm z_j\right)$
and its counterpart with $k$ in place of $i$. Since the residuals are
orthogonal to $\bm u$ and $\norm{\bm u}=1$,
\[
\begin{aligned}
 0
 &=\ip{\bm c_i-\bm c_j}{\bm c_k-\bm c_j}\\
 &=\left(\alpha_i-\alpha_j\right)\left(\alpha_k-\alpha_j\right)
       +\ip{\bm z_i-\bm z_j}{\bm z_k-\bm z_j},\\
 |\left(\alpha_i-\alpha_j\right)\left(\alpha_k-\alpha_j\right)|
 &=|\ip{\bm z_i-\bm z_j}{\bm z_k-\bm z_j}|\\
 &\le\norm{\bm z_i-\bm z_j}\norm{\bm z_k-\bm z_j}\\
 &\le\left(\norm{\bm z_i}+\norm{\bm z_j}\right)
       \left(\norm{\bm z_k}+\norm{\bm z_j}\right)\\
 &\le\left(2\rho_\ell\right)\left(2\rho_\ell\right)=4\rho_\ell^2.
\end{aligned}
\]
The inequalities are Cauchy--Schwarz, the triangle inequality, and
membership in $\mathcal I$, respectively.
Consequently no three projections from $\mathcal I$ can be pairwise more than $2\rho_\ell$ apart: ordering their center indices would contradict this bound.

\item \textbf{Cover the projections by at most two short intervals.}
If $\mathcal I$ is empty, no center ball meets the line. Otherwise let $a_1$ be the smallest projection, and cover all projections in $[a_1,a_1+2\rho_\ell]$. If any remain, let $a_2$ be the smallest remaining projection. Every remaining projection lies in $[a_2,a_2+2\rho_\ell]$; a projection larger than $a_2+2\rho_\ell$ would form three pairwise separated projections with $a_1,a_2$, contrary to Step~1.

\item \textbf{Bound the length within distance $\rho_\ell$ of a center.}
If $\norm{\bm x+t\bm u-\bm c_i}<\rho_\ell$, then $i\in\mathcal I$ and $|t-\alpha_i|<\rho_\ell$. Each interval from Step~2 therefore contributes only within
\[
 [a_b-\rho_\ell,a_b+3\rho_\ell],
 \qquad\text{an interval of length }4\rho_\ell.
\]
There are at most two such intervals, so the part of $[\bm x,\bm y]$ at distance less than $\rho_\ell$ from $\mathcal C$ has length at most $8\rho_\ell$. Its complement is measurable and has length at least
\[
 \ell-8\rho_\ell=\frac{\ell}{2}.
\]
Every point in this complement has distance at least $\rho_\ell=\frac{\ell}{16}$ from $\mathcal C$, as required.
\end{enumerate}
\end{itemize}
\end{proof}

%% file: figures/segment-length-figure.tex
\begin{figure}[!htbp]
\centering
\includegraphics[width=.485\linewidth]{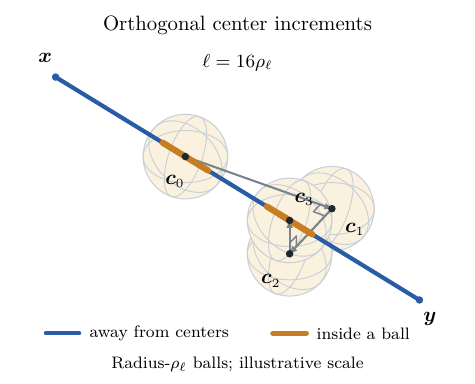}\hfill
\includegraphics[width=.485\linewidth]{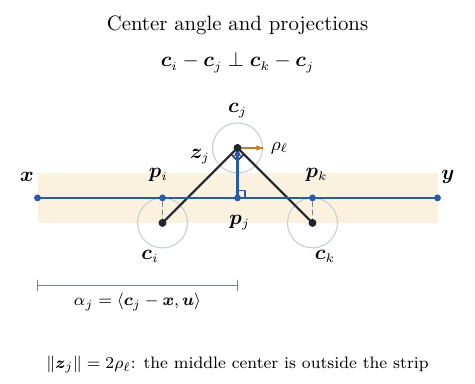}
\caption{Geometry for Lemma~\ref{base:porosity}, with $\rho_\ell=\frac{\ell}{16}$.
Left: orthogonal center increments at illustrative scales and
radius-$\rho_\ell$ balls, with gold showing the parts of the segment inside
their union. Right: line projections $\bm p_i=\bm x+\alpha_i\bm u$ and
residual $\bm z_j=\bm c_j-\bm p_j$. The three projections are pairwise more
than $2\rho_\ell$ apart, while the middle center lies outside the
radius-$\rho_\ell$ strip. These fixed-radius balls differ from the
replacement neighborhoods $S_j$.}
\label{fig:segment-length-geometry}
\end{figure}

%% file: appendix-phases.tex
\section{Supporting Proofs for the Phase Functions}\label{app:paper-phases}

\subsection{Constructing and verifying the convex phase functions}
\label{collar:section}

\subsubsection{Requirements for one phase}
\label{collar:gradient-requirements}

The uniformly convex base from Appx.~\ref{app:paper-base} gives the
geometry that the phase functions must preserve.
Lemma~\resultpart{base:nesting-lemma}{base:nesting-lemma:part:recurrences} gives
\begin{equation}\label{collar:base-gradient-increment}
 \bm b_j:=\D B\left(\bm c_{j+1}\right)-\D B\left(\bm c_j\right)
 =\mu R_j^{q-2}\left(\bm c_{j+1}-\bm c_j\right).
\end{equation}
For $F=B+\ell_{\mathrm{init}}+\sum_k A_k$, with
$\ell_{\mathrm{init}}\left(\bm x\right)=-s_0\ip{\bm v_0}{\bm x}$, the phase conditions are:
\begin{enumerate}[label=\textbf{R\arabic*.},ref=R\arabic*]
\item\label{collar:req:gradient}\textbf{Center gradients.}
The target gradients are
\begin{equation}\label{collar:required-gradients}
 \D F\left(\bm c_j\right)=-s_j\bm v_j,\qquad
 \D F\left(\bm c_{j+1}\right)=-s_{j+1}\bm v_{j+1}
\end{equation}
with $\D A_j\left(\bm c_j\right)=0$ and $A_j=I_j$ near $\bm c_{j+1}$, where
$I_j\left(\bm x\right)=\ip{\bm g_j}{\bm x-\bm c_{j+1}}$ and
\begin{equation}\label{collar:required-auxiliary-gradient}
 \bm g_j=s_j\bm v_j-s_{j+1}\bm v_{j+1}-\bm b_j.
\end{equation}
Lemma~\resultpart{design:balance-lemma}{design:balance-lemma:part:slope} in Appx.~\ref{app:paper-assembly}
derives this cancellation under \ref{collar:req:other-phases}.
We use $s_L\bm v_L=0$ in the terminal phase.
\item\label{collar:req:gap}\textbf{Optimality gap on chain prefixes.}
For the subspaces $W_{j,b}$ containing early evaluation points,
\begin{equation}\label{collar:required-residual}
 \begin{aligned}
 r_j:=A_j-I_j&=0&&\text{near }\bm c_{j+1},\\
 r_j\left(\bm x\right)&\ge e_j:=7\delta_j>0&& \left(\bm x\in W_{j,b}\right).
 \end{aligned}
\end{equation}
Appx.~\ref{global:gradient-and-gap} uses this excess and the orthogonal
center gradient to prove $F\left(\bm x\right)-\min F\ge7\delta_j$.
\item\label{collar:req:convexity}\textbf{Convexity:}
$A_j$ is convex, so Lemma~\ref{paper:bregman-additions} ensures that
adding it preserves the base's Bregman lower bound.
\item\label{collar:req:oracle}\textbf{Oracle consistency:}
future base changes miss $W_{j,b}$, and the current phase's oracle output uses at most one new direction.
\item\label{collar:req:holder}\textbf{Smoothness:}
$A_j\in C^{p,\nu}$ and
$\norm{\D^pA_j\left(\bm x\right)-\D^pA_j\left(\bm y\right)}
\le C_A\lambda\norm{\bm x-\bm y}^\nu$ for all $\bm x,\bm y\in\R^d$,
with $C_A$ independent of $m$.
\item\label{collar:req:other-phases}\textbf{Other phases:}
$\D A_k\left(\bm c_{j+1}\right)=\D A_k\left(\bm c_j\right)$ for $k\ne j$.
Appx.~\ref{global:verification} uses $E_j\perp E_k$ for $j\ne k$,
$\bm h_j\in E_j$, and $\bm v_{k+1}\in E_{k+1}$ when the next direction is present.
For $k=j-1$, the gradient remains $\bm g_{j-1}$ along the completed-center segment
by Lemma~\resultpart{paper:phase}{paper:phase:completed}.
\end{enumerate}
Throughout, $p\ge2$, $0<\nu\le1$, $q=p+\nu$, and
$q^*:=\frac{q}{q-1}$. Define $\norm{\bm a}_r=\left(\sum_i|a_i|^r\right)^{\frac{1}{r}}$;
unqualified norms are Euclidean and derivative norms are operator norms.

\subsubsection{Adapt the classical power chain to one phase}
\label{collar:chain}\label{sources:chain-section}

The phase chain follows the construction of
\cite[Sec.~5.1, Eq.~(5.4)]{GrapigliaNesterov2020}, whose definition,
matrix estimate, and minimizer formula are restated in
Lemma~\ref{sources:power-chain} below. The same chain structure appears in
\cite[Sec.~5.1]{BaiBullins2025}.\footnote{For the classical quadratic and cubic
constructions, see \cite[Secs.~2.1.2 and~4.3.1]{Nesterov2018}; for higher-order
power chains, see \cite[Sec.~3]{ArjevaniShamirShiff2019} and
\cite[Sec.~4]{Nesterov2021Tensor}. These citations identify the construction's
lineage; the source facts used below are stated in Lemma~\ref{sources:power-chain}.}

\begin{lemma}[Power chain {\cite[Lemmas~5.1--5.2 and their proofs]{GrapigliaNesterov2020}}]
\label{sources:power-chain}
Let $p\ge2$ be an integer, $0<\nu\le1$, $q=p+\nu$, and let
$2\le k\le d$ be integers. Define
\[
 f_k\left(\bm a\right):=\frac{1}{q}\left[
 \sum_{i=1}^{k-1}|a_i-a_{i+1}|^q+\sum_{i=k}^d|a_i|^q
 \right]-a_1\qquad\left(\bm a\in\R^d\right).
\]
Let $\mathsf U_k$ have diagonal entries $1$, superdiagonal entries $-1$,
and all other entries zero. Then
\begin{equation}\label{sources:chain-matrix-bound}
 \mathsf A_k=\operatorname{diag}\left(\mathsf U_k,I_{d-k}\right),\qquad
 f_k\left(\bm a\right)=\frac{1}{q}\norm{\mathsf A_k\bm a}_q^q-a_1,\qquad
 \norm{\mathsf A_k}_{2\to2}\le2.
\end{equation}
For $k=d$, the second diagonal block is absent.
The function $f_k$ has a unique global minimizer $\bm a_k^*$ and
\[
 \min_{\bm a\in\R^d}f_k\left(\bm a\right)=-\frac{q-1}{q}k,
 \qquad
 \norm{\bm a_k^*}<\frac{\left(k+1\right)^{\frac{3}{2}}}{\sqrt3}.
\]
The minimizer satisfies
\begin{equation}\label{sources:chain-ramp}
 \left(\bm a_k^*\right)_i=\left(k-i+1\right)_+\quad\left(1\le i\le d\right),\qquad
 \mathsf A_k\bm a_k^*
 =\left(\underbrace{1,\ldots,1}_{k},0,\ldots,0\right)^\top,
\end{equation}
where $\left(a\right)_+:=\max\{a,0\}$.
\end{lemma}

\begin{definition}[$p^{th}$-order zero-chain {\cite[Sec.~3.1, Definition~3]{CarmonEtAl2020}}]
\label{sources:zero-chain}
For a tensor $Q$, let $\operatorname{supp}Q$ contain every coordinate index
appearing in a nonzero entry of $Q$; for a vector this is its usual support.
For an integer $p\ge1$, a function $f:\R^d\to\R$ with derivatives
through order $p$ is a
$p^{th}$-order zero-chain if, for every $\bm x\in\R^d$ and $1\le i\le d$,
\[
 \operatorname{supp}\bm x\subseteq\{1,\ldots,i-1\}
 \quad\Longrightarrow\quad
 \bigcup_{r=1}^p\operatorname{supp}\D^r f\left(\bm x\right)
 \subseteq\{1,\ldots,i\}.
\]
If this holds for every positive integer $p$, it is called a zero-chain.
\end{definition}

\begingroup
\theoremstyle{definition}
\newtheorem*{phasechaindefinition}{Definition}
\begin{phasechaindefinition}[Phase chain and endpoint]
Fix $0\le j<L$, an integer $m\ge1$, and $\lambda,\tau_j>0$.
Let $\bm v_j,\bm u_{j,1},\ldots,\bm u_{j,m}$ be orthonormal, and put
\[
 E_j:=\operatorname{span}\{\bm v_j,\bm u_{j,1},\ldots,\bm u_{j,m}\},
 \qquad \bm w=t\bm v_j+\sum_{i=1}^m x_i\bm u_{j,i}\in E_j.
\]
Matrix products use coordinates in the ordered basis
$\left(\bm v_j,\bm u_{j,1},\ldots,\bm u_{j,m}\right)$. Thus
\[
 D_m\bm w:=D_m\left(t,x_1,\ldots,x_m\right)^\top,\qquad
 \bm e_0=\left(1,0,\ldots,0\right)^\top,
\]
where $\bm e_0$ represents $\bm v_j$.
Define the residual matrix by
\begin{equation}
D_m=
\begin{pmatrix}
m^{-\frac{1}{2}}&-1&0&\cdots&0\\
0&1&-1&\ddots&\vdots\\
\vdots&\ddots&\ddots&\ddots&0\\
0&\cdots&0&1&-1\\
0&\cdots&\cdots&0&1
\end{pmatrix},
\qquad
D_m\begin{pmatrix}t\\x_1\\\vdots\\x_m\end{pmatrix}
=\begin{pmatrix}
\frac{t}{\sqrt{m}}-x_1\\x_1-x_2\\\vdots\\x_{m-1}-x_m\\x_m
\end{pmatrix}.
\label{collar:map}
\end{equation}
For $m=1$, $D_m=\left(\begin{smallmatrix}1&-1\\0&1\end{smallmatrix}\right)$.
With $\phi\left(a\right)=\frac{|a|^q}{q}$, define the chain on the ambient space by
\[
 Q_j\left(\bm x\right):=\frac{\lambda}{q}\norm{D_mP_{E_j}\bm x}_q^q
 \qquad\left(\bm x\in\R^d\right).
\]
Thus $Q_j\left(\bm w\right)=\lambda\sum_{i=0}^m\phi\left(\left(D_m\bm w\right)_i\right)$ on $E_j$;
rows are indexed $0,\ldots,m$. Define the endpoint by
\begin{equation}
 t_j^\star:=\sqrt m\left(m+1\right)\tau_j,\qquad
 \bm h_j:=t_j^\star\bm v_j+\sum_{i=1}^m\left(m-i+1\right)\tau_j\bm u_{j,i},
 \qquad R_j:=\norm{\bm h_j}.
\label{collar:endpoint}
\end{equation}
\end{phasechaindefinition}
\endgroup

\begin{lemma}[The chain matrix and its endpoint]
\label{collar:chain-lemma}
For the chain $Q_j$ and endpoint defined above:
\begin{claims}
\item\label{collar:chain-lemma:part:residuals}\textbf{Residuals and gradient.}
$D_m$ is invertible, and
\begin{equation}
 D_m\begin{pmatrix}t_j^\star\\m\tau_j\\\vdots\\\tau_j\end{pmatrix}=\tau_j\bm1,
 \qquad D_m^\top\bm1=\frac{\bm e_0}{\sqrt m}.
\label{collar:residual-identities}
\end{equation}
\[
 \D Q_j\left(\bm h_j\right)=\frac{\lambda\tau_j^{q-1}}{\sqrt m}\bm v_j.
\]
\item\label{collar:chain-lemma:part:norm}\textbf{Operator norm.}
\begin{equation}
 \norm{D_m}_{2\to2}\le2,
 \qquad \norm{D_m\bm w}_r\le2\norm{\bm w}\quad\left(\bm w\in E_j,\ r\ge2\right).
\label{collar:map-norm}
\end{equation}
\item\label{collar:chain-lemma:part:geometry}\textbf{Endpoint lengths.}
\begin{equation}
 R_j^2=\tau_j^2\frac{m\left(m+1\right)\left(8m+7\right)}{6},
\label{collar:R}
\end{equation}
\begin{equation}
 \tau_j m^{\frac{3}{2}}\le R_j\le\sqrt5\,\tau_j m^{\frac{3}{2}},\qquad
 \frac{\sqrt3}{2}R_j\le t_j^\star\le R_j,\qquad
 \sum_{i=1}^m\left(m-i+1\right)\tau_j\le m^2\tau_j.
\label{collar:endpoint-bounds}
\end{equation}
\[
 \frac{R_j}{\sqrt5}\le\norm{\bm h_j-t_j^\star\bm v_j}\le\frac{R_j}{2}.
\]
\item\label{collar:chain-lemma:part:prefix}\textbf{One additional direction.}
For every integer $0\le k<m$ and
$\bm w\in\operatorname{span}\{\bm v_j,\bm u_{j,1},\ldots,\bm u_{j,k}\}$,
\[
 Q_j\left(\bm w\right),\ \D Q_j\left(\bm w\right),\ldots,\D^pQ_j\left(\bm w\right)
 \quad\text{depend only on}\quad
 \bm v_j,\bm u_{j,1},\ldots,\bm u_{j,k+1},m,\lambda,q,\bm w.
\]
These are full ambient derivative tensors.
\end{claims}
\end{lemma}

\begin{proof}
\begin{itemize}
\item \textbf{Proof of \ref{collar:chain-lemma:part:residuals}.}
Define $\upsilon_{j,i}=\left(D_m\bm w\right)_i$. Consecutive residuals contain
$x_i$ with opposite signs, so
\[
 \begin{aligned}
 \frac{\partial Q_j}{\partial x_i}\left(\bm w\right)
 &=\lambda[\phi'\left(\upsilon_{j,i}\right)-\phi'\left(\upsilon_{j,i-1}\right)]
 &&\left(1\le i\le m\right),\\
 \frac{\partial Q_j}{\partial t}\left(\bm w\right)
 &=\frac{\lambda}{\sqrt m}\phi'\left(\upsilon_{j,0}\right).
 \end{aligned}
\]
Equal residuals therefore cancel every link-gradient component. Solving
for a common value $\tau_j$ gives
\[
 \begin{aligned}
 x_m=\tau_j,\quad x_i-x_{i+1}=\tau_j\ \left(1\le i<m\right)
 &\quad\Longrightarrow\quad x_i=\left(m-i+1\right)\tau_j\ \left(1\le i\le m\right),\\
 \frac{t}{\sqrt{m}}-x_1=\tau_j
 &\quad\Longrightarrow\quad t=\sqrt m\left(m+1\right)\tau_j=t_j^\star.
 \end{aligned}
\]
This is exactly $\bm h_j$ in Eq.~\eqref{collar:endpoint}, so
\[
 D_m\bm h_j=\tau_j\bm1,\qquad
 \D Q_j\left(\bm h_j\right)=\frac{\lambda}{\sqrt m}\phi'\left(\tau_j\right)\bm v_j
 =\frac{\lambda\tau_j^{q-1}}{\sqrt m}\bm v_j.
\]
To represent full ambient derivatives, define the row vectors by
\[
 \bm a_{j,0}=\frac{\bm v_j}{\sqrt{m}}-\bm u_{j,1},\quad
 \bm a_{j,i}=\bm u_{j,i}-\bm u_{j,i+1}\ \left(1\le i<m\right),\quad
 \bm a_{j,m}=\bm u_{j,m}.
\]
Then $\left(D_mP_{E_j}\bm x\right)_i=\ip{\bm a_{j,i}}{\bm x}$, and the link vectors
cancel consecutively:
\[
\begin{aligned}
 \sum_{i=0}^m\bm a_{j,i}
 &=\frac{\bm v_j}{\sqrt m}-\bm u_{j,1}
   +\sum_{i=1}^{m-1}\left(\bm u_{j,i}-\bm u_{j,i+1}\right)
   +\bm u_{j,m}\\
 &=\frac{\bm v_j}{\sqrt m}-\bm u_{j,1}
   +\left(\bm u_{j,1}-\bm u_{j,m}\right)+\bm u_{j,m}
 =\frac{\bm v_j}{\sqrt m}.
\end{aligned}
\]
For $m=1$, the middle sum is empty and the same cancellation holds.
This proves the second identity
in Eq.~\eqref{collar:residual-identities}. The triangular matrix $D_m$ has
nonzero diagonal, so it is invertible.

\item \textbf{Proof of \ref{collar:chain-lemma:part:norm}.}
The matrix factorization below reduces the norm bound to
Eq.~\eqref{sources:chain-matrix-bound} applied to $\mathsf U_{m+1}$.
Put $\mathsf S_m:=\operatorname{diag}\left(m^{-\frac{1}{2}},1,\ldots,1\right)$. Then
\[
 \begin{aligned}
 D_m&=\mathsf U_{m+1}\mathsf S_m,\qquad
 \norm{\mathsf S_m}_{2\to2}=1,\\
 \norm{D_m}_{2\to2}
 &\le\norm{\mathsf U_{m+1}}_{2\to2}\norm{\mathsf S_m}_{2\to2}\le2,\\
 \norm{D_m\bm w}_r&\le\norm{D_m\bm w}_2\le2\norm{\bm w}\quad\left(r\ge2\right).
 \end{aligned}
\]
The last line uses the comparison of $\ell_r$ and Euclidean norms.

\item \textbf{Proof of \ref{collar:chain-lemma:part:geometry}.}
The two orthogonal components of Eq.~\eqref{collar:endpoint} have squared lengths
\[
 \begin{aligned}
 \left(t_j^\star\right)^2&=m\left(m+1\right)^2\tau_j^2,\\
 \norm{\bm h_j-t_j^\star\bm v_j}^2
 &=\tau_j^2\sum_{a=1}^m a^2
 =\frac{m\left(m+1\right)\left(2m+1\right)}{6}\tau_j^2.
 \end{aligned}
\]
Hence
\[
 \begin{aligned}
 R_j^2
 &=\left(t_j^\star\right)^2+\norm{\bm h_j-t_j^\star\bm v_j}^2
  =\tau_j^2\left[m\left(m+1\right)^2+\frac{m\left(m+1\right)\left(2m+1\right)}{6}\right]\\
  &=\frac{\tau_j^2m\left(m+1\right)}{6}\left[6\left(m+1\right)+\left(2m+1\right)\right]\\
  &=\frac{m\left(m+1\right)\left(8m+7\right)}{6}\tau_j^2,\\
 \frac{R_j^2}{\tau_j^2m^3}
  &=\frac{\left(m+1\right)\left(8m+7\right)}{6m^2}
   =\frac{8m^2+15m+7}{6m^2}\\
  &=\frac{4}{3}+\frac{5}{2m}+\frac{7}{6m^2}\in[1,5],\\
 \frac{\norm{\bm h_j-t_j^\star\bm v_j}^2}{\left(t_j^\star\right)^2}
 &=\frac{2m+1}{6\left(m+1\right)}\in\left[\frac{1}{4},\frac{1}{3}\right].
 \end{aligned}
\]
Define $r=\frac{\norm{\bm h_j-t_j^\star\bm v_j}^2}{\left(t_j^\star\right)^2}\in\left[\frac{1}{4},\frac{1}{3}\right]$.
Since $R_j^2=\left(t_j^\star\right)^2\left(1+r\right)$, the two fractions are
\[
\begin{aligned}
 \frac{\left(t_j^\star\right)^2}{R_j^2}
 &=\frac{1}{1+r}\in\left[\frac{3}{4},\frac{4}{5}\right],\\
 \frac{\norm{\bm h_j-t_j^\star\bm v_j}^2}{R_j^2}
 &=\frac{r}{1+r}\in\left[\frac{1}{5},\frac{1}{4}\right].
\end{aligned}
\]
Taking square roots proves the length bounds. Also
$\sum_{i=1}^m\left(m-i+1\right)\tau_j=\frac{m\left(m+1\right)\tau_j}{2}\le m^2\tau_j$.

This explains the coefficient $m^{-\frac{1}{2}}$ in the first residual:
using $\beta t-x_1$ with $\beta>0$ instead would give $t=\frac{\left(m+1\right)\tau_j}{\beta}$ at equal
residuals; $\beta=m^{-\frac{1}{2}}$ makes the first coordinate and link vector
comparable in length. The first coordinate supplies $s_jt_j^\star$ in
Lemma~\ref{global:exact-drop-lemma} in Appx.~\ref{app:paper-assembly}; the unexposed links give the separation
proved in Lemma~\resultpart{collar:separation-lemma}{collar:separation-lemma:part:distance}.

\item \textbf{Proof of \ref{collar:chain-lemma:part:prefix}.}
For the stated $\bm w$, its coordinates satisfy
$x_{k+1}=\cdots=x_m=0$. Thus
\[
 \left(D_m\bm w\right)_i=0\quad\left(i>k\right),\qquad
 \phi^{\left(r\right)}\left(0\right)=0\quad\left(0\le r\le p\right),
\]
where $\phi^{\left(0\right)}=\phi$ and the second equality follows from
Lemma~\resultpart{elem:power}{elem:power:part:flatness} in Appx.~\ref{app:paper-elementary}.
Using the row vectors from part~\ref{collar:chain-lemma:part:residuals},
for $1\le r\le p$ and arbitrary $\bm d_1,\ldots,\bm d_r\in\R^d$ we obtain
\[
 \begin{aligned}
 Q_j\left(\bm w\right)&=\lambda\sum_{i=0}^k\phi\left(\left(D_m\bm w\right)_i\right),\\
 \D^rQ_j\left(\bm w\right)[\bm d_1,\ldots,\bm d_r]
 &=\lambda\sum_{i=0}^k\phi^{\left(r\right)}\left(\left(D_m\bm w\right)_i\right)
   \prod_{\ell=1}^r\ip{\bm a_{j,i}}{\bm d_\ell}.
 \end{aligned}
\]
Rows $0,\ldots,k$ use only
$\bm v_j,\bm u_{j,1},\ldots,\bm u_{j,k+1}$, proving the claim.
At $\bm w=0$ all residuals vanish, so $\D^rQ_j\left(0\right)=0$ for $1\le r\le p$.
Together these statements give the $p^{th}$-order zero-chain property in the
ordered coordinates of $E_j$, in the sense of
\cite[Definition~3]{CarmonEtAl2020}, restated as
Definition~\ref{sources:zero-chain}.
\end{itemize}
\end{proof}

\subsubsection{Joining the affine and chain formulas}
\label{collar:affine}
The affine formula supplies gradient $\bm g_j$ after completing the phase,
whereas the chain formula gives positive excess at early evaluation points.
We use the smooth maximum to retain each formula on its required region.

\begin{definition}[Affine and chain formulas]\label{collar:phase-definition}
For the phase-$j$ chain and endpoint from Appx.~\ref{collar:chain}, let
\[
 \lambda,\alpha,\delta_j>0,\qquad G_j,s_{j+1},t_{j+1}^\star\ge0.
\]
The local origin is $\bm c_j$ and $\bm h_j=\bm c_{j+1}-\bm c_j$.
For $j<L-1$, choose a unit $\bm v_{j+1}\perp E_j$; local inputs are
\[
 \left(\bm w,t_{j+1}\right)
 =\left(P_{E_j}\left(\bm x-\bm c_j\right),
        \ip{\bm v_{j+1}}{\bm x-\bm c_j}\right)
 \in E_j\times\R.
\]
The same symbols denote the ambient functions obtained by this substitution.
Here $s_{j+1},t_{j+1}^\star$ are the next phase's starting-gradient magnitude
and first endpoint coordinate. In the terminal phase set $s_L=t_L^\star=0$,
omit $\bm v_L$, and restrict to $t_{j+1}=0$.
Set
\[
 \bm x_j^*:=\left(\bm h_j,0\right),\qquad
 \Gamma_j:=\frac{\lambda\tau_j^{q-1}}{\sqrt m},\qquad
 \bm b_j:=\frac{G_j\bm h_j}{R_j},\qquad \mathcal E_j:=\lambda m\tau_j^q,
\]
\begin{equation}\label{collar:s}
 s_j:=\Gamma_j+\frac{G_jt_j^\star}{R_j},\qquad
 \bm g_j:=s_j\bm v_j-s_{j+1}\bm v_{j+1}-\bm b_j.
\end{equation}
Define
\begin{equation}\label{collar:P-definition}
 \begin{aligned}
 P_j\left(\bm w,t_{j+1}\right)&:=Q_j\left(\bm w\right)+\frac{\alpha}{q}|t_{j+1}|^q
 =\frac{\lambda}{q}\norm{D_m\bm w}_q^q+\frac{\alpha}{q}|t_{j+1}|^q,\\
 P_j\left(\bm h_j,0\right)&=\frac{\lambda}{q}\left(m+1\right)\tau_j^q,
 \end{aligned}
\end{equation}
\begin{equation}\label{collar:inner}
 I_j\left(\bm w,t_{j+1}\right):=\ip{\bm g_j}{\bm w-\bm h_j+t_{j+1}\bm v_{j+1}},
\end{equation}
\begin{equation}\label{collar:outer}
 O_j:=P_j-P_j\left(\bm h_j,0\right)-4\delta_j,\qquad \zeta_j:=O_j-I_j,
\end{equation}
\begin{equation}\label{collar:Q}
 A_j:=M_{\delta_j}\left(I_j,O_j\right)=I_j+\psi_{\delta_j}\left(\zeta_j\right).
\end{equation}
The base construction will supply $G_j=\mu R_j^{q-1}$.
\end{definition}
The term $\frac{\alpha|t_{j+1}|^q}{q}$ controls the linear term $s_{j+1}t_{j+1}$ in $\zeta_j$;
Corollary~\ref{collar:next-coordinate-lemma} quantifies this bound.

\begin{lemma}[The affine gradient in residual coordinates]
\label{collar:affine-lemma}
With Definition~\ref{collar:phase-definition}, put
$u_j:=\lambda\tau_j^{q-1}=\sqrt m\,\Gamma_j$.
\begin{claims}
\item\label{collar:affine-lemma:part:mismatch}\textbf{Gradient difference.}
\begin{equation}
 \bm e_j:=\bm g_j-\D P_j\left(\bm x_j^*\right)
 =-\frac{G_j}{R_j}\sum_{i=1}^m x_{j,i}^*\bm u_{j,i}-s_{j+1}\bm v_{j+1},
 \qquad x_{j,i}^*=\left(m-i+1\right)\tau_j.
 \label{collar:mismatch}
\end{equation}
\[
 I_j\left(\bm h_j,0\right)=0,\qquad
 I_j\left(\bm h_j,t_{j+1}^\star\right)=-s_{j+1}t_{j+1}^\star,\qquad
 \D I_j=\bm g_j.
\]
\item\label{collar:affine-lemma:part:residuals}\textbf{Residual coordinates.}
The coefficients
\begin{equation}
 \pi_{j,0}=u_j,\qquad
 \pi_{j,i}=u_j-\frac{G_j}{R_j}\sum_{k=1}^i x_{j,k}^*\quad\left(1\le i\le m\right)
 \label{collar:dual-coordinates}
\end{equation}
satisfy $D_m^\top\bm\pi_j=P_{E_j}\bm g_j$ in the coordinates of
Appx.~\ref{collar:chain}, and
\begin{equation}
 \begin{aligned}
 \zeta_j\left(\bm w,t_{j+1}\right)&=P_j\left(\bm w,t_{j+1}\right)-\ip{\bm\pi_j}{D_m\bm w}+s_{j+1}t_{j+1}+k_{j,\rm aff},\\
 k_{j,\rm aff}&:=\tau_j\sum_{i=0}^m\pi_{j,i}-\frac{\lambda}{q}\left(m+1\right)\tau_j^q-4\delta_j.
 \end{aligned}
 \label{collar:dual-identity}
\end{equation}
\item\label{collar:affine-lemma:part:bound}\textbf{Coefficient bounds.}
For any $0<c_0<1$, if $G_j\le c_0\Gamma_j$, then
\begin{equation}
 \begin{aligned}
 \left(1-c_0\right)u_j&\le\pi_{j,i}\le u_j\quad\left(0\le i\le m\right),\\
 k_{j,\rm aff}&\ge\left(m+1\right)\left(\frac{1}{q^*}-c_0\right)\lambda\tau_j^q-4\delta_j,\\
 k_{j,\rm aff}&>0\quad\text{under the scalar choice Eq.~\eqref{collar:constants} below}.
 \end{aligned}
 \label{collar:dual-bounds}
\end{equation}
\end{claims}
\end{lemma}
\begin{proof}
\begin{itemize}
\item \textbf{Proof of \ref{collar:affine-lemma:part:mismatch}.}
Lemma~\resultpart{collar:chain-lemma}{collar:chain-lemma:part:residuals}
gives $\D P_j\left(\bm x_j^*\right)=\Gamma_j\bm v_j$. Subtract it from $\bm g_j$:
\[
 \begin{aligned}
 \bm g_j-\D P_j\left(\bm x_j^*\right)
 &=\underbrace{\left(s_j-\Gamma_j-\frac{G_jt_j^\star}{R_j}\right)}_{0\text{ by Eq.~\eqref{collar:s}}}\bm v_j\\
 &\quad-\frac{G_j}{R_j}\sum_{i=1}^m x_{j,i}^*\bm u_{j,i}
       -s_{j+1}\bm v_{j+1}.
 \end{aligned}
\]
Thus the choice of $s_j$ removes the $\bm v_j$-component, proving
Eq.~\eqref{collar:mismatch}. For $j<L-1$, the affine definition
Eq.~\eqref{collar:inner} also gives
\[
\begin{aligned}
 I_j\left(\bm h_j,t_{j+1}\right)
 &=\ip{\bm g_j}{\bm h_j-\bm h_j+t_{j+1}\bm v_{j+1}}\\
 &=t_{j+1}\ip{\bm g_j}{\bm v_{j+1}}=-s_{j+1}t_{j+1},\\
 I_j\left(\bm h_j,0\right)&=0,\qquad
 I_j\left(\bm h_j,t_{j+1}^\star\right)=-s_{j+1}t_{j+1}^\star,\qquad
 \D I_j=\bm g_j.
\end{aligned}
\]
The third equality uses $\bm v_{j+1}\perp E_j$ and
$\bm b_j,\bm v_j\in E_j$.
In the terminal phase $t_{j+1}=0$ and the next-direction term is omitted;
the affine definition directly gives $I_j\left(\bm h_j,0\right)=0$ and
$\D I_j=\bm g_j$.
\item \textbf{Proof of \ref{collar:affine-lemma:part:residuals}.}
Using the columns of $D_m$ and the gradient in
part~\ref{collar:affine-lemma:part:mismatch}, the equation
$D_m^\top\bm\pi_j=P_{E_j}\bm g_j$ is
\[
 \frac{\pi_{j,0}}{\sqrt m}=\Gamma_j,\qquad
 \pi_{j,i}-\pi_{j,i-1}=-\frac{G_jx_{j,i}^*}{R_j}\quad\left(1\le i\le m\right).
\]
The first equation gives $\pi_{j,0}=\sqrt m\,\Gamma_j=u_j$.
Telescoping the next $i$ equations gives
\[
\begin{aligned}
 \pi_{j,i}-\pi_{j,0}
 &=\sum_{k=1}^i\left(\pi_{j,k}-\pi_{j,k-1}\right)
 =-\frac{G_j}{R_j}\sum_{k=1}^i x_{j,k}^*,\\
 \pi_{j,i}
 &=\pi_{j,0}-\frac{G_j}{R_j}\sum_{k=1}^i x_{j,k}^*
 =u_j-\frac{G_j}{R_j}\sum_{k=1}^i x_{j,k}^*.
\end{aligned}
\]
This is Eq.~\eqref{collar:dual-coordinates}.
Since $D_m\bm h_j=\tau_j\bm1$, the affine and chain formulas become
\[
\begin{aligned}
 I_j\left(\bm w,t_{j+1}\right)
 &=\ip{\bm\pi_j}{D_m\left(\bm w-\bm h_j\right)}-s_{j+1}t_{j+1}\\
 &=\ip{\bm\pi_j}{D_m\bm w}-\tau_j\sum_{i=0}^m\pi_{j,i}-s_{j+1}t_{j+1},\\
 O_j\left(\bm w,t_{j+1}\right)&=P_j\left(\bm w,t_{j+1}\right)-\frac{\lambda}{q}\left(m+1\right)\tau_j^q-4\delta_j.
\end{aligned}
\]
Subtracting the affine formula from the chain formula and distributing the negative sign gives
\[
\begin{aligned}
 \zeta_j\left(\bm w,t_{j+1}\right)
 &=P_j\left(\bm w,t_{j+1}\right)-\frac{\lambda}{q}\left(m+1\right)\tau_j^q-4\delta_j\\
 &\quad-\left(
      \ip{\bm\pi_j}{D_m\bm w}
      -\tau_j\sum_{i=0}^m\pi_{j,i}-s_{j+1}t_{j+1}\right)\\
 &=P_j\left(\bm w,t_{j+1}\right)-\ip{\bm\pi_j}{D_m\bm w}+s_{j+1}t_{j+1}\\
 &\quad+\tau_j\sum_{i=0}^m\pi_{j,i}
      -\frac{\lambda}{q}\left(m+1\right)\tau_j^q-4\delta_j\\
 &=P_j\left(\bm w,t_{j+1}\right)-\ip{\bm\pi_j}{D_m\bm w}
      +s_{j+1}t_{j+1}+k_{j,\rm aff}.
\end{aligned}
\]
The last two constant terms and the coefficient sum give precisely
$k_{j,\rm aff}$ in Eq.~\eqref{collar:dual-identity}.
\item \textbf{Proof of \ref{collar:affine-lemma:part:bound}.}
Lemma~\resultpart{collar:chain-lemma}{collar:chain-lemma:part:geometry} gives
$\sum_kx_{j,k}^*\le\tau_j m^2$ and $R_j\ge\tau_j m^{\frac{3}{2}}$. Hence
\[
\begin{aligned}
 0\le u_j-\pi_{j,i}
 &=\frac{G_j}{R_j}\sum_{k=1}^i x_{j,k}^*\\
 &\le\frac{G_j}{\tau_jm^{\frac{3}{2}}}\tau_jm^2
 =G_j\sqrt m\\
 &\le c_0\Gamma_j\sqrt m=c_0u_j.
\end{aligned}
\]
Insert $\left(1-c_0\right)u_j\le\pi_{j,i}\le u_j$ into the constant formula:
\[
\begin{aligned}
 k_{j,\rm aff}
 &=\tau_j\sum_{i=0}^m\pi_{j,i}
      -\frac{\lambda}{q}\left(m+1\right)\tau_j^q-4\delta_j\\
 &\ge\left(m+1\right)\left(1-c_0\right)\tau_ju_j
      -\frac{\lambda}{q}\left(m+1\right)\tau_j^q-4\delta_j\\
 &=\left(m+1\right)\left(1-c_0-\frac{1}{q}\right)\lambda\tau_j^q-4\delta_j\\
 &=\left(m+1\right)\left(\frac{1}{q^*}-c_0\right)\lambda\tau_j^q-4\delta_j.
\end{aligned}
\]
The third line uses $u_j=\lambda\tau_j^{q-1}$, and the last uses
$1-\frac{1}{q}=\frac{1}{q^*}$.
Under the later choice Eq.~\eqref{collar:constants}, $c_0=c_1<\frac{1}{64}$,
$\delta_j=c_1\mathcal E_j$, and $\frac{1}{q^*}>\frac{1}{2}$. Consequently
\[
\begin{aligned}
 k_{j,\rm aff}
 &\ge\left(m+1\right)\left(\frac{1}{q^*}-c_0\right)\lambda\tau_j^q-4c_1\mathcal E_j\\
 &\ge\left(\frac{1}{q^*}-c_0\right)\mathcal E_j-4c_1\mathcal E_j\\
 &=\mathcal E_j\left(\frac{1}{q^*}-c_0-4c_1\right)\\
 &>\mathcal E_j\left(\frac{1}{2}-\frac{5}{64}\right)>0.
\end{aligned}
\]
Here $\frac{1}{q^*}-c_0>0$ permits $m+1\ge m$, and the last line uses
$c_0=c_1<\frac{1}{64}$.
\end{itemize}
\end{proof}

\subsubsection{Value differences and the next-coordinate bound}
\label{collar:result}\label{sources:conjugacy}
\begin{lemma}[Fenchel--Young inequality {\cite[Proposition~2.3.1]{BorweinVanderwerff2010}}]
\label{sources:fenchel-young}
Let $f:\R^n\to(-\infty,+\infty]$, let
$\bm a\in\operatorname{dom}f$, and let $\bm b\in\R^n$.
Define
\[
 f^*\left(\bm b\right)=\sup_{\bm z\in\R^n}
             \{\ip{\bm b}{\bm z}-f\left(\bm z\right)\}.
\]
Then
\[
 f\left(\bm a\right)+f^*\left(\bm b\right)\ge\ip{\bm b}{\bm a},
 \qquad
 f\left(\bm a\right)+f^*\left(\bm b\right)=\ip{\bm b}{\bm a}
 \quad\Longleftrightarrow\quad\bm b\in\partial f\left(\bm a\right).
\]
Here
\[
 \operatorname{dom}f=\{\bm z:f\left(\bm z\right)<+\infty\},\qquad
 \partial f\left(\bm a\right)=
 \{\bm b:f\left(\bm z\right)\ge f\left(\bm a\right)+\ip{\bm b}{\bm z-\bm a}
                  \text{ for every }\bm z\in\R^n\}.
\]
The proposition does not require convexity or lower semicontinuity of $f$.
\end{lemma}

\begin{lemma}[Power conjugates and positive scaling {\cite[Tables~2.1--2.2, pp.~45--46]{BorweinVanderwerff2010}}]
\label{sources:power-conjugates}
For $q>1$ and $q^*=\frac{q}{q-1}$, the conjugate pair is
\[
 \phi\left(a\right)=\frac{|a|^q}{q},\qquad
 \phi^*\left(b\right)=\frac{|b|^{q^*}}{q^*},
 \qquad \operatorname{dom}\phi=\operatorname{dom}\phi^*=\R.
\]
For $h:\R\to(-\infty,+\infty]$ with conjugate
$h^*\left(b\right)=\sup_{a\in\R}\{ba-h\left(a\right)\}$ and $\ell>0$,
\[
 \left(\ell h\right)^*\left(b\right)=\ell h^*\left(\frac{b}{\ell}\right).
\]
\end{lemma}

\begin{lemma}[Scaled scalar and separable power inequalities]
\label{elem:power-young}
Let $q>1$, $q^*=\frac{q}{q-1}$, and $\ell>0$. Then
\begin{equation}\label{sources:scalar-power-young}
 ba\le\frac{\ell}{q}|a|^q
       +\frac{\ell^{-\frac{1}{q-1}}}{q^*}|b|^{q^*}
 \qquad\left(a,b\in\R\right).
\end{equation}
Consequently, for any finite collection of real $a_i,b_i$ and
positive $\ell_i$,
\begin{equation}\label{sources:separable-power-young}
 \sum_i b_i a_i\le\sum_i\frac{\ell_i}{q}|a_i|^q
 +\frac{1}{q^*}\sum_i\ell_i^{-\frac{1}{q-1}}|b_i|^{q^*}.
\end{equation}
\end{lemma}
\begin{proof}
Put $\phi\left(a\right)=\frac{|a|^q}{q}$ and $h=\ell\phi$.
The conjugate pair and positive-scaling rule in
Lemma~\ref{sources:power-conjugates} give
\[
 \begin{aligned}
 h^*\left(b\right)
 &=\ell\phi^*\left(\frac{b}{\ell}\right)
 &&\text{by the positive-scaling rule},\\
  &=\frac{\ell}{q^*}\left|\frac{b}{\ell}\right|^{q^*}
 &&\text{by the conjugate of }\phi,\\
 &=\frac{\ell^{1-q^*}}{q^*}|b|^{q^*}\\
 &=\frac{\ell^{-\frac{1}{q-1}}}{q^*}|b|^{q^*}
 &&\text{since }1-q^*=-\frac{1}{q-1}.
 \end{aligned}
\]
Apply Lemma~\ref{sources:fenchel-young} with $f=h$:
\[
 ba\le h\left(a\right)+h^*\left(b\right)
 =\frac{\ell}{q}|a|^q+
   \frac{\ell^{-\frac{1}{q-1}}}{q^*}|b|^{q^*}.
\]
This is Eq.~\eqref{sources:scalar-power-young}.
Applying it with $\left(a,b,\ell\right)=\left(a_i,b_i,\ell_i\right)$ for each index
and summing gives Eq.~\eqref{sources:separable-power-young}.
\end{proof}

To determine which function equals $A_j$, Lemma~\ref{paper:smooth-maximum-properties} gives
\[
 \zeta_j>\delta_j\ \Longrightarrow\ A_j=O_j,\qquad
 \zeta_j<-\delta_j\ \Longrightarrow\ A_j=I_j.
\]

\begin{lemma}[Constants giving the required value differences]
\label{collar:margin-lemma}
Use Definition~\ref{collar:phase-definition}. Choose
\begin{equation}\label{collar:constants}
 c_0=c_1=\frac{2^{-q-3}}{q},\qquad \delta_j=c_1\mathcal E_j,
\end{equation}
and assume
\begin{equation}\label{collar:G-assumption}
 G_j\le c_0\Gamma_j,
\end{equation}
\begin{equation}\label{collar:exit-port}
 s_{j+1}t_{j+1}^\star+\frac{\alpha}{q}\left(t_{j+1}^\star\right)^q\le c_0\mathcal E_j.
\end{equation}
Use the chain-prefix subspace $U_{j,b}$ from Step~2.1, where the integer
$b$ counts links and the coordinate $\ip{\bm v_j}{\bm w}$ is unrestricted.
\begin{claims}
\item\label{collar:margin-lemma:part:prefix}\textbf{Value difference on chain prefixes.}
\[
 \zeta_j\left(\bm w,0\right)\ge7\delta_j
 \qquad\left(\bm w\in U_{j,b},\ 0\le b\le\lfloor \frac{m}{4}\rfloor\right).
\]
\item\label{collar:margin-lemma:part:endpoints}\textbf{Value differences at the endpoints.}
\[
 \zeta_j\left(\bm h_j,0\right)=-4\delta_j,\qquad
 \zeta_j\left(\bm h_j,t_{j+1}^\star\right)\le-3\delta_j.
\]
\end{claims}
\end{lemma}
\begin{proof}
\begin{itemize}
\item \textbf{Proof of \ref{collar:margin-lemma:part:prefix}.}
Lemma~\ref{sources:power-conjugates} gives the scalar conjugate formula
\begin{equation}
 \sup_{a\in\R}\left\{va-\frac{\ell}{q}|a|^q\right\}
 =\frac{\ell^{-\frac{1}{q-1}}}{q^*}|v|^{q^*}
 \qquad\left(\ell>0,\ v\in\R\right).
 \label{collar:power-conjugate}
\end{equation}
The chosen constants satisfy
\[
 2c_0+11c_1=13c_1<\frac{13}{64}<\frac{3}{8}<\frac{3}{4q^*},
\]
since $c_1<\frac{1}{64}$ and $q^*<2$. Fix $\bm w\in U_{j,b}$, $0\le b\le\lfloor \frac{m}{4}\rfloor$, and put $\bm a=D_m\bm w$.
Then $a_i=0$ for $i>b$. For $0\le i\le b$, applying
Eq.~\eqref{sources:scalar-power-young} with coefficient $\pi_{j,i}$
and weight $\lambda$ gives
\[
 \begin{aligned}
 \frac{\lambda}{q}|a_i|^q-\pi_{j,i}a_i
 &\ge-\frac{\lambda^{-\frac{1}{q-1}}}{q^*}\pi_{j,i}^{q^*}\\
 &\ge-\frac{\lambda^{-\frac{1}{q-1}}}{q^*}
        \left(\lambda\tau_j^{q-1}\right)^{q^*}
 =-\frac{\lambda\tau_j^q}{q^*}.
 \end{aligned}
\]
Here $0<\pi_{j,i}\le\lambda\tau_j^{q-1}$ by
Lemma~\resultpart{collar:affine-lemma}{collar:affine-lemma:part:bound}.
The power simplification on the last line is
\[
 \lambda^{-\frac{1}{q-1}}\left(\lambda\tau_j^{q-1}\right)^{q^*}
 =\lambda^{q^*-\frac{1}{q-1}}\tau_j^{\left(q-1\right)q^*}
 =\lambda\tau_j^q,
\]
since $q^*-\frac{1}{q-1}=1$ and $\left(q-1\right)q^*=q$.
Thus Eq.~\eqref{collar:dual-identity} and Eq.~\eqref{collar:dual-bounds} give
\begin{align*}
 \zeta_j\left(\bm w,0\right)
 &=k_{j,\rm aff}+\sum_{i=0}^b
     \left(\frac{\lambda}{q}|a_i|^q-\pi_{j,i}a_i\right)\\
 &\ge k_{j,\rm aff}-\frac{b+1}{q^*}\lambda\tau_j^q\\
 &\ge\left[
    \underbrace{\frac{m+1}{q^*}}_{\text{all rows}}
   -\underbrace{\frac{b+1}{q^*}}_{\text{rows }0,\ldots,b}
   -\underbrace{c_0\left(m+1\right)}_{\text{base error}}
   \right]\lambda\tau_j^q-4\delta_j\\
 &=\left(\frac{m-b}{q^*}-c_0\left(m+1\right)\right)\lambda\tau_j^q-4\delta_j\\
 &\ge\left(\frac{3}{4q^*}-2c_0-4c_1\right)\mathcal E_j
 \qquad\left(\substack{b\le \frac{m}{4},\ m+1\le2m\\\delta_j=c_1\mathcal E_j}\right)\\
 &\ge7c_1\mathcal E_j=7\delta_j
 \qquad\left(2c_0+11c_1\le\frac{3}{4q^*}\right).
\end{align*}
Hence
\begin{equation}
 \zeta_j\left(\bm w,0\right)\ge7\delta_j
 \qquad\left(\bm w\in U_{j,b},\ b\le\lfloor \frac{m}{4}\rfloor\right).
 \label{collar:prefix-margin}
\end{equation}

\item \textbf{Proof of \ref{collar:margin-lemma:part:endpoints}.}
By Eq.~\eqref{collar:outer} and
Lemma~\resultpart{collar:affine-lemma}{collar:affine-lemma:part:mismatch},
for $t_{j+1}\ge0$ in the nonterminal phase,
\[
\begin{aligned}
 O_j\left(\bm h_j,t_{j+1}\right)
 &=P_j\left(\bm h_j,0\right)+\frac{\alpha}{q}t_{j+1}^q
      -P_j\left(\bm h_j,0\right)-4\delta_j
 =\frac{\alpha}{q}t_{j+1}^q-4\delta_j,\\
 I_j\left(\bm h_j,t_{j+1}\right)&=-s_{j+1}t_{j+1},\\
 \zeta_j\left(\bm h_j,t_{j+1}\right)
 &=O_j\left(\bm h_j,t_{j+1}\right)-I_j\left(\bm h_j,t_{j+1}\right)\\
 &=\frac{\alpha}{q}t_{j+1}^q+s_{j+1}t_{j+1}-4\delta_j.
\end{aligned}
\]
At $t_{j+1}=0$, the same identity holds in the terminal phase with the
next-coordinate terms omitted.
The definitions and $c_0\mathcal E_j=\delta_j$ therefore give
\[
 \begin{aligned}
 \zeta_j\left(\bm h_j,0\right)&=-4\delta_j,\\
 \zeta_j\left(\bm h_j,t_{j+1}^\star\right)
 &=s_{j+1}t_{j+1}^\star+\frac{\alpha}{q}\left(t_{j+1}^\star\right)^q-4\delta_j
 \le-3\delta_j,
 \end{aligned}
\]
using Eq.~\eqref{collar:exit-port}. The subtraction $4\delta_j$ thus permits a
$\delta_j$ increase while keeping $\zeta_j<-\delta_j$.
For $j<L-1$, the assembly in Appx.~\ref{app:paper-assembly}
identifies the second input with the following center:
\[
 \begin{aligned}
 P_{E_j}\left(\bm c_{j+2}-\bm c_j\right)
 &=P_{E_j}\left(\bm h_j+\bm h_{j+1}\right)=\bm h_j,\\
 \ip{\bm v_{j+1}}{\bm c_{j+2}-\bm c_j}
 &=\ip{\bm v_{j+1}}{\bm h_{j+1}}=t_{j+1}^\star.
 \end{aligned}
\]
Keeping the same affine formula there leaves the completed phase's gradient unchanged.
\end{itemize}
\end{proof}

\begin{corollary}[Fenchel--Young bound for the next coordinate]
\label{collar:next-coordinate-lemma}
For $\alpha>0$ and $s_{j+1}\ge0$,
\[
 \inf_{t_{j+1}\in\R}\left\{\frac{\alpha}{q}|t_{j+1}|^q+s_{j+1}t_{j+1}\right\}
 =-\frac{\alpha^{-\frac{1}{q-1}}s_{j+1}^{q^*}}{q^*}.
\]
Consequently, if
\begin{equation}\label{collar:exit-dual}
 \alpha^{-\frac{1}{q-1}}s_{j+1}^{q^*}\le c_0\mathcal E_j,
\end{equation}
then
\[
 \frac{\alpha}{q}|t_{j+1}|^q+s_{j+1}t_{j+1}\ge-\frac{c_0\mathcal E_j}{q^*}
 \qquad\left(t_{j+1}\in\R\right).
\]
\end{corollary}
\begin{proof}
Apply the power-conjugate formula Eq.~\eqref{collar:power-conjugate} with
$\left(a,v,\ell\right)=\left(t_{j+1},-s_{j+1},\alpha\right)$:
\[
 \inf_{t_{j+1}}\left\{\frac{\alpha}{q}|t_{j+1}|^q+s_{j+1}t_{j+1}\right\}
 =-\sup_{t_{j+1}}\left\{-s_{j+1}t_{j+1}-\frac{\alpha}{q}|t_{j+1}|^q\right\}
 =-\frac{\alpha^{-\frac{1}{q-1}}s_{j+1}^{q^*}}{q^*}.
\]
For each $t_{j+1}\in\R$, Eq.~\eqref{collar:exit-dual} therefore implies
\[
\begin{aligned}
 \frac{\alpha}{q}|t_{j+1}|^q+s_{j+1}t_{j+1}
 &\ge\inf_{a\in\R}\left\{\frac{\alpha}{q}|a|^q+s_{j+1}a\right\}\\
 &=-\frac{\alpha^{-\frac{1}{q-1}}s_{j+1}^{q^*}}{q^*}
 \ge-\frac{c_0\mathcal E_j}{q^*}.
\end{aligned}
\]
\end{proof}

Set
\begin{equation}\label{collar:alpha}
 \alpha:=\frac{\mathcal E_j}{\left(\tau_jm^{\frac{3}{2}}\right)^q}
 =\lambda m^{1-\frac{3q}{2}}.
\end{equation}
Lemma~\ref{collar:exit-lemma} in Appx.~\ref{app:paper-assembly}
verifies both next-coordinate bounds with this choice and the base's ratio $\eta$.

\subsubsection{Verify the phase function}
\begin{theorem}[Phase function with $A_j-I_j\ge7\delta_j$ on short chain prefixes]
\label{collar:theorem}
Fix a phase $0\le j<L$, $p\ge2$, $0<\nu\le1$, $q=p+\nu$, and an orthonormal block
as in Eq.~\eqref{collar:map}--\eqref{collar:endpoint}. Take
\[
 m\ge1,\qquad \lambda,\tau_j>0,\qquad G_j,s_{j+1},t_{j+1}^\star\ge0,
\]
with $m,p$ integers. For $j<L-1$, take a unit $\bm v_{j+1}\perp E_j$;
part~\ref{collar:theorem:part:terminal} specifies the terminal convention.
Use $c_0,c_1$ from Eq.~\eqref{collar:constants}, and set
\[
\begin{gathered}
 \Gamma_j=\frac{\lambda\tau_j^{q-1}}{\sqrt{m}},\qquad
 \bm b_j=\frac{G_j\bm h_j}{R_j},\qquad s_j=\Gamma_j+\frac{G_jt_j^\star}{R_j},\\
 \alpha=\lambda m^{1-\frac{3q}{2}},\qquad \delta_j=c_1\mathcal E_j,
 \qquad\mathcal E_j=\lambda m\tau_j^q,\\
 \bm g_j=s_j\bm v_j-s_{j+1}\bm v_{j+1}-\bm b_j.
\end{gathered}
\]
Assume
\[
 \begin{aligned}
 G_j&\le c_0\Gamma_j,\\
 s_{j+1}t_{j+1}^\star+\frac{\alpha\left(t_{j+1}^\star\right)^q}{q}&\le c_0\mathcal E_j,
 &&\text{Eq.~\eqref{collar:exit-port}},\\
 \alpha^{-\frac{1}{q-1}}s_{j+1}^{q^*}&\le c_0\mathcal E_j,
 &&\text{Eq.~\eqref{collar:exit-dual}}.
 \end{aligned}
\]
For $I_j,O_j,A_j$ from Eq.~\eqref{collar:inner}, Eq.~\eqref{collar:outer}, Eq.~\eqref{collar:Q},
and $U_{j,b}=\operatorname{span}\{\bm v_j,\bm u_{j,1},\ldots,\bm u_{j,b}\}$:
\begin{claims}
\item\label{collar:theorem:part:convexity}\textbf{Convexity.}
$A_j$ is convex and $A_j\ge\max\{I_j,O_j\}$.
\item\label{collar:theorem:part:prefix}\textbf{Chain-prefix excess and oracle output.}
For every integer $0\le b\le\lfloor \frac{m}{4}\rfloor$ and every $\bm w\in U_{j,b}$,
without a bound on $\norm{\bm w}$,
\[
 O_j\left(\bm w,0\right)-I_j\left(\bm w,0\right)\ge7\delta_j,\qquad
 A_j=O_j\quad\text{near }\left(\bm w,0\right).
\]
The value and full ambient derivatives through order $p$, given by
Eq.~\eqref{collar:prefix-derivatives}--\eqref{collar:prefix-value}, depend only on
$\bm v_j,\bm u_{j,1},\ldots,\bm u_{j,b+1}$ and fixed scalars.
\item\label{collar:theorem:part:endpoints}\textbf{Affine region.}
$A_j=I_j$ near $\{\left(\bm h_j,t_{j+1}\right):0\le t_{j+1}\le t_{j+1}^\star\}$, and
\begin{align}
A_j\left(\bm h_j,0\right)&=0,& A_j\left(\bm h_j,t_{j+1}^\star\right)&=-s_{j+1}t_{j+1}^\star,\label{collar:port-values}\\
\D A_j\left(\bm h_j,0\right)&=\bm g_j,& \D A_j\left(\bm h_j,t_{j+1}^\star\right)&=\bm g_j,\label{collar:port-gradients}\\
\D^k A_j\left(\bm h_j,0\right)&=0,&
\D^k A_j\left(\bm h_j,t_{j+1}^\star\right)&=0\quad\left(2\le k\le p\right).\label{collar:port-higher}
\end{align}
\item\label{collar:theorem:part:origin}\textbf{Starting derivatives and gradient change.}
\begin{equation}
A_j\left(0,0\right)=-\frac{\lambda}{q}\left(m+1\right)\tau_j^q-4\delta_j,
\qquad \D^k A_j\left(0,0\right)=0\quad\left(1\le k\le p\right).
\label{collar:origin}
\end{equation}
Consequently,
\begin{equation}
\bm b_j+\left(\D A_j\left(\bm h_j,0\right)-\D A_j\left(0,0\right)\right)
=\bm b_j+\bm g_j=s_j\bm v_j-s_{j+1}\bm v_{j+1}.
\label{collar:source}
\end{equation}
\item\label{collar:theorem:part:regularity}\textbf{Smoothness.}
$A_j\in C^{p,\nu}$ and, for all $\bm x,\bm y\in\R^d$,
\[
 \norm{\D^pA_j\left(\bm x\right)-\D^pA_j\left(\bm y\right)}
 \le C_A\lambda\norm{\bm x-\bm y}^\nu,
\]
where $C_A$ in Eq.~\eqref{collar:Q-holder} depends only on $p,\nu$ and the fixed smoothing function.
\item\label{collar:theorem:part:terminal}\textbf{Terminal phase.}
For $j=L-1$, parts~\ref{collar:theorem:part:convexity}--\ref{collar:theorem:part:regularity}
hold with $s_L=t_L^\star=0$, $\bm v_L$ omitted, and the functions restricted to $t_{j+1}=0$.
\end{claims}
\end{theorem}
\begin{proof}
Define $\mathcal E_j=\lambda m\tau_j^q$ and $\zeta_j=O_j-I_j$.
\begin{itemize}
\item \textbf{Proof of \ref{collar:theorem:part:convexity}.}
The scalar function $\phi\left(a\right)=\frac{|a|^q}{q}$ is convex and $C^p$ by
Lemma~\resultpart{elem:power}{elem:power:part:flatness}, \ref{elem:power:part:hessian} in Appx.~\ref{app:paper-elementary}.
Sums and affine substitutions preserve convexity
\cite[Secs.~3.2.1--3.2.2]{BoydVandenberghe2004}
(Lemma~\resultpart{sources:convex-rules}{sources:convex-rules:sum}--\ref{sources:convex-rules:affine}),
so $P_j,O_j,\zeta_j$ are convex and $C^p$.
Since $I_j$ is affine, Lemma~\resultpart{elem:composition-lemma}{elem:composition-lemma:part:convexity}
gives $A_j=M_{\delta_j}\left(I_j,O_j\right)$ convex and $C^p$;
Lemma~\ref{paper:smooth-maximum-properties} gives $A_j\ge I_j,O_j$.

\item \textbf{Proof of \ref{collar:theorem:part:prefix}.}
The affine gradient $\bm g_j$ can involve unrevealed directions,
so the oracle output at every $\bm w\in U_{j,b}$ with
$b\le\lfloor\frac{m}{4}\rfloor$ must come from $O_j$.
Fix $\bm w\in U_{j,b}$, $0\le b\le\lfloor \frac{m}{4}\rfloor$.
The theorem's choices of $c_0,c_1,\delta_j$ and its base and endpoint assumptions
are precisely those of Lemma~\ref{collar:margin-lemma}. Its
part~\ref{collar:margin-lemma:part:prefix} gives
\[
 \zeta_j\left(\bm w,0\right)\ge7\delta_j
 \qquad\left(\bm w\in U_{j,b},\ 0\le b\le\lfloor \frac{m}{4}\rfloor\right).
\]
Continuity and $7\delta_j>\delta_j$ give $A_j=O_j$ near $\left(\bm w,0\right)$ by
Lemma~\resultpart{elem:composition-lemma}{elem:composition-lemma:part:exact};
hence $A_j\left(\bm w,0\right)-I_j\left(\bm w,0\right)\ge7\delta_j$.

\textbf{The full oracle output.}
Apply Lemma~\resultpart{collar:chain-lemma}{collar:chain-lemma:part:prefix}
to $Q_j$ with prefix index $k=b$, which is permitted because
$b\le\lfloor \frac{m}{4}\rfloor<m$. The added term $\alpha\phi\left(t_{j+1}\right)$ and its
derivatives through order $p$ vanish at $t_{j+1}=0$ by
Lemma~\resultpart{elem:power}{elem:power:part:flatness}.
Thus $A_j=O_j$ near the evaluation point and the formulas in
Lemma~\resultpart{collar:chain-lemma}{collar:chain-lemma:part:prefix} give,
for $1\le k\le p$ and arbitrary ambient $\bm d_1,\ldots,\bm d_k$,
\begin{equation}
\D^kA_j\left(\bm w,0\right)[\bm d_1,\ldots,\bm d_k]
=\lambda\sum_{i=0}^{b}
\phi^{\left(k\right)}\left(\left(D_m\bm w\right)_i\right)
\prod_{\ell=1}^k\ip{\bm a_{j,i}}{\bm d_\ell}.
\label{collar:prefix-derivatives}
\end{equation}
The value is
\begin{equation}
A_j\left(\bm w,0\right)=\frac{\lambda}{q}\sum_{i=0}^b|\left(D_m\bm w\right)_i|^q
-\frac{\lambda}{q}\left(m+1\right)\tau_j^q-4\delta_j.
\label{collar:prefix-value}
\end{equation}
The row vectors $\bm a_{j,0},\ldots,\bm a_{j,b}$ involve only
$\bm v_j,\bm u_{j,1},\ldots,\bm u_{j,b+1}$, so they belong to
$U_{j,b+1}$. For the ambient functions in
Eq.~\eqref{paper:ambient-phase}, the point $\bm c_j+\bm w$
has local inputs $\left(\bm w,0\right)$. Orthogonal projection
therefore leaves every inner product in
Eq.~\eqref{collar:prefix-derivatives} unchanged, giving
\[
 \begin{aligned}
 \D^kA_j\left(\bm c_j+\bm w\right)[\bm d_1,\ldots,\bm d_k]
 &=\lambda\sum_{i=0}^{b}
   \phi^{\left(k\right)}\left(\left(D_m\bm w\right)_i\right)
   \prod_{\ell=1}^k
       \ip{\bm a_{j,i}}{P_{U_{j,b+1}}\bm d_\ell}\\
 &=\D^kQ_j\left(\bm w\right)
   [P_{U_{j,b+1}}\bm d_1,\ldots,P_{U_{j,b+1}}\bm d_k]
 &&\text{by Lemma~\resultpart{collar:chain-lemma}{collar:chain-lemma:part:prefix}}.
 \end{aligned}
\]
Thus the value and full ambient tensors agree for all remaining completions.
This is pointwise agreement at the evaluation point; completions need not agree nearby.

\item \textbf{Proof of \ref{collar:theorem:part:endpoints}.}
After the chain is complete, $A_j$ must equal $I_j$ near the
endpoint segment to supply the fixed gradient $\bm g_j$.
Lemma~\resultpart{collar:margin-lemma}{collar:margin-lemma:part:endpoints} gives
$\zeta_j\left(\bm h_j,0\right)=-4\delta_j$ and
$\zeta_j\left(\bm h_j,t_{j+1}^\star\right)\le-3\delta_j$.
The function $\zeta_j=P_j-I_j-P_j\left(\bm h_j,0\right)-4\delta_j$
is convex by Lemma~\resultpart{sources:convex-rules}{sources:convex-rules:sum}.
Its value at an intermediate point is therefore at most the
weighted average of the endpoint values, hence at most
$-3\delta_j$, as proved in
Lemma~\resultpart{collar:branches-lemma}{collar:branches-lemma:part:endpoints}.
The segment therefore lies in the open set $\{\zeta_j<-\delta_j\}$, where
Lemma~\resultpart{elem:composition-lemma}{elem:composition-lemma:part:exact}
gives $A_j=I_j$. Evaluating and differentiating
Eq.~\eqref{collar:inner} on this region gives
\[
 I_j\left(\bm h_j,t_{j+1}\right)=-s_{j+1}t_{j+1},\qquad \D I_j=\bm g_j,\qquad
 \D^kI_j=0\quad\left(2\le k\le p\right),
\]
which gives Eq.~\eqref{collar:port-values}--\eqref{collar:port-higher}.

\item \textbf{Proof of \ref{collar:theorem:part:origin}.}
At $\left(0,0\right)$, part~\ref{collar:theorem:part:prefix}
gives $A_j=O_j$ on a neighborhood. Evaluating
Eq.~\eqref{collar:outer} with the power terms and their derivatives
zero gives the starting data. More explicitly,
$D_m\bm h_j=\tau_j\bm1$, so
\[
 \begin{aligned}
 A_j\left(0,0\right)
 &=P_j\left(0,0\right)-P_j\left(\bm h_j,0\right)-4\delta_j\\
 &=-\frac{\lambda}{q}\left(m+1\right)\tau_j^q-4\delta_j,\\
 \D^kA_j\left(0,0\right)
 &=\D^kP_j\left(0,0\right)=0
 &&\text{by Lemma~\resultpart{elem:power}{elem:power:part:flatness}, for }1\le k\le p.
 \end{aligned}
\]
Using the endpoint gradient from part~\ref{collar:theorem:part:endpoints},
\[
 \begin{aligned}
 \bm b_j+\left(\D A_j\left(\bm h_j,0\right)-\D A_j\left(0,0\right)\right)
 &=\bm b_j+\bm g_j\\
 &=\bm b_j+\left(s_j\bm v_j-s_{j+1}\bm v_{j+1}-\bm b_j\right)\\
 &=s_j\bm v_j-s_{j+1}\bm v_{j+1}.
 \end{aligned}
\]

\item \textbf{Proof of \ref{collar:theorem:part:regularity}.}
Part~\ref{collar:theorem:part:convexity} gives convexity and
$C^p$ regularity. It remains to bound the H\"older parameter
independently of $m$ and the phase scale.
The additional correction $A_j-O_j$ is zero where
$O_j-I_j>\delta_j$. We therefore define
\[
 A_j=O_j+\mathcal K_j,\qquad
 \mathcal K_j:=\psi_{\delta_j}\left(\zeta_j\right)-\zeta_j,\qquad
 \mathcal T_j:=\left\{\bm x:\zeta_j\left(\bm x\right)\le\delta_j\right\}.
\]
Lemma~\resultpart{elem:max}{elem:max:part:psi-exact} gives
$\mathcal K_j=0$ on $\mathcal T_j^c$.
On the remaining region, Lemma~\ref{collar:region-lemma}
controls the values of $P_j$ by a fixed multiple of $\delta_j$.
To keep the smoothing correction's H\"older parameter independent
of the phase scale, we normalize its lengths by
\[
 \rho_j^A:=\left(\frac{\delta_j}{\lambda}\right)^{\frac{1}{q}},
 \qquad \lambda\left(\rho_j^A\right)^q=\delta_j,
\]
and values by $\delta_j$.
The composition estimate needs these bounds on a neighborhood,
so Lemma~\resultpart{collar:region-lemma}{collar:region-lemma:part:neighborhood}
also controls the points within distance $\rho_j^A$ of
$\mathcal T_j$. On this neighborhood,
\[
 \begin{aligned}
 qP_j\left(\bm w,t_{j+1}\right)
 &=\lambda\norm{D_m\bm w}_q^q+\alpha|t_{j+1}|^q
 \le C_N\delta_j
 &&\text{by Lemma~\resultpart{collar:region-lemma}{collar:region-lemma:part:neighborhood}},\\
 \norm{\D^k\zeta_j}&\le M_k\lambda\left(\rho_j^A\right)^{q-k},\quad1\le k\le p
 &&\text{by Lemma~\ref{collar:derivative-lemma}}.
 \end{aligned}
\]
Here Eq.~\eqref{collar:exit-dual} controls the linear term
$s_{j+1}t_{j+1}$ in $\zeta_j$, so the constants are independent of
the next phase's direction and scale.
With the normalized value and derivative bounds,
Lemma~\resultpart{collar:composition-lemma}{collar:composition-lemma:part:correction}
gives, for all $\bm x,\bm y\in\R^d$,
\[
 \begin{aligned}
 \norm{\D^p\mathcal K_j\left(\bm x\right)
                   -\D^p\mathcal K_j\left(\bm y\right)}
 &\le C_{\mathcal K}\frac{\delta_j}{\left(\rho_j^A\right)^{p+\nu}}
       \norm{\bm x-\bm y}^\nu\\
 &=C_{\mathcal K}\frac{\delta_j}{\left(\rho_j^A\right)^q}
       \norm{\bm x-\bm y}^\nu
 &&\text{since }q=p+\nu,\\
 &=C_{\mathcal K}\lambda\norm{\bm x-\bm y}^\nu
 &&\text{by the definition of }\rho_j^A.
 \end{aligned}
\]
It remains to add $P_j$'s global bound.
Since $O_j-P_j$ is constant and $p\ge2$,
$\D^pO_j=\D^pP_j$.
With $\alpha=\lambda m^{1-\frac{3q}{2}}\le\lambda$,
Lemma~\resultpart{collar:P-lemma}{collar:P-lemma:part:residual}
gives H\"older parameter $L_P\lambda$ for $P_j$.
The decomposition $A_j=O_j+\mathcal K_j$ therefore gives
\[
\begin{aligned}
 &\norm{\D^pA_j\left(\bm x\right)-\D^pA_j\left(\bm y\right)}\\
 &\quad\le\norm{\D^pO_j\left(\bm x\right)-\D^pO_j\left(\bm y\right)}
        +\norm{\D^p\mathcal K_j\left(\bm x\right)-\D^p\mathcal K_j\left(\bm y\right)}\\
 &\quad\le L_P\lambda\norm{\bm x-\bm y}^{\nu}
        +C_{\mathcal K}\lambda\norm{\bm x-\bm y}^{\nu}\\
 &\quad=\left(L_P+C_{\mathcal K}\right)\lambda\norm{\bm x-\bm y}^{\nu}
 =C_A\lambda\norm{\bm x-\bm y}^{\nu}.
\end{aligned}
\]
The constants $L_P,C_{\mathcal K},C_A$ depend only on
$p,\nu$ and the fixed smooth maximum, giving the claimed global
bound independently of $m$, dimension, and phase scale.

\item \textbf{Proof of \ref{collar:theorem:part:terminal}.}
For $s_{j+1}=t_{j+1}^\star=0$, both next-coordinate bounds reduce to
$0\le c_0\mathcal E_j$. Restrict the construction on $E_j\times\R$
by the linear isometry $J_j\bm w=\left(\bm w,0\right)$. Then
\[
 \D^k\left(A_j\circ J_j\right)\left(\bm w\right)[\bm d_1,\ldots,\bm d_k]
 =\D^kA_j\left(J_j\bm w\right)[J_j\bm d_1,\ldots,J_j\bm d_k],
\]
so convexity and the H\"older bound restrict to $E_j$. The formulas above
restrict by setting $t_{j+1}=0$, and $\bm g_j=s_j\bm v_j-\bm b_j\in E_j$.
No additional ambient direction $\bm v_{j+1}$ is needed.
\end{itemize}
\end{proof}

\subsubsection{Prove preservation of the oracle outputs}
\label{collar:base-separation}\label{app:phase-separation}

Theorem~\resultpart{collar:theorem}{collar:theorem:part:prefix}
proves \ref{collar:req:oracle} for the phase function. For the base, we must
separate $W_{j,b}$ from the next modification ball.

Apply the endpoint construction in mutually orthogonal blocks
\(E_j=\operatorname{span}\left\{\bm v_j,\bm u_{j,1},\ldots,\bm u_{j,m}\right\}\):
\begin{equation}
t_j^\star=\left(m+1\right)\sqrt m\,\tau_j,\qquad
\bm h_j=t_j^\star\bm v_j+\sum_{i=1}^m\left(m-i+1\right)\tau_j\bm u_{j,i}.
\label{base:blockendpoint}
\end{equation}
Set \(\bm c_0=0\), \(\bm c_{j+1}=\bm c_j+\bm h_j\), and
\(R_j=\norm{\bm h_j}\). In particular,
\begin{equation}
R_j^2=\tau_j^2\left[m\left(m+1\right)^2+\sum_{k=1}^m k^2\right]
\le5\tau_j^2m^3.
\label{base:endpointnorm}
\end{equation}
Thus \(\tau_j=\eta^j\tau_0\) gives \(R_j=\eta^jR_0\).

\begin{lemma}[Future base changes miss every early prefix]
\label{collar:separation-lemma}
Let \(W_{j,b}\) be the span of all earlier blocks,
\(\bm v_j,\bm u_{j,1},\ldots,\bm u_{j,b}\), and any additional directions
orthogonal to all hard blocks. Suppose \(0\le b\le\left\lfloor \frac{m}{2}\right\rfloor\)
and \(\theta\le\frac{1}{16}\).
\begin{claims}
\item\label{collar:separation-lemma:part:distance}\textbf{Separation.}
\begin{equation}
\operatorname{dist}\left(\bm c_{j+1},W_{j,b}\right)
\ge \frac{R_j}{\sqrt{120}}>\theta R_j.
\label{base:supportexclusion}
\end{equation}
\item\label{collar:separation-lemma:part:agreement}\textbf{Base agreement.}
If every modification from $B_j$ onward lies in $S_j=B\left(\bm c_{j+1},\theta R_j\right)$, then
\[
 B=B_j\quad\text{near every }\bm x\in W_{j,b}.
\]
\end{claims}
\end{lemma}
\begin{proof}
\begin{itemize}
\item \textbf{Proof of \ref{collar:separation-lemma:part:distance}.}
The previous center $\bm c_j$ belongs to $W_{j,b}$. From Eq.~\eqref{base:blockendpoint}, the component of $\bm c_{j+1}=\bm c_j+\bm h_j$ perpendicular to $W_{j,b}$ is
\[
 \sum_{i=b+1}^m\left(m-i+1\right)\tau_j\bm u_{j,i}.
\]
Orthogonality and $b\le \frac{m}{2}$ therefore give
\begin{align*}
  \operatorname{dist}\left(\bm c_{j+1},W_{j,b}\right)^2
  &=\tau_j^2\sum_{i=b+1}^m\left(m-i+1\right)^2\\
  &=\tau_j^2\sum_{a=1}^{m-b}a^2\\
  &=\tau_j^2\frac{\left(m-b\right)\left(m-b+1\right)\left(2\left(m-b\right)+1\right)}{6}\\
  &=\frac{\tau_j^2}{6}
      \left(2\left(m-b\right)^3+3\left(m-b\right)^2+\left(m-b\right)\right)\\
  &\ge\frac{\tau_j^2\left(m-b\right)^3}{3}\\
  &\ge\frac{\tau_j^2\left(\frac{m}{2}\right)^3}{3}
   =\frac{\tau_j^2m^3}{24}.
\end{align*}
Apply Lemma~\resultpart{collar:chain-lemma}{collar:chain-lemma:part:geometry} to the $j^{th}$ endpoint. Its bound $R_j^2\le5\tau_j^2m^3$ then yields
\[
\begin{aligned}
 \operatorname{dist}\left(\bm c_{j+1},W_{j,b}\right)^2
 &\ge\frac{\tau_j^2m^3}{24}
 \ge\frac{R_j^2}{5\cdot24}
 =\frac{R_j^2}{120},\\
 \operatorname{dist}\left(\bm c_{j+1},W_{j,b}\right)
 &\ge\frac{R_j}{\sqrt{120}}>\theta R_j,
\end{aligned}
\]
because $\theta\le\frac{1}{16}<\frac{1}{\sqrt{120}}$. This is Eq.~\eqref{base:supportexclusion}.
\item \textbf{Proof of \ref{collar:separation-lemma:part:agreement}.}
For any $\bm x\in W_{j,b}$, part~\ref{collar:separation-lemma:part:distance} gives
\[
 \operatorname{dist}\left(\bm x,\overline S_j\right)
 \ge\operatorname{dist}\left(\bm c_{j+1},W_{j,b}\right)-\theta R_j>0.
\]
The assumed containment of every later modification in $S_j$ and this separation verify the hypotheses of Lemma~\resultpart{base:prefix-lemma}{base:prefix-lemma:part:neighborhood} in Appx.~\ref{app:paper-base}. Applying it with $W=W_{j,b}$ gives
\[
 B=B_j\quad\text{on a neighborhood of every }\bm x\in W_{j,b}.
\]
The displayed distance bound is independent of $\norm{\bm x}$.
\end{itemize}
\end{proof}

Since $B_j$ uses only selected centers and Euclidean distances, this
neighborhood agreement makes the base's oracle output independent of remaining directions.
A fixed fraction of missing links is essential: with one link left,
$\operatorname{dist}\left(\bm c_{j+1},W_{j,m-1}\right)=\tau_j$, whereas
$\theta R_j\ge\theta\tau_jm^{\frac{3}{2}}$.

The scalar assumptions of Theorem~\ref{collar:theorem} are checked in
Appx.~\ref{app:paper-assembly}. In particular,
Corollary~\ref{collar:coefficient-corollary}
derives the coefficient
\(\lambda=K\mu m^{\frac{3q-2}{2}}\) needed to control the base increment.

\subsection{Technical lemmas for the phase construction}
\label{app:phase-support}

Fix $0\le j<L$. Use $\bm x=\left(\bm w,t_{j+1}\right)\in E_j\times\R$ and the ordered basis
$\left(\bm v_j,\bm u_{j,1},\ldots,\bm u_{j,m}\right)$ for $D_m\bm w$.
The terminal phase omits $t_{j+1}$ and all next-phase terms.

\subsubsection{Value and derivative formulas under weaker assumptions}
\label{collar:gradient-and-gap}
\label{collar:scalar}
\label{collar:width}
\label{collar:branches}

\begin{lemma}[Convexity, the residual gap, and the endpoint derivatives]
\label{collar:branches-lemma}
Under the assumptions of Lemma~\ref{collar:margin-lemma}:
\begin{claims}
\item\label{collar:branches-lemma:part:convexity}\textbf{Convexity.}
$A_j$ is convex and $A_j\ge I_j,O_j$.
\item\label{collar:branches-lemma:part:prefix}\textbf{Prefix.}
For $\bm w\in U_{j,b}$ and $0\le b\le\lfloor \frac{m}{4}\rfloor$,
\[
 A_j=O_j\ \text{near }\left(\bm w,0\right),\qquad
 A_j\left(\bm w,0\right)-I_j\left(\bm w,0\right)\ge7\delta_j.
\]
\item\label{collar:branches-lemma:part:origin}\textbf{Origin.}
As in Eq.~\eqref{collar:origin},
\[
 A_j\left(0,0\right)=-\frac{\lambda}{q}\left(m+1\right)\tau_j^q-4\delta_j,\qquad
 \D^kA_j\left(0,0\right)=0\quad\left(1\le k\le p\right).
\]
\item\label{collar:branches-lemma:part:endpoints}\textbf{Affine region.}
$A_j=I_j$ near $[\left(\bm h_j,0\right),\left(\bm h_j,t_{j+1}^\star\right)]$.
The endpoint data are Eq.~\eqref{collar:port-values}--\eqref{collar:port-higher};
the gradient increment is Eq.~\eqref{collar:source}.
\end{claims}
\end{lemma}
\begin{proof}
\begin{itemize}
\item \textbf{Proof of \ref{collar:branches-lemma:part:convexity}.}
The positive weighted sum defining $P_j$ is convex, since
$\phi\left(a\right)=\frac{|a|^q}{q}$ is convex by
Lemma~\resultpart{elem:power}{elem:power:part:hessian} in
Appx.~\ref{app:paper-elementary}. Hence $O_j$ is convex and $I_j$ is
affine. Lemma~\resultpart{elem:composition-lemma}{elem:composition-lemma:part:convexity}
gives convexity of $A_j$, and Lemma~\ref{paper:smooth-maximum-properties}
gives $A_j\ge I_j,O_j$.
\item \textbf{Proof of \ref{collar:branches-lemma:part:prefix}.}
Lemma~\resultpart{collar:margin-lemma}{collar:margin-lemma:part:prefix} gives
$\zeta_j\left(\bm w,0\right)\ge7\delta_j>\delta_j$.
Continuity preserves the strict last inequality on a neighborhood of
$\left(\bm w,0\right)$. Lemma~\resultpart{elem:composition-lemma}{elem:composition-lemma:part:exact} then gives
\[
 A_j\left(\bm w,0\right)-I_j\left(\bm w,0\right)
 =O_j\left(\bm w,0\right)-I_j\left(\bm w,0\right)
 =\zeta_j\left(\bm w,0\right)\ge7\delta_j.
\]
\item \textbf{Proof of \ref{collar:branches-lemma:part:origin}.}
Apply the preceding part with $\bm w=0\in U_{j,0}$.
The power value and derivatives through order $p$ vanish at zero by
Lemma~\resultpart{elem:power}{elem:power:part:flatness}. Thus
\[
\begin{aligned}
 A_j\left(0,0\right)
 &=O_j\left(0,0\right)=P_j\left(0,0\right)-P_j\left(\bm h_j,0\right)-4\delta_j\\
 &=-\frac{\lambda}{q}\left(m+1\right)\tau_j^q-4\delta_j,\\
 \D^kA_j\left(0,0\right)
 &=\D^kO_j\left(0,0\right)=\D^kP_j\left(0,0\right)=0
 \qquad\left(1\le k\le p\right).
\end{aligned}
\]
The derivative equalities follow from neighborhood agreement, rather
than only the equality of values at the origin.
\item \textbf{Proof of \ref{collar:branches-lemma:part:endpoints}.}
Lemma~\resultpart{collar:margin-lemma}{collar:margin-lemma:part:endpoints}
gives $\zeta_j\left(\bm h_j,0\right)=-4\delta_j$ and
$\zeta_j\left(\bm h_j,t_{j+1}^\star\right)\le-3\delta_j$.
Since $P_j$ is convex and $I_j$ is affine,
$\zeta_j=P_j-I_j-P_j\left(\bm h_j,0\right)-4\delta_j$ is convex.
For $0\le r\le1$, its value at an intermediate point is at
most the weighted average of these endpoint values, i.e.,
\[
\begin{aligned}
 \zeta_j\left(\bm h_j,rt_{j+1}^\star\right)
 &\le\left(1-r\right)\zeta_j\left(\bm h_j,0\right)
             +r\zeta_j\left(\bm h_j,t_{j+1}^\star\right)
 &&\text{by convexity of }\zeta_j,\\
 &\le-4\left(1-r\right)\delta_j-3r\delta_j
 &&\text{by Lemma~\resultpart{collar:margin-lemma}{collar:margin-lemma:part:endpoints}},\\
 &=\left(-4+r\right)\delta_j\le-3\delta_j<-\delta_j
 &&\text{since }0\le r\le1.
\end{aligned}
\]
Thus the entire segment lies in the open set
$\{\zeta_j<-\delta_j\}$.
Lemma~\resultpart{elem:composition-lemma}{elem:composition-lemma:part:exact}
gives $A_j=I_j$ on this set, including equality of derivatives.
Evaluating and differentiating Eq.~\eqref{collar:inner} gives
\[
 \begin{aligned}
 A_j\left(\bm h_j,t_{j+1}\right)
 &=\ip{\bm g_j}{t_{j+1}\bm v_{j+1}}=-s_{j+1}t_{j+1}
 &&\left(0\le t_{j+1}\le t_{j+1}^\star\right),\\
 \D A_j&=\D I_j=\bm g_j,\qquad
 \D^kA_j=\D^kI_j=0
 &&\left(2\le k\le p\right).
 \end{aligned}
\]
These identities give Eq.~\eqref{collar:port-values}--\eqref{collar:port-higher}.
Combining the fixed endpoint gradient with the zero starting
gradient in part~\ref{collar:branches-lemma:part:origin} gives
\[
\begin{aligned}
 \bm b_j+\D A_j\left(\bm h_j,0\right)-\D A_j\left(0,0\right)
 &=\bm b_j+\bm g_j\\
 &=\bm b_j+s_j\bm v_j-s_{j+1}\bm v_{j+1}-\bm b_j\\
 &=s_j\bm v_j-s_{j+1}\bm v_{j+1}.
\end{aligned}
\]
None of these calculations uses Eq.~\eqref{collar:exit-dual}
or the particular choice of $\alpha$.
\end{itemize}
\end{proof}

\subsubsection{Details for preservation of the oracle outputs}
\label{app:phase-oracle}

\Needspace{8\baselineskip}
\begin{lemma}[Only one additional direction enters a phase's oracle output]
\label{collar:prefix-oracle-lemma}
Under the assumptions of Lemma~\ref{collar:margin-lemma}:
\begin{claims}
\item\label{collar:prefix-oracle-lemma:part:dependence}
For $\bm w\in U_{j,b}$ and $0\le b\le\lfloor \frac{m}{4}\rfloor$, the full tensors
\[
 \D^kA_j\left(\bm w,0\right),\quad 0\le k\le p\qquad\left(\D^0A_j=A_j\right),
\]
given by Eq.~\eqref{collar:prefix-derivatives}--\eqref{collar:prefix-value},
depend only on $\bm v_j,\bm u_{j,1},\ldots,\bm u_{j,b+1}$ and fixed scalars.
\end{claims}
\end{lemma}
\begin{proof}
Lemma~\resultpart{collar:branches-lemma}{collar:branches-lemma:part:prefix}
gives $A_j=O_j$ near the evaluation point, so every derivative through order $p$
there agrees with that of $O_j$. Apply
Lemma~\resultpart{collar:chain-lemma}{collar:chain-lemma:part:prefix}
with $k=b<m$. For $1\le r\le p$ and arbitrary ambient vectors
$\bm d_1,\ldots,\bm d_r$, it gives
\[
\begin{aligned}
 \D^rA_j\left(\bm w,0\right)[\bm d_1,\ldots,\bm d_r]
 &=\D^rP_j\left(\bm w,0\right)[\bm d_1,\ldots,\bm d_r]\\
 &=\lambda\sum_{i=0}^{b}\phi^{\left(r\right)}\left(\left(D_m\bm w\right)_i\right)
         \prod_{\ell=1}^{r}\ip{\bm a_{j,i}}{\bm d_\ell}.
\end{aligned}
\]
The term $\alpha\phi\left(t_{j+1}\right)$ contributes zero because
$\phi^{\left(r\right)}\left(0\right)=0$ by
Lemma~\resultpart{elem:power}{elem:power:part:flatness} in
Appx.~\ref{app:paper-elementary}. The same statement at order zero,
and the constant in $O_j$, give Eq.~\eqref{collar:prefix-value}.
Rows $0,\ldots,b$ contain only
$\bm v_j,\bm u_{j,1},\ldots,\bm u_{j,b+1}$, so these are full ambient
tensor identities independent of all remaining directions.
They recover Eq.~\eqref{collar:prefix-derivatives}--\eqref{collar:prefix-value}
in Theorem~\resultpart{collar:theorem}{collar:theorem:part:prefix}
under the weaker assumptions of Lemma~\ref{collar:margin-lemma}.
\end{proof}

\subsubsection{H\"older estimates}
\label{collar:power-regularity}
\label{collar:correction-region}
\label{collar:local-derivatives}
\label{collar:composition}

Use Eq.~\eqref{collar:s}, Eq.~\eqref{collar:constants}, and Eq.~\eqref{collar:alpha},
and assume Eq.~\eqref{collar:G-assumption}, Eq.~\eqref{collar:exit-port}, and
Eq.~\eqref{collar:exit-dual}; in particular, $0<\alpha\le\lambda$.

\Needspace{10\baselineskip}

\begin{lemma}[Derivative bounds for the residual expression]
\label{collar:P-lemma}
Define
\[
 \mathcal R_j\bm x:=\left(D_m\bm w,\left(\frac{\alpha}{\lambda}\right)^{\frac{1}{q}}t_{j+1}\right),
 \qquad P_j\left(\bm x\right)=\frac{\lambda}{q}\norm{\mathcal R_j\bm x}_q^q,
 \qquad \norm{\mathcal R_j}_{2\to2}\le2.
\]
Omit the last coordinate in the terminal phase. Then:
\begin{claims}
\item\label{collar:P-lemma:part:residual}\textbf{Global H\"older bound.}
For all $\bm x,\bm y\in E_j\times\R$,
\begin{equation}
 \begin{aligned}
 \norm{\D^pP_j\left(\bm x\right)-\D^pP_j\left(\bm y\right)}
 &\le L_P\lambda\norm{\bm x-\bm y}^\nu,\\
 L_P&=L_\phi\left(2^{p+\nu}+1\right).
 \end{aligned}
\label{collar:P-holder}
\end{equation}
\item\label{collar:P-lemma:part:difference}\textbf{Difference between the formulas.}
For all $\bm x,\bm y\in E_j\times\R$,
\[
 \norm{\D^p\zeta_j\left(\bm x\right)-\D^p\zeta_j\left(\bm y\right)}
 \le L_P\lambda\norm{\bm x-\bm y}^\nu.
\]
\item\label{collar:P-lemma:part:orders}\textbf{Orders $1$ through $p$.}
\[
 \norm{\D^kP_j\left(\bm x\right)}
 \le2^kS_{q,k}\lambda\norm{\mathcal R_j\bm x}_q^{q-k}
 \qquad\left(\bm x\in E_j\times\R,\ 1\le k\le p\right).
\]
\end{claims}
\end{lemma}
\begin{proof}
By definition,
\[
 \lambda\norm{\mathcal R_j\bm x}_q^q
 =\lambda\norm{D_m\bm w}_q^q
   +\lambda\left(\frac{\alpha}{\lambda}\right)|t_{j+1}|^q
 =\lambda\norm{D_m\bm w}_q^q+\alpha|t_{j+1}|^q.
\]
By $\alpha\le\lambda$ and Lemma~\resultpart{collar:chain-lemma}{collar:chain-lemma:part:norm},
\[
 \norm{\mathcal R_j\bm x}_2^2
 =\norm{D_m\bm w}_2^2+\left(\frac{\alpha}{\lambda}\right)^{\frac{2}{q}}|t_{j+1}|^2
 \le4\norm{\bm w}^2+|t_{j+1}|^2
 \le4\norm{\bm x}^2.
\]
Hence $\norm{\mathcal R_j}_{2\to2}\le2$.

Put $\Phi\left(\bm a\right)=\sum_i\phi\left(a_i\right)$, so $P_j=\lambda\Phi\circ\mathcal R_j$.
For $2\le k\le p$, the $k^{th}$ derivative is diagonal:
\[
 \D^k\Phi\left(\bm a\right)[\bm z_1,\ldots,\bm z_k]
 =\sum_i\phi^{\left(k\right)}\left(a_i\right)\prod_{\ell=1}^k z_{\ell,i}.
\]
Lemma~\ref{sources:orthogonal-forms} in Appx.~\ref{app:paper-elementary}, applied to the one-dimensional
orthogonal coordinate spaces for every $k\ge2$, gives
\[
 \begin{aligned}
 \norm{\D^k\Phi\left(\bm a\right)}&=\max_i|\phi^{\left(k\right)}\left(a_i\right)|,\\
 \norm{\D^p\Phi\left(\bm a\right)-\D^p\Phi\left(\bm b\right)}
 &=\max_i|\phi^{\left(p\right)}\left(a_i\right)-\phi^{\left(p\right)}\left(b_i\right)|.
 \end{aligned}
\]
There is no factor $m$. For $1\le k\le p$, the linear chain rule gives
\[
 \D^kP_j\left(\bm x\right)[\bm z_1,\ldots,\bm z_k]
 =\lambda\D^k\Phi\left(\mathcal R_j\bm x\right)
   [\mathcal R_j\bm z_1,\ldots,\mathcal R_j\bm z_k].
\]
\begin{itemize}
\item \textbf{Proof of \ref{collar:P-lemma:part:residual}.}
Apply the preceding identity and Lemma~\resultpart{elem:power}{elem:power:part:holder}:
\[
\begin{aligned}
 \norm{\D^pP_j\left(\bm x\right)-\D^pP_j\left(\bm x'\right)}
 &\le\lambda\norm{\mathcal R_j}_{2\to2}^p
   \max_i\left|\phi^{\left(p\right)}\left(\left(\mathcal R_j\bm x\right)_i\right)
                 -\phi^{\left(p\right)}\left(\left(\mathcal R_j\bm x'\right)_i\right)\right|\\
 &\le2^p L_\phi\lambda
       \norm{\mathcal R_j\left(\bm x-\bm x'\right)}_\infty^\nu\\
 &\le2^p L_\phi\lambda
       \left(2\norm{\bm x-\bm x'}\right)^\nu\\
 &=2^{p+\nu}L_\phi\lambda\norm{\bm x-\bm x'}^\nu.
\end{aligned}
\]
Here $\norm{\mathcal R_j\bm z}_\infty\le\norm{\mathcal R_j\bm z}_2\le2\norm{\bm z}$
and $2^{p+\nu}L_\phi\le L_P$. The scalar bound includes $\nu=1$.

\item \textbf{Proof of \ref{collar:P-lemma:part:difference}.}
Since $\zeta_j-P_j$ is affine and $p\ge2$,
$\D^p\zeta_j=\D^pP_j$; apply part~\ref{collar:P-lemma:part:residual}.

\item \textbf{Proof of \ref{collar:P-lemma:part:orders}.}
For $2\le k\le p$, use the chain rule, diagonal norm identity, and
Lemma~\resultpart{elem:power}{elem:power:part:derivatives}:
\[
\begin{aligned}
 \norm{\D^kP_j\left(\bm x\right)}
 &\le\lambda\norm{\mathcal R_j}_{2\to2}^k
       \max_i|\phi^{\left(k\right)}\left(\left(\mathcal R_j\bm x\right)_i\right)|\\
 &\le2^kS_{q,k}\lambda\norm{\mathcal R_j\bm x}_\infty^{q-k}
 \le2^kS_{q,k}\lambda\norm{\mathcal R_j\bm x}_q^{q-k}.
\end{aligned}
\]
Here $q-k>0$, with zero arguments covered by
Lemma~\resultpart{elem:power}{elem:power:part:flatness}.
For $k=1$, set $\bm a=\mathcal R_j\bm x$ and use
$\D P_j=\lambda\mathcal R_j^\top\left(|\bm a|^{q-2}\bm a\right)$ (componentwise powers):
\[
\begin{aligned}
 \norm{\D P_j\left(\bm x\right)}
 &\le2\lambda\left(\sum_i|a_i|^{2\left(q-1\right)}\right)^{\frac{1}{2}}
 =2\lambda\norm{\bm a}_{2\left(q-1\right)}^{q-1}\\
 &\le2\lambda\norm{\bm a}_q^{q-1}
 =2S_{q,1}\lambda\norm{\mathcal R_j\bm x}_q^{q-1}.
\end{aligned}
\]
Here $2\left(q-1\right)\ge q$ and $S_{q,1}=1$.
\end{itemize}
\end{proof}

Set $\chi_\delta\left(a\right):=\psi_\delta\left(a\right)-a$ for $\delta>0$ and
$\mathcal K_j:=\chi_{\delta_j}\circ\zeta_j=A_j-O_j$.
Define
\begin{equation}
\rho_j^A=\left(\frac{\delta_j}{\lambda}\right)^{\frac{1}{q}},\qquad
\mathcal T_j=\{\bm x:\zeta_j\left(\bm x\right)\le\delta_j\}.
\label{collar:rho}
\end{equation}
Let $\mathcal T_{j,\rho_j^A}:=\{\bm x:\operatorname{dist}\left(\bm x,\mathcal T_j\right)\le\rho_j^A\}$
be its closed Euclidean $\rho_j^A$-neighborhood; $\rho_j^A$ differs from
the base support radius $\rho_j=\theta R_j$ in Eq.~\eqref{base:scaledconstants}.

\begin{lemma}[Nearest point in a nonempty closed set {\cite[Sec.~8.1, p.~397]{BoydVandenberghe2004}}]
\label{sources:nearest-point}\label{sources:projection}
If $S\subset\R^d$ is nonempty and closed, then for every $\bm x\in\R^d$
there exists $\bm a\in S$ such that
\[
 \norm{\bm x-\bm a}=\operatorname{dist}\left(\bm x,S\right).
\]
Such a point is called a projection; there can be more than one.
If $S$ is also convex, the Euclidean projection is unique.
\end{lemma}

\begin{lemma}[Residual bounds on the enlarged region]
\label{collar:region-lemma}
Under the assumptions stated before Lemma~\ref{collar:P-lemma}, set
\[
 C_{\rm dual}:=\frac{2^{\frac{1}{q-1}}\left(2+c_0\right)}{q^*c_1},\qquad
 C_{\mathcal T}:=2q\left(1+C_{\rm dual}\right),\qquad
 C_N:=2^{q-1}\left(C_{\mathcal T}+2^q+1\right).
\]
These constants depend only on $q$.
\begin{claims}
\item\label{collar:region-lemma:part:sublevel}\textbf{Sublevel set.}
\[
 \begin{aligned}
  \zeta_j\left(\bm x\right)&\ge-C_{\rm dual}\delta_j
      &&\left(\bm x\in E_j\times\R\right),\\
  \lambda\norm{D_m\bm w}_q^q+\alpha|t_{j+1}|^q&\le C_{\mathcal T}\delta_j
      &&\left(\left(\bm w,t_{j+1}\right)\in\mathcal T_j\right).
 \end{aligned}
\]
\item\label{collar:region-lemma:part:neighborhood}\textbf{Enlarged region.}
\begin{equation}
 \lambda\norm{D_m\bm w}_q^q+\alpha|t_{j+1}|^q
 \le C_N\delta_j\qquad\left(\left(\bm w,t_{j+1}\right)\in\mathcal T_{j,\rho_j^A}\right).
 \label{collar:neighborhood-energy}
\end{equation}
\end{claims}
\end{lemma}
\begin{proof}
\begin{itemize}
\item \textbf{Proof of \ref{collar:region-lemma:part:sublevel}.}
To bound the residuals where $\zeta_j$ is small, we first bound its
linear terms by $\frac{P_j}{2}$ and a constant. The remaining
$\frac{P_j}{2}$ then bounds the residual energy.
In Eq.~\eqref{collar:power-conjugate}, substitute
\[
 \left(a,b,\ell\right)=\left(\left(D_m\bm w\right)_i,\pi_{j,i},\frac{\lambda}{2}\right)
 \qquad\left(0\le i\le m\right).
\]
Each chain row gives
\[
 \pi_{j,i}\left(D_m\bm w\right)_i
 \le\frac{\lambda}{2q}|\left(D_m\bm w\right)_i|^q
   +\frac{2^{\frac{1}{q-1}}}{q^*}\lambda^{-\frac{1}{q-1}}\pi_{j,i}^{q^*},
\]
where $\left(\frac{\lambda}{2}\right)^{-\frac{1}{q-1}}=2^{\frac{1}{q-1}}\lambda^{-\frac{1}{q-1}}$.
For the next coordinate, apply Corollary~\ref{collar:next-coordinate-lemma}
with $\frac{\alpha}{2}$ in place of $\alpha$. Its infimum identity gives, for every $t_{j+1}\in\R$,
\[
 -s_{j+1}t_{j+1}
 \le\frac{\alpha}{2q}|t_{j+1}|^q
   +\frac{2^{\frac{1}{q-1}}}{q^*}\alpha^{-\frac{1}{q-1}}s_{j+1}^{q^*}.
\]
Summing these inequalities gives
\begin{align*}
 \ip{\bm\pi_j}{D_m\bm w}-s_{j+1}t_{j+1}
 &\le\frac{1}{2}P_j\left(\bm w,t_{j+1}\right)\\
 &\quad+\frac{2^{\frac{1}{q-1}}}{q^*}
  \left(\lambda^{-\frac{1}{q-1}}\sum_{i=0}^m\pi_{j,i}^{q^*}
       +\alpha^{-\frac{1}{q-1}}s_{j+1}^{q^*}\right).
\end{align*}
By Eq.~\eqref{collar:dual-bounds} and $\left(q-1\right)q^*=q$,
\[
\begin{aligned}
 \lambda^{-\frac{1}{q-1}}\sum_{i=0}^m\pi_{j,i}^{q^*}
 &\le\left(m+1\right)\lambda^{-\frac{1}{q-1}}
                  \left(\lambda\tau_j^{q-1}\right)^{q^*}\\
 &=\left(m+1\right)\lambda^{q^*-\frac{1}{q-1}}\tau_j^{\left(q-1\right)q^*}
 =\left(m+1\right)\lambda\tau_j^q\le2\mathcal E_j.
\end{aligned}
\]
Use Eq.~\eqref{collar:exit-dual} for the other term and $\mathcal E_j=\frac{\delta_j}{c_1}$:
the total remainder is at most
\[
 \frac{2^{\frac{1}{q-1}}}{q^*}\left(2+c_0\right)\mathcal E_j
 =\frac{2^{\frac{1}{q-1}}\left(2+c_0\right)}{q^*c_1}\delta_j
 =C_{\rm dual}\delta_j.
\]
Substitute in Eq.~\eqref{collar:dual-identity}, using $k_{j,\rm aff}>0$ from Eq.~\eqref{collar:dual-bounds}:
\[
\begin{aligned}
 \zeta_j
 &=P_j-\left(\ip{\bm\pi_j}{D_m\bm w}-s_{j+1}t_{j+1}\right)+k_{j,\rm aff}\\
 &\ge\frac{1}{2}P_j-C_{\rm dual}\delta_j+k_{j,\rm aff}
 \ge\frac{1}{2}P_j-C_{\rm dual}\delta_j
 \ge-C_{\rm dual}\delta_j.
\end{aligned}
\]
For $\bm x\in\mathcal T_j$, $\zeta_j\left(\bm x\right)\le\delta_j$ therefore gives
\[
 \frac{1}{2}P_j\left(\bm x\right)\le\left(1+C_{\rm dual}\right)\delta_j,
 \qquad
 \lambda\norm{D_m\bm w}_q^q+\alpha|t_{j+1}|^q
 =qP_j\left(\bm x\right)\le C_{\mathcal T}\delta_j.
\]
Omit the $t_{j+1}$-terms in the terminal phase.

\item \textbf{Proof of \ref{collar:region-lemma:part:neighborhood}.}
Continuity of $\zeta_j$ and $\zeta_j\left(\bm x_j^*\right)=-4\delta_j$ make
$\mathcal T_j$ closed and nonempty. For $\bm x\in\mathcal T_{j,\rho_j^A}$,
Lemma~\ref{sources:nearest-point} supplies
\[
 \bm x=\bm a+\bm d,\quad \norm{\bm d}\le\rho_j^A,
 \qquad
 \bm a=\left(\bm w_a,a_{j+1}\right),\quad \bm d=\left(\bm d_E,d_{j+1}\right).
\]
Here $\bm a\in\mathcal T_j$,
$\bm w_a=P_{E_j}\bm a$, and $\bm d_E=P_{E_j}\bm d$.
For $j<L-1$, we define
$a_{j+1}=\ip{\bm v_{j+1}}{\bm a}$ and
$d_{j+1}=\ip{\bm v_{j+1}}{\bm d}$.
Thus $\bm w=\bm w_a+\bm d_E$ and $t_{j+1}=a_{j+1}+d_{j+1}$.
The terminal phase omits the scalar components.
Apply $|a+b|^q\le2^{q-1}\left(|a|^q+|b|^q\right)$ coordinatewise:
\[
\begin{aligned}
 \lambda\norm{D_m\bm w}_q^q+\alpha|t_{j+1}|^q
 &\le2^{q-1}\left[
       \lambda\norm{D_m\bm w_a}_q^q+\alpha|a_{j+1}|^q
       +\lambda\norm{D_m\bm d_E}_q^q+\alpha|d_{j+1}|^q\right]\\
 &\le2^{q-1}\left[C_{\mathcal T}\delta_j
       +2^q\lambda\norm{\bm d_E}^q+\lambda|d_{j+1}|^q\right]\\
 &\le2^{q-1}\left[C_{\mathcal T}\delta_j
       +\left(2^q+1\right)\lambda\left(\rho_j^A\right)^q\right]\\
 &=2^{q-1}\left(C_{\mathcal T}+2^q+1\right)\delta_j=C_N\delta_j.
\end{aligned}
\]
The bounds used are part~\ref{collar:region-lemma:part:sublevel} at $\bm a$,
$\norm{D_m\bm d_E}_q\le2\norm{\bm d_E}$
(Lemma~\resultpart{collar:chain-lemma}{collar:chain-lemma:part:norm}),
$\alpha\le\lambda$, $\norm{\bm d_E},|d_{j+1}|\le\rho_j^A$, and
$\lambda\left(\rho_j^A\right)^q=\delta_j$.
\end{itemize}
\end{proof}

\begin{lemma}[Bounds for each derivative order]
\label{collar:derivative-lemma}
Assume Lemma~\ref{collar:region-lemma}'s hypotheses. The constants
$M_k$ in Eq.~\eqref{collar:Mj} and Eq.~\eqref{collar:M1} depend only on $p,\nu$;
$k$ is the derivative order and $j$ the phase index.
\begin{claims}
\item\label{collar:derivative-lemma:part:higher}\textbf{Orders $2$ through $p$.}
\begin{equation}
 \norm{\D^k\zeta_j\left(\bm x\right)}_{\mathrm{op}}
 \le M_k\lambda\left(\rho_j^A\right)^{q-k}
 \quad\left(\bm x\in\mathcal T_{j,\rho_j^A},\ 2\le k\le p\right).
\label{collar:local-derivative-bound}
\end{equation}
\item\label{collar:derivative-lemma:part:gradient}\textbf{Gradient.}
\[
 \norm{\D\zeta_j\left(\bm x\right)}\le M_1\lambda\left(\rho_j^A\right)^{q-1}
 \qquad\left(\bm x\in\mathcal T_{j,\rho_j^A}\right).
\]
\end{claims}
\end{lemma}
\begin{proof}
For $\bm x\in\mathcal T_{j,\rho_j^A}$, Lemma~\resultpart{collar:region-lemma}{collar:region-lemma:part:neighborhood} gives
\[
 \lambda\norm{\mathcal R_j\bm x}_q^q
 =\lambda\norm{D_m\bm w}_q^q+\alpha|t_{j+1}|^q
 \le C_N\delta_j=C_N\lambda\left(\rho_j^A\right)^q.
\]
Taking the $q^{th}$ root after division by $\lambda>0$, we obtain
\[
 \norm{\mathcal R_j\bm x}_q\le V\rho_j^A,\qquad V:=C_N^{\frac{1}{q}}.
\]
Lemma~\resultpart{collar:P-lemma}{collar:P-lemma:part:orders} then gives
\[
 \norm{\D^kP_j\left(\bm x\right)}
 \le2^kS_{q,k}\lambda\left(V\rho_j^A\right)^{q-k}
 =2^kS_{q,k}V^{q-k}\lambda\left(\rho_j^A\right)^{q-k}
 \quad\left(1\le k\le p\right).
\]
\begin{itemize}
\item \textbf{Proof of \ref{collar:derivative-lemma:part:higher}.}
For $2\le k\le p$, $\D^k\zeta_j=\D^kP_j$ because $\zeta_j-P_j$ is affine.
Thus Eq.~\eqref{collar:local-derivative-bound} holds with
\begin{equation}
 M_k:=2^kS_{q,k}V^{q-k}\qquad\left(2\le k\le p\right).
\label{collar:Mj}
\end{equation}

\item \textbf{Proof of \ref{collar:derivative-lemma:part:gradient}.}
At order one, $\D\zeta_j=\D P_j-\bm g_j$, and Eq.~\eqref{collar:mismatch} gives
\[
 \bm g_j=\Gamma_j\bm v_j
    -\frac{G_j}{R_j}\sum_{i=1}^m x_{j,i}^*\bm u_{j,i}
    -s_{j+1}\bm v_{j+1}.
\]
The link projection has norm at most $R_j$. Hence, using $G_j\le c_0\Gamma_j$,
\[
\begin{aligned}
 \norm{\bm g_j}
 &\le\Gamma_j\norm{\bm v_j}
       +\frac{G_j}{R_j}
             \left\|\sum_{i=1}^m x_{j,i}^*\bm u_{j,i}\right\|
       +s_{j+1}\\
 &\le\Gamma_j+\frac{G_j}{R_j}R_j+s_{j+1}\\
 &=\Gamma_j+G_j+s_{j+1}\\
 &\le\left(1+c_0\right)\Gamma_j+s_{j+1}.
\end{aligned}
\]
Since $\rho_j^A=\left(\frac{\delta_j}{\lambda}\right)^{\frac{1}{q}}=\left(c_1m\right)^{\frac{1}{q}}\tau_j$ and $m\ge1$,
\[
 \frac{\Gamma_j}{\lambda\left(\rho_j^A\right)^{q-1}}
 =\frac{\lambda\tau_j^{q-1}m^{-\frac{1}{2}}}{\lambda\left(c_1m\right)^{\frac{q-1}{q}}\tau_j^{q-1}}
 =c_1^{-\frac{q-1}{q}}m^{-\frac{1}{2}-\frac{q-1}{q}}
 \le c_1^{-\frac{q-1}{q}}.
\]
From Eq.~\eqref{collar:exit-dual} and $c_0=c_1$,
$\alpha^{-\frac{1}{q-1}}s_{j+1}^{q^*}\le\delta_j$. Its $\left(q^*\right)^{th}$ root gives
\[
\begin{aligned}
 s_{j+1}
 &\le\alpha^{\frac{1}{q}}\delta_j^{\frac{1}{q^*}}
 =\alpha^{\frac{1}{q}}\lambda^{\frac{q-1}{q}}\left(\rho_j^A\right)^{q-1}\\
 &\le\lambda\left(\rho_j^A\right)^{q-1},
\end{aligned}
\]
using $\left(q-1\right)q^*=q$ and $\alpha\le\lambda$. Add the $k=1$ bound for $P_j$:
\[
\begin{aligned}
 \norm{\D\zeta_j}
 &\le\norm{\D P_j}+\norm{\bm g_j}\\
  &\le2V^{q-1}\lambda\left(\rho_j^A\right)^{q-1}
       +\left(1+c_0\right)\Gamma_j+s_{j+1}\\
  &\le2V^{q-1}\lambda\left(\rho_j^A\right)^{q-1}
       +\left(1+c_0\right)c_1^{-\frac{q-1}{q}}\lambda\left(\rho_j^A\right)^{q-1}
       +\lambda\left(\rho_j^A\right)^{q-1}\\
 &\le\left[2V^{q-1}+\left(1+c_0\right)c_1^{-\frac{q-1}{q}}+1\right]
       \lambda\left(\rho_j^A\right)^{q-1}.
\end{aligned}
\]
Take
\begin{equation}
 M_1:=2V^{q-1}+\left(1+c_0\right)c_1^{-\frac{q-1}{q}}+1.
\label{collar:M1}
\end{equation}
All constants depend only on $p,\nu$; the terminal case has $s_{j+1}=0$.
\end{itemize}
\end{proof}

\Needspace{6\baselineskip}
\begin{lemma}[The smoothing correction has a uniform H\"older bound]
\label{collar:composition-lemma}
Under the assumptions of Lemma~\ref{collar:region-lemma}:
\begin{claims}
\item\label{collar:composition-lemma:part:correction}\textbf{Correction.}
$\mathcal K_j=\chi_{\delta_j}\circ\zeta_j\in C^{p,\nu}$ and, for all $\bm x,\bm y\in\R^d$,
\[
 \norm{\D^p\mathcal K_j\left(\bm x\right)-\D^p\mathcal K_j\left(\bm y\right)}
 \le C_{\mathcal K}\lambda\norm{\bm x-\bm y}^\nu.
\]
\item\label{collar:composition-lemma:part:auxiliary}\textbf{Phase function.}
$A_j\in C^{p,\nu}$ and, for all $\bm x,\bm y\in\R^d$,
\[
 \norm{\D^pA_j\left(\bm x\right)-\D^pA_j\left(\bm y\right)}
 \le C_A\lambda\norm{\bm x-\bm y}^\nu,
\]
with $C_A$ given by Eq.~\eqref{collar:Q-holder}.
\end{claims}
The constants depend only on $p,\nu$ and the fixed smoothing function,
independently of $m$, dimension, and phase scales.
No $\D^{p+1}P_j$ is required, including at $\nu=1$.
\end{lemma}
\begin{proof}
\begin{itemize}
\item \textbf{Proof of \ref{collar:composition-lemma:part:correction}.}
Since $\zeta_j\in C^p$ and $\chi_{\delta_j}\in C^\infty$, $\mathcal K_j\in C^p$.
Lemma~\ref{elem:holder-composition} in Appx.~\ref{app:paper-elementary} requires a convex open domain and bounds
on the value and derivatives of the inner function.

\emph{Domain and value bound.}
Since $\zeta_j$ is continuous and convex, $\mathcal T_j$ is closed and convex;
$\zeta_j\left(\bm x_j^*\right)=-4\delta_j$ makes it nonempty. Set
\[
 \mathcal U_j:=\{\bm x:\operatorname{dist}\left(\bm x,\mathcal T_j\right)<\rho_j^A\}
 =\mathcal T_j+B\left(0,\rho_j^A\right).
\]
The sum is convex and open, with
$\mathcal T_j\subset\mathcal U_j\subset\mathcal T_{j,\rho_j^A}$.
For $\bm x\in\mathcal U_j$, choose a nearest $\bm a\in\mathcal T_j$;
$\norm{\bm x-\bm a}<\rho_j^A$. Set
\[
 M:=M_1\lambda\left(\rho_j^A\right)^{q-1},\qquad
 v\left(t\right):=\zeta_j\left(\bm a+t\left(\bm x-\bm a\right)\right),\qquad
 g\left(t\right):=Mt\norm{\bm x-\bm a}\quad\left(0\le t\le1\right).
\]
The segment lies in the open convex set $\mathcal U_j$.
Lemma~\resultpart{collar:derivative-lemma}{collar:derivative-lemma:part:gradient}
gives
\[
\begin{aligned}
 |v'\left(t\right)|
 &\le\norm{\D\zeta_j\left(\bm a+t\left(\bm x-\bm a\right)\right)}
       \norm{\bm x-\bm a}\\
 &\le M\norm{\bm x-\bm a}=g'\left(t\right).
\end{aligned}
\]
Theorem~\ref{sources:mean-value} therefore implies
$|\zeta_j\left(\bm x\right)-\zeta_j\left(\bm a\right)|\le M\norm{\bm x-\bm a}$, so
\[
\begin{aligned}
 \zeta_j\left(\bm x\right)
 &\le\zeta_j\left(\bm a\right)+M_1\lambda\left(\rho_j^A\right)^{q-1}\norm{\bm x-\bm a}\\
 &\le\delta_j+M_1\lambda\left(\rho_j^A\right)^q=\left(1+M_1\right)\delta_j.
\end{aligned}
\]
Here $\zeta_j\left(\bm a\right)\le\delta_j$ and $\lambda\left(\rho_j^A\right)^q=\delta_j$.
Together with $\zeta_j\ge-C_{\rm dual}\delta_j$
(Lemma~\resultpart{collar:region-lemma}{collar:region-lemma:part:sublevel}),
this gives $|\zeta_j|\le M_0\delta_j$ on $\mathcal U_j$, where
$M_0:=\max\{C_{\rm dual},1+M_1\}$.

\emph{Normalized bounds.}
Define
\[
 \widetilde{\mathcal U}_j:=\left(\rho_j^A\right)^{-1}\mathcal U_j,
 \qquad
 \widetilde\zeta_j\left(\bm u\right):=\frac{\zeta_j\left(\rho_j^A\bm u\right)}{\delta_j}.
\]
The domain remains convex and open, and $\norm{\widetilde\zeta_j}_\infty\le M_0$.
For $1\le k\le p$, the chain rule and Lemma~\ref{collar:derivative-lemma} give
\[
\begin{aligned}
 \D^k\widetilde\zeta_j\left(\bm u\right)
 &=\frac{\left(\rho_j^A\right)^k}{\delta_j}\D^k\zeta_j\left(\rho_j^A\bm u\right),\\
 \norm{\D^k\widetilde\zeta_j}_{\infty;\widetilde{\mathcal U}_j}
 &\le\frac{\left(\rho_j^A\right)^k}{\delta_j}
       M_k\lambda\left(\rho_j^A\right)^{q-k}
 =M_k.
\end{aligned}
\]
Here $\delta_j=\lambda\left(\rho_j^A\right)^q$. For distinct
$\bm u,\bm v\in\widetilde{\mathcal U}_j$,
Lemma~\resultpart{collar:P-lemma}{collar:P-lemma:part:difference} gives
\[
\begin{aligned}
 \frac{\norm{\D^p\widetilde\zeta_j\left(\bm u\right)-\D^p\widetilde\zeta_j\left(\bm v\right)}}{\norm{\bm u-\bm v}^\nu}
 &\le\frac{\left(\rho_j^A\right)^p}{\delta_j}
       L_P\lambda\frac{\norm{\rho_j^A\left(\bm u-\bm v\right)}^\nu}{\norm{\bm u-\bm v}^\nu}\\
 &=L_P\frac{\lambda\left(\rho_j^A\right)^{p+\nu}}{\delta_j}=L_P.
\end{aligned}
\]
Since $p+\nu=q$, the bounded norm therefore satisfies
\[
 \norm{\widetilde\zeta_j}_{p,\nu;\widetilde{\mathcal U}_j}
 \le C_\zeta:=\max\{1,M_0,M_1,\ldots,M_p,L_P\}.
\]
Choose $J:=\left(-M_0-1,M_0+1\right)$, so that
$\widetilde\zeta_j\left(\widetilde{\mathcal U}_j\right)\subset[-M_0,M_0]\subset J$.
Smoothness of $\chi_1$ on $\overline J$ makes
$C_{p+1}:=\norm{\chi_1^{\left(p+1\right)}}_{\infty;J}$ finite.
For any distinct $a,b\in J$, the mean-value theorem on the interval
between them gives
\[
\begin{aligned}
 |\chi_1^{\left(p\right)}\left(a\right)-\chi_1^{\left(p\right)}\left(b\right)|
 &\le C_{p+1}|a-b|,\\
 \frac{|\chi_1^{\left(p\right)}\left(a\right)-\chi_1^{\left(p\right)}\left(b\right)|}{|a-b|^\nu}
 &\le C_{p+1}|a-b|^{1-\nu}\\
 &\le C_{p+1}|J|^{1-\nu}.
\end{aligned}
\]
The last line uses $|a-b|\le|J|$ and $1-\nu\ge0$, including $\nu=1$.
Thus, for all $a,b\in J$,
\[
 |\chi_1^{\left(p\right)}\left(a\right)-\chi_1^{\left(p\right)}\left(b\right)|
 \le\norm{\chi_1^{\left(p+1\right)}}_{\infty;J}|J|^{1-\nu}|a-b|^\nu,
\]
where $|J|=2M_0+2$; thus $C_\chi:=\norm{\chi_1}_{p,\nu;J}<\infty$.
Apply Lemma~\ref{elem:holder-composition} with
\[
 U=\widetilde{\mathcal U}_j,\qquad
 f=\widetilde\zeta_j,\qquad h=\chi_1,\qquad J=\left(-M_0-1,M_0+1\right).
\]
The bounds $\norm{f}_{p,\nu;U}\le C_\zeta$ and $\norm{h}_{p,\nu;J}=C_\chi$ yield
\[
 \norm{\chi_1\circ\widetilde\zeta_j}_{p,\nu;\widetilde{\mathcal U}_j}
 \le C_q^{\rm comp}C_\chi\max\{1,C_\zeta^q\}
 =C_{\rm comp},\qquad
 C_{\rm comp}:=C_q^{\rm comp}C_\chi C_\zeta^q.
\]
At $\nu=1$ the lemma uses the Lipschitz $\D^p\widetilde\zeta_j$;
only the smooth outer function $\chi_1$ uses a $\left(p+1\right)^{th}$ derivative.

\emph{Return to the original scale.}
By Eq.~\eqref{elem:psi-M-def},
\[
\begin{aligned}
 \chi_{\delta_j}\left(a\right)
 &=\psi_{\delta_j}\left(a\right)-a
 =\frac{a+\delta_j\Theta\left(\frac{a}{\delta_j}\right)}{2}-a\\
 &=\delta_j\left(
       \frac{\frac{a}{\delta_j}+\Theta\left(\frac{a}{\delta_j}\right)}{2}-\frac{a}{\delta_j}\right)
 =\delta_j\chi_1\left(\frac{a}{\delta_j}\right).
\end{aligned}
\]
Put $H_j=\chi_1\circ\widetilde\zeta_j$. The normalized definition then gives
\[
\begin{aligned}
 \mathcal K_j\left(\rho_j^A\bm u\right)
 &=\chi_{\delta_j}\left(\zeta_j\left(\rho_j^A\bm u\right)\right)\\
 &=\delta_j\chi_1\left(\frac{\zeta_j\left(\rho_j^A\bm u\right)}{\delta_j}\right)
 =\delta_jH_j\left(\bm u\right).
\end{aligned}
\]
Differentiating with respect to $\bm u$ gives
$\left(\rho_j^A\right)^p\D^p\mathcal K_j\left(\rho_j^A\bm u\right)=\delta_j\D^pH_j\left(\bm u\right)$.
Thus, for distinct $\bm u,\bm v\in\widetilde{\mathcal U}_j$,
\[
\begin{aligned}
 &\frac{\norm{\D^p\mathcal K_j\left(\rho_j^A\bm u\right)
                   -\D^p\mathcal K_j\left(\rho_j^A\bm v\right)}}{\norm{\rho_j^A\bm u-\rho_j^A\bm v}^{\nu}}\\
 &\quad=\frac{\delta_j}{\left(\rho_j^A\right)^p}
       \frac{\norm{\D^pH_j\left(\bm u\right)-\D^pH_j\left(\bm v\right)}}{\left(\rho_j^A\right)^\nu\norm{\bm u-\bm v}^\nu}\\
 &\quad=\frac{\delta_j}{\left(\rho_j^A\right)^{p+\nu}}
       \frac{\norm{\D^pH_j\left(\bm u\right)-\D^pH_j\left(\bm v\right)}}{\norm{\bm u-\bm v}^\nu}.
\end{aligned}
\]
The composition bound and the bijection $\bm u\mapsto\rho_j^A\bm u$
therefore give, for all $\bm x,\bm y\in\mathcal U_j$,
\[
\begin{aligned}
 \norm{\D^p\mathcal K_j\left(\bm x\right)-\D^p\mathcal K_j\left(\bm y\right)}
 &\le C_{\rm comp}\frac{\delta_j}{\left(\rho_j^A\right)^{p+\nu}}
       \norm{\bm x-\bm y}^\nu\\
 &=C_{\rm comp}\lambda\norm{\bm x-\bm y}^\nu
 &&\left(p+\nu=q\right).
\end{aligned}
\]

\emph{Global extension.}
Lemma~\resultpart{elem:max}{elem:max:part:psi-exact} and smoothness give
$\chi_{\delta_j}^{\left(k\right)}=0$ on $[\delta_j,\infty)$ for $0\le k\le p$;
the chain rule gives $\D^p\mathcal K_j=0$ on $\{\zeta_j\ge\delta_j\}$.
Pairs in $\mathcal T_j$ satisfy the bound on $\mathcal U_j$;
pairs outside have zero tensors. For $\bm x\in\mathcal T_j$ and
$\bm y\notin\mathcal T_j$, continuity supplies $\bm z\in[\bm x,\bm y]$
with $\zeta_j\left(\bm z\right)=\delta_j$. Thus
$\bm x,\bm z\in\mathcal T_j\subset\mathcal U_j$,
$\D^p\mathcal K_j\left(\bm z\right)=\D^p\mathcal K_j\left(\bm y\right)=0$, and
\[
\begin{aligned}
 \norm{\D^p\mathcal K_j\left(\bm x\right)-\D^p\mathcal K_j\left(\bm y\right)}
 &=\norm{\D^p\mathcal K_j\left(\bm x\right)-\D^p\mathcal K_j\left(\bm z\right)}\\
 &\le C_{\rm comp}\lambda\norm{\bm x-\bm z}^\nu
 \le C_{\rm comp}\lambda\norm{\bm x-\bm y}^\nu.
\end{aligned}
\]
Hence $C_{\mathcal K}:=C_{\rm comp}$ is independent of $m$, dimension, and phase scale.

\item \textbf{Proof of \ref{collar:composition-lemma:part:auxiliary}.}
Use $A_j=O_j+\mathcal K_j$, $\D^pO_j=\D^pP_j$, and
Lemma~\resultpart{collar:P-lemma}{collar:P-lemma:part:residual}.
For all $\bm x,\bm y\in\R^d$,
\begin{equation}
\begin{aligned}
 &\norm{\D^pA_j\left(\bm x\right)-\D^pA_j\left(\bm y\right)}\\
 &\quad\le\norm{\D^pP_j\left(\bm x\right)-\D^pP_j\left(\bm y\right)}
         +\norm{\D^p\mathcal K_j\left(\bm x\right)-\D^p\mathcal K_j\left(\bm y\right)}\\
 &\quad\le\left(L_P+C_{\rm comp}\right)\lambda\norm{\bm x-\bm y}^\nu\\
 &\quad=C_A\lambda\norm{\bm x-\bm y}^\nu,\qquad C_A:=L_P+C_{\rm comp}.
\end{aligned}
\label{collar:Q-holder}
\end{equation}
\end{itemize}
\end{proof}

%% file: appendix-assembly.tex
\section{Compatible Parameters and Assembly}\label{app:paper-assembly}

\subsection{Parameter choices and the complete function}\label{global:assembly}

\subsubsection{Phase coordinates and scales}\label{global:compatible-parameters}

\begin{definition}[Phase coordinates and scales]\label{global:phase-data}
Choose $m,L\in\N$, $\mu,\tau_0>0$, and mutually orthogonal subspaces
\[
 E_j=\operatorname{span}\{\bm v_j,\bm u_{j,1},\ldots,\bm u_{j,m}\},\quad 0\le j<L,
\]
with the displayed orthonormal bases. For a contraction $\eta\in\left(0,1\right)$ to be chosen later, put
\[
 \tau_j=\eta^j\tau_0,\quad
 t_j^\star=\sqrt m\left(m+1\right)\tau_j,\quad
 \bm h_j=t_j^\star\bm v_j+\sum_{i=1}^m\left(m-i+1\right)\tau_j\bm u_{j,i},
 \quad \bm c_j=\sum_{k<j}\bm h_k.
\]
Thus $\bm h_j=\bm c_{j+1}-\bm c_j$. Lemma~\resultpart{collar:chain-lemma}{collar:chain-lemma:part:geometry} in Appx.~\ref{app:paper-phases} gives
\[
 R_j:=\norm{\bm h_j}=\beta_m\tau_jm^{\frac{3}{2}},\qquad
 \beta_m^2=\frac{4}{3}+\frac{5}{2m}+\frac{7}{6m^2},\qquad 1\le\beta_m\le\sqrt5.
\]
In particular $R_j=\eta^jR_0$.
\end{definition}
For $\eta$ satisfying the base nesting inequality Eq.~\eqref{base:etaconstraint}, Lemma~\resultpart{base:nesting-lemma}{base:nesting-lemma:part:recurrences} in Appx.~\ref{app:paper-base} gives
\[
 \bm b_j:=\D B\left(\bm c_{j+1}\right)-\D B\left(\bm c_j\right)
 =\mu R_j^{q-2}\bm h_j,\qquad G_j:=\norm{\bm b_j}=\mu R_j^{q-1}.
\]
To match the chain's endpoint gradient with its affine continuation,
the starting slope must compensate the base increment's projection
$\frac{G_jt_j^\star}{R_j}$ along $\bm v_j$. For the chain coefficient
$\lambda>0$, we therefore define
\[
 \Gamma_j:=\frac{\lambda\tau_j^{q-1}}{\sqrt{m}},\qquad
 s_j:=\Gamma_j+\frac{G_jt_j^\star}{R_j}.
\]

\subsubsection{Center-gradient conditions}
\label{design:gradient-target}\label{app:design-support}
For the mutually orthogonal phase subspaces $E_0,\ldots,E_{L-1}$ and unit vectors $\bm v_j\in E_j$, the centers satisfy
\[
 \bm c_0=0,\qquad \bm c_{j+1}-\bm c_j\in E_j.
\]
Let the phase-$j$ prefix be a subspace $W_{j,b}\subset E_0\oplus\cdots\oplus E_j$. For $j<L-1$,
\[
 \bm v_{j+1}\perp W_{j,b},\qquad \bm v_{j+1}\perp\bm c_{j+1}.
\]
Consequently the choice $\D F\left(\bm c_{j+1}\right)=-s_{j+1}\bm v_{j+1}$ gives $\ip{\D F\left(\bm c_{j+1}\right)}{\bm x-\bm c_{j+1}}=0$ for every $\bm x\in W_{j,b}$, as required by Lemma~\resultpart{design:gap}{design:gap:part:gap}, regardless of $\norm{\bm x}$.

\begin{condition}[Center gradients for the gap proof]\label{design:gradient-condition}
Choose positive $s_j$ and require
\begin{equation}\label{design:full-gradient-target}
 \D F\left(\bm c_j\right)=-s_j\bm v_j\quad\left(0\le j<L\right),\qquad
 \D F\left(\bm c_L\right)=0.
\end{equation}
The terminal equality also identifies the minimizer once uniform convexity is proved. We use $s_L\bm v_L=0$ as a convention; no additional direction $\bm v_L$ is needed.
\end{condition}
\subsubsection{Phase decomposition}
Fix a phase $j$ and $e_j,\delta_j>0$. Let $O_j$ be convex and $C^p$,
let $I_j$ be affine, and suppose
\[
 O_j-I_j\ge e_j\text{ on }W_{j,b},\qquad
 O_j-I_j< -\delta_j\text{ near }\bm c_{j+1}.
\]
The smooth positive part from Appx.~\ref{elem:smooth-max} satisfies
\[
 \psi_{\delta_j}\left(t\right)\ge t_+,\qquad
 \psi_{\delta_j}\left(t\right)=0\ \left(t\le-\delta_j\right),\qquad
 \psi_{\delta_j}\left(t\right)=t\ \left(t\ge\delta_j\right).
\]
For a differentiable function $\Psi_j$, define
\begin{equation}\label{design:phase-decomposition}
 \begin{aligned}
 r_j&:=\psi_{\delta_j}\left(O_j-I_j\right),\\
 A_j&:=I_j+r_j=M_{\delta_j}\left(I_j,O_j\right),\\
 F&:=\Psi_j+r_j=\Psi_j+A_j-I_j.
 \end{aligned}
\end{equation}
Then $r_j\ge e_j$ on $W_{j,b}$ and $r_j=0$ near $\bm c_{j+1}$. Lemma~\resultpart{elem:composition-lemma}{elem:composition-lemma:part:convexity} in Appx.~\ref{app:paper-elementary} gives convexity of $A_j$.

\Needspace{8\baselineskip}
\begin{condition}[Two choices that isolate the current gradient change]\label{design:two-choices}
Let $\bm g_j:=\D I_j$. In addition to Eq.~\eqref{design:phase-decomposition}, assume:
\begin{claims}
\item\label{design:two-choices:part:zero-start} \textbf{Zero starting gradient for the current phase function.}
\begin{equation}\label{design:zero-start}
 \D A_j\left(\bm c_j\right)=0.
\end{equation}
\item\label{design:two-choices:part:unchanged-rest} \textbf{The remaining terms contribute no gradient change.}
\begin{equation}\label{design:unchanged-rest}
 \D\left(\Psi_j-B\right)\left(\bm c_{j+1}\right)=\D\left(\Psi_j-B\right)\left(\bm c_j\right).
\end{equation}
\end{claims}
Part~\ref{design:two-choices:part:zero-start} preserves the starting gradient of $\Psi_j-I_j$ when $A_j$ is added. Part~\ref{design:two-choices:part:unchanged-rest} leaves the gradient change of $\Psi_j$ equal to that of $B$.
\end{condition}
For $\bm b_j:=\D B\left(\bm c_{j+1}\right)-\D B\left(\bm c_j\right)$, the target gradient relation is
\begin{equation}\label{design:gradient-cancellation}
 \begin{aligned}
 -s_j\bm v_j+\bm b_j+\bm g_j&=-s_{j+1}\bm v_{j+1},\\
 \text{equivalently}\qquad\bm g_j&=s_j\bm v_j-s_{j+1}\bm v_{j+1}-\bm b_j.
 \end{aligned}
\end{equation}

\subsubsection{The affine-gradient identity}
\begin{lemma}[The two choices determine the affine gradient]\label{design:balance-lemma}
Assume Condition~\ref{design:two-choices}, the decomposition Eq.~\eqref{design:phase-decomposition}, $A_j=I_j$ near $\bm c_{j+1}$, and $\D F\left(\bm c_j\right)=-s_j\bm v_j$. Keep the definition
$\bm b_j:=\D B\left(\bm c_{j+1}\right)-\D B\left(\bm c_j\right)$. Then:
\begin{claims}
\item\label{design:balance-lemma:part:increments} \textbf{The actual gradient changes are}
\[
 \D\Psi_j\left(\bm c_{j+1}\right)-\D\Psi_j\left(\bm c_j\right)=\bm b_j,\qquad
 \D r_j\left(\bm c_{j+1}\right)-\D r_j\left(\bm c_j\right)=\bm g_j.
\]
\item\label{design:balance-lemma:part:slope} \textbf{The next-gradient target determines the affine slope.}
The target $\D F\left(\bm c_{j+1}\right)=-s_{j+1}\bm v_{j+1}$ holds if and only if the affine gradient satisfies Eq.~\eqref{design:gradient-cancellation}.
\end{claims}
\end{lemma}
\begin{proof}
\begin{itemize}
\item \textbf{Proof of \ref{design:balance-lemma:part:increments}.}
For the first identity, split off the base term:
\[
 \begin{aligned}
 \D\Psi_j\left(\bm c_{j+1}\right)-\D\Psi_j\left(\bm c_j\right)
 &=\D B\left(\bm c_{j+1}\right)-\D B\left(\bm c_j\right)\\
 &\quad+\D\left(\Psi_j-B\right)\left(\bm c_{j+1}\right)-\D\left(\Psi_j-B\right)\left(\bm c_j\right)\\
 &=\D B\left(\bm c_{j+1}\right)-\D B\left(\bm c_j\right)
 &&\text{by Condition~\resultpart{design:two-choices}{design:two-choices:part:unchanged-rest}},\\
 &=\bm b_j.
 \end{aligned}
\]
For the second, use $r_j=A_j-I_j$ from
Eq.~\eqref{design:phase-decomposition}, $\D I_j=\bm g_j$,
Condition~\resultpart{design:two-choices}{design:two-choices:part:zero-start},
and the assumed neighborhood equality $A_j=I_j$:
\[
 \begin{aligned}
 \D r_j\left(\bm c_j\right)
 &=\D A_j\left(\bm c_j\right)-\D I_j\left(\bm c_j\right)=0-\bm g_j,\\
 \D r_j\left(\bm c_{j+1}\right)
 &=\D A_j\left(\bm c_{j+1}\right)-\D I_j\left(\bm c_{j+1}\right)
 =\bm g_j-\bm g_j=0,\\
 \D r_j\left(\bm c_{j+1}\right)-\D r_j\left(\bm c_j\right)
 &=0-\left(-\bm g_j\right)=\bm g_j.
 \end{aligned}
\]
\item \textbf{Proof of \ref{design:balance-lemma:part:slope}.}
Differentiate $F=\Psi_j+r_j$ at the two centers. Part~\ref{design:balance-lemma:part:increments} and the assumed starting gradient yield
\[
 \D F\left(\bm c_{j+1}\right)
 =\D F\left(\bm c_j\right)+\bm b_j+\bm g_j
 =-s_j\bm v_j+\bm b_j+\bm g_j.
\]
Thus the next-center target in Condition~\ref{design:gradient-condition}
holds exactly when
$-s_j\bm v_j+\bm b_j+\bm g_j=-s_{j+1}\bm v_{j+1}$.
Solving for $\bm g_j$ gives Eq.~\eqref{design:gradient-cancellation}.
\end{itemize}
\end{proof}

The assembled phase functions therefore use
\begin{equation}\label{global:cancellation}
 \bm g_j:=s_j\bm v_j-s_{j+1}\bm v_{j+1}-\bm b_j.
\end{equation}
In the terminal phase, we use $s_L=t_L^\star=0$ and omit $\bm v_L$.
The identity $\ip{\bm g_j}{\bm v_j}=\Gamma_j$ follows from
Lemma~\resultpart{collar:affine-lemma}{collar:affine-lemma:part:mismatch}.

\subsubsection{One parameter choice satisfying all requirements}

\begin{proposition}[Parameters satisfying all phase conditions]\label{global:parameters}
Let $c_0,c_1$ be the constants in Theorem~\ref{collar:theorem} in Appx.~\ref{app:paper-phases}, and set
\[
 a_q:=\frac{3q-2}{2},\qquad
 K:=\frac{5^{\frac{q-1}{2}}}{c_0},\qquad
 \eta:=\frac{\gamma}{2\left(1+\theta\right)}.
\]
These constants depend only on $p,\nu$, and $0<\eta<1$. For every $m,L\in\N$ and $\mu,\tau_0>0$,
use Definition~\ref{global:phase-data} and choose
\[
 \lambda:=K\mu m^{a_q},\qquad
 \alpha:=\lambda m^{-a_q}=K\mu,\qquad
 \delta_j:=c_1\lambda m\tau_j^q.
\]
Parts~\ref{global:parameters:part:base}--\ref{global:parameters:part:next} hold for every $0\le j<L$:
\begin{claims}
\item\label{global:parameters:part:base} \textbf{Base gradient.}
\begin{equation}\label{global:condition-base}
 G_j\le c_0\Gamma_j,\qquad \Gamma_j\le s_j\le\left(1+c_0\right)\Gamma_j.
\end{equation}
\item\label{global:parameters:part:next} \textbf{Next-coordinate bounds.}
\begin{align}
 \alpha^{-\frac{1}{q-1}}s_{j+1}^{\frac{q}{q-1}}
 &\le c_0\lambda m\tau_j^q,\label{global:condition-next-gradient}\\
 s_{j+1}t_{j+1}^\star+\frac{\alpha}{q}\left(t_{j+1}^\star\right)^q
 &\le c_0\lambda m\tau_j^q.\label{global:condition-next-value}
\end{align}
\item\label{global:parameters:part:nesting} \textbf{Nested base supports.}
\[
 0<\eta=\frac{\gamma}{2\left(1+\theta\right)},\qquad
 \overline S_{j+1}\subseteq
 \overline B\left(\bm c_{j+1},\frac{\gamma R_j}{2}\right)\subset U_j
 \quad\left(0\le j<L-1\right).
\]
\end{claims}
\end{proposition}

\subsubsection{The base-gradient estimate}

\begin{corollary}[The coefficient required for a block of length \(m\)]
\label{collar:coefficient-corollary}
Suppose \(G_j=\mu R_j^{q-1}\), and choose
\[
\lambda=K\mu m^{\frac{3q-2}{2}},\qquad
K\ge\frac{5^{\frac{q-1}{2}}}{c_0}.
\]
Then the base increment satisfies \(G_j\le c_0\Gamma_j\).
\end{corollary}
\begin{proof}
Substitute $R_j=\beta_m\tau_jm^{\frac{3}{2}}$ from
Definition~\ref{global:phase-data}:
\[
 \begin{aligned}
 \frac{G_j}{\Gamma_j}
 &=\frac{\mu\left(\beta_m\tau_jm^{\frac{3}{2}}\right)^{q-1}}{\lambda\tau_j^{q-1}m^{-\frac{1}{2}}}
 &&\text{by Definition~\ref{global:phase-data}},\\
 &=\frac{\mu}{\lambda}\beta_m^{q-1}
       \frac{\tau_j^{q-1}}{\tau_j^{q-1}}
       m^{\frac{3\left(q-1\right)}{2}+\frac{1}{2}}\\
 &=\frac{\mu}{\lambda}\beta_m^{q-1}m^{\frac{3q-2}{2}}
 &&\text{since }3\left(q-1\right)+1=3q-2.
 \end{aligned}
\]
With $\lambda=K\mu m^{\frac{3q-2}{2}}$ and
$\beta_m\le\sqrt5$ from Lemma~\resultpart{collar:chain-lemma}{collar:chain-lemma:part:geometry} in Appx.~\ref{app:paper-phases},
\[
 \frac{G_j}{\Gamma_j}\le\frac{\mu}{\lambda}5^{\frac{q-1}{2}}m^{\frac{3q-2}{2}}
 =\frac{5^{\frac{q-1}{2}}}{K}\le c_0.
\]
\end{proof}

\subsubsection{The two next-coordinate estimates}
\label{collar:exit}

\begin{lemma}[A common coefficient for the two next-center bounds]
\label{collar:exit-lemma}
Fix a nonterminal phase $0\le j<L-1$. Assume
\[
 \begin{gathered}
 G_j\le c_0\Gamma_j,\qquad
 s_j=\Gamma_j+\frac{G_jt_j^\star}{R_j},\qquad \alpha=\lambda m^{1-\frac{3q}{2}},\\
 s_{j+1}=\eta^{q-1}s_j,\qquad t_{j+1}^\star=\eta t_j^\star,\\
 \eta^q\left(1+c_0\right)^{q^*}\le c_0,\qquad
 \eta^q\left(2\left(1+c_0\right)+\frac{2^q}{q}\right)\le c_0.
 \end{gathered}
\]
The restrictions on $\eta$ depend only on $q$. Then:
\begin{claims}
\item\label{collar:exit-lemma:part:dual}
\textbf{Dual bound Eq.~\eqref{collar:exit-dual}.}
\[
 \alpha^{-\frac{1}{q-1}}s_{j+1}^{q^*}\le c_0\mathcal E_j.
\]
\item\label{collar:exit-lemma:part:endpoint}
\textbf{Endpoint bound Eq.~\eqref{collar:exit-port}.}
\[
 s_{j+1}t_{j+1}^\star+\frac{\alpha}{q}\left(t_{j+1}^\star\right)^q\le c_0\mathcal E_j.
\]
\end{claims}
\end{lemma}
\begin{proof}
\begin{itemize}
\item \textbf{Proof of \ref{collar:exit-lemma:part:dual}.}
Substitute $\alpha=\lambda m^{1-\frac{3q}{2}}$ from Eq.~\eqref{collar:alpha}
and $\Gamma_j=\frac{\lambda\tau_j^{q-1}}{\sqrt{m}}$:
\begin{align}
\alpha^{-\frac{1}{q-1}}\Gamma_j^{q^*}
&=\left(\lambda m^{1-\frac{3q}{2}}\right)^{-\frac{1}{q-1}}
\left(\lambda\tau_j^{q-1}m^{-\frac{1}{2}}\right)^{\frac{q}{q-1}}\notag\\
&=\lambda^{-\frac{1}{q-1}+\frac{q}{q-1}}
 \tau_j^{\left(q-1\right)\frac{q}{q-1}}
 m^{\frac{\frac{3q}{2}-1}{q-1}-\frac{q}{2\left(q-1\right)}}\notag\\
&=\lambda\tau_j^q m^{\frac{q-1}{q-1}}\notag\\
&=\lambda m\tau_j^q=\mathcal E_j.
\label{collar:balanced-exit}
\end{align}
The third line uses $-1+q=q-1$ for the power of $\lambda$ and
$\frac{3q}{2}-1-\frac{q}{2}=q-1$ for the power of $m$.
By Lemma~\resultpart{collar:chain-lemma}{collar:chain-lemma:part:geometry} in Appx.~\ref{app:paper-phases},
$0<\frac{t_j^\star}{R_j}\le1$, so
\[
 \Gamma_j\le s_j=\Gamma_j+\frac{G_jt_j^\star}{R_j}
 \le\Gamma_j+G_j\le\left(1+c_0\right)\Gamma_j.
\]
Using \(s_{j+1}=\eta^{q-1}s_j\) and
Eq.~\eqref{collar:balanced-exit},
\[
 \begin{aligned}
 \alpha^{-\frac{1}{q-1}}s_{j+1}^{q^*}
 &=\eta^{\left(q-1\right)q^*}\alpha^{-\frac{1}{q-1}}s_j^{q^*}\\
 &\le\eta^q\left(1+c_0\right)^{q^*}
       \alpha^{-\frac{1}{q-1}}\Gamma_j^{q^*}
 &&\text{since }\left(q-1\right)q^*=q,\\
 &=\eta^q\left(1+c_0\right)^{q^*}\mathcal E_j
 &&\text{by Eq.~\eqref{collar:balanced-exit}},\\
 &\le c_0\mathcal E_j
 &&\text{by the first bound on }\eta.
 \end{aligned}
\]

\item \textbf{Proof of \ref{collar:exit-lemma:part:endpoint}.}
The bound \(s_j\le\left(1+c_0\right)\Gamma_j\) proved in
part~\ref{collar:exit-lemma:part:dual}, and the definitions of
\(\Gamma_j,t_j^\star,\alpha,\mathcal E_j\), give
\[
 \begin{aligned}
 \frac{\Gamma_j t_j^\star}{\mathcal E_j}
 &=\frac{\left(\frac{\lambda\tau_j^{q-1}}{\sqrt{m}}\right)
         \left(\sqrt m\left(m+1\right)\tau_j\right)}{\lambda m\tau_j^q}\\
 &=\frac{m+1}{m}\le2,\\
 \frac{\alpha\left(t_j^\star\right)^q}{\mathcal E_j}
 &=\frac{\lambda m^{1-\frac{3q}{2}}
         m^{\frac{q}{2}}\left(m+1\right)^q\tau_j^q}{\lambda m\tau_j^q}\\
 &=m^{-q}\left(m+1\right)^q
 =\left(\frac{m+1}{m}\right)^q\le2^q.
 \end{aligned}
\]
Since \(t_{j+1}^\star=\eta t_j^\star\), it follows that
\[
 \begin{aligned}
 s_{j+1}t_{j+1}^\star+\frac{\alpha}{q}\left(t_{j+1}^\star\right)^q
 &=\eta^q\left(s_jt_j^\star+\frac{\alpha}{q}\left(t_j^\star\right)^q\right)\\
 &\le\eta^q\left(\left(1+c_0\right)\Gamma_jt_j^\star
                   +\frac{\alpha}{q}\left(t_j^\star\right)^q\right)\\
 &\le\eta^q\left(2\left(1+c_0\right)+\frac{2^q}{q}\right)\mathcal E_j\\
 &\le c_0\mathcal E_j.
 \end{aligned}
\]
The second line uses the bound on $s_j$, the third uses the two
ratios just proved, and the last uses the second hypothesis on $\eta$.
\end{itemize}
\end{proof}

\subsubsection{Verify the common parameter choice}
\begin{proof}[Proof of Proposition~\ref{global:parameters}]
\begin{itemize}
\item \textbf{Proof of \ref{global:parameters:part:base}.}
Corollary~\ref{collar:coefficient-corollary} gives $G_j\le c_0\Gamma_j$. Since $0<\frac{t_j^\star}{R_j}\le1$,
\[
 \Gamma_j\le s_j=\Gamma_j+\frac{G_jt_j^\star}{R_j}
 \le\Gamma_j+G_j\le\left(1+c_0\right)\Gamma_j.
\]

\item \textbf{Proof of \ref{global:parameters:part:next}.}
The choice $\alpha=\lambda m^{-a_q}=\lambda m^{1-\frac{3q}{2}}$ is exactly
Eq.~\eqref{collar:alpha}.
For $j<L-1$,
Definition~\ref{global:phase-data} gives
\[
 G_{j+1}=\eta^{q-1}G_j,\qquad \Gamma_{j+1}=\eta^{q-1}\Gamma_j,
 \qquad t_{j+1}^\star=\eta t_j^\star,\qquad R_{j+1}=\eta R_j.
\]
Substituting all four relations into the definition of $s_{j+1}$ gives
\[
 \begin{aligned}
 s_{j+1}
 &=\Gamma_{j+1}+G_{j+1}\frac{t_{j+1}^\star}{R_{j+1}}\\
 &=\eta^{q-1}\Gamma_j+
   \eta^{q-1}G_j\frac{\eta t_j^\star}{\eta R_j}\\
 &=\eta^{q-1}\left(\Gamma_j+G_j\frac{t_j^\star}{R_j}\right)
 =\eta^{q-1}s_j.
 \end{aligned}
\]
This verifies the scaling assumptions of Lemma~\ref{collar:exit-lemma}.
It remains to check its two bounds on $\eta$.
Lemma~\ref{base:reset-lemma} in Appx.~\ref{app:paper-base} gives $\gamma<\frac{\theta}{2}$ and $\theta\le\frac{1}{16}$,
so
\[
 \eta=\frac{\gamma}{2\left(1+\theta\right)}
 <\frac{\gamma}{2}<\frac{\theta}{4}\le\frac{1}{64}<\frac{1}{32}.
\]
Also $c_0=\frac{2^{-q-3}}{q}<1$ by
Eq.~\eqref{collar:constants}, and $q^*<2$. Hence
\[
 \left(1+c_0\right)^{q^*}\le4\le2^{q+1},\qquad
 2\left(1+c_0\right)+\frac{2^q}{q}\le4+\frac{2^q}{q}\le2^{q+1}.
\]
Both hypotheses of Lemma~\ref{collar:exit-lemma} follow from
$\eta<2^{-5}$, $c_0^{-1}=q2^{q+3}$, and
\[
 \begin{aligned}
 &\frac{\eta^q}{c_0}
 \max\left\{\left(1+c_0\right)^{q^*},\ 2\left(1+c_0\right)+\frac{2^q}{q}\right\}\\
 &\qquad\le 2^{-5q}\,q2^{q+3}\,2^{q+1}
 =q2^{4-3q}\le2^{3-2q}<1.
 \end{aligned}
\]
The penultimate inequality uses $q\le2^{q-1}$ for $q\ge2$.
Apply Lemma~\resultpart{collar:exit-lemma}{collar:exit-lemma:part:dual}--\ref{collar:exit-lemma:part:endpoint} to obtain
Eq.~\eqref{global:condition-next-gradient}--\eqref{global:condition-next-value}.
For $j=L-1$, $s_L=t_L^\star=0$, so the left sides of
Eq.~\eqref{global:condition-next-gradient} and \eqref{global:condition-next-value}
are zero.

\item \textbf{Proof of \ref{global:parameters:part:nesting}.}
Definition~\ref{global:phase-data} gives orthogonal $\bm h_j$ and $R_j=\eta^jR_0$.
The chosen $\eta$ satisfies Eq.~\eqref{base:etaconstraint} with equality;
Lemma~\resultpart{base:nesting-lemma}{base:nesting-lemma:part:construction} and
\eqref{base:containment} give the stated containment.
\end{itemize}
\end{proof}

\subsubsection{Define the complete objective}\label{global:define-objective}

\begin{definition}[The complete objective]\label{global:function-definition}
With the parameters of Proposition~\ref{global:parameters}, define
\begin{equation}\label{global:objective}
 \ell_{\mathrm{init}}\left(\bm x\right):=-s_0\ip{\bm v_0}{\bm x},\qquad
 F\left(\bm x\right):=B\left(\bm x\right)+\ell_{\mathrm{init}}\left(\bm x\right)+\sum_{j=0}^{L-1}A_j\left(\bm x\right).
\end{equation}
Here $B=B_L$, and each $A_j$ is the phase function of Appx.~\ref{collar:affine}
evaluated at the local coordinates
\[
 \bm w=P_{E_j}\left(\bm x-\bm c_j\right)
       =t_j\left(\bm x\right)\bm v_j+\sum_{i=1}^m x_{j,i}\left(\bm x\right)\bm u_{j,i},
\]
where, for $\bm x\in\R^d$,
\[
 \begin{aligned}
 t_j\left(\bm x\right)&:=\ip{\bm v_j}{\bm x-\bm c_j},\\
 x_{j,i}\left(\bm x\right)&:=\ip{\bm u_{j,i}}{\bm x-\bm c_j}\quad\left(1\le i\le m\right).
 \end{aligned}
\]
For $j<L-1$, orthogonality of $\bm h_j$ and $\bm v_{j+1}$ gives
$t_{j+1}\left(\bm x\right)
=\ip{\bm v_{j+1}}{\bm x-\bm c_{j+1}}
=\ip{\bm v_{j+1}}{\bm x-\bm c_j}$.
We retain $P_j,I_j,O_j,A_j$ for the ambient functions obtained from
the local inputs $\left(\bm w,t_{j+1}\left(\bm x\right)\right)$
in Eq.~\eqref{paper:ambient-phase}. Their formulas are
\begin{align*}
 P_j\left(\bm x\right)&:=\frac{\lambda}{q}\left(\left|\frac{t_j\left(\bm x\right)}{\sqrt m}-x_{j,1}\left(\bm x\right)\right|^q\right.\\
           &\qquad\left.+\sum_{i=1}^{m-1}|x_{j,i}\left(\bm x\right)-x_{j,i+1}\left(\bm x\right)|^q
           +|x_{j,m}\left(\bm x\right)|^q\right)+\frac{\alpha}{q}|t_{j+1}\left(\bm x\right)|^q,\\
 I_j\left(\bm x\right)&:=\ip{\bm g_j}{\bm x-\bm c_{j+1}},\\
 O_j\left(\bm x\right)&:=P_j\left(\bm x\right)-P_j\left(\bm c_{j+1}\right)-4\delta_j,\\
 A_j\left(\bm x\right)&:=I_j\left(\bm x\right)+\psi_{\delta_j}\left(O_j\left(\bm x\right)-I_j\left(\bm x\right)\right).
\end{align*}
Here $P_j\left(\bm c_{j+1}\right)=\frac{\lambda\left(m+1\right)\tau_j^q}{q}$ and
$A_j=M_{\delta_j}\left(I_j,O_j\right)$. Since
$\bm g_j\in E_j\oplus\operatorname{span}\{\bm v_{j+1}\}$,
with the last summand omitted in the terminal phase, the displayed $I_j$ agrees
with the local affine function under the same coordinate map.
The last phase omits $t_L$ and its term; when $m=1$, the sum over $i$ is empty.

Orthogonality gives the following local inputs:
\[
\begin{array}{c|cc}
 \text{Ambient point}&\bm w=P_{E_j}\left(\bm x-\bm c_j\right)&t_{j+1}\left(\bm x\right)\\\hline
 \bm c_j&0&0\\
 \bm c_{j+1}&\bm h_j&0\\
 \bm c_{j+2}\ \left(j<L-1\right)&\bm h_j&t_{j+1}^\star
\end{array}
\]
For the terminal phase the last coordinate and last row are omitted.
The endpoint value of the first coordinate is $t_j\left(\bm c_{j+1}\right)=\ip{\bm v_j}{\bm h_j}=t_j^\star$.
\end{definition}

\subsection{Verification of the complete objective}\label{global:verification}

\subsubsection{Retain uniform convexity and smoothness}\label{global:class-membership}

\begin{lemma}[Uniform convexity of the complete objective]\label{global:uc-lemma}
The complete objective $F$ in Eq.~\eqref{global:objective} satisfies
\[
 \Breg F{\bm y}{\bm x}\ge c_B\mu\norm{\bm y-\bm x}^q,
 \qquad c_B=\frac{1}{2q\cdot16^{q-2}}.
\]
\end{lemma}
\begin{proof}
Theorem~\resultpart{collar:theorem}{collar:theorem:part:convexity}
gives convexity in the local phase coordinates.
Since the coordinate map is affine,
Lemma~\resultpart{sources:convex-rules}{sources:convex-rules:affine}
in Appx.~\ref{app:paper-elementary} makes each ambient $A_j$ convex.
Lemma~\ref{paper:bregman-additions} applies to the convex sum
$G=\sum_jA_j$ and affine term $\ell_{\mathrm{init}}$.
Together with Lemma~\resultpart{base:uc-lemma}{base:uc-lemma:part:remainder}
in Appx.~\ref{app:paper-base}, it gives
\[
 \begin{aligned}
 \Breg F{\bm y}{\bm x}
 &=\Breg B{\bm y}{\bm x}
   +\Breg{\ell_{\mathrm{init}}}{\bm y}{\bm x}
   +\sum_j\Breg{A_j}{\bm y}{\bm x}
 &&\text{by Eq.~\eqref{global:objective} and the Bregman definition},\\
 &=\Breg B{\bm y}{\bm x}+\sum_j\Breg{A_j}{\bm y}{\bm x}
 &&\text{since }\ell_{\mathrm{init}}\text{ is affine},\\
 &\ge\Breg B{\bm y}{\bm x}
 &&\text{by Lemma~\ref{paper:bregman-additions}},\\
 &\ge c_B\mu\norm{\bm y-\bm x}^q
 &&\text{by Lemma~\resultpart{base:uc-lemma}{base:uc-lemma:part:remainder}}.
 \end{aligned}
\]
The modulus is therefore at least $qc_B\mu$ under the
uniform-convexity convention of Theorem~\ref{paper:main-theorem}.
\end{proof}

\begin{lemma}[Smoothness independent of the number of phases]\label{global:smoothness-lemma}
For $F$ in Eq.~\eqref{global:objective}, we have
$F\in C^{p,\nu}\left(\R^d\right)$ and, for all $\bm x,\bm y\in\R^d$,
\[
 \norm{\D^pF\left(\bm x\right)-\D^pF\left(\bm y\right)}
 \le C_H\mu m^{\frac{3q-2}{2}}\norm{\bm x-\bm y}^{\nu}.
\]
The constant $C_H:=C_B+2C_AK$ depends only on $p,\nu$.
\end{lemma}
\begin{proof}
All terms are $C^p$ by Lemma~\resultpart{base:nesting-lemma}{base:nesting-lemma:part:construction} in Appx.~\ref{app:paper-base}
and Theorem~\resultpart{collar:theorem}{collar:theorem:part:regularity} in Appx.~\ref{app:paper-phases}.
Set
\[
 V_j=E_j\oplus\operatorname{span}\{\bm v_{j+1}\}\ \left(j<L-1\right),\qquad
 V_{L-1}=E_{L-1},\qquad \Pi_j=P_{V_j}.
\]
Then $|j-k|\ge2\Longrightarrow V_j\perp V_k$. Thus each parity class
$J=\{j:0\le j<L,\ j\equiv a\pmod2\}$, $a\in\{0,1\}$, is orthogonal.

Fix $\bm x,\bm y$ and a nonempty parity class $J$. Define the tensor on $V_j$ by
\[
 \mathcal Q_j:=\left.\left(\D^pA_j\left(\bm x\right)-\D^pA_j\left(\bm y\right)\right)\right|_{V_j^p}.
\]
The linear chain rule and the phase bound in Theorem~\resultpart{collar:theorem}{collar:theorem:part:regularity} give
\[
 \begin{aligned}
 \left(\D^pA_j\left(\bm x\right)-\D^pA_j\left(\bm y\right)\right)[\bm z_1,\ldots,\bm z_p]
 &=\mathcal Q_j[\Pi_j\bm z_1,\ldots,\Pi_j\bm z_p],\\
 \norm{\mathcal Q_j}
 &\le C_A\lambda\norm{\Pi_j\left(\bm x-\bm y\right)}^\nu.
 \end{aligned}
\]
Each tensor acts only on $V_j$, and these subspaces are orthogonal within $J$.
Lemma~\ref{elem:orthogonal-sum} in Appx.~\ref{app:paper-elementary}
therefore controls their sum by the largest phase contribution, avoiding an $L$ factor.
Applied with order $p$ and tensors $\mathcal Q_j$ on $V_j^p$, it gives
\[
 \begin{aligned}
 \left\|\sum_{j\in J}\left(\D^pA_j\left(\bm x\right)-\D^pA_j\left(\bm y\right)\right)\right\|
 &=\max_{j\in J}\norm{\mathcal Q_j}\\
 &\le C_A\lambda\max_{j\in J}\norm{\Pi_j\left(\bm x-\bm y\right)}^\nu\\
 &\le C_A\lambda\norm{\bm x-\bm y}^\nu.
 \end{aligned}
\]
For $J=\varnothing$, the sum is $0$ and satisfies the same bound.
Adding the two parity bounds and the base
bound from Lemma~\resultpart{base:global-holder-lemma}{base:global-holder-lemma:part:holder}, with $\D^p\ell_{\mathrm{init}}=0$, yields
\begin{equation}\label{global:holder}
 \begin{aligned}
 \norm{\D^pF\left(\bm x\right)-\D^pF\left(\bm y\right)}
 &\le\left(C_B\mu+2C_A\lambda\right)\norm{\bm x-\bm y}^{\nu}\\
 &=\left(C_B\mu+2C_AK\mu m^{a_q}\right)\norm{\bm x-\bm y}^{\nu}\\
 &\le\left(C_B+2C_AK\right)\mu m^{a_q}\norm{\bm x-\bm y}^{\nu}
 &&\left(m^{a_q}\ge1\right),\\
 &=C_H\mu m^{a_q}\norm{\bm x-\bm y}^{\nu}.
 \end{aligned}
\end{equation}
\end{proof}

\subsubsection{Verify the center gradients and identify the minimizer}
\label{global:center-verification}
The affine slopes telescope in the center gradients. These gradients also
provide the orthogonality used in the prefix-gap estimate.

\begin{lemma}[All center gradients and the minimizer]\label{global:center-lemma}
For the objective in Definition~\ref{global:function-definition} and
$0\le k<L$, $0\le j\le L$:
\begin{claims}
\item\label{global:center-lemma:part:phase-data} \textbf{Phase derivatives.}
\[
 \D A_k\left(\bm c_j\right)=
 \begin{cases}\bm g_k,&k<j,\\0,&k\ge j,\end{cases}
 \qquad
 \D^\ell A_k\left(\bm c_j\right)=0\quad\left(2\le\ell\le p\right).
\]
\item\label{global:center-lemma:part:gradients} \textbf{Full gradients.}
\[
 \D F\left(\bm c_j\right)=\begin{cases}-s_j\bm v_j,&j<L,\\0,&j=L.\end{cases}
\]
\item\label{global:center-lemma:part:minimizer} \textbf{Minimizer.}
\[
 \operatorname*{argmin}_{\bm x\in\R^d}F\left(\bm x\right)=\{\bm c_L\}.
\]
\end{claims}
\end{lemma}
\begin{proof}
\begin{itemize}
\item \textbf{Proof of \ref{global:center-lemma:part:phase-data}.}
Definition~\ref{global:phase-data} and orthogonality give three cases for the arguments of $A_k$ at $\bm c_j$:
\[
\begin{array}{c|c|cc}
 \text{Case}&k&P_{E_k}\left(\bm c_j-\bm c_k\right)&\ip{\bm v_{k+1}}{\bm c_j-\bm c_k}\\\hline
 1&k\ge j&0&0\\
 2&k=j-1&\bm h_k&0\\
 3&k\le j-2&\bm h_k&t_{k+1}^\star
\end{array}
\]
Omit the second coordinate for $k=L-1$, as in Theorem~\resultpart{collar:theorem}{collar:theorem:part:terminal} in Appx.~\ref{app:paper-phases}.
Apply its origin formula (part~\ref{collar:theorem:part:origin}) and affine endpoint formula (part~\ref{collar:theorem:part:endpoints}):
\[
 \begin{aligned}
 k\ge j&:\quad \D^\ell A_k\left(\bm c_j\right)=0 &&\left(1\le\ell\le p\right),\\
 k<j&:\quad A_k=I_k\text{ near }\bm c_j,\\
 &\qquad\D A_k\left(\bm c_j\right)=\bm g_k,\quad
       \D^\ell A_k\left(\bm c_j\right)=0 &&\left(2\le\ell\le p\right).
 \end{aligned}
\]
\item \textbf{Proof of \ref{global:center-lemma:part:gradients}.}
For $1\le j\le L$, use
Lemma~\resultpart{base:nesting-lemma}{base:nesting-lemma:part:center-formulas}
in Appx.~\ref{app:paper-base},
part~\ref{global:center-lemma:part:phase-data}, and
Eq.~\eqref{global:cancellation}:
\[
\begin{aligned}
 \D F\left(\bm c_j\right)
 &=\sum_{k=0}^{j-1}\bm b_k-s_0\bm v_0+\sum_{k=0}^{j-1}\bm g_k\\
 &=-s_0\bm v_0+\sum_{k=0}^{j-1}
       \left(\bm b_k+s_k\bm v_k-s_{k+1}\bm v_{k+1}-\bm b_k\right)\\
 &=-s_0\bm v_0+\sum_{k=0}^{j-1}s_k\bm v_k
                    -\sum_{k=0}^{j-1}s_{k+1}\bm v_{k+1}\\
 &=-s_0\bm v_0+
       \left(s_0\bm v_0+\sum_{k=1}^{j-1}s_k\bm v_k\right)
       -\left(\sum_{k=1}^{j-1}s_k\bm v_k+s_j\bm v_j\right)\\
 &=-s_j\bm v_j.
\end{aligned}
\]
In the fourth line, the second sum is reindexed from $k+1$ to $k$.
At $j=L$, the convention $s_L\bm v_L=0$ gives $\D F\left(\bm c_L\right)=0$.
At $j=0$, both sums in the first line are empty and
$\D B\left(\bm c_0\right)=0$, giving $\D F\left(\bm c_0\right)=-s_0\bm v_0$ directly.
\item \textbf{Proof of \ref{global:center-lemma:part:minimizer}.}
Lemma~\ref{global:uc-lemma} at $\bm c_L$, together with the zero gradient in part~\ref{global:center-lemma:part:gradients}, gives
\[
 F\left(\bm x\right)\ge F\left(\bm c_L\right)+c_B\mu\norm{\bm x-\bm c_L}^q.
\]
The right side is strictly larger than $F\left(\bm c_L\right)$ when $\bm x\ne\bm c_L$.
\end{itemize}
\end{proof}

\subsubsection{A convex gap comparison}
\begin{lemma}[A convex comparison proves an optimality gap]\label{design:gap}
Let $W\subset\R^d$, $e>0$, and $\bm c\in\R^d$. Suppose $F=\Psi+r$, $\Psi$ is convex and differentiable, and $r=0$ on a neighborhood of $\bm c$. Assume also that $r\left(\bm x\right)\ge e$ and $\ip{\D F\left(\bm c\right)}{\bm x-\bm c}=0$ for every $\bm x\in W$. Then:
\begin{claims}
\item\label{design:gap:part:contact} $\Psi\left(\bm c\right)=F\left(\bm c\right)$ and $\D\Psi\left(\bm c\right)=\D F\left(\bm c\right)$.
\item\label{design:gap:part:gap} $F\left(\bm x\right)\ge F\left(\bm c\right)+e\ge\inf F+e$ for every $\bm x\in W$.
\end{claims}
\end{lemma}
\begin{proof}
\begin{itemize}
\item \textbf{Proof of \ref{design:gap:part:contact}.}
The neighborhood identity $F=\Psi$ gives equality of values and gradients at $\bm c$.
\item \textbf{Proof of \ref{design:gap:part:gap}.}
By the first-order convexity inequality in
Lemma~\resultpart{sources:convex-rules}{sources:convex-rules:first},
for $\bm x\in W$, part~\ref{design:gap:part:contact} and the two assumptions on $W$ give
\[
 \begin{aligned}
 F\left(\bm x\right)
 &=\Psi\left(\bm x\right)+r\left(\bm x\right)
 &&\text{by the decomposition},\\
 &\ge\Psi\left(\bm c\right)+\ip{\D\Psi\left(\bm c\right)}{\bm x-\bm c}+r\left(\bm x\right)
 &&\text{by Lemma~\resultpart{sources:convex-rules}{sources:convex-rules:first}},\\
 &=F\left(\bm c\right)+\ip{\D F\left(\bm c\right)}{\bm x-\bm c}+r\left(\bm x\right)
 &&\text{by part~\ref{design:gap:part:contact}},\\
 &\ge F\left(\bm c\right)+e
 &&\text{by the two assumptions on }W,\\
 &\ge\inf F+e.
 \end{aligned}
\]
\end{itemize}
\end{proof}

\subsubsection{Optimality gaps on short prefixes}\label{global:gradient-and-gap}

For each $0\le j<L$ and $0\le b\le\lfloor \frac{m}{4}\rfloor$, define the subspace containing all earlier phases and the selected part of phase $j$:
\[
 W_{j,b}:=\left(\bigoplus_{k<j}E_k\right)
 \oplus\operatorname{span}\{\bm v_j,\bm u_{j,1},\ldots,\bm u_{j,b}\}.
\]
Using $A_j$ and $I_j$ from Definition~\ref{global:function-definition}, define
\[
 r_j:=A_j-I_j,\qquad
 \Psi_j:=F-A_j+I_j
 =B+\ell_{\mathrm{init}}+\sum_{k\ne j}A_k+I_j.
\]

\subsubsection{Convex comparison for the complete objective}\label{app:assembly-support}

\begin{lemma}[Convex comparison with local equality]\label{global:comparison-lemma}
For $0\le j<L$, use $F$, $r_j=A_j-I_j$, and $\Psi_j=F-A_j+I_j$ as defined in Secs.~\ref{global:define-objective} and \ref{global:gradient-and-gap}. Then
\[
 \begin{gathered}
 \Psi_j\text{ is convex and }C^p,\qquad r_j=0\text{ near }\bm c_{j+1},\\
 \Psi_j\left(\bm c_{j+1}\right)=F\left(\bm c_{j+1}\right),\qquad
 \D\Psi_j\left(\bm c_{j+1}\right)=\D F\left(\bm c_{j+1}\right).
 \end{gathered}
\]
\end{lemma}
\begin{proof}
Lemma~\resultpart{base:nesting-lemma}{base:nesting-lemma:part:construction} in Appx.~\ref{app:paper-base} and
Theorem~\resultpart{collar:theorem}{collar:theorem:part:convexity}, \ref{collar:theorem:part:regularity} in Appx.~\ref{app:paper-phases}
give convexity and $C^p$ regularity of $B$ and the $A_k$.
The convex-sum rule in
Lemma~\resultpart{sources:convex-rules}{sources:convex-rules:sum}
applies to
\[
 \Psi_j=B+\sum_{k\ne j}A_k+\underbrace{\ell_{\mathrm{init}}+I_j}_{\text{affine}},
\]
so $\Psi_j$ is convex and $C^p$.
Theorem~\resultpart{collar:theorem}{collar:theorem:part:endpoints} gives
$A_j=I_j$ near $\bm c_{j+1}$. On that neighborhood,
\[
 r_j=A_j-I_j=0,\qquad \Psi_j=F-r_j=F.
\]
Evaluating and differentiating the second equality proves the center identities.
\end{proof}
\Needspace{8\baselineskip}
\begin{lemma}[Gap on every short prefix]\label{global:gap-lemma}
For the objective $F$ in Definition~\ref{global:function-definition}:
\begin{claims}
\item\label{global:gap-lemma:part:prefix} \textbf{Gap for every point in an early prefix.}
\begin{equation}\label{global:prefix-gap}
 \bm x\in W_{j,b},\quad b\le\lfloor \frac{m}{4}\rfloor
 \quad\Longrightarrow\quad F\left(\bm x\right)-\min F\ge7\delta_j.
\end{equation}
\end{claims}
\end{lemma}
\begin{proof}
\begin{itemize}
\item \textbf{Proof of \ref{global:gap-lemma:part:prefix}.}
Fix $\bm x\in W_{j,b}$. Lemma~\ref{global:comparison-lemma}
gives the comparison hypotheses:
\[
 \begin{gathered}
 F=\Psi_j+r_j,\qquad \Psi_j\text{ is convex and }C^1,\qquad
 r_j=0\text{ near }\bm c_{j+1},\\
 \Psi_j\left(\bm c_{j+1}\right)=F\left(\bm c_{j+1}\right),\qquad
 \D\Psi_j\left(\bm c_{j+1}\right)=\D F\left(\bm c_{j+1}\right).
 \end{gathered}
\]
The phase coordinates of $\bm x$ satisfy
\[
 P_{E_j}\left(\bm x-\bm c_j\right)\in
 \operatorname{span}\{\bm v_j,\bm u_{j,1},\ldots,\bm u_{j,b}\},
 \qquad t_{j+1}\left(\bm x\right)=0,
\]
with the next coordinate omitted in the terminal phase.
Thus Theorem~\resultpart{collar:theorem}{collar:theorem:part:prefix}
gives $r_j\left(\bm x\right)\ge7\delta_j$.
For $j<L-1$, we have $\bm v_{j+1}\in E_{j+1}$, whereas
$\bm x,\bm c_{j+1}\in\bigoplus_{k\le j}E_k$. Thus
$\ip{\bm v_{j+1}}{\bm x}=\ip{\bm v_{j+1}}{\bm c_{j+1}}=0$.
Lemma~\resultpart{global:center-lemma}{global:center-lemma:part:gradients}
then gives
\[
 \ip{\D F\left(\bm c_{j+1}\right)}{\bm x-\bm c_{j+1}}
 =\begin{cases}
 -s_{j+1}\ip{\bm v_{j+1}}{\bm x-\bm c_{j+1}}=0,&j<L-1,\\
 0,&j=L-1.
 \end{cases}
\]
In the first line $\bm v_{j+1}$ is orthogonal to both points; in the second
$\D F\left(\bm c_L\right)=0$.
The first-order convexity inequality for $\Psi_j$ therefore has zero
linear term on $W_{j,b}$. Lemma~\resultpart{design:gap}{design:gap:part:gap}
applies with
\[
 \left(\Psi,r,W,\bm c,e\right)=\left(\Psi_j,r_j,W_{j,b},\bm c_{j+1},7\delta_j\right).
\]
Its conclusion is
\[
 F\left(\bm x\right)\ge F\left(\bm c_{j+1}\right)+7\delta_j\ge\min F+7\delta_j.
\]
\end{itemize}
\end{proof}

\subsubsection{Compute the initial gap for normalization}\label{global:initial-gap-calibration}

To normalize the objective to the prescribed initial gap, we first
compute its decrease between consecutive centers. The base and phase
value changes give an exact formula for each decrease. Their geometric
decay then determines $D_{m,L}$.

\begin{lemma}[Exact decrease between consecutive centers]\label{global:exact-drop-lemma}
For $0\le j<L$, use the decrease
$\chi_j:=F\left(\bm c_j\right)-F\left(\bm c_{j+1}\right)$ from
Lemma~\resultpart{paper:assembly}{paper:assembly:initial}. Then:
\begin{claims}
\item\label{global:value-lemma:part:exact-drop} \textbf{Exact decrease.}
\begin{equation}\label{global:exact-drop}
 \chi_j=s_jt_j^\star-\beta_q\mu R_j^q
 -\frac{\lambda\left(m+1\right)\tau_j^q}{q}-4\delta_j,
\end{equation}
where $\beta_q=\frac{1}{q}+a_{\mathrm{lift}}\theta^2$.
\end{claims}
\end{lemma}
\begin{proof}
We compare the base, the phase functions, and the initial affine term
at $\bm c_j$ and $\bm c_{j+1}=\bm c_j+\bm h_j$.
The base gradient at $\bm c_j$ lies in the completed phase subspaces,
so its linear contribution along $\bm h_j$ vanishes. Indeed,
Lemma~\resultpart{base:nesting-lemma}{base:nesting-lemma:part:center-formulas}
in Appx.~\ref{app:paper-base} gives
\[
 \begin{aligned}
 \ip{\D B\left(\bm c_j\right)}{\bm h_j}
 &=\left\langle\sum_{k<j}\bm b_k,\bm h_j\right\rangle\\
 &=\sum_{k<j}\mu R_k^{q-2}\ip{\bm h_k}{\bm h_j}
 &&\text{since }\bm b_k=\mu R_k^{q-2}\bm h_k,\\
 &=0
 &&\text{since }E_k\perp E_j\text{ for }k<j.
 \end{aligned}
\]
The center-value recurrence in
Lemma~\resultpart{base:nesting-lemma}{base:nesting-lemma:part:recurrences}
therefore identifies the base contribution as
\[
 \begin{aligned}
 B\left(\bm c_j\right)-B\left(\bm c_{j+1}\right)
 &=-\ip{\D B\left(\bm c_j\right)}{\bm h_j}-\beta_q\mu R_j^q\\
 &=-\beta_q\mu R_j^q
 &&\text{by the zero inner product above}.
 \end{aligned}
\]
For the phase functions, moving from $\bm c_j$ to $\bm c_{j+1}$
changes only the inputs of $A_{j-1}$, when $j\ge1$, and $A_j$.
The other completed phases and the later phases have the same inputs
at both centers. The center evaluations in
Lemma~\resultpart{global:center-lemma}{global:center-lemma:part:phase-data},
using Theorem~\ref{collar:theorem}\,\ref{collar:theorem:part:origin},
\ref{collar:theorem:part:endpoints} in Appx.~\ref{app:paper-phases}, give
\[
 A_k\left(\bm c_j\right)-A_k\left(\bm c_{j+1}\right)=
 \begin{cases}
 0,&k\le j-2\text{ or }k\ge j+1\quad\text{(same inputs)},\\
 0-\left(-s_jt_j^\star\right),&k=j-1,\ j\ge1,\\
 \left(-\frac{\lambda\left(m+1\right)\tau_j^q}{q}-4\delta_j\right)-0,&k=j.
 \end{cases}
\]
When $j=0$, there is no preceding phase, and the initial affine term
supplies the same starting-slope contribution, since
\[
 \ell_{\mathrm{init}}\left(\bm c_j\right)-\ell_{\mathrm{init}}\left(\bm c_{j+1}\right)
 =s_0\ip{\bm v_0}{\bm h_j}
 =\begin{cases}s_0t_0^\star,&j=0,\\0,&j\ge1.\end{cases}
\]
Adding these value changes and then arranging them in the stated order gives
\[
 \begin{aligned}
 \chi_j
 &=\left(B\left(\bm c_j\right)-B\left(\bm c_{j+1}\right)\right)
   +\left(\ell_{\mathrm{init}}\left(\bm c_j\right)
          -\ell_{\mathrm{init}}\left(\bm c_{j+1}\right)\right)\\
 &\qquad+\sum_{k=0}^{L-1}
      \left(A_k\left(\bm c_j\right)-A_k\left(\bm c_{j+1}\right)\right)\\
 &=-\beta_q\mu R_j^q+s_jt_j^\star
   -\frac{\lambda\left(m+1\right)\tau_j^q}{q}-4\delta_j
 &&\text{by the value changes above},\\
 &=s_jt_j^\star-\beta_q\mu R_j^q
   -\frac{\lambda\left(m+1\right)\tau_j^q}{q}-4\delta_j.
 \end{aligned}
\]
The last line is Eq.~\eqref{global:exact-drop}.
\end{proof}

\begin{lemma}[Bounds and exact normalization of the initial gap]\label{global:value-lemma}
For the decreases $\chi_j$ in Lemma~\ref{global:exact-drop-lemma}:
\begin{claims}
\item\label{global:value-lemma:part:bounds} \textbf{Decrease bounds.}
\[
 0<c_B\mu R_j^q\le\chi_j\le2\left(1+c_0\right)\lambda m\tau_j^q
 \qquad\left(0\le j<L\right).
\]
\item\label{global:value-lemma:part:exact-gap} \textbf{Initial gap.}
\[
 \chi_j=\eta^{qj}\chi_0,\qquad
 D_{m,L}:=F\left(0\right)-\min F
 =\chi_0\frac{1-\eta^{qL}}{1-\eta^q}>0.
\]
For every orthogonal matrix $Q$, replacing all $\bm v_j,\bm u_{j,i}$ by
$Q\bm v_j,Q\bm u_{j,i}$ leaves $D_{m,L}$ unchanged.
\end{claims}
\end{lemma}
\begin{proof}
\begin{itemize}
\item \textbf{Proof of \ref{global:value-lemma:part:bounds}.}
Uniform convexity gives a positive decrease between consecutive centers,
while convexity bounds that decrease by the starting slope.
Lemma~\resultpart{global:center-lemma}{global:center-lemma:part:gradients}
and orthogonality give $\ip{\D F\left(\bm c_{j+1}\right)}{\bm h_j}=0$,
including the terminal phase. Hence
\[
 \begin{aligned}
 \Breg F{\bm c_j}{\bm c_{j+1}}
 &=F\left(\bm c_j\right)-F\left(\bm c_{j+1}\right)
   -\ip{\D F\left(\bm c_{j+1}\right)}{\bm c_j-\bm c_{j+1}}\\
 &=\chi_j+\ip{\D F\left(\bm c_{j+1}\right)}{\bm h_j}=\chi_j.
 \end{aligned}
\]
Lemma~\ref{global:uc-lemma} now gives
\[
 \chi_j=\Breg F{\bm c_j}{\bm c_{j+1}}
 \ge c_B\mu\norm{\bm c_j-\bm c_{j+1}}^q
 =c_B\mu R_j^q>0.
\]
The first-order convexity inequality in
Lemma~\resultpart{sources:convex-rules}{sources:convex-rules:first},
now based at $\bm c_j$, gives
\[
 \begin{aligned}
 F\left(\bm c_{j+1}\right)
 &\ge F\left(\bm c_j\right)+\ip{\D F\left(\bm c_j\right)}{\bm h_j},\\
 \chi_j
 &\le-\ip{\D F\left(\bm c_j\right)}{\bm h_j}\\
 &=s_j\ip{\bm v_j}{\bm h_j}
 =s_jt_j^\star.
 \end{aligned}
\]
By Proposition~\resultpart{global:parameters}{global:parameters:part:base},
\[
 s_jt_j^\star\le\left(1+c_0\right)\Gamma_jt_j^\star
 =\left(1+c_0\right)\lambda\left(m+1\right)\tau_j^q\le2\left(1+c_0\right)\lambda m\tau_j^q.
\]
Therefore
\begin{equation}\label{global:drop}
 0<c_B\mu R_j^q\le F\left(\bm c_j\right)-F\left(\bm c_{j+1}\right)
 \le s_jt_j^\star\le2\left(1+c_0\right)\lambda m\tau_j^q.
\end{equation}
\item \textbf{Proof of \ref{global:value-lemma:part:exact-gap}.}
The phase parameters have the same scaling in every phase.
Definition~\ref{global:phase-data} and the formulas for $s_j$ and $\delta_j$ give
\[
 s_j=\eta^{\left(q-1\right)j}s_0,\quad t_j^\star=\eta^jt_0^\star,\quad
 \tau_j=\eta^j\tau_0,\quad R_j=\eta^jR_0,\quad \delta_j=\eta^{qj}\delta_0.
\]
Each term in Eq.~\eqref{global:exact-drop} therefore scales by $\eta^{qj}$, giving
\[
 \begin{aligned}
 \chi_j
 &=\left(\eta^{\left(q-1\right)j}s_0\right)\left(\eta^jt_0^\star\right)
   -\frac{\lambda}{q}\left(m+1\right)\left(\eta^j\tau_0\right)^q
   -4\eta^{qj}\delta_0\\
 &\qquad-\beta_q\mu\left(\eta^jR_0\right)^q\\
 &=\eta^{qj}\left(s_0t_0^\star
   -\frac{\lambda}{q}\left(m+1\right)\tau_0^q-4\delta_0
   -\beta_q\mu R_0^q\right)\\
 &=\eta^{qj}\chi_0
 &&\text{by Eq.~\eqref{global:exact-drop} at }j=0.
 \end{aligned}
\]
Since $\bm c_0=\bm0$ and $\bm c_L$ is the minimizer by
Lemma~\resultpart{global:center-lemma}{global:center-lemma:part:minimizer},
summing these decreases gives the exact initial gap
\begin{equation}\label{global:initial}
 \begin{aligned}
 D_{m,L}:=F\left(0\right)-F\left(\bm c_L\right)
 &=\sum_{j=0}^{L-1}\left(F\left(\bm c_j\right)-F\left(\bm c_{j+1}\right)\right)
 &&\text{since }\bm c_0=0,\\
 &=\sum_{j=0}^{L-1}\chi_j
 =\chi_0\sum_{j=0}^{L-1}\eta^{qj}\\
 &=\chi_0\frac{1-\eta^{qL}}{1-\eta^q}
 &&\text{by the finite geometric sum},\\
 &>0.
 \end{aligned}
\end{equation}
Positivity follows from part~\ref{global:value-lemma:part:bounds}, $L\ge1$, and $0<\eta<1$.
The values $s_0,t_0^\star,R_0,\tau_0,\delta_0$ in the formula for
$\chi_0$ depend only on the fixed scalar parameters, by
Definition~\ref{global:phase-data}. Consequently, the formula for
$D_{m,L}$ is independent of the orthonormal directions, proving
the stated invariance under $Q$.
\end{itemize}
\end{proof}

\subsection{Proof of Lemma~\ref{paper:assembly}}
\label{app:assembly-full-proof}

\begin{proof}[Proof of Lemma~\ref{paper:assembly}]
\textbf{Proof of \ref{paper:assembly:class}.}
Since the added phase functions are convex,
Lemma~\ref{paper:bregman-additions} with
$G=\sum_{j=0}^{L-1}A_j$ and $\ell=\ell_{\mathrm{init}}$ gives
$D_F\ge D_B$. The bound for the base in
Lemma~\resultpart{paper:base}{paper:base:class} therefore yields
\[
 \Breg F{\bm y}{\bm x}
 \ge\Breg B{\bm y}{\bm x}
 \ge c_B\mu\norm{\bm y-\bm x}^q.
\]
Summing phase smoothness bounds would cost a factor $L$.
Instead, each $A_j$ depends only on
$V_j=E_j\oplus\operatorname{span}\left\{\bm v_{j+1}\right\}$,
with $V_{L-1}=E_{L-1}$. These subspaces are mutually orthogonal
within each parity group
$J_a=\left\{0\le j<L:j\equiv a\pmod2\right\}$
for $a\in\{0,1\}$, allowing a maximum bound on the derivative differences
$\mathcal Q_j=\D^pA_j\left(\bm x\right)-\D^pA_j\left(\bm y\right)$.
Each $\mathcal Q_j$ acts only on the projections of its arguments
onto $V_j$, as verified in Lemma~\ref{global:smoothness-lemma}.
For nonempty $J_a$,
\[
 \begin{aligned}
 \left\|\sum_{j\in J_a}\mathcal Q_j\right\|
 &=\max_{j\in J_a}\norm{\mathcal Q_j}
 &&\text{by Lemma~\ref{elem:orthogonal-sum}},\\
 &\le C_A\lambda\norm{\bm x-\bm y}^\nu
 &&\text{by Lemma~\resultpart{paper:phase}{paper:phase:class}}.
 \end{aligned}
\]
The same bound holds with a zero sum when $J_a=\varnothing$.
Adding both groups and the base bound gives
\[
 \begin{aligned}
 \norm{\D^pF\left(\bm x\right)-\D^pF\left(\bm y\right)}
 &\le\left(C_B\mu+2C_A\lambda\right)\norm{\bm x-\bm y}^\nu
 &&\text{by Lemma~\resultpart{paper:base}{paper:base:class}},\\
 &\le C_H\mu m^{a_q}\norm{\bm x-\bm y}^\nu
 &&\text{since }\lambda=K\mu m^{a_q},\ m\ge1.
 \end{aligned}
\]
All terms of $F$ are $C^p$ by Lemmas~\ref{paper:base} and
\ref{paper:phase}, so this proves both stated regularity bounds.

\textbf{Proof of \ref{paper:assembly:centers}.}
By Lemma~\ref{paper:phase}\,\ref{paper:phase:completed}--\ref{paper:phase:starting},
completed phases contribute $\bm g_k$ and the others $0$ at
$\bm c_j$. Adding the base increments gives
\[
 \begin{aligned}
 \D F\left(\bm c_j\right)
 &=-s_0\bm v_0+\sum_{k<j}\left(\bm b_k+\bm g_k\right)
 &&\text{by Lemma~\resultpart{paper:base}{paper:base:centers}},\\
 &=-s_0\bm v_0+\sum_{k<j}
               \left(s_k\bm v_k-s_{k+1}\bm v_{k+1}\right)
 &&\text{by Eq.~\eqref{paper:affine-slope}},\\
 &=-s_j\bm v_j
 &&\text{by telescoping}\qquad\left(0\le j<L\right).
 \end{aligned}
\]
For $j=L$, the same telescope ends with the omitted terminal
term $s_L\bm v_L=0$, giving $\D F\left(\bm c_L\right)=0$.
Convexity makes $\bm c_L$ a minimizer, and the positive
uniform-convexity bound from part~\ref{paper:assembly:class}
makes it unique. Lemma~\ref{global:center-lemma} verifies the
individual phase evaluations and this terminal case.

\textbf{Proof of \ref{paper:assembly:prefix}.}
To prevent other components from offsetting $A_j-I_j\ge7\delta_j$,
we define $\Psi_j:=F-A_j+I_j$ and
$r_j\left(\bm x\right):=A_j\left(\bm x\right)-I_j\left(\bm x\right)$.
The function
$\Psi_j=B+\ell_{\mathrm{init}}+\sum_{k\ne j}A_k+I_j$ is convex
by Lemma~\resultpart{sources:convex-rules}{sources:convex-rules:sum}.
Lemma~\resultpart{paper:phase}{paper:phase:completed} gives $\Psi_j=F$
near $\bm c_{j+1}$, hence equality of their values and gradients there.
For $\bm x\in W_{j,b}$ with $b\le\left\lfloor\frac{m}{4}\right\rfloor$,
define $\bm w=P_{E_j}\left(\bm x-\bm c_j\right)\in U_{j,b}$.
The next coordinate $\ip{\bm v_{j+1}}{\bm x-\bm c_j}$ is $0$
when $j<L-1$. Thus Eq.~\eqref{paper:ambient-phase} and
Lemma~\resultpart{paper:phase}{paper:phase:prefix} give
\[
 r_j\left(\bm x\right)
 =r_j\left(\bm c_j+\bm w\right)\ge7\delta_j.
\]
For $j<L-1$, $\bm x-\bm c_{j+1}\in\bigoplus_{k\le j}E_k$
is perpendicular to $\bm v_{j+1}$. By part~\ref{paper:assembly:centers},
$\ip{\D F\left(\bm c_{j+1}\right)}{\bm x-\bm c_{j+1}}=0$.
The same equality holds for $j=L-1$ because
$\D F\left(\bm c_L\right)=0$. The first-order convexity inequality
therefore yields
\[
 \begin{aligned}
 \Psi_j\left(\bm x\right)
 &\ge\Psi_j\left(\bm c_{j+1}\right)
       +\ip{\D\Psi_j\left(\bm c_{j+1}\right)}{\bm x-\bm c_{j+1}}
 &&\text{by Lemma~\resultpart{sources:convex-rules}{sources:convex-rules:first}},\\
 &=F\left(\bm c_{j+1}\right)
 &&\text{by }\Psi_j=F\text{ nearby and the zero inner product},\\
 F\left(\bm x\right)
 &=\Psi_j\left(\bm x\right)+r_j\left(\bm x\right)
 &&\text{by the definition of }\Psi_j,\\
 &\ge F\left(\bm c_{j+1}\right)+7\delta_j
 &&\text{by the preceding bounds},\\
 &\ge F\left(\bm c_L\right)+7\delta_j
 &&\text{since }\bm c_L\text{ is the minimizer}.
 \end{aligned}
\]
This proves the prefix-gap conclusion for all $\bm x\in W_{j,b}$.

\textbf{Proof of \ref{paper:assembly:initial}.}
Lemmas~\ref{global:exact-drop-lemma} and \ref{global:value-lemma}
give the base and phase value changes, their geometric decay,
and $0<\chi_0\le2\left(1+c_0\right)\lambda m\tau_0^q$.
Summing these decreases gives
\[
 \begin{aligned}
 D_{m,L}
 &=F\left(\bm c_0\right)-F\left(\bm c_L\right)
 =\sum_{j=0}^{L-1}\chi_j
 &&\text{by the definition of }\chi_j,\\
 &=\chi_0\sum_{j=0}^{L-1}\eta^{qj}
 &&\text{by Lemma~\resultpart{global:value-lemma}{global:value-lemma:part:exact-gap}},\\
 &=\chi_0\frac{1-\eta^{qL}}{1-\eta^q}
 &&\text{by the finite geometric sum}.
 \end{aligned}
\]
The value-change formula in Lemma~\ref{global:exact-drop-lemma}
depends only on $s_j,t_j^\star,\tau_j,\delta_j,R_j$ and the fixed
scalar coefficients. These quantities are independent of the
orthonormal phase directions, so both $\chi_0$ and $D_{m,L}$
have the claimed independence.
\end{proof}

%% file: appendix-oracle.tex
\section{The Oracle Transcript and One Fixed Completion}\label{app:paper-oracle}

\subsection{Preservation of earlier oracle outputs}\label{design:consistency}
To apply Lemma~\resultpart{global:gap-lemma}{global:gap-lemma:part:prefix}
to a deterministic algorithm $\mathcal A$ as in Definition~\ref{def:oracle},
we must place its
candidate solutions in short prefixes $W_{j,b}$. We therefore choose the
unfixed chain directions as the algorithm runs, while fixing the
dimension and all scalar parameters beforehand.
For oracle call $s$, denote the evaluation point by $\bm x_s$,
the returned value and derivative tensors by $a_{s,k}$, and the
candidate solution after this call by $\widehat{\bm x}_s$.
We impose no norm bound on these points, and a candidate solution
need not be an evaluation point.
We collect the oracle outputs and available information as
\begin{equation}\label{paper:oracle-data}
 \mathcal O_s=\left(a_{s,0},\ldots,a_{s,p}\right),\qquad
 \mathcal H_s=\left(\left(\bm x_\ell,\mathcal O_\ell\right)\right)_{\ell=1}^{s},
 \qquad \mathcal H_0=\varnothing.
\end{equation}
With fixed parameters and initialization, determinism gives
\begin{equation}\label{paper:oracle-trajectory}
 \bm x_{s+1}=\mathcal A_{s+1}\left(\mathcal H_s\right),\qquad
 \widehat{\bm x}_s=\widehat{\mathcal A}_s\left(\mathcal H_s\right).
\end{equation}
These rules define the evaluation trajectory
$\bm x_1,\ldots,\bm x_{n_mL}$ and candidate solutions
$\widehat{\bm x}_0,\ldots,\widehat{\bm x}_{n_mL}$.
To obtain one fixed objective that reproduces this adaptive trajectory,
we must preserve every recorded value and derivative tensor when
choosing the remaining directions.
For the nonempty completion family $\mathcal F_s$ in
Eq.~\eqref{paper:oracle-consistency}, preservation of the recorded
oracle outputs requires
\begin{equation}\label{design:answer-consistency}
 \D^kF\left(\bm x_\ell\right)=a_{\ell,k}
 \quad\left(F\in\mathcal F_s,\ 1\le\ell\le s,\ 0\le k\le p\right),
\end{equation}
where $\D^0F=F$ and $a_{\ell,0}$ is the returned value.

\subsubsection{Replay and preservation of earlier oracle outputs}
\begin{lemma}[Completing directions gives one fixed oracle instance]\label{design:replay}
Suppose the dimension and all scalar parameters are fixed before directions are chosen adaptively. After each oracle call, suppose the nonempty family $\mathcal F_s$ of final completions satisfies the oracle consistency requirement Eq.~\eqref{design:answer-consistency}. Running the algorithm on any final completion has both properties:
\begin{claims}
\item\label{design:replay:part:transcript} It generates the recorded evaluation points and oracle outputs.
\item\label{design:replay:part:outputs} It generates the recorded deterministic candidate solutions, including the zero-call candidate solution.
\end{claims}
\end{lemma}
\begin{proof}
Fix any final completion $F$. Its oracle agrees with every recorded
oracle output by Eq.~\eqref{design:answer-consistency}. The deterministic
rules in Eq.~\eqref{paper:oracle-trajectory} then reproduce both the
evaluation trajectory and the candidate solutions, as follows.
\begin{itemize}
\item \textbf{Proof of \ref{design:replay:part:transcript}.}
Eq.~\eqref{paper:oracle-trajectory} gives
$\bm x_1=\mathcal A_1\left(\mathcal H_0\right)$ from the fixed input.
If the first $\ell-1$ oracle outputs agree, their transcript
$\mathcal H_{\ell-1}$ is the same, so
$\bm x_\ell=\mathcal A_\ell\left(\mathcal H_{\ell-1}\right)$ is unchanged.
Eq.~\eqref{design:answer-consistency} fixes every component
$a_{\ell,k}$ of $\mathcal O_\ell$. Induction therefore reproduces
the recorded transcript $\mathcal H_s$.
\item \textbf{Proof of \ref{design:replay:part:outputs}.}
Part~\ref{design:replay:part:transcript} fixes each $\mathcal H_s$.
Eq.~\eqref{paper:oracle-trajectory} then gives the same candidate solution
$\widehat{\bm x}_s=\widehat{\mathcal A}_s\left(\mathcal H_s\right)$,
including $s=0$, when $\mathcal H_0=\varnothing$.
\end{itemize}
\end{proof}

\begin{lemma}[A local change preserves earlier oracle outputs of the base]\label{design:allowed-changes}
Suppose all partial bases belong to $C^p$, $B_{j+1}=B_j$ outside an open ball $S_j$, all later changes up to the finite final base are contained in $S_j$, and $\operatorname{dist}\left(\bm x_s,\overline S_j\right)>0$. Then:
\begin{claims}
\item\label{design:allowed-changes:part:neighborhood} The final base $B$ equals $B_j$ on a neighborhood of $\bm x_s$.
\item\label{design:allowed-changes:part:oracle} $\D^kB\left(\bm x_s\right)=\D^kB_j\left(\bm x_s\right)$ for $0\le k\le p$.
\end{claims}
\end{lemma}
\begin{proof}
\begin{itemize}
\item \textbf{Proof of \ref{design:allowed-changes:part:neighborhood}.}
Define $d_s:=\operatorname{dist}\left(\bm x_s,\overline S_j\right)>0$.
For $\bm z\in B\left(\bm x_s,\frac{d_s}{2}\right)$ and $\bm w\in\overline S_j$,
\[
 \begin{aligned}
 \norm{\bm z-\bm w}
 &\ge\norm{\bm x_s-\bm w}-\norm{\bm z-\bm x_s}
 &&\text{by the reverse triangle inequality},\\
 &\ge d_s-\norm{\bm z-\bm x_s}>\frac{d_s}{2}.
 \end{aligned}
\]
Thus this ball is disjoint from every later modification. If $B=B_L$
is the final base, then at every point $\bm z$ of the ball,
\[
 \begin{aligned}
 B\left(\bm z\right)-B_j\left(\bm z\right)
 &=\sum_{k=j}^{L-1}\left(B_{k+1}\left(\bm z\right)-B_k\left(\bm z\right)\right)\\
 &=0,
 \end{aligned}
\]
where the second line uses the assumed exterior agreement at each step.
This proves equality throughout the neighborhood.
\item \textbf{Proof of \ref{design:allowed-changes:part:oracle}.}
Differentiate the neighborhood equality in part~\ref{design:allowed-changes:part:neighborhood}; the case $k=0$ is the equality itself.
\end{itemize}
\end{proof}

\subsection{One fixed instance for deterministic algorithms}\label{oracle:section}

\subsubsection{Two facts for selecting directions}\label{oracle:independence}\label{oracle:arbitrary-points}

\begin{lemma}[One direction includes an arbitrary point]\label{oracle:point}
Let $V\subseteq\R^d$ be a linear subspace and $\bm x\in\R^d$.
At most one additional unit direction perpendicular to $V$ is needed to
include $\bm x$ in the enlarged span.
\end{lemma}
\begin{proof}
To include $\bm x$ in the enlarged span, we isolate its component
perpendicular to $V$. Using the orthogonal projection from
Sec.~\ref{paper:settings}, define $\bm\omega=\bm x-P_V\bm x$.
The projection identity $P_V^2=P_V$ gives
\[
 P_V\bm\omega=P_V\bm x-P_V^2\bm x=0,
\]
so $\bm\omega\perp V$. If $\bm\omega=0$, then $\bm x=P_V\bm x\in V$.
Otherwise $\bm u=\frac{\bm\omega}{\norm{\bm\omega}}$ is a unit vector perpendicular to $V$, and
\[
 \bm x=P_V\bm x+\norm{\bm\omega}\bm u
 \in V\oplus\operatorname{span}\{\bm u\}.
\]
In particular, for $V=W_{j,b}$, the additional chain direction is
chosen by
\begin{equation}\label{paper:point-inclusion}
 \bm\omega=\bm x-P_{W_{j,b}}\bm x,\qquad
 \bm u_{j,b+1}=\frac{\bm\omega}{\norm{\bm\omega}}
 \quad\left(\bm\omega\ne\bm0\right).
\end{equation}
Then $\bm x\in W_{j,b+1}$, so we increase $b$ by one.
If $\bm\omega=\bm0$, no direction is needed.
\end{proof}

\begin{lemma}[One additional direction specifies a prefix oracle output]\label{oracle:answer}
Use the scalar parameters of Definition~\ref{global:function-definition} in Appx.~\ref{app:paper-assembly}. Fix indices $0\le j<L$ and
$0\le b\le\left\lfloor \frac{m}{4}\right\rfloor$, the ordered orthonormal bases of
$E_0,\ldots,E_{j-1}$, and $\bm v_j,\bm u_{j,1},\ldots,\bm u_{j,b+1}$.
The objectives $F,\widetilde F$ arising from any two remaining orthonormal
completions satisfy
\[
 \D^kF\left(\bm x\right)=\D^k\widetilde F\left(\bm x\right)
 \qquad\left(\bm x\in W_{j,b},\ 0\le k\le p\right),\qquad \D^0F:=F.
\]
\end{lemma}
\begin{proof}
Fix $\bm x\in W_{j,b}$. To show that its oracle output is unchanged
by the remaining completion, we check each component of
$F=B+\ell_{\mathrm{init}}+\sum_k A_k$ from
Definition~\ref{global:function-definition}.
The terminal convention is supplied by
Theorem~\resultpart{collar:theorem}{collar:theorem:part:terminal}.
\begin{enumerate}[label=\textbf{Case \arabic*.}]
\item \textbf{Base term $B$.}
The change from $B_j$ to $B_{j+1}$ is confined to $S_j$,
and all subsequent replacement balls lie inside $S_j$ by
Lemma~\resultpart{base:nesting-lemma}{base:nesting-lemma:part:construction}.
Lemma~\resultpart{collar:separation-lemma}{collar:separation-lemma:part:distance}
gives $\operatorname{dist}\left(\bm c_{j+1},W_{j,b}\right)>\theta R_j$,
so $\bm x$ lies outside $\overline S_j$.
Lemma~\resultpart{base:prefix-lemma}{base:prefix-lemma:part:neighborhood}
therefore gives $B=B_j$ near $\bm x$, and hence
\[
 \D^kB\left(\bm x\right)=\D^kB_j\left(\bm x\right)
 \qquad\left(0\le k\le p\right).
\]
These derivatives use only the fixed centers
$\bm c_0,\ldots,\bm c_j$ and Euclidean distances, so they are
independent of all remaining directions.
\item \textbf{Current phase $A_j$.}
The evaluation point has phase coordinates
\[
 \bm w=P_{E_j}\left(\bm x-\bm c_j\right)
 \in\operatorname{span}\{\bm v_j,\bm u_{j,1},\ldots,\bm u_{j,b}\},
 \qquad
 \ip{\bm v_{j+1}}{\bm x-\bm c_j}=0,
\]
with the second coordinate omitted in the terminal phase.
The coordinate definition in Eq.~\eqref{paper:ambient-phase} therefore gives
\[
 \D^kA_j\left(\bm x\right)=\D^kA_j\left(\bm c_j+\bm w\right)
 \qquad\left(0\le k\le p\right).
\]
By Lemma~\ref{collar:prefix-oracle-lemma}, these tensors and the value
depend only on the prefix and its next direction
$\bm u_{j,b+1}$. Fixing this direction determines the full
ambient tensors, including their action on test vectors outside
the selected span.
\item \textbf{Future phases $A_k$, $k>j$.}
Both $\bm x\in W_{j,b}$ and $\bm c_k=\sum_{\ell<k}\bm h_\ell$
belong to $\bigoplus_{\ell<k}E_\ell$. Hence
\[
 P_{E_k}\left(\bm x-\bm c_k\right)=0,\qquad
 \ip{\bm v_{k+1}}{\bm x-\bm c_k}=0,
\]
with the second coordinate omitted for $k=L-1$.
Theorem~\resultpart{collar:theorem}{collar:theorem:part:origin} at these
zero inputs gives
\[
 A_k\left(\bm x\right)
 =-\frac{\lambda\left(m+1\right)\tau_k^q}{q}-4\delta_k,\qquad
 \D^\ell A_k\left(\bm x\right)=0\quad\left(1\le\ell\le p\right).
\]
Thus unchosen future directions enter neither the derivatives nor
the value returned now.
\item \textbf{Completed phases and $\ell_{\mathrm{init}}$.}
For $k<j$, the coordinates of $A_k$ use only the fixed $E_k$ and $\bm v_{k+1}$. The initial affine term uses only the fixed $\bm v_0$.
\end{enumerate}
Adding all components gives the same value and full derivative tensors
for every remaining completion and every $0\le k\le p$.
The base agrees on a neighborhood; for the phase functions, equality
of the value and derivatives at the evaluation point is sufficient for
Eq.~\eqref{design:answer-consistency}.
\end{proof}

\subsubsection{Construct the transcript and complete one fixed instance}

\begin{proposition}[Direction selection and fixed-function realization]\label{oracle:realization}
Before the first oracle call, fix the scalar parameters of Definition~\ref{global:function-definition} in Appx.~\ref{app:paper-assembly}, with
\[
 m,L\in\N,\qquad m\ge32,\qquad L\ge1,\qquad
 d\ge L\left(m+1\right),\qquad n_m:=\left\lfloor \frac{m}{16}\right\rfloor.
\]
For every deterministic algorithm, the ordered bases can be chosen so that:
\begin{claims}
\item\label{oracle:realization:part:prefix} \textbf{Prefix placement.}
$n_m\ge \frac{m}{32}$. For $0\le j<L$ and $jn_m<s\le\left(j+1\right)n_m$, there is
an index $0\le b\le\left\lfloor \frac{m}{4}\right\rfloor$ with
\[
 \bm x_s,\widehat{\bm x}_s\in W_{j,b}.
\]
Also $\widehat{\bm x}_0\in W_{0,b}$ for some $0\le b\le1$.
\item\label{oracle:realization:part:fixed} \textbf{One fixed objective.}
The complete bases define one $F$ by Eq.~\eqref{global:objective}.
Running the algorithm with its oracle reproduces all recorded
evaluation points, oracle outputs, and candidate solutions through
oracle call $n_mL$.
\item\label{oracle:realization:part:gap} \textbf{Gaps at candidate solutions for that same $F$.}
\[
 \begin{aligned}
 F\left(\widehat{\bm x}_s\right)-\min F&\ge7\delta_j
 &&\left(0\le j<L,\ jn_m<s\le\left(j+1\right)n_m\right),\\
 F\left(\widehat{\bm x}_0\right)-\min F&\ge7\delta_0.
 \end{aligned}
\]
\end{claims}
\end{proposition}
\begin{proof}
We prove Proposition~\ref{paper:oracle} together with the prefix-placement
claim. Steps~5.1--5.2 establish its completion-family invariant, and
Step~5.3 gives the fixed objective and candidate gaps.
\begin{itemize}
\item \textbf{Proof of \ref{oracle:realization:part:prefix}.}
Work at initialization zero, with the normalization in
Appx.~\ref{final:rescaling-section} fixed before selecting directions.
The initial gap $D_{m,L}$ is independent of the orthonormal directions
by Lemma~\resultpart{global:value-lemma}{global:value-lemma:part:exact-gap},
so this normalization does not depend on later choices.
When $d>L\left(m+1\right)$, the blocks need not span $\R^d$.
The base starts from $B_0\left(\bm x\right)=\frac{\mu}{q}\norm{\bm x}^q$,
and its replacements use the full ambient Euclidean norm.
Thus its uniform-convexity bound controls the unused directions as well,
by Lemma~\resultpart{paper:base}{paper:base:class}.
\textbf{Step 5.1: Available directions.}\quad
Group $j$ consists of oracle calls $jn_m+1,\ldots,\left(j+1\right)n_m$,
and $b$ counts the chosen current-block chain directions.
To keep every candidate solution within a prefix where the gap bound
applies, we leave room for the zero-call candidate solution and at most
three new chain directions per oracle call. For
$n_m=\left\lfloor\frac{m}{16}\right\rfloor$,
\[
 \begin{aligned}
 1+3n_m
 &<3n_m+2
 \le\frac{3m}{16}+2
 &&\text{since }n_m=\left\lfloor \frac{m}{16}\right\rfloor,\\
 &\le\frac{3m}{16}+\frac{m}{16}=\frac{m}{4}
 &&\text{since }m\ge32.
 \end{aligned}
\]
Since $1+3n_m$ is an integer, this gives
$1+3n_m\le\left\lfloor\frac{m}{4}\right\rfloor$.
Also $\frac{m}{16}\ge2$, so the inequality $\left\lfloor x\right\rfloor\ge \frac{x}{2}$
for $x\ge2$ gives
\[
 n_m=\left\lfloor \frac{m}{16}\right\rfloor\ge\frac{\frac{m}{16}}{2}=\frac{m}{32}.
\]

All choices below are orthonormally completable: whenever a new direction
is required in group $j$,
\[
 \begin{aligned}
 \dim W_{j,b}^{\perp}
 &=d-j\left(m+1\right)-\left(b+1\right)\\
 &\ge\left(L-j\right)\left(m+1\right)-\left(b+1\right)
 &&\text{since }d\ge L\left(m+1\right),\\
 &\ge m-b
 &&\text{since }j\le L-1,\\
 &\ge m-\left(1+3n_m\right)\\
 &\ge\frac{13m}{16}-1>0
 &&\text{since }n_m\le \frac{m}{16},\ m\ge32.
 \end{aligned}
\]
Consequently all link selections admit an orthonormal completion.

\textbf{Step 5.2: Identical oracle outputs.}\quad
\emph{Induction invariant.}
Before oracle call $s+1$, where $0\le s<n_mL$, the current group is
\[
 j=\left\lfloor\frac{s}{n_m}\right\rfloor.
\]
The earlier blocks are complete, and the current prefix is $W_{j,b}$.
Let $\mathcal F_s$ be the family of objectives obtained by completing
the directions selected after $s$ oracle calls, including the direction for the candidate solution
and any preparation for the next group.
Equality of the recorded oracle outputs across these completions
will make Eq.~\eqref{paper:oracle-trajectory} reproduce the trajectory,
allowing us to apply the prefix gaps to its candidate solutions.
We therefore maintain the invariant in
Eq.~\eqref{paper:oracle-consistency}:
\[
 \mathcal F_s\ne\varnothing,\qquad
 \D^kF\left(\bm x_\ell\right)=a_{\ell,k}
 \quad\left(F\in\mathcal F_s,\ 1\le\ell\le s,\ 0\le k\le p\right).
\]
\textbf{Base case: $s=0$.}
Initially, we choose $\bm v_0$, set $b=0$, and use
Eq.~\eqref{paper:point-inclusion} to include $\widehat{\bm x}_0$.
Lemma~\ref{oracle:point} gives
$\widehat{\bm x}_0\in W_{0,b}$ with $b\le1$.
Orthonormal completion gives $\mathcal F_0\ne\varnothing$;
there are no oracle output equalities to check.

\textbf{Induction step: from $s$ oracle calls to $s+1$ oracle calls.}
Assume the invariant for $\mathcal F_s$ and use
$j=\left\lfloor \frac{s}{n_m}\right\rfloor$ as above.
The procedure consists of the following three selections, with $b$
updated after each selection. Denote the families after the
evaluation point and oracle output selections by $\mathcal F^{\rm q},\mathcal F^{\rm a}$.
\begin{enumerate}
\item \textbf{Include the evaluation point.}
The common transcript $\mathcal H_s$ determines
$\bm x_{s+1}=\mathcal A_{s+1}\left(\mathcal H_s\right)$ by
Eq.~\eqref{paper:oracle-trajectory}.
Set $\bm\omega=\bm x_{s+1}-P_{W_{j,b}}\bm x_{s+1}$.
Lemma~\ref{oracle:point} gives
\[
 \begin{cases}
 \mathcal F^{\rm q}:=\mathcal F_s,\quad b\text{ unchanged},&\bm\omega=0,\\
 \mathcal F^{\rm q}:=\{F\in\mathcal F_s:
       \bm u_{j,b+1}=\frac{\bm\omega}{\norm{\bm\omega}}\},\quad b\leftarrow b+1,&\bm\omega\ne0.
 \end{cases}
\]
Thus $\bm x_{s+1}\in W_{j,b}$.

\item \textbf{Determine and return the oracle output.}
The evaluation point already lies in $W_{j,b}$, but the current
chain can involve $\bm u_{j,b+1}$ in its derivative tensors.
We therefore fix one further unit direction $\bm e\perp W_{j,b}$
before returning the oracle output, restricting to
\[
 \mathcal F^{\rm a}:=\{F\in\mathcal F^{\rm q}:\bm u_{j,b+1}=\bm e\}.
\]
Lemma~\ref{oracle:answer}, with the current $b$, gives
\[
 \D^kF\left(\bm x_{s+1}\right)=\D^k\widetilde F\left(\bm x_{s+1}\right)=:a_{s+1,k}
 \quad\left(F,\widetilde F\in\mathcal F^{\rm a},\ 0\le k\le p\right).
\]
Return these values and tensors as $\mathcal O_{s+1}$ and append
$\left(\bm x_{s+1},\mathcal O_{s+1}\right)$ to $\mathcal H_s$,
giving $\mathcal H_{s+1}$ as in
Eq.~\eqref{paper:oracle-data}. Then set $b\leftarrow b+1$.

\item \textbf{Include the candidate solution.}
The enlarged transcript gives
$\widehat{\bm x}_{s+1}=\widehat{\mathcal A}_{s+1}\left(\mathcal H_{s+1}\right)$
by Eq.~\eqref{paper:oracle-trajectory}.
Set $\bm\omega=\widehat{\bm x}_{s+1}-P_{W_{j,b}}\widehat{\bm x}_{s+1}$.
Lemma~\ref{oracle:point} gives
\[
 \begin{cases}
 \mathcal F_{s+1}:=\mathcal F^{\rm a},\quad b\text{ unchanged},&\bm\omega=0,\\
 \mathcal F_{s+1}:=\{F\in\mathcal F^{\rm a}:
       \bm u_{j,b+1}=\frac{\bm\omega}{\norm{\bm\omega}}\},\quad b\leftarrow b+1,&\bm\omega\ne0.
 \end{cases}
\]
Hence $\bm x_{s+1},\widehat{\bm x}_{s+1}\in W_{j,b}$.
Each oracle call adds at most three links. Including the zero-call candidate solution when
$j=0$, the final prefix length is at most $1+3n_m\le\left\lfloor \frac{m}{4}\right\rfloor$.
\end{enumerate}
Further direction choices restrict the completion family, so the
induction hypothesis preserves all earlier oracle outputs.
The second selection supplies the new oracle output, giving
\[
 \begin{gathered}
 \varnothing\ne\mathcal F_{s+1}\subseteq\mathcal F^{\rm a}
 \subseteq\mathcal F^{\rm q}\subseteq\mathcal F_s,\\
 \D^kF\left(\bm x_\ell\right)=a_{\ell,k}
 \quad\left(F\in\mathcal F_{s+1},\ 1\le\ell\le s+1,\ 0\le k\le p\right).
 \end{gathered}
\]

\textbf{Prepare the next oracle call.}
If $s+1<\left(j+1\right)n_m$, retain the current group and prefix.
If $s+1=\left(j+1\right)n_m$, complete $E_j$; when $j<L-1$, choose
$\bm v_{j+1}$ perpendicular to the completed blocks and start group
$j+1$ with $b=0$. There are at least $\left(L-j-1\right)\left(m+1\right)$ dimensions for the remaining
blocks. Restrict $\mathcal F_{s+1}$ to these choices: it remains nonempty,
and the displayed oracle output equalities persist. The recorded evaluation point and
candidate solution remain in their group-$j$ prefix, whose directions are now fixed.
This proves the invariant at $s+1$, including group boundaries.
Induction covers all $n_mL$ oracle calls and completes all bases.

\item \textbf{Proof of \ref{oracle:realization:part:fixed}.}
\textbf{Step 5.3: One fixed objective.}\quad
The complete ordered bases define $F\in\mathcal F_{n_mL}$ by
Eq.~\eqref{global:objective}. The invariant gives
\[
 \D^kF\left(\bm x_s\right)=a_{s,k}
 \qquad\left(1\le s\le n_mL,\ 0\le k\le p\right).
\]
Lemma~\resultpart{design:replay}{design:replay:part:transcript} reproduces
the recorded evaluation points and oracle outputs on this fixed $F$;
part~\ref{design:replay:part:outputs} reproduces every candidate solution,
including $\widehat{\bm x}_0$.
If the algorithm stops early, extend its rule arbitrarily after stopping
to finish the construction; its recorded stopping candidate solution is unchanged.

\item \textbf{Proof of \ref{oracle:realization:part:gap}.}
Parts~\ref{oracle:realization:part:prefix}--\ref{oracle:realization:part:fixed}
and Lemma~\resultpart{global:gap-lemma}{global:gap-lemma:part:prefix} in Appx.~\ref{app:paper-assembly} give,
for each candidate solution in group $j$,
\[
 \widehat{\bm x}_s\in W_{j,b},\quad b\le\left\lfloor \frac{m}{4}\right\rfloor
 \quad\Longrightarrow\quad F\left(\widehat{\bm x}_s\right)-\min F\ge7\delta_j.
\]
For $s=0$, use $j=0$ and $b\le1$.
\end{itemize}
\end{proof}

\Needspace{9\baselineskip}

%% file: appendix-rate.tex
\section{Detailed Parameter Calibration and Critical Rescaling}\label{app:paper-calibration}

\subsection{Integer calibration of the chain length and phase count}
\label{final:parameters}\label{design:counts}

\begin{lemma}[Calibrating the oracle lower bound]\label{final:calibration}
Let $a_q=\frac{3q-2}{2}$ and $\kappa=\frac{H}{\sigma}$.
With the constants of Lemma~\ref{paper:assembly}, define
\[
 \begin{aligned}
 c_m&=\left(\frac{qc_B}{C_H}\right)^{\frac{1}{a_q}},&
 \kappa_{p,\nu}&=\left(\frac{32}{c_m}\right)^{a_q},\\
 c_*&=\frac{7c_1\left(1-\eta^q\right)}{2\left(1+c_0\right)},&
 c_L&=\frac{1}{q\log\left(\frac{1}{\eta}\right)+\log\left(\frac{2}{c_*}\right)}.
 \end{aligned}
\]
Suppose $\kappa\ge\kappa_{p,\nu}$ and
$0<\epsilon\le\frac{c_*\Delta_0}{2}$.
Take $\tau_0=1$ and define
\[
 \mu=\frac{\sigma}{qc_B},\qquad
 m=\left\lfloor c_m\kappa^{\frac{1}{a_q}}\right\rfloor,\qquad
 L=1+\left\lfloor
 \frac{\log\left(\frac{c_*\Delta_0}{2\epsilon}\right)}
 {q\log\left(\frac{1}{\eta}\right)}
 \right\rfloor.
\]
Then, with $n_m=\left\lfloor\frac{m}{16}\right\rfloor$,
\begin{claims}
\item\label{final:calibration:part:chain}
\textbf{Smoothness, uniform convexity, and per-phase cost.}
$m\ge32$, $qc_B\mu=\sigma$, $C_H\mu m^{a_q}\le H$, and
$n_m\ge\frac{c_m}{32}\kappa^{\frac{1}{a_q}}$.
\item\label{final:calibration:part:phases}
\textbf{Normalized phase gaps and phase count.}
For $0\le j<L$,
\[
 \begin{gathered}
 \frac{7\delta_j}{D_{m,L}}\ge c_*\eta^{qj},\qquad
 c_*\Delta_0\eta^{q\left(L-1\right)}\ge2\epsilon,\\
 L\asymp_{p,\nu}\log\left(\frac{\Delta_0}{\epsilon}\right).
 \end{gathered}
\]
\item\label{final:calibration:part:dimension}
\textbf{Ambient dimension.}
For $K_{p,\nu}:=2c_m\left(1+\frac{1}{q\log\left(\frac{1}{\eta}\right)}\right)$,
\[
 L\left(m+1\right)
 \le K_{p,\nu}\left(\frac{H}{\sigma}\right)^{\frac{2}{3q-2}}
 \log\left(\frac{\Delta_0}{\epsilon}\right).
\]
\end{claims}
\end{lemma}
\begin{proof}

\textbf{Proof of \ref{final:calibration:part:chain}.}
The initial scale $\tau_0$ will be adjusted by rescaling, so we
first take $\tau_0=1$ and define $\kappa=\frac{H}{\sigma}$.
By Lemma~\resultpart{paper:assembly}{paper:assembly:class},
the theorem's two conditions hold provided
$qc_B\mu\ge\sigma$ and $C_H\mu m^{a_q}\le H$.
Equating the first coefficient to $\sigma$ allows the largest
chain length under the second inequality. Indeed, for
$\mu=\frac{\sigma}{qc_B}$, the smoothness requirement becomes
\[
 \begin{aligned}
 C_H\mu m^{a_q}\le H
 &\quad\Longleftrightarrow\quad
 m^{a_q}\le\frac{H}{C_H\mu}\\
 &\quad\Longleftrightarrow\quad
 m^{a_q}\le\frac{qc_B}{C_H}\frac{H}{\sigma}\\
 &\quad\Longleftrightarrow\quad
 m\le c_m\kappa^{\frac{1}{a_q}}
 &&\text{by the definitions of }c_m,\kappa.
 \end{aligned}
\]
Our choice
$m=\left\lfloor c_m\kappa^{\frac{1}{a_q}}\right\rfloor$
is the largest integer allowed by this bound.
Proposition~\ref{global:parameters} in Appx.~\ref{app:paper-assembly}
supplies compatible phase coefficients for this choice.
For $\kappa\ge\kappa_{p,\nu}$, define
$x:=c_m\kappa^{\frac{1}{a_q}}\ge32$.
Then $m=\left\lfloor x\right\rfloor\ge32$ and, for all $\bm u,\bm v\in\R^d$,
\[
 \begin{aligned}
 \norm{\D^pF\left(\bm u\right)-\D^pF\left(\bm v\right)}
 &\le C_H\mu m^{a_q}\norm{\bm u-\bm v}^{\nu}\\
 &\le H\norm{\bm u-\bm v}^{\nu}
 &&\text{by the chosen }m,\\
 \Breg F{\bm v}{\bm u}
 &\ge c_B\mu\norm{\bm v-\bm u}^q
 =\frac{\sigma}{q}\norm{\bm v-\bm u}^q
 &&\text{since }\mu=\frac{\sigma}{qc_B}.
 \end{aligned}
\]
Thus Proposition~\ref{paper:oracle} applies, with the per-phase
number of calls bounded by
\[
 \begin{aligned}
 n_m
 &=\left\lfloor\frac{\left\lfloor x\right\rfloor}{16}\right\rfloor\\
 &=\left\lfloor\frac{x}{16}\right\rfloor\\
 &\ge\frac{x}{32}
 &&\text{since }\frac{x}{16}\ge2,\ \left\lfloor t\right\rfloor\ge \frac{t}{2},\\
 &=\frac{c_m}{32}\kappa^{\frac{2}{3q-2}}
 &&\text{since }\frac{1}{a_q}=\frac{2}{3q-2}.
 \end{aligned}
\]
For the second equality, define $n=\left\lfloor \frac{x}{16}\right\rfloor$.
Then $16n\le x<16\left(n+1\right)$, and the integer endpoints imply
$16n\le\left\lfloor x\right\rfloor<16\left(n+1\right)$.
Dividing by $16$ and taking floors gives
$\left\lfloor\frac{\left\lfloor x\right\rfloor}{16}\right\rfloor=n$.

\textbf{Proof of \ref{final:calibration:part:phases}.}
The choices above meet the smoothness and uniform-convexity bounds
and fix the number of calls required in each phase, independently
of $L$. We can therefore increase the total oracle cost by taking
as many phases as possible while keeping their candidate gaps
above $\epsilon$ after normalization.
Rescaling to $\Delta_0$ multiplies each candidate gap by
$\frac{\Delta_0}{D_{m,L}}$, so we compare $7\delta_j$ from
Proposition~\ref{paper:oracle} to $D_{m,L}$.

\emph{Compare the phase gaps with the initial gap.}
Lemma~\resultpart{global:value-lemma}{global:value-lemma:part:bounds} in Appx.~\ref{app:paper-assembly}
and part~\ref{global:value-lemma:part:exact-gap} give, for every $L\ge1$,
\[
 \begin{aligned}
 D_{m,L}
 &=\sum_{j=0}^{L-1}\chi_j
 &&\text{by Eq.~\eqref{global:initial}},\\
 &\le2\left(1+c_0\right)\lambda m\sum_{j=0}^{L-1}\tau_j^q
 &&\text{by Eq.~\eqref{global:drop}},\\
 &=2\left(1+c_0\right)\lambda m\tau_0^q\sum_{j=0}^{L-1}\eta^{qj}\\
 &\le\frac{2\left(1+c_0\right)}{1-\eta^q}\lambda m\tau_0^q
 &&\text{since }0<\eta^q<1,\\
 \frac{7\delta_j}{D_{m,L}}
 &=\frac{7c_1\lambda m\tau_0^q\eta^{qj}}{D_{m,L}}
 &&\text{by the definition of }\delta_j,\\
 &\ge\frac{7c_1\left(1-\eta^q\right)}{2\left(1+c_0\right)}\eta^{qj}
 =c_*\eta^{qj}
 &&\left(0\le j<L\right).
 \end{aligned}
\]
Thus $c_*$ is independent of $L$, and $0<c_*<\frac{7}{128}<1$ because $c_1<\frac{1}{64}$.
Lemma~\ref{final:rescaling} verifies the stated effect of rescaling.
Since $0<\eta<1$, a last-phase gap above $2\epsilon$
also keeps all earlier gaps above $2\epsilon$, as
\[
 \frac{7\delta_j}{D_{m,L}}\Delta_0
 \ge c_*\Delta_0\eta^{qj}
 \ge c_*\Delta_0\eta^{q\left(L-1\right)}\qquad\left(0\le j<L\right).
\]

\emph{Choose the largest phase count keeping the last gap above $2\epsilon$.}
For $\epsilon\le\frac{c_*\Delta_0}{2}$, solving
$c_*\Delta_0\eta^{q\left(L-1\right)}\ge2\epsilon$
for the largest integer $L$ gives
\[
 \begin{aligned}
 c_*\Delta_0\eta^{q\left(L-1\right)}\ge2\epsilon
 &\quad\Longleftrightarrow\quad
 q\left(L-1\right)\log\eta
 \ge\log\left(\frac{2\epsilon}{c_*\Delta_0}\right)\\
 &\quad\Longleftrightarrow\quad
 L-1\le
 \frac{\log\left(\frac{c_*\Delta_0}{2\epsilon}\right)}
      {q\log\left(\frac{1}{\eta}\right)}
 &&\text{since }\log\eta<0.
 \end{aligned}
\]
Denoting the upper bound by $u$, this choice is
\[
 u:=\frac{\log\left(\frac{c_*\Delta_0}{2\epsilon}\right)}{q\log\left(\frac{1}{\eta}\right)}\ge0,
 \qquad L:=1+\left\lfloor u\right\rfloor\ge1.
\]
Indeed, $L-1\le u$ gives
\[
 \eta^{q\left(L-1\right)}\ge\eta^{qu}
 =\exp\left(-\log\left(\frac{c_*\Delta_0}{2\epsilon}\right)\right)
 =\frac{2\epsilon}{c_*\Delta_0},
\]
which proves the lower bound on the last normalized phase gap.

\emph{Verify the logarithmic phase count.}
This choice also gives $L\ge c_L\log\left(\frac{\Delta_0}{\epsilon}\right)$.
To verify it, define $z:=\log\left(\frac{\Delta_0}{\epsilon}\right)$
and $b_*:=\log\left(\frac{2}{c_*}\right)>0$.
Then $z\ge b_*$, $u=\frac{z-b_*}{q\log\left(\frac{1}{\eta}\right)}<L$, and $L\ge1$ imply
\[
 \begin{aligned}
 z&=q\log\left(\frac{1}{\eta}\right)u+b_*\\
  &<q\log\left(\frac{1}{\eta}\right)L+b_*\\
  &\le\left(q\log\left(\frac{1}{\eta}\right)+b_*\right)L=\frac{L}{c_L}.
 \end{aligned}
\]
Hence $L\ge c_L\log\left(\frac{\Delta_0}{\epsilon}\right)$, proving the phase-count bound.
For the corresponding upper bound, $u\le\frac{z}{q\log\left(\frac{1}{\eta}\right)}$
and $z\ge b_*$ give
\[
 L=1+\left\lfloor u\right\rfloor
 \le1+u
 \le\left(\frac{1}{b_*}+
 \frac{1}{q\log\left(\frac{1}{\eta}\right)}\right)z.
\]
Thus $L\asymp_{p,\nu}\log\left(\frac{\Delta_0}{\epsilon}\right)$,
as used in Step~6.

\textbf{Proof of \ref{final:calibration:part:dimension}.}
The direction-selection argument requires $L\left(m+1\right)$
dimensions. With $x=c_m\kappa^{\frac{1}{a_q}}\ge32$ and
$z=\log\left(\frac{\Delta_0}{\epsilon}\right)$ as above,
the accuracy restriction gives
$z\ge\log\left(\frac{2}{c_*}\right)>1$ because $c_*<\frac{7}{128}$.
Consequently,
\[
 \begin{aligned}
 m+1
 &\le x+1\le2x
 &&\text{since }m=\left\lfloor x\right\rfloor,\ x\ge32,\\
 L
 &\le1+\frac{z}{q\log\left(\frac{1}{\eta}\right)}
 \le\left(1+\frac{1}{q\log\left(\frac{1}{\eta}\right)}\right)z
 &&\text{since }z>1,\\
 L\left(m+1\right)
 &\le2c_m\left(1+\frac{1}{q\log\left(\frac{1}{\eta}\right)}\right)
 \kappa^{\frac{1}{a_q}}z\\
 &=K_{p,\nu}\left(\frac{H}{\sigma}\right)^{\frac{2}{3q-2}}
 \log\left(\frac{\Delta_0}{\epsilon}\right)
 &&\text{by }\frac{1}{a_q}=\frac{2}{3q-2}.
 \end{aligned}
\]
\end{proof}

To use one constant in the accuracy range and oracle lower bound of
Theorem~\ref{paper:main-theorem}, choose
$c_{p,\nu}:=\min\left\{\frac{c_*}{2},\frac{c_m c_L}{32}\right\}$.
The first term ensures $\epsilon\le\frac{c_*\Delta_0}{2}$;
the second is the coefficient in the product of the two count bounds.

\subsection{Critical rescaling}\label{final:rescaling-section}

The choice of $L$ compares the phase gaps with the target
$\epsilon$ for an initial gap of $\Delta_0$, whereas the current
objective has initial gap $D_{m,L}$. We now rescale the objective
to make this comparison apply.
Eq.~\eqref{global:initial} determines $D_{m,L}$ independently
of the chain directions, so we can fix
\[
 a=\left(\frac{\Delta_0}{D_{m,L}}\right)^{\frac{1}{q}}
\]
before applying the trajectory construction.
The following lemma verifies that this normalization adjusts the
initial gap while preserving both parameters because $q=p+\nu$.

\begin{lemma}[Rescaling preserves the critical class]\label{final:rescaling}
Let $F$ be the objective in Definition~\ref{global:function-definition}.
For $a>0$ and an initialization $\bm x_0$, define
\[
 F_a\left(\bm x\right):=a^qF\left(\frac{\bm x}{a}\right),\qquad
 \widehat F\left(\bm x\right):=F_a\left(\bm x-\bm x_0\right).
\]
For $0\le k\le p$, with $\D^0F=F$,
\[
 \D^k\widehat F\left(\bm x\right)
 =a^{q-k}\D^kF\left(\frac{\bm x-\bm x_0}{a}\right).
\]
Every $p^{th}$-order H\"older smoothness parameter $H$ of $F$ is preserved:
for all $\bm x,\bm y\in\R^d$,
\begin{equation}\label{paper:critical-normalization}
\begin{aligned}
 \norm{\D^p\widehat F\left(\bm x\right)-\D^p\widehat F\left(\bm y\right)}
 &\le H\norm{\bm x-\bm y}^{\nu},\\
 \Breg{\widehat F}{\bm y}{\bm x}
 &=a^q\Breg F{\frac{\bm y-\bm x_0}{a}}{\frac{\bm x-\bm x_0}{a}}.
\end{aligned}
\end{equation}
Every admissible degree-$q$ uniform-convexity modulus is preserved, and
\[
 \operatorname*{argmin}_{\bm x\in\R^d}\widehat F\left(\bm x\right)=\{\bm x_0+a\bm c_L\},
 \qquad \widehat F\left(\bm x_0\right)-\min\widehat F=a^qD_{m,L}.
\]
In particular, $a=\left(\frac{\Delta_0}{D_{m,L}}\right)^{\frac{1}{q}}$
gives $\widehat F\left(\bm x_0\right)-\min\widehat F=\Delta_0$.
\end{lemma}
\begin{proof}
Define $\bm u=\frac{\bm x-\bm x_0}{a}$ and
$\bm v=\frac{\bm y-\bm x_0}{a}$.
For $1\le k\le p$ and arbitrary test vectors $\bm z_1,\ldots,\bm z_k$,
\[
 \begin{aligned}
 \D^k\widehat F\left(\bm x\right)\left[\bm z_1,\ldots,\bm z_k\right]
 &=a^q\D^kF\left(\bm u\right)\left[\frac{\bm z_1}{a},\ldots,\frac{\bm z_k}{a}\right]
 &&\text{by the linear chain rule},\\
 &=a^{q-k}\D^kF\left(\bm u\right)\left[\bm z_1,\ldots,\bm z_k\right]
 &&\text{by multilinearity}.
 \end{aligned}
\]
This proves the full derivative identity; $k=0$ is the definition of
$\widehat F$. For any smoothness parameter $H$ of $F$,
\begin{align*}
 \norm{\D^p\widehat F\left(\bm x\right)-\D^p\widehat F\left(\bm y\right)}
 &=a^{q-p}\norm{\D^pF\left(\bm u\right)-\D^pF\left(\bm v\right)}
 &&\text{by the derivative identity above},\\
 &\le a^{q-p}H\norm{\bm u-\bm v}^{\nu}
 &&\text{by the H\"older smoothness assumption on }F,\\
 &=a^{q-p-\nu}H\norm{\bm x-\bm y}^{\nu}
 &&\text{since }\bm u-\bm v=\frac{\bm x-\bm y}{a},\\
 &=H\norm{\bm x-\bm y}^{\nu}
 &&\text{since }q=p+\nu.
\end{align*}
For uniform convexity, we first expand the Bregman divergence.
Using $\bm y-\bm x=a\left(\bm v-\bm u\right)$ gives
\[
 \begin{aligned}
 \Breg{\widehat F}{\bm y}{\bm x}
 &=a^q\left(F\left(\bm v\right)-F\left(\bm u\right)\right)
   -a^{q-1}\ip{\D F\left(\bm u\right)}{\bm y-\bm x}
 &&\text{by the definition and chain rule},\\
 &=a^q\left(F\left(\bm v\right)-F\left(\bm u\right)
       -\ip{\D F\left(\bm u\right)}{\bm v-\bm u}\right)
 &&\text{since }\bm y-\bm x=a\left(\bm v-\bm u\right),\\
 &=a^q\Breg F{\bm v}{\bm u}
 &&\text{by the Bregman definition}.
 \end{aligned}
\]
If $F$ has modulus $\sigma$, this identity gives
\[
 \begin{aligned}
 \Breg{\widehat F}{\bm y}{\bm x}
 &=a^q\Breg F{\bm v}{\bm u}
 &&\text{by the identity above},\\
 &\ge a^q\frac{\sigma}{q}\norm{\bm v-\bm u}^q
 &&\text{by uniform convexity of }F,\\
 &=\frac{a^q}{a^q}\frac{\sigma}{q}\norm{\bm y-\bm x}^q
 &&\text{since }\bm v-\bm u=\frac{\bm y-\bm x}{a},\\
 &=\frac{\sigma}{q}\norm{\bm y-\bm x}^q.
 \end{aligned}
\]
The affine change of variables is bijective and $a^q>0$.
Lemma~\resultpart{global:center-lemma}{global:center-lemma:part:minimizer} therefore gives
\[
 \begin{aligned}
 \operatorname*{argmin}_{\bm x}\widehat F\left(\bm x\right)
 &=\left\{\bm x_0+a\bm u:
       \bm u\in\operatorname*{argmin}_{\bm u}F\left(\bm u\right)\right\}\\
 &=\{\bm x_0+a\bm c_L\}.
 \end{aligned}
\]
Using this minimizer and the definition of $D_{m,L}$,
\[
 \begin{aligned}
 \widehat F\left(\bm x_0\right)-\min\widehat F
 &=a^qF\left(\bm0\right)-a^qF\left(\bm c_L\right)\\
 &=a^q\left(F\left(\bm0\right)-F\left(\bm c_L\right)\right)
 &&\text{by factoring }a^q,\\
 &=a^qD_{m,L}
 &&\text{by Eq.~\eqref{global:initial}}.
 \end{aligned}
\]
For the normalizing scale fixed before the trajectory construction,
the last expression is
\[
 \begin{aligned}
 a^qD_{m,L}
 &=\left(\frac{\Delta_0}{D_{m,L}}\right)D_{m,L}
 &&\text{by the choice of }a,\\
 &=\Delta_0
 &&\text{since }D_{m,L}>0
   \text{ by Lemma~\ref{global:value-lemma}\,\ref{global:value-lemma:part:bounds}}.
 \end{aligned}
\]
\end{proof}